\documentclass{CSML}

\def\dOi{11(3:24)2015}
\lmcsheading%
{\dOi}
{1--76}
{}
{}
{Dec.~10, 2014}
{Sep.~30, 2015}
{}

\ACMCCS{[{\bf Theory of computation}]: Models of
computation---Probabilistic computation /
Quantum computation theory; Semantics and reasoning---Program
semantics—Categorical semantics}

\newif\ifignore 
\ignorefalse

\newcommand{\auxproof}[1]{
\ifignore\mbox{}\newline
\textit{Proof:} \dotfill\newline
{\it #1}\mbox{}\newline
\qed
\fi}

\usepackage[english]{babel}
\usepackage{amsmath}
\usepackage{amssymb}
\usepackage{xspace}
\usepackage{stmaryrd}
\usepackage{cmll}
\usepackage{url,hyperref}
\usepackage{fancybox}
\usepackage{xypic}
\usepackage{xy}
\xyoption{v2}
\xyoption{curve}
\xyoption{2cell}
\UseTwocells

\usepackage{url}
\usepackage{amssymb}
\usepackage{amstext}
\usepackage{xspace}

\newdimen\proofrulebreadth \proofrulebreadth=.05em
\newdimen\proofdotseparation \proofdotseparation=1.25ex
\newdimen\proofrulebaseline \proofrulebaseline=2ex
\newcount\proofdotnumber \proofdotnumber=3
\let\then\relax
\def\hfi{\hskip0pt plus.0001fil}
\mathchardef\squigto="3A3B
\newif\ifinsideprooftree\insideprooftreefalse
\newif\ifonleftofproofrule\onleftofproofrulefalse
\newif\ifproofdots\proofdotsfalse
\newif\ifdoubleproof\doubleprooffalse
\let\wereinproofbit\relax
\newdimen\shortenproofleft
\newdimen\shortenproofright
\newdimen\proofbelowshift
\newbox\proofabove
\newbox\proofbelow
\newbox\proofrulename
\def\shiftproofbelow{\let\next\relax\afterassignment\setshiftproofbelow\dimen0 }
\def\shiftproofbelowneg{\def\next{\multiply\dimen0 by-1 }%
\afterassignment\setshiftproofbelow\dimen0 }
\def\setshiftproofbelow{\next\proofbelowshift=\dimen0 }
\def\setproofrulebreadth{\proofrulebreadth}

\def\prooftree{
\ifnum  \lastpenalty=1
\then   \unpenalty
\else   \onleftofproofrulefalse
\fi
\ifonleftofproofrule
\else   \ifinsideprooftree
        \then   \hskip.5em plus1fil
        \fi
\fi
\bgroup
\setbox\proofbelow=\hbox{}\setbox\proofrulename=\hbox{}%
\let\justifies\proofover\let\leadsto\proofoverdots\let\Justifies\proofoverdbl
\let\using\proofusing\let\[\prooftree
\ifinsideprooftree\let\]\endprooftree\fi
\proofdotsfalse\doubleprooffalse
\let\thickness\setproofrulebreadth
\let\shiftright\shiftproofbelow \let\shift\shiftproofbelow
\let\shiftleft\shiftproofbelowneg
\let\ifwasinsideprooftree\ifinsideprooftree
\insideprooftreetrue
\setbox\proofabove=\hbox\bgroup$\displaystyle 
\let\wereinproofbit\prooftree
\shortenproofleft=0pt \shortenproofright=0pt \proofbelowshift=0pt
\onleftofproofruletrue\penalty1
}

\def\eproofbit{
\ifx    \wereinproofbit\prooftree
\then   \ifcase \lastpenalty
        \then   \shortenproofright=0pt  
        \or     \unpenalty\hfil         
        \or     \unpenalty\unskip       
        \else   \shortenproofright=0pt  
        \fi
\fi
\global\dimen0=\shortenproofleft
\global\dimen1=\shortenproofright
\global\dimen2=\proofrulebreadth
\global\dimen3=\proofbelowshift
\global\dimen4=\proofdotseparation
\global\count255=\proofdotnumber
$\egroup  
\shortenproofleft=\dimen0
\shortenproofright=\dimen1
\proofrulebreadth=\dimen2
\proofbelowshift=\dimen3
\proofdotseparation=\dimen4
\proofdotnumber=\count255
}

\def\proofover{
\eproofbit 
\setbox\proofbelow=\hbox\bgroup 
\let\wereinproofbit\proofover
$\displaystyle
}%
\def\proofoverdbl{
\eproofbit 
\doubleprooftrue
\setbox\proofbelow=\hbox\bgroup 
\let\wereinproofbit\proofoverdbl
$\displaystyle
}%
\def\proofoverdots{
\eproofbit 
\proofdotstrue
\setbox\proofbelow=\hbox\bgroup 
\let\wereinproofbit\proofoverdots
$\displaystyle
}%
\def\proofusing{
\eproofbit 
\setbox\proofrulename=\hbox\bgroup 
\let\wereinproofbit\proofusing
\kern0.3em$
}

\def\endprooftree{
\eproofbit 
  \dimen5 =0pt
\dimen0=\wd\proofabove \advance\dimen0-\shortenproofleft
\advance\dimen0-\shortenproofright
\dimen1=.5\dimen0 \advance\dimen1-.5\wd\proofbelow
\dimen4=\dimen1
\advance\dimen1\proofbelowshift \advance\dimen4-\proofbelowshift
\ifdim  \dimen1<0pt
\then   \advance\shortenproofleft\dimen1
        \advance\dimen0-\dimen1
        \dimen1=0pt
        \ifdim  \shortenproofleft<0pt
        \then   \setbox\proofabove=\hbox{%
                        \kern-\shortenproofleft\unhbox\proofabove}%
                \shortenproofleft=0pt
        \fi
\fi
\ifdim  \dimen4<0pt
\then   \advance\shortenproofright\dimen4
        \advance\dimen0-\dimen4
        \dimen4=0pt
\fi
\ifdim  \shortenproofright<\wd\proofrulename
\then   \shortenproofright=\wd\proofrulename
\fi
\dimen2=\shortenproofleft \advance\dimen2 by\dimen1
\dimen3=\shortenproofright\advance\dimen3 by\dimen4
\ifproofdots
\then
        \dimen6=\shortenproofleft \advance\dimen6 .5\dimen0
        \setbox1=\vbox to\proofdotseparation{\vss\hbox{$\cdot$}\vss}%
        \setbox0=\hbox{%
                \advance\dimen6-.5\wd1
                \kern\dimen6
                $\vcenter to\proofdotnumber\proofdotseparation
                        {\leaders\box1\vfill}$%
                \unhbox\proofrulename}%
\else   \dimen6=\fontdimen22\the\textfont2 
        \dimen7=\dimen6
        \advance\dimen6by.5\proofrulebreadth
        \advance\dimen7by-.5\proofrulebreadth
        \setbox0=\hbox{%
                \kern\shortenproofleft
                \ifdoubleproof
                \then   \hbox to\dimen0{%
                        $\mathsurround0pt\mathord=\mkern-6mu%
                        \cleaders\hbox{$\mkern-2mu=\mkern-2mu$}\hfill
                        \mkern-6mu\mathord=$}%
                \else   \vrule height\dimen6 depth-\dimen7 width\dimen0
                \fi
                \unhbox\proofrulename}%
        \ht0=\dimen6 \dp0=-\dimen7
\fi
\let\doll\relax
\ifwasinsideprooftree
\then   \let\VBOX\vbox
\else   \ifmmode\else$\let\doll=$\fi
        \let\VBOX\vcenter
\fi
\VBOX   {\baselineskip\proofrulebaseline \lineskip.2ex
        \expandafter\lineskiplimit\ifproofdots0ex\else-0.6ex\fi
        \hbox   spread\dimen5   {\hfi\unhbox\proofabove\hfi}%
        \hbox{\box0}%
        \hbox   {\kern\dimen2 \box\proofbelow}}\doll%
\global\dimen2=\dimen2
\global\dimen3=\dimen3
\egroup 
\ifonleftofproofrule
\then   \shortenproofleft=\dimen2
\fi
\shortenproofright=\dimen3
\onleftofproofrulefalse
\ifinsideprooftree
\then   \hskip.5em plus 1fil \penalty2
\fi
}

\renewcommand{\arraystretch}{1.3}
\newtheorem{theorem}{Theorem}

\newtheorem{assumption}{Assumption}[]

\newenvironment{myproof}[1][Proof]%
   { \begin{trivlist}%
     \item[\hskip \labelsep {\it #1.}]%
   }%
   { \end{trivlist}%
   }

\newcommand{\after}{\mathrel{\circ}}
\newcommand{\set}[2]{\{#1\;|\;#2\}}
\newcommand{\setin}[3]{\{#1\in#2\;|\;#3\}}
\newcommand{\conjun}{\mathrel{\wedge}}
\newcommand{\disjun}{\mathrel{\vee}}
\newcommand{\andthen}[2]{\ensuremath{\tuple{#1?}(#2)}}
\renewcommand{\implies}[2]{\ensuremath{[#1?](#2)}}

\newcommand{\all}[2]{\forall{#1}.\,#2}
\newcommand{\allin}[3]{\forall{#1\in#2}.\,#3}

\newcommand{\ex}[2]{\exists{#1}.\,#2}

\newcommand{\lam}[2]{\lambda#1.\,#2}
\newcommand{\lamin}[3]{\lambda#1\in#2.\,#3}
\newcommand{\tuple}[1]{\langle#1\rangle}

\newcommand{\downset}{\mathop{\downarrow}}
\newcommand{\upset}{\mathop{\uparrow}\!}
\newcommand{\NNO}{\mathbb{N}}
\newcommand{\C}{\mathbb{C}}
\newcommand{\R}{\mathbb{R}}
\newcommand{\Z}{\mathbb{Z}}
\newcommand{\ket}[1]{\ensuremath{|{\kern.1em}#1{\kern.1em}\rangle}}
\newcommand{\bra}[1]{\langle\,#1\,|}
\newcommand{\upstate}{\ensuremath{\mathop{\uparrow}}}
\newcommand{\downstate}{\ensuremath{\mathop{\downarrow}}}
\newcommand{\newstate}{\ensuremath{\textsl{new}}}
\newcommand{\measure}{\ensuremath{\textsl{measure}}}
\newcommand{\B}{\mathcal{B}}

\newcommand{\bigovee}{\mathop{\vphantom{\sum}\mathchoice%
        {\vcenter{\hbox{\huge $\ovee$}}}%
        {\vcenter{\hbox{\Large $\ovee$}}}%
        {\ovee}{\ovee}}\displaylimits}

\newcommand{\Dst}{\mathcal{D}}

\newcommand{\Giry}{\mathcal{G}}
\newcommand{\Rad}{\mathcal{R}}
\newcommand{\Exp}{\mathcal{E}}
\newcommand{\Lift}{\mathcal{L}}

\newcommand{\Pow}{\mathcal{P}}
\newcommand{\Powfin}{\mathcal{P}_{\mathit{fin}}}

\newcommand{\idmap}[1][]{\ensuremath{\mathrm{id}_{#1}}}

\newcommand{\op}[1]{#1^{\textrm{op}}}
\newcommand{\orthogonal}{\mathrel{\bot}}
\newcommand{\supp}{\textsl{supp}}

\newcommand{\ev}{\textsl{ev}}
\newcommand{\Hom}{\textsl{Hom}}
\newcommand{\st}{\ensuremath{\textsl{st}}}

\newcommand{\tr}{\ensuremath{\textsl{tr}}}
\newcommand{\isup}{\ensuremath{\textsl{isup}}}
\newcommand{\charac}{\ensuremath{\textsl{char}}}
\newcommand{\indic}[1]{\mathbf{1}_{#1}}
\newcommand{\intd}{{\kern.2em}\mathrm{d}{\kern.03em}}

\newcommand{\instr}{\ensuremath{\textsl{instr}}}
\renewcommand{\swap}{\ensuremath{\textsl{swap}}}
\newcommand{\FstAnd}{\ensuremath{\textsl{FstAnd}}}
\newcommand{\FstThen}{\ensuremath{\textsl{FstThen}}}
\newcommand{\SndAnd}{\ensuremath{\textsl{SndAnd}}}
\newcommand{\SndThen}{\ensuremath{\textsl{SndThen}}}
\newcommand{\La}[1]{\ensuremath{[#1, 0]}}
\newcommand{\Ra}[1]{\ensuremath{[#1, 1]}}
\renewcommand{\wp}{\ensuremath{\textsl{wp}}}

\newcommand{\Mat}{\ensuremath{\textsl{M}}}

\newcommand{\dst}{\ensuremath{\textsl{dst}}}
\newcommand{\sotimes}{\mathrel{\raisebox{.05pc}{$\scriptstyle \otimes$}}}
\newcommand{\sodot}{\mathrel{\raisebox{.05pc}{$\scriptstyle \odot$}}}

\newcommand{\inprod}[2]{\ensuremath{\langle #1\,|\,#2 \rangle}}

\newcommand{\cat}[1]{\ensuremath{\textbf{#1}}\xspace}
\newcommand{\Sets}{\ensuremath{\textbf{Sets}}\xspace}

\newcommand{\CH}{\ensuremath{\textbf{CH}}\xspace}
\newcommand{\EA}{\ensuremath{\textbf{EA}}\xspace}
\newcommand{\BA}{\ensuremath{\textbf{BA}}\xspace}
\newcommand{\DL}{\ensuremath{\textbf{DL}}\xspace}

\newcommand{\EMod}{\ensuremath{\textbf{EMod}}\xspace}
\newcommand{\Mon}{\ensuremath{\textbf{Mon}}\xspace}
\newcommand{\CMon}{\ensuremath{\textbf{CMon}}\xspace}

\newcommand{\Ab}{\ensuremath{\textbf{Ab}}\xspace}
\newcommand{\Conv}{\ensuremath{\textbf{Conv}}\xspace}
\newcommand{\Meas}{\ensuremath{\textbf{Meas}}\xspace}

\newcommand{\CL}{\ensuremath{\textbf{CL}}\xspace}
\newcommand{\Hilb}{\ensuremath{\textbf{Hilb}}\xspace}
\newcommand{\FdHilb}{\ensuremath{\textbf{FdHilb}}\xspace}
\newcommand{\Vect}{\ensuremath{\textbf{Vect}}\xspace}

\newcommand{\Rng}{\ensuremath{\textbf{Rng}}\xspace}

\newcommand{\CRng}{\ensuremath{\textbf{CRng}}\xspace}

\newcommand{\Mod}{\ensuremath{\textbf{Mod}}\xspace}

\newcommand{\Cliq}{\ensuremath{\textsl{Clq}}\xspace}

\newcommand{\Act}{\ensuremath{\textbf{Act}}\xspace}

\newcommand{\Cstar}{\ensuremath{\textbf{Cstar}}\xspace}
\newcommand{\CstarMap}[1]{\ensuremath{\Cstar_{#1}}\xspace}
\newcommand{\CstarPU}{\CstarMap{\textrm{PU}}}
\newcommand{\CstarCPU}{\CstarMap{\textrm{CPU}}}
\newcommand{\CstarMIU}{\CstarMap{\textrm{MIU}}}

\newcommand{\CCstarMap}[1]{\ensuremath{\textbf{CCstar}_{#1}}\xspace}
\newcommand{\CCstarPU}{\CCstarMap{\textrm{PU}}}
\newcommand{\CCstarMIU}{\CCstarMap{\textrm{MIU}}}

\newcommand{\Coh}{\ensuremath{\textbf{Coh}}\xspace}
\newcommand{\CohMap}[1]{\ensuremath{\Coh_{#1}}\xspace}
\newcommand{\CohL}{\CohMap{\textrm{L}}}

\newcommand{\Pred}{\ensuremath{\textsl{Pred}}\xspace}

\newcommand{\Stat}{\ensuremath{\textsl{Stat}}\xspace}

\newcommand{\scalar}{\mathrel{\bullet}}

\newcommand{\inv}{\mathop{\rlap{\raisebox{.3ex}{${\kern.6ex}\cdot$}}-}}
\newcommand{\simop}{\mathop{\sim}}
\newcommand{\polar}{\ensuremath{\mathrel{\smash{\stackrel{\raisebox{-.7em}{$|$}}{\raisebox{-.3em}{$\sim$}}}}}}

\renewcommand{\H}{\ensuremath{\mathcal{H}}}
\newcommand{\K}{\ensuremath{\mathcal{K}}}

\newcommand{\Kl}{\mathcal{K}{\kern-.2ex}\ell}

\newcommand{\KlL}[1]{{#1}_{+1}}
\newcommand{\EM}{\mathcal{E}{\kern-.2ex}\mathcal{M}}

\newcommand{\Ef}{\ensuremath{\mathcal{E}{\kern-.5ex}f}}
\newcommand{\DM}{\ensuremath{\mathcal{D}{\kern-.85ex}\mathcal{M}}}
\newcommand{\klafter}{\mathrel{\raisebox{.15em}{$\scriptscriptstyle\circledcirc$}}}

\newcommand{\QEDbox}{\square}
\newcommand{\QED}{\hspace*{\fill}$\QEDbox$}

\newcommand{\conglongrightarrow}{\mathrel{\smash{\stackrel{
           \raisebox{.5ex}{$\scriptstyle\cong$}}{
           \raisebox{0ex}[0ex][0ex]{$\longrightarrow$}}}}}

\newcommand{\pullback}[1][dr]{\save*!/#1-1.2pc/#1:(-1,1)@^{|-}\restore}
\makeatother
 \newdir{ >}{{}*!/-7.5pt/@{>}}
 \newdir{|>}{!/4.5pt/@{|}*:(1,-.2)@^{>}*:(1,+.2)@_{>}}
 \newdir{ |>}{{}*!/-3pt/@{|}*!/-7.5pt/:(1,-.2)@^{>}*!/-7.5pt/:(1,+.2)@_{>}}
\newcommand{\xyline}[2][]{\ensuremath{\smash{\xymatrix@1#1{#2}}}}
\newcommand{\xyinline}[2][]{\ensuremath{\smash{\xymatrix@1#1{#2}}}^{\rule[8.5pt]{0pt}{0pt}}}
\makeatletter

\newcommand{\IV}{\raisebox{-.2em}{$\xy 
(0,2)*{\cdot}; 
(1,2)*{\cdot};
(2,2)*{\cdot};
(0,0)*{\cdot};
(1,0)*{\cdot};
{\xline(0,2);(0,0)};
{\xline(1,2); (1,0)};
{\xline(2,2); (1,0)};
\endxy$}}

\newcommand{\XI}{\raisebox{-.2em}{$\xy
(0,2)*{\cdot};
(1,2)*{\cdot};
(2,2)*{\cdot};
(0,0)*{\cdot};
(1,0)*{\cdot};
{\xline(0,2); (1,0)};
{\xline(1,2); (0,0)};
{\xline(2,2); (1,0)};
\endxy$}}

\begin{document}

\title[New Directions in Categorical Logic]{New Directions in Categorical Logic, \\
   for Classical, Probabilistic and Quantum Logic}

\author[B.~Jacobs]{Bart Jacobs}
\address{
Institute for Computing and Information Sciences, 
Radboud University Nijmegen, The Netherlands. }
\urladdr{www.cs.ru.nl/B.Jacobs} 
\email{bart@cs.ru.nl}

\date{}

\maketitle

\begin{abstract}
Intuitionistic logic, in which the double negation law $\neg\neg P =
P$ fails, is dominant in categorical logic, notably in topos
theory. This paper follows a different direction in which double
negation does hold, especially in quantitative logics for
probabilistic and quantum systems. The algebraic notions of effect
algebra and effect module that emerged in theoretical physics form the
cornerstone. It is shown that under mild conditions on a category, its
maps of the form $X \rightarrow 1+1$ carry such effect module
structure, and can be used as predicates. Maps of this form $X
\rightarrow 1+1$ are identified in many different situations, and
capture for instance ordinary subsets, fuzzy predicates in a
probabilistic setting, idempotents in a ring, and effects (positive
elements below the unit) in a $C^*$-algebra or Hilbert space.

In quantum foundations the duality between states and effects
(predicates) plays an important role. This duality appears in the form
of an adjunction in our categorical setting, where we use maps $1
\rightarrow X$ as states. For such a state $\omega$ and a predicate
$p$, the validity probability $\omega\models p$ is defined, as an
abstract Born rule. It captures many forms of (Boolean or
probabilistic) validity known from the literature.

Measurement from quantum mechanics is formalised categorically in
terms of `instruments', using L{\"u}ders rule in the quantum
case. These instruments are special maps associated with predicates
(more generally, with tests), which perform the act of measurement and
may have a side-effect that disturbs the system under
observation. This abstract description of side-effects is one of the
main achievements of the current approach. It is shown that in the
special case of $C^*$-algebras, side-effects appear exclusively in the
non-commutative (properly quantum) case. Also, these instruments are
used for test operators in a dynamic logic that can be used for
reasoning about quantum programs/protocols.

The paper describes four successive assumptions, towards a categorical
axiomatisation of quantitative logic for probabilistic and quantum
systems, in which the above mentioned elements occur.
\end{abstract}

\section{Introduction}\label{IntroSec}

Mathematical logic started in the 19th century with George Boole's
\emph{Laws of Thought}. Since then `Boolean algebras' form the
(algebraic) basis of formal logic. In the 20th century these Boolean
algebras have been generalised, mainly in two directions, see
Figure~\ref{BAFig}. Of these two directions the one on the left is
most familiar and dominant. It has given rise to the edifice of topos
theory, where the logic naturally associated with geometric structures
is intuitionistic --- in the sense that the double-negation law
$\neg\neg P = P$ fails. Also, it has given rise to Martin-L\"of's
intuitionistic type theory~\cite{MartinLof84,NordstromPS90}, which
forms the foundation for the theorem prover Coq (and other such
computerised formal systems). In this tradition it is claimed that the
natural logic of computation is intuitionistic.

\begin{figure}
\label{BAFig}
$$\vcenter{\xymatrix@C-2.5pc@R-.5pc{
& & 
\ovalbox{\bf \begin{tabular}{c} categorical \\[-.2em] 
           quantum logic?  \end{tabular}}
    \\
\ovalbox{\bf\begin{tabular}{c} topos theory \\[-.2em] 
          type theory \end{tabular}}
& & 
\ovalbox{\bf \begin{tabular}{c} effect algebras \\[-.2em] 
           logic \& probability \end{tabular}}\ar[u]
    \\
\ovalbox{\bf\begin{tabular}{c} Intuitionistic logic \\[-.2em] 
          Heyting algebra \end{tabular}}\ar[u]
& & 
\ovalbox{\bf \begin{tabular}{c} Quantum logic \\[-.2em] 
           Orthomodular lattice \end{tabular}}\ar[u]
    \\ 
& \ovalbox{\bf \begin{tabular}{c} Boolean \\[-.2em] 
                     logic/algebra \end{tabular}}
   \ar[ur]_(0.6){\begin{array}{c} 
        \;\mbox{\footnotesize drop distributivity} \\[-.3em]
        \mbox{\footnotesize keep double negation} \end{array}}
   \ar[ul]^(0.6){\begin{array}{c} 
        \;\mbox{\footnotesize drop double negation} \\[-.3em]
        \mbox{\footnotesize keep distributivity} \end{array}}
}}$$
\caption{Generalisations of Boolean logic}
\end{figure}

However, there are generalisations of Boolean logic where
double-negation does hold, like linear logic or quantum logic. In this
tradition negation is often written as superscript $P^{\perp}$ and
called orthocomplement. The logic of the quantum world has first been
formalised by von Neumann and Birkhoff in terms of orthomodular
lattices~\cite{Neumann32,BirkhoffN36} (see also the
monograph~\cite{Kalmbach83}). But these orthomodular lattices are
rather awkward mathematical structures that have not led to a formal
calculus of propositions (or predicates) capturing relevant logical
phenomena in physics.

A second generalisation, started in the 1990s, involves so-called
effect algebras~\cite{FoulisB94} and effect modules, introduced in
various forms (\textit{e.g.}~as D-posets~\cite{ChovanecK07} or as weak
orthoalgebras in~\cite{GiuntiniG94}) and by various people (see
\cite{DvurecenskijP00}), notably by theoretical physicists Foulis,
Bennet, Gudder and Pulmannov{\'a}, who published their work in the
physics literature, largely outside the range of perception of
traditional logicians. These effect algebras are partial commutative
monoids with an orthocomplement $(-)^\perp$. Examples come from
classical logic (Boolean algebras), quantum logic in a Hilbert space
or $C^*$-algebra, probability theory (the unit interval $[0,1]$) and
measure theory (step functions, taking finitely many values). Effect
algebras form a common generalisation of structures from logic and
probability theory, making it possible to see both negation $P \mapsto
\neg P$ in logic and opposite probability $r\mapsto 1-r$ as instances
of the same orthocomplement operation.  Partial commutative monoids
have appeared in programming language semantics early on, notably in
the partially additive structures of~\cite{ArbibM86}, and also in
models of linear logic~\cite{Haghverdi00b}. But orthocomplement is
typical of the more restrictive quantum context. 

As an aside: the name `effect' in `effect algebra' refers to the
observer effect\footnote{To be clear: we are concerned with the
  \emph{observer} effect, on the system that is being measured, and
  not with the possible effect on the measurement apparatus.} that may
arise in a quantum setting, where the act of observation may disturb
the system under observation,
see~\cite{Kraus83,BuschLM96,Busch09}. There is no historical
connection with the use of the word `effect' in the context of
computer programming with monads (going back
to~\cite{Moggi91a}). However, there is a similarity in meaning. Both
usages of the word `effect' --- in physics and in computer science ---
refer to a `side-effect', that is, to a change in an underlying state
space.

Effect modules are effect algebras with an additional scalar
multiplication, where the scalars are typically taken from the unit
interval $[0,1]$. Intuitively, these effect modules can be seen as
vector spaces, defined not over the real or complex numbers, but over
$[0,1]$. They are the algebraic counterparts of quantitative logics,
where predicates can be multiplied with a scaling factor. This
includes fuzzy predicates and effects in Hilbert spaces or
$C^*$-algebras. In fact, we present a general construction for effect
modules of predicates in a category, see
Proposition~\ref{CoprodEModProp}. This forms an important part of our
development of the right leg in Figure~\ref{BAFig}.

The basic theory of effect algebras and modules has been developed by
several physicists who (understandably) lacked a background in modern
logic, (program) semantics, and categorical logic. Especially the
unfamiliarity with the modern methods of categorical logic seriously
hampered, or even blocked, progress in this field. As a result these
founders failed to identify certain fundamental relationships and
similarities with programming semantics and logic (notably the
categorical dualities between convex sets and effect
modules~\cite{JacobsM12a}, described in quantum theoretic terms as the
duality between states and effects, or between Schr{\"o}dinger's and
Heisenberg's view, see~\cite[Chapter~2]{HeinosaariZ12}). Here this
duality is described in terms of an adjunction $\op{\EMod}
\rightleftarrows \Conv$ between categories of effect modules and
convex sets, see Figure~\ref{TriangleFig} and
Proposition~\ref{ConvEModAdjProp}.

One of the great contributions in the left leg of Figure~\ref{BAFig}
is the notion of a topos (see \textit{e.g.}~\cite{MacLaneM92}), which
can be seen as a categorical formulation of intuitionistic set theory
and logic. One wonders: is there a comparably important (categorical)
notion in the right leg? Is there an analogue of a topos for the
quantum world, incorporating a quantitative logic? This question forms
the main motivation for the research presented in this paper. We do
not claim to have the definitive answer at this stage, but we do think
that we identify some key properties of such an analogue, belonging to
the right leg in Figure~\ref{BAFig}.

In topos theory predicates on an object $X$ can equivalently be
described as subobjects of $X$ and as maps $X \rightarrow \Omega$,
where $\Omega$ is a special object called the subobject classifier.
This $\Omega$ carries the structure of a Heyting algebra. Here we
shall use maps of the form $p\colon X \rightarrow 1+1$ as predicates,
where $1$ is the final object and $+$ is coproduct. Such predicates
have a classical flavour, since there is an associated orthocomplement
$p^{\perp}$ obtained by swapping the output options. We write this as
$p^{\perp} = [\kappa_{2}, \kappa_{1}] \after p$, where the $\kappa_{i}
\colon 1 \rightarrow 1+1$ are coprojections and $[-,-]$ is cotupling.
Double negation $p^{\perp\perp} = p$ holds by construction. 

In the ordinary set-theoretic world there is a one-to-one
correspondence between subsets $S\subseteq X$ of a set $X$ and
characteristic functions $X \rightarrow 1+1$ on $X$. Both these
notions can be used as interpretation of predicates. This one-to-one
correspondence fails in many other settings. One can say that the
emphasis in the left leg in Figure~\ref{BAFig} lies on the `spatial'
interpretation of predicates as subsets. The path on the right in
Figure~\ref{BAFig} rests, so we claim, on the view of predicates as
characteristic maps.

The two main discoveries wrt.\ these predicates as characteristic maps
are:
\begin{enumerate}
\item in general, these predicates $X \rightarrow 1+1$ do not form a
  Boolean algebra, as one might have expected, but an effect algebra
  --- or more specifically, an effect module;

\item these predicates $X \rightarrow 1+1$ look rather simple, but,
  interpreted in different categories, they naturally capture various
  descriptions in a common framework, including for instance ordinary
  set-theoretic predicates $X \rightarrow \{0,1\}$, fuzzy predicates
  $X \rightarrow [0,1]$, effects in a Hilbert space or $C^*$-algebra,
  idempotents in a ring, complementable elements in a distributive
  lattice, \textit{etc.}
\end{enumerate}

\noindent For instance, for a ring $R$ we show that predicates $R
\rightarrow 1+1$, in the opposite of the category of rings, correspond
to idempotents in $R$, and form an effect algebra. When the ring $R$
is commutative, this effect algebra is actually a Boolean algebra.
The fact that idempotents in a \emph{commutative} ring form a Boolean
algebra is well-known (see~\cite{Johnstone82}). But the fact that
these idempotents form an effect algebra in the non-commutative case
is new (as far as we know).

There are four categories that will be used as leading examples. We
briefly mention them now, but refer to the Appendix for more details.
\begin{enumerate}
\item The category $\Sets$ of sets and functions, as model of
  deterministic computation with Boolean logic.

\item The Kleisli category $\Kl(\Dst)$ of the distribution monad
  $\Dst$ on the category $\Sets$, as model for \emph{discrete}
  probabilistic computation and quantitative (fuzzy) logic.

\item The Kleisli category $\Kl(\Giry)$ of the Giry monad $\Giry$ on
  the category $\Meas$ of measurable spaces; this category
  $\Kl(\Giry)$ is used for \emph{continuous} probabilistic computation
  and quantitative logic.

\item The opposites $\op{(\CstarPU)}$ and $\op{(\CstarCPU)}$ of the
  categories of $C^*$-algebras and (completely) positive unital maps,
  as models for quantum computation and logic. The (common)
  subcategory of \emph{commutative} $C^*$-algebras and positive unital
  maps captures the probabilistic case.
\end{enumerate}

\begin{figure}
\begin{center}
\hspace*{-.5em}\begin{tabular}{|c||c|c|c|}
\hline
\quad\textbf{category }\qquad & 
   \quad\textbf{predicate } $X \stackrel{p}{\rightarrow} 1+1$ \quad & 
   \quad\textbf{state } $1 \stackrel{\omega}{\rightarrow} X$ \quad & 
   \quad\textbf{validity } $\omega\models p$ 
   \vrule height5mm depth3mm width0mm \\
\hline\hline
$\Sets$ & 
   subset $p\subseteq X$ &
   element $\omega\in X$ &
   $\omega\in p$ 
   \vrule height5mm depth3mm width0mm \\
\hline
$\Kl(\Dst)$ & 
   fuzzy $p\colon X \rightarrow [0,1]$ &
   distribution $\omega\in \Dst(X)$ &
   $\sum_{x\in X} p(x)\cdot\omega(x)$ 
   \vrule height5mm depth3mm width0mm \\
\hline
$\Kl(\Giry)$ & 
   measurable $p\colon X \rightarrow [0,1]$ &
 $\begin{array}{c}
   \mbox{probability measure} \\[-.4em]
   \omega\in \Giry(X) \end{array}$ &
   $\int p \intd \omega$ 
   \vrule height5mm depth3mm width0mm \\
\hline
$\op{\DL}$ & 
   complementable $p\in X$ &
   prime filter $\omega\subseteq X$ &
   $p\in\omega$ 
   \vrule height5mm depth3mm width0mm \\
\hline
$\op{\BA}$ & 
   $p\in X$ &
   ultrafilter $\omega\subseteq X$ &
   $p\in\omega$ 
   \vrule height5mm depth3mm width0mm \\
\hline
$\op{\Rng}$ & 
   idempotent $p\in X$ &
   $\omega \colon X \rightarrow \Z$ &
   $\omega(p)$ 
   \vrule height5mm depth3mm width0mm \\
\hline
$\op{(\CstarPU)}$ & 
   $p\in X$ with $0 \leq p \leq 1$ &
   $\omega\colon X \rightarrow \C$ &
   $\omega(p)$ 
   \vrule height5mm depth3mm width0mm \\
\cline{2-4}
$\begin{array}{c} \mbox{special case } X = \\[-.4em]
   \B(\H)\in\op{(\CstarPU)}, \\[-.4em]
   \mbox{for } \H\in\FdHilb \end{array}$  & 
$\begin{array}{c}
   \H\stackrel{p}{\rightarrow} \H \\[-.4em]
   \mbox{with } 0 \leq p \leq \idmap \end{array}$ &
$\begin{array}{c} \mbox{density matrix} \\[-.4em] 
    \omega\in\DM(\H) \end{array}$ &
   $\tr(\omega p)$ 
   \vrule height5mm depth3mm width0mm \\
\hline
\end{tabular}
\end{center}
\caption{Examples of predicates $p$, states $\omega$, and validity
  $(\omega\models p) = p \after \omega$ in different categories, where
  $\DL$, $\BA$, $\Rng$, and $\FdHilb$ are the categories of
  distributive lattices, Boolean algebras, rings, and of
  finite-dimensional Hilbert spaces respectively. Validity
  $\omega\models p$ is either Boolean, in $\{0,1\}$, or probabilistic,
  in $[0,1]$, depending on what the scalars $1 \rightarrow 1+1$ are in
  the category at hand.}
\label{ValidityFig}
\end{figure}

\noindent In our axiomatisation of key aspects of these (and a few other)
categories we start from predicates as maps of the form $X \rightarrow
1+1$. We then presents four successive assumptions on a category that
are important, so we think, for the right leg in
Figure~\ref{BAFig}. These assumptions on a category $\cat{B}$ can be
summarised as follows.
{
\begin{enumerate}
\item The category $\cat{B}$ is an `effectus', that is, it has a final
  object $1$ and finite coproducts $0,+$ satisfying some mild
  properties. Then we prove that the predicates $X \rightarrow 1+1$
  form an effect module. These predicates are related, via an
  adjunction, to states, which are maps of the form $1\rightarrow X$
  in $\cat{B}$. Moreover, for these predicates $p$ and states $\omega$
  we can define validity $\omega\models p$ via an abstract Born rule
  $p\after \omega$, producing a scalar. This scalar may for instance
  live in $\{0,1\}$ or in $[0,1]$, corresponding to Boolean and
  probabilistic validity respectively. The table in
  Figure~\ref{ValidityFig} gives an overview of the different forms
  that this uniform definition of validity $\models$ takes in
  different categories. The duality between states and effects is an
  intrinsic part of our framework. The situation can be summarised in
  a ``state-and-effect'' triangle as in Figure~\ref{TriangleFig}.
  Many more such triangles are described in~\cite{Jacobs15b}.

\item For each $n$-test, consisting of $n$ predicates $p_{i} \colon X
  \rightarrow 1+1$ which add up to 1, there is a (chosen)
  \emph{instrument} map $\instr_{p} \colon X \rightarrow n\cdot X = X
  + \cdots + X$ in $\cat{B}$ that performs measurement. These maps are
  a categorical formalisation of (discrete) instruments in quantum
  theory, see~\cite{DaviesL70,Ozawa84,HeinosaariZ12}.  Intuitively, an
  instrument map sends an element in $X$ to a suitable spread over the
  $n$ coproduct options in the output $X+\cdots+X$, determined by the
  predicates $p_i$. Such a map $\instr_{p}$ may have a side-effect
  (aka.\ observer effect) on the system under observation, describing
  the state change resulting form measurement/observation. Such state
  changes are typical for the quantum world.  We show that these
  side-effects do not occur in set-theoretic or probabilistic models,
  including commutative $C^*$-algebras. However, side-effects do arise
  in proper quantum models, given by non-commutative $C^*$-algebras.

Via these instruments $\instr_{p}$ a guarded command for a programming
language and test operators for a dynamic logic can be defined. There
are two of these test operations, labelled as `test-andthen' and as
`test-then'.  Both of them may have side-effects (in quantum
models). Thus, side-effects have a prominent role in our framework.

\item The category $\cat{B}$ has tensors $\otimes$ which interact
  appropriately with the previous two points. In particular, the unit
  of this tensor is the final object $1$, giving a tensor with
  projections $X \leftarrow X\otimes Y \rightarrow Y$. These
  projections are used for weakening on predicates and for marginal /
  partial trace computations on states. They can be used to pinpoint
  dependence in probabilistic and entanglement in quantum models.

\item Finally, an object $Q$ is assumed together with two states and
  an appropriate predicate that capture quantum bits (qubits). This is
  rather standard, but is required to model quantum computations.
\end{enumerate}}

\begin{figure}
$$\vcenter{\xymatrix{
\op{(\EMod_{M})}\ar@/^1.5ex/[rr]^-{\Hom(-,M)} & \top & 
   \Conv_{M}\ar@/^1.5ex/[ll]^-{\Hom(-,M)} \\
& \cat{B}\ar[ul]^{\Hom(-,1+1)=\Pred\quad}\ar[ur]_{\;\Stat=\Hom(1,-)} &
}}$$
\caption{State-and-effect triangle for an effectus $\cat{B}$, with
  predicate and state functors, and with scalars $M = \Pred(1) =
  \Stat(1+1) = \Hom(1,1+1)$ in $\cat{B}$.}
\label{TriangleFig}
\end{figure}

\noindent The first of these three points apply to set-theoretic,
probabilistic, and to quantum models. Only the last, fourth point
holds exclusively for quantum models (esp.\ $C^*$-algebras).

Thus, one of the achievements of this paper is the uniform description
of these set-theoretic, probabilistic and quantum models. Quantum
mechanics is often described in terms of Hilbert spaces, see textbooks
like~\cite{HeinosaariZ12}, but
also~\cite{AbramskyC04,AbramskyC09}. However, it is the `algebraic'
formulation of quantum mechanics in term of $C^*$-algebras that fits
in the current uniform framework. This algebraic approach stems from
Heisenberg's matrix mechanics picture, see
\textit{e.g.}~\cite{Landsman09b,Strocchi05}.

This switch to $C^*$-algebras was one of three steps that were crucial
in the development of this work. We list these steps explicitly.
\begin{enumerate}
\item Not to use Hilbert spaces but $C^*$-algebras as (categorical)
  model for quantum computation.

\item Not to use $C^*$-algebras with standard *-homomorphisms
  (preserving multiplication, involution and unit), but with
  (completely) positive unital maps. The latter preserve the unit and
  positive elements, and, if they are completely positive, preserve
  positive elements under every extension.

\item Not to use the standard direction of maps between $C^*$-algebra,
  but to work in the opposite category, and understand computations
  (maps) between $C^*$-algebras as predicate transformers, working
  backwards.
\end{enumerate}

\noindent These three points come together in the simple but crucial
observation: maps of the form $A \rightarrow 1+1$,
\textit{i.e.}~predicates, in the \emph{opposite} of the category of
$C^*$-algebras with \emph{(completely) positive unital} maps are in
one-to-one correspondence with elements $a\in A$ satisfying $0 \leq a
\leq 1$; these elements are called effects, and form an effect
module. Also, instruments for $C^*$-algebras are given by L{\"u}ders
rule, see~\eqref{MeasurementCstarEqn}, and form a completely positive
map $A \rightarrow A + \cdots + A$, again in the opposite category.

From duality theory (see~\cite{Johnstone82}) it is already known that
some categories occur most naturally in dual form. A well-known
example is the category $\cat{cHA}$ of complete Heyting algebras
(aka.\ frames), whose opposite category $\op{\cat{cHA}}$ even has a
separate name, namely the category $\cat{Loc}$ of locales. Similarly,
categories of $C^*$-algebras occur naturally in opposite form, as
categories of non-commutative geometries. This is in line with
Heisenberg's picture where computations act in opposite
direction. Here we describe this phenomenon more generally as: maps
between `algebraic models of logics' can be understood as predicate
transformers acting in the opposite direction. This applies for
instance to categories of Boolean algebras, distributive lattices, or
rings (like in Figure~\ref{ValidityFig}).

This paper is organised as follows. Background information about
effect algebras and modules is provided in Section~\ref{PrelimSec},
and about our running probabilistic and quantum example categories in
the Appendix. Predicates of the form $X \rightarrow 1+1$ in various
categories are described in Section~\ref{PredSec}. Their effect
structure is investigated subsequently in Section~\ref{EAPredSec}.
States, as maps $1 \rightarrow X$, are added to the picture in
Section~\ref{StatesSec}. At that stage we can give our semantic
framework involving the satisfaction relation $\models$ (see the
overview in Figure~\ref{ValidityFig}), predicate and state
transformers, and the state-and-effect triangle in
Figure~\ref{TriangleFig}. Section~\ref{PredicateCoprodSec} explains a
side-topic, namely how predicates $X+Y\rightarrow 1+1$ on coproducts
$X+Y$ correspond (via an isomorphism of effect modules) to pairs of
predicates on $X$ and on $Y$. This general result includes familiar
isomorphisms, like $\Pow(X+Y) \cong \Pow(X) \times \Pow(Y)$.

Section~\ref{MeasurementSec} introduces the next major contribution,
namely a categorical formalisation of the notion of measurement
instrument. Via these instruments we express what it means when a
predicate (or more generally, a test) has a side-effect. This is
studied more systematically in Section~\ref{SefSec}, where
side-effect-freeness is related to commutativity in $C^*$-algebras.
This gives mathematical expression to the idea that observations may
disturb, which is typical for the quantum setting, as modelled by
non-commutative $C^*$-algebras. These instruments are used to perform
the action in a dynamic logic, involving `andthen' and `then' test
operators $\andthen{p}{q}$ and $\implies{p}{q}$, for predicates
$p,q$. These operators are defined in the general setting provided by
our axiomatisation. Instantiations in our running examples give
familiar logical constructions, like conjunction and implication,
multiplication and Reichenbach implication, or the sequential
operations from~\cite{GudderG02,GudderN02}. These logical test
operators $\andthen{p}{q}$ and $\implies{p}{q}$ satisfy some of the
standard properties from dynamic logic, but they represent a new
approach in that they capture side-effects of the action $p$, via the
associated instrument that is used in their definition.

In Sections~\ref{TensorSec} and~\ref{QuantumStateSec} the third and
fourth assumptions are discussed, involving tensor products $\otimes$
and a special object $Q$ for qubits. The tensor product $\otimes$
comes with projections, which lead to weakening as predicate
transformer and to marginals as state transformer. The composition
operation on scalars turns out to be commutative in the presence of
tensors --- via a more-or-less standard, Eckmann-Hilton style argument
--- and allow us to pinpoint dependence and entanglement. We conclude
in Section~\ref{QuantumStateSec} by elaborating the superdense coding
example in our newly developed setting.

\section{Preliminaries on effect algebras, effect modules, and convex 
   sets}\label{PrelimSec}

This paper uses categorical language to organise its material. It
assumes a basic level of familiarity with category theory including
(co)product and monoidal structure, adjunctions, and monads (including
Kleisli and Eilenberg-Moore categories), see
\textit{e.g.}~\cite{Awodey06,MacLane71}. Also, it uses examples
involving the distribution and Giry monad $\Dst$ and $\Giry$, Hilbert
spaces, and $C^*$-algebras. Background information about these
examples is relegated to the appendix. It may be consulted on a
call-by-need basis.

This preliminary section does describe the basics of effect algebras,
effect modules, and convex sets. This material is less standard, and
an essential part of the approach presented here.

In order to define an effect algebra, we need the notion of partial
commutative monoid (PCM). Before reading the definition of PCM, think
of the unit interval $[0,1]$ with addition $+$. This $+$ is obviously
only a partial operation, which is commutative and associative in a
suitable sense. This will be formalised next.

\begin{defi}
\label{PCMDef}
A \emph{partial commutative monoid} (PCM) consists of a set $M$ with a zero
element $0\in M$ and a partial binary operation $\ovee\colon M\times
M\rightarrow M$ satisfying the three requirements below. They involve
the notation $x\orthogonal y$ for: $x\ovee y$ is defined; in that case
$x,y$ are called \emph{orthogonal}.
\begin{enumerate}
\item Commutativity: $x\orthogonal y$ implies $y\orthogonal x$ and
$x\ovee y = y\ovee x$;

\item Associativity: $y\orthogonal z$ and $x \orthogonal (y\ovee z)$
implies $x\orthogonal y$ and $(x\ovee y) \orthogonal z$ and also
$x \ovee (y\ovee z) = (x\ovee y)\ovee z$;

\item Zero: $0\orthogonal x$ and $0\ovee x = x$.
\end{enumerate}
\end{defi}

Later on we shall also use `joint orthogonality' of elements $x_{1},
\ldots, x_{n}$. This means that the sum $x_{1} \ovee \cdots \ovee
x_{n}$ is defined. The associativity property means that we can write
such expressions $x_{1} \ovee \cdots \ovee x_{n}$ without brackets.

The notion of effect algebra extends a PCM with an orthocomplement
$(-)^{\perp}$. It is due to~\cite{FoulisB94}, and in slightly
different form also to~\cite{ChovanecK07} and~\cite{GiuntiniG94}. The
reference~\cite{DvurecenskijP00} gives an overview.

\begin{defi}
\label{EffAlgDef}
An \emph{effect algebra} is a PCM $(E, 0, \ovee)$ with an
orthocomplement, that is, with a (total) unary operation $(-)^{\perp}
\colon E\rightarrow E$, such that:
\begin{itemize}
\item $x^{\perp}\in E$ is the unique element in $E$ with $x\ovee
  x^{\perp} = 1$, where $1 = 0^\perp$;

\item $x\orthogonal 1 \Rightarrow x=0$.
\end{itemize}

\noindent For such an effect algebra one defines,
$$\begin{array}{rclcrcl}
x \leq y
& \Longleftrightarrow &
\ex{z}{x\ovee z = y}
& \qquad\mbox{and}\qquad &
x\ominus y
& = &
(x^{\perp} \ovee y)^{\perp}, \mbox{ for } x \geq y.
\end{array}$$

\noindent We shall shortly see, in
Lemma~\ref{EffAlgBasicsLem}~\eqref{EffAlgBasicsLemCanc}, that there is
at most one element $z$ with $x\ovee z = y$.

A homomorphism $E\rightarrow D$ of effect algebras is given by a
function $f\colon E\rightarrow D$ between the underlying sets
satisfying $f(1) = 1$, and if $x\orthogonal x'$ in $E$ then both $f(x)
\orthogonal f(x')$ in $D$ and $f(x\ovee x') = f(x) \ovee f(x')$.
Effect algebras and their homomorphisms form a category, written as
\EA.
\end{defi}

The standard example of an effect algebra is the unit interval.  In
$[0,1]$ one has $r \orthogonal s$ iff $r + s \leq 1$, and in that case
$r\ovee s = r+s$. The orthocomplement is $r^{\perp} = 1-r$, obviously
with $r \orthogonal (1-r)$ and $r \ovee (1-r) = 1$. But the 2-element
set $2 = \{0,1\}$ is also an example of an effect algebra --- it is in
fact the initial one. Hence both Booleans and probabilities form
instances of the notion of effect algebra.

When writing $x\ovee y$ we shall implicitly assume that $x\ovee y$ is
defined, \textit{i.e.}~that $x,y$ are orthogonal: $x\orthogonal
y$. Similarly for $\owedge$. The notation $\ovee$ might suggest that
this operation is idempotent, but this is not the case (and not
intended): for instance, in the unit interval $[0,1]$ one has
$\frac{1}{2} \ovee \frac{1}{2} = 1 \neq \frac{1}{2}$. But more
generally, $x\ovee x = x$ only holds for $x = 0$, since $x \ovee x = x
= x \ovee 0$ yields $x = 0$ by cancellativity, see
Lemma~\ref{EffAlgBasicsLem}~\eqref{EffAlgBasicsLemCanc} below.

Aside: an alternative axiomatisation of the unit interval exists in
terms of MV-algebras \cite{Chang58,Mundici11}, in which the sum
operation is forced to be total, via truncation: $\min(r+s, 1)$. The
partiality of $\ovee$ in effect algebras may look strange at first,
but turns out to be quite natural\footnote{Interestingly, in his
  \emph{An Investigation of the Laws of Thought} from 1854 George
  Boole himself considered disjunction to be a partial operation. He
  writes, on p.66 in the original edition: ``Now those laws have been
  determined from the study of instances, in all of which it has been
  a necessary condition, that the classes or things added together in
  thought should be mutually exclusive. The expression $x+y$ seems
  indeed uninterpretable, unless it be assumed that the things
  represented by $x$ and the things represented by $y$ are entirely
  separate; that they embrace no individuals in common.''}, occurring
in many categories (see Proposition~\ref{CoprodEAProp} below).


We shall see many examples of effect algebras later on. Here we just
mention that each orthomodular lattice is an effect algebra, in which
elements $x,y$ are orthogonal iff $x \leq y^{\perp}$ iff $y^{\perp}
\leq x$, and also that each Boolean algebra is an effect algebra with
$x \orthogonal y$ iff $x \wedge y = 0$.

\auxproof{
A separate class of examples has a join as sum $\ovee$.
Let $(L, \vee, 0, (-)^{\perp})$ be an ortholattice: $\vee, 0$ are
finite joins and complementation $(-)^{\perp}$ satisfies $x\leq y
\Rightarrow y^{\perp} \leq x^{\perp}$, $x^{\perp\perp} = x$ and
$x\disjun x^{\perp} = 1 = 0^{\perp}$. This $L$ is called an
orthomodular lattice if $x \leq y$ implies $y = x \disjun (x^{\perp}
\conjun y)$. Such an orthomodular lattice forms an effect algebra in
which $x\ovee y$ is defined if and only if $x\orthogonal y$
(\textit{i.e.}~$x\leq y^{\perp}$, or equivalently, $y \leq
x^{\perp}$); and in that case $x\ovee y = x\disjun y$. This
restriction of $\vee$ is needed for the validity of requirements~(1)
and~(2) in Definition~\ref{EffAlgDef}: 
\begin{itemize}
\item suppose $x\ovee y = 1$, where $x\orthogonal y$,
  \textit{i.e.}~$x\leq y^{\perp}$. Then, by the orthomodularity
  property, we get uniqueness of orthocomplements:
$$\begin{array}{rcccccccccl}
y^{\perp}
& = &
x \disjun (x^{\perp} \conjun y^{\perp})
& = &
x\disjun (x\disjun y)^{\perp}
& = &
x \disjun 1^{\perp}
& = &
x\disjun 0
& = &
x.
\end{array}$$

\item $x\orthogonal 1$ means $x \leq 1^{\perp} = 0$, so that $x=0$.
\end{itemize}
}

We summarise some basic properties, without proof.

\begin{lem}
\label{EffAlgBasicsLem}
In an effect algebra one has:
\begin{enumerate}
\item \label{EffAlgBasicsLemBotbot} Orthocomplementing is an
  involution: $x^{\perp\perp} = x$;

\item Top and bottom are each others orthocomplements: $1^{\perp} = 0$ and 
$0^{\perp} = 1$;

\item \label{EffAlgBasicsLemCanc} Cancellation holds: $x\ovee y =
  x\ovee z$ implies $y=z$;

\item \label{EffAlgBasicsLemPos} Positivity (or zerosumfreeness)
  holds: $x\ovee y = 0$ implies $x=y=0$;

\item \label{EffAlgBasicsLemPoSet} $\leq$ is a partial order with $1$ as
  top and $0$ as bottom element;

\item \label{EffAlgBasicsLemBotRev} $x\leq y$ implies $y^{\perp} \leq
  x^{\perp}$;

\item \label{EffAlgBasicsLemOrtho} $x\ovee y$ is defined iff
  $x\orthogonal y$ iff $y\leq x^{\perp}$ iff $x\leq y^{\perp}$;

\item \label{EffAlgBasicsLemSumMon} $x\leq y$ and $y \orthogonal z$
  implies $x \orthogonal z$ and $x\ovee z \leq y\ovee z$;

\item \label{EffAlgBasicsLemOminus} $x\ovee y = z$ iff $y = z\ominus
  x$. \QED

\end{enumerate}
\end{lem}

\auxproof{
\begin{myproof} For point~\eqref{EffAlgBasicsLemBotbot} we use $x\ovee
    x^{\perp} = 1 = x^{\perp\perp} \ovee x^{\perp}$, so that $x =
    x^{\perp\perp}$ since the orthocomplement of $x^{\perp}$ is
    unique.  Uniqueness also yields $1^{\perp} = 0$ from $1\ovee
    1^{\perp} = 1 = 1 \ovee 0$. Hence $0^{\perp} = 1^{\perp\perp} =
    1$.

For cancellation in~\eqref{EffAlgBasicsLemCanc}, assume $x\ovee y =
x\ovee z$. The orthocomplement $u = (x\ovee y)^{\perp} = (x\ovee
z)^{\perp}$ satisfies $x\ovee y \ovee u = 1 = x\ovee z \ovee u$, so
that $y = (x\ovee u)^{\perp} = z$ by uniqueness of
orthocomplements. If $x\ovee y = 0$ in~\eqref{EffAlgBasicsLemPos},
then $y^{\perp} = 0 \ovee y^{\perp} = x \ovee y \ovee y^{\perp} =
x\ovee 1$. Hence $x\orthogonal 1$ and thus $x=0$. Similarly $y=0$.

For~\eqref{EffAlgBasicsLemPoSet} we obtain reflexivity $x\leq x$ from
$x\ovee 0 = x$. Transitivity is easy: if $x\leq y$ and $y\leq z$, say
via $y = x\ovee u$ and $z = y\ovee v$, then $z = y\ovee v = x \ovee
u\ovee v$, so that $x\leq z$. Finally, if $x\leq y$ and $y\leq x$, say
via $y = x\ovee u$ and $x = y\ovee v$, then $x\ovee 0 = x = x\ovee u
\ovee v$.  By cancelling $x$ we get $0 = u\ovee v$, so that $u = v =
0$ by positivity, and thus $x=y$. We get $0\leq x$ since $0\ovee x =
x$, and $x\leq 1$ since $x\ovee x^{\perp} = 1$.

For~\eqref{EffAlgBasicsLemBotRev}, let $x \leq y$, say via $x \ovee u
= y$. From $1 = y \ovee y^{\perp} = x\ovee u \ovee y^{\perp}$, we get
$x^{\perp} = u\ovee y^{\perp}$, showing that $y^{\perp} \leq
x^{\perp}$.

For~\eqref{EffAlgBasicsLemOrtho} first assume $x \orthogonal y$. Then
$x\ovee y \ovee (x\ovee y)^{\perp} = 1$, so that $x^{\perp} = y \ovee
(x\ovee y)^{\perp}$, and thus $y \leq x^{\perp}$. Conversely, if
$x^{\perp} = y\ovee u$, then $x \ovee y \ovee u = 1$, so that $x\ovee
y$ exists.

For~\eqref{EffAlgBasicsLemSumMon} assume $x\leq y$, say via $y =
x\ovee u$. If $y\orthogonal z$, then $y \leq z^{\bot}$, so also $x
\leq z^{\bot}$ and so $x \orthogonal z$. Then $x\ovee z\ovee u =
x\ovee u \ovee z = y\ovee z$, so that $x\ovee z\leq y\ovee z$.


For~\eqref{EffAlgBasicsLemOminus}, first assume $x\ovee y = z$. Then
$(z^{\bot} \ovee x) \ovee y = z^{\bot} \ovee z = 1$. Hence $y =
(z^{\bot} \ovee x)^{\bot} = z \ominus x$. In the other direction,
$z^{\bot} \ovee x \ovee (z \ominus x) = z^{\bot} \ovee x \ovee
(z^{\bot} \ovee x)^{\bot} = 1$, so that $x \ovee (z \ominus x) =
z$. \QED
\end{myproof}

The minus (or difference) operation $\ominus$ introduced in
Definition~\ref{EffAlgDef} satisfies many expected properties. We only
list them here and leave the proof to the interested reader.

\begin{lem}
\label{MinusLem}
In an effect algebra:
\begin{enumerate}
\item $x\ominus y$ is defined if and only if $y\leq x$;

\item $x\ominus 0 = x$ and $1\ominus x = x^{\perp}$;

\item $x \ominus y = (x^{\perp} \ovee y)^{\perp} = x \owedge y^{\perp}$;

\item $x\ovee y \leq z$ iff $x \leq z\ominus y$;

\item $z\ominus x = z\ominus y$ implies $x=y$;

\item $(x\ovee y) \ominus y = x$;

\item $x\leq y$ implies $y\ominus x \leq y$ and $y \ominus (y\ominus x) = x$;

\item $x\leq y\leq z$ implies $z\ominus y \leq z\ominus x$ and
$(z\ominus x) \ominus (z\ominus y) = y\ominus x$;

\item $x\ominus y \leq z$ iff $x\leq z\ovee y$. \QED
\end{enumerate}
\end{lem}

From points~(iv) and~(ix) we may conclude that $-\ominus y$,
considered as functor $\upset y\rightarrow \downset y^{\perp}$ between
posets, has $-\ovee y \colon \downset y^{\perp} \rightarrow \upset y$
both as left and right adjoint: $-\ovee y \dashv -\ominus y \dashv
-\ovee y$. Hence $-\ominus y$ and $-\ovee y$ preserve all meets and
joins that exist in $\upset y$ and $\downset y^{\perp}$, respectively.

Note that $x\ominus y$ is defined for $x\geq y$. The result
satisfies $x\ominus y = x\owedge y^{\perp} \leq y^{\perp}$, so lands
in $\downset y^{\perp}$. In the other direction, $x\ovee y$ is
defined for $x\orthogonal y$, \textit{i.e.}~$x\leq y^{\perp}$.
Obviously $x\ovee y \geq y$.

Notice that in the context of orthomodular lattices, for $y\leq x$ one
has $x \ominus y = x \conjun y^{\perp}$. Since $x\ominus y \leq
y^{\perp}$ one has $x\ominus y \orthogonal y$, so that $(x\ominus y)
\ovee y = (x \conjun y^{\perp}) \disjun y = x$, the latter by
orthomodularity.

The adjunction:
$$x\ominus y^{\perp} = x\conjun y \leq z
\Longrightarrow
x \leq y^{\perp} \ovee z$$

\noindent does not give a right adjoint to conjunction: the
restriction in this case is $y^{\perp} \leq$ which means that
this boils down to the andthen-sasaki adjunction since
$x\andthen y = y \conjun (y^{\perp} \disjun x) = y \conjun x$.

The proof of Lemma~\ref{MinusLem} goes as follows.

For~(i) we first observe that if $x\ominus y$ is defined, then $x =
y\ovee (x\ominus z)$, showing that $y\leq x$. Conversely, if $y\leq
x$, say via $x = y\ovee z$, then $x\ominus y = z$. In point~(ii) $x
\ominus 0 = x$ follows directly from $0\ovee x = x$, and $1 \ominus x =
x^{\perp}$ from $x\ovee (1\ominus x) = 1$.

In~(iii) we use $x^{\perp} \ovee y \ovee (x^{\perp} \ovee y)^{\perp} =
1$ to deduce that $x = x^{\perp\perp} = y \ovee (x^{\perp} \ovee
y)^{\perp}$, by uniqueness of orthocomplements. This means $(x^{\perp}
\ovee y)^{\perp} = x \ominus y$.

For~(iv) first assume $x\ovee y \leq z$, say via $x\ovee y \ovee u =
z$.  This means $x\ovee u = z\ominus y$, and thus $x\leq z\ominus
y$. Conversely, if $x\leq z\ominus y$, say via $x\ovee v = z\ominus y$,
then $x\ovee u \ovee y = z$, so that $x\ovee y \leq z$.

In~(v) assume $z\ominus x = z\ominus y$. Then $(z\ominus x) \ovee x = z
= (z\ominus y) \ovee y = (z\ominus x) \ovee y$, so that $x=y$ by
cancelling $z$ as in Lemma~\ref{EffAlgBasicsLem}~(3).

Point~(vi) holds by definition, since $(x\ovee y)\ominus y = z$ follows
because $y\ovee z = x\ovee y$ holds for $z=x$. 

For~(vii) assume $x\leq y$, say via $y = x\ovee u$. Then $u = y\ominus x
\leq y$. Moreover, $y \ominus (y\ominus x) = x$ holds because
$(y\ominus x) \ovee x = y$. 

For~(viii) we assume $x\leq y\leq z$, say with $x\ovee u = y$. Then
$z\ominus y \leq z\ominus x$, since $(z\ominus y) \ovee w = z\ominus
x$, namely for $w = u$: $x\ovee ((z\ominus y)\ovee u) = y \ovee
(z\ominus y) = z$. Further, $(z\ominus x) \ominus (z\ominus y) =
y\ominus x$ since $w = y\ominus x$ satisfies $(z\ominus y) \ovee w =
(z\ominus x)$, because $x \ovee ((z\ominus y) \ovee w) = x \ovee
(y\ominus x) \ovee (z\ominus y) = y \ovee (z\ominus y) = z$.

Finally, in~(ix) we assume $y\leq x$. In one direction, if $x\ominus y
\leq z$, then $x = (x\ominus y) \ovee y \leq z\ovee y$. Conversely, we
have $y\leq x \leq z\ovee y$, so that the previous point yields
$(z\ovee y)\ominus x \leq (z\ovee y)\ominus y = z$ and $x\ominus y =
z\ominus ((z\ovee y)\ominus x) \leq z$. \QED
}

Homomorphisms of effect algebras preserve the sums $\ovee$ that exist.
This is like for a probability measure $\mu \colon \Sigma \rightarrow
[0,1]$ satisfying $\mu(U\cup V) = \mu(U)+\mu(V)$ but only if the
measurable subsets $U,V\in\Sigma$ are disjoint, \textit{i.e.}~if the
sum $U\ovee V$ is defined in the Boolean algebra $\Sigma$ of
measurable subsets. Such a $\mu$ is thus a map of effect algebras. We
mention two properties.

\begin{lem}
\label{EffAlgHomLem}
Let $f\colon E\rightarrow D$ be a homomorphism of effect algebras. 
\begin{enumerate}
\item Then:
$$\begin{array}{rclcrclcrcl}
f(x^{\perp})
& = &
f(x)^{\perp}
& \qquad &
f(0)
& = &
0
& \qquad &
x\leq x' 
& \Longrightarrow &
f(x) \leq f(x').
\end{array}$$

\item If $E,D$ happen to be Boolean algebras, then $f$ is also a map
  of Boolean algebras: the inclusion functor $\BA \hookrightarrow \EA$
  is full and faithful --- where $\BA$ is the category of Boolean
  algebras.
\end{enumerate}
\end{lem}

\begin{myproof}
For the first point notice that $1 = f(1) = f(x \ovee x^{\perp}) =
f(x) \ovee f(x^{\perp})$, so $f(x^{\perp}) = f(x)^{\perp}$ by
uniqueness of orthocomplements. In particular, $f(0) = f(1^{\perp}) =
f(1)^{\perp} = 1^{\perp} = 0$. Monotonicity is trivial. For the second
point one uses that a join $x\vee y$ in a Boolean algebra can
equivalently be expressed as a disjoint join $(x\wedge \neg y) \vee
(x\wedge y) \vee (y\wedge \neg x)$. Therefor it is preserved by a map
of effect algebras. \QED

\auxproof{
The first observation is that in a Boolean algebra $x\orthogonal y$
iff $x\conjun y = 0$, since:
\begin{itemize}
\item if $x\orthogonal y$, \textit{i.e.}~$x\leq \neg y$, one gets
$x\conjun y \leq \neg y \conjun y = 0$.

\item if $x\conjun y = 0$, then by distributivity: $x = x\conjun 1 = x
  \conjun \neg(x\conjun y) = x \conjun (\neg x \disjun \neg y) =
  (x\conjun \neg x) \disjun (x\conjun \neg y) = 0 \disjun (x\conjun
  \neg y) = x\conjun \neg y$. Hence $x \leq \neg y$.
\end{itemize}

\noindent Next notice that the join of two elements $x,y$ can be written
jointly as:
\begin{equation}
\begin{array}{rcl}
x\disjun y
& = &
(x\conjun \neg y) \disjun (x\conjun y) \disjun (y\conjun \neg x).
\end{array}
\end{equation}

Disjointness of the three parts is obvious, so we only check that
the join is $x\disjun y$.
$$\begin{array}{rcl}
\lefteqn{(x\conjun \neg y) \disjun (x\conjun y) \disjun (y\conjun \neg x)} \\
& = &
[x \disjun (x\conjun y) \disjun (y\conjun \neg x)] \conjun
   [\neg y \disjun (x\conjun y) \disjun (y\conjun \neg x)] \\
& = &
[x \disjun x \disjun (y\conjun \neg x)] \conjun
[x \disjun y \disjun (y\conjun \neg x)] \conjun \\
& & \qquad
   [\neg y \disjun x \disjun (y\conjun \neg x)] \conjun
   [\neg y \disjun y \disjun (y\conjun \neg x)] \\
& = &
[x\disjun y] \conjun
[x \disjun \neg x] \conjun
[x \disjun y \disjun y] \conjun 
[x \disjun y \disjun \neg x] \conjun \\
& & \qquad
   [\neg y \disjun x \disjun y] \conjun
   [\neg y \disjun x \disjun \neg x] \conjun
   [\neg y \disjun y \disjun y] \conjun
   [\neg y \disjun y \disjun \neg x] \\
& = &
[x\disjun y] \conjun
0 \conjun
[x \disjun y] \conjun 
0 \conjun \\
& & \qquad
   0 \conjun
   0 \conjun
   0 \conjun
   0 \\
& = &
x\disjun y.
\end{array}$$

\noindent If $f\colon B\rightarrow C$ is a morphism of effect
algebras, it preserves disjoint joins, and by monotonicity it
satisfies $f(x) \disjun f(y) \leq f(x\disjun y)$ and $f(x\conjun y)
\leq f(x) \conjun f(y)$. Therefor, using the disjoint join expression:
$$\begin{array}{rcl}
f(x\disjun y)
& = &
f\big((x\conjun \neg y) \disjun (x\conjun y) \disjun (y\conjun \neg x)\big) \\
& = &
f(x\conjun \neg y) \disjun f(x\conjun y) \disjun f(y\conjun \neg x) \\
& \leq &
\big(f(x)\conjun \neg f(y)\big) \disjun \big(f(x)\conjun f(y)\big) 
   \disjun \big(f(y)\conjun \neg f(x)\big) \\
& = &
f(x) \disjun f(y).
\end{array}$$

\noindent Hence $f$ preserves joins --- and also meets, since
$x\conjun y = \neg (\neg x \disjun \neg y)$ --- and is thus morphism
of Boolean algebras. 
}
\end{myproof}

A \emph{test}, or more precisely, an \emph{$n$-test} in an effect
algebra $E$ is given by $n$ elements $e_{1}, \ldots, e_{n} \in E$ with
$e_{1} \ovee \cdots \ovee e_{n} = 1$. Such tests will be used as the
basis for measurement instruments in Section~\ref{MeasurementSec}. An
$n$-test in $E$ can be identified with a map of effect algebras
$\Pow(n) \rightarrow E$. If $E$ is the set of effects $\Ef(\H)$ in a
Hilbert space $\H$ --- consisting of the positive maps $\H \rightarrow
\H$ below the identity --- then an $n$-test in $\Ef(\H)$ is an
observable on the discrete measure space $n$, in the sense
of~\cite{HeinosaariZ12}\footnote{An observable on an arbitrary measure
  space $(X,\Sigma)$ can be described as a map of $\sigma$-effect
  algebras $\Sigma \rightarrow \Ef(\H)$, where the $\sigma$ refers to
  the preservation of countable joins of pairwise orthogonal elements,
  see~\cite[Defn.~3.5]{HeinosaariZ12}.}.

Two extensions of the notion of effect algebra will be used, namely
extension of an effect algebra $E$ with:
\begin{itemize}
\item a multiplication $x\cdot y$ of its elements $x,y\in E$; we call
  such a structure an \emph{effect monoid};

\item a scalar multiplication $r\scalar x$, where $x\in E$ and $r$ is
  a scalar, belonging for instance to $[0,1]$, or more generally to a
  given effect monoid; such a structure will be called an \emph{effect
    module}.
\end{itemize}

\noindent In an \emph{effect monoid} we require that the
multiplication operation $\cdot$ is associative, preserves $0,\ovee$
in each argument separately, and satisfies $1\cdot x = x = x\cdot
1$. An effect monoid is called commutative if its multiplication is
commutative. The unit interval $[0,1]$ is an example of a commutative
effect monoid, with the usual multiplication of
probabilities. Similarly, the Booleans $\{0,1\}$ form a commutative
effect monoid with multiplication (conjunction).

Given an effect monoid $M$, an effect algebra $E$ is an \emph{effect
  module} over $M$ if there is a scalar multiplication $\scalar \colon
M \times E \rightarrow E$ which preserves $0,\ovee$ in each argument
separately and additionally satisfies: $1\scalar x = x$ and $r\scalar
(s\scalar x) = (r\cdot s) \scalar x$. A map of effect modules is a map
of effect algebras that commutes with scalar multiplication. We write
$\EMod_{M}$ for the resulting category. There is an obvious forgetful
functor $\EMod_{M} \rightarrow \EA$, comparable to the forgetful
functor $\Vect_{K} \rightarrow \Ab$ of vector spaces over a field $K$
to Abelian groups.

Here is an example of an effect module --- many more will be given
below. The set $[0,1]^{X}$ contains functions from a set $X$ to
$[0,1]$ that can be understood as fuzzy predicates. This set
$[0,1]^{X}$ is an effect algebra, via pointwise operations, with $p
\orthogonal q$ if $p(x) + q(x) \leq 1$ for all $x\in X$. It is also an
effect module over $[0,1]$ with scalar multiplication $r\scalar p \in
[0,1]^{X}$, for $r\in [0,1]$, given by $(r\scalar p)(x) = r\cdot
p(x)$. Thus, $[0,1]^{X} \in \EMod_{[0,1]}$, for each set $X$.

In a trivial manner each effect algebra $E$ is an effect module over
$\{0,1\}$, with obvious scalar multiplication $\{0,1\} \times E
\rightarrow E$. This is like: each Abelian group is a trivial module
over the (initial) ring $\Z$ of integers.

\begin{rem}
\label{EffectMonoidalRem}
The notions of effect monoid and module can be described more
abstractly: an effect monoid is a monoid in the category of effect
algebras, just like a semiring is a monoid in the category of
commutative monoids. In~\cite{JacobsM12a} it is shown that the
category \EA of effect algebras is symmetric monoidal with the
two-element (initial) effect algebra $2 = \{0,1\}$ as tensor
unit. Then one can consider, in a standard way, the categories
$\Mon(\EA)$ and $\CMon(\EA)$ of (commutative) monoids $2 \rightarrow M
\leftarrow M\otimes M$ in the category $\EA$ of effect algebras.

For an effect monoid $M\in\Mon(\EA)$ the category $\EMod_{M}$ of
effect modules over $M$ is the the category $\Act_{M}(\EA)$ of
$M$-actions $M\otimes E \rightarrow E$, see
\cite[VII,\S4]{MacLane71}. This resembles the category $\Mod_{S}$ of
modules over a semiring $S$ which may be described as the category
$\Act_{S}(\CMon)$ of commutative monoids with $S$-scalar
multiplication. Effect modules over $[0,1]$ have appeared under the
name `convex effect algebras', see~\cite{PulmannovaG98}. Via suitable
totalisations of the partial operation $\ovee$ (see~\cite{JacobsM12a})
it is shown in~\cite[Prop.~9]{JacobsM12b} that the category
$\EMod_{[0,1]}$ is equivalent to the category of order unit spaces:
(real) partially ordered vector spaces with an order unit 1
satisfying: for each vector $v$ there is an $n\in\NNO$ with $-n\cdot 1
\leq v \leq n\cdot 1$.
\end{rem}

The Appendix~\ref{DiscProbSubsec} describes the distribution monad
$\Dst\colon\Sets\rightarrow\Sets$ where $\Dst(X)$ consists of formal
sums $\sum_{i}r_{i}\ket{x_{i}}$ with $r_{i}\in [0,1]$ satisfying
$\sum_{i}r_{i} = 1$. In~\cite{Jacobs11c} it is shown that such a
distribution monad can be defined more generally wrt.~a (commutative)
effect monoid $M$ in place of $[0,1]$. This works as follows.
$$\begin{array}{rcl}
\Dst_{M}(X)
& = &
\set{m_{1}\ket{x_{1}} + \cdots + m_{k}\ket{x_{k}}}{x_{i}\in X, m_{i}\in M 
   \mbox{ with } \bigovee_{i}m_{i}=1\in M}.
\end{array}$$

\noindent Implicitly in this formulation we assume that the finite sum
$\ovee_{i}\,m_{i}$ exists. It can be shown that $\Dst_M$ is a monad on
\Sets, just like $\Dst = \Dst_{[0,1]}$ is a monad. For the trivial
effect monoid $\{0,1\}$ the associated monad $\Dst_{\{0,1\}}$ is the
identity functor.

In a next step we can form the category $\Conv_{M} = \EM(\Dst_{M})$ of
Eilenberg-Moore algebras of this monad $\Dst_M$. Its objects are
`convex sets over $M$', that is, sets $X$ in which for each formal
convex sum $\sum_{i}m_{i}\ket{x_{i}}\in\Dst_{M}(X)$ an element
$\bigovee_{i}\,m_{i}x_{i}\in X$ exists. Such convex sets can also be
described in terms of `weighted sums' $x +_{r} y$, interpreted as $rx
+ (1-r)y$, see
\textit{e.g.}~\cite{Stone49,Swirszcz74,Jacobs10e}. Morphisms in
$\EM(\Dst_{M})$ are \emph{affine} maps; they preserve such convex sums
$\ovee_{i}\,m_{i}x_{i}$.

In~\cite[Prop.~5]{JacobsM12b} (based on~\cite{Jacobs10e}) a (dual)
adjunction between convex sets over $[0,1]$ and effect modules over
$[0,1]$ is described. This adjunction exists in fact for an arbitrary
effect monoid $M$ --- instead of $[0,1]$ --- by using $M$ as dualising
object. It formalises the duality between effects and states in
quantum foundations, and also the dual relations between Heisenberg's
and Schr{\"o}dinger's view on quantum computation.

\begin{prop}
\label{ConvEModAdjProp}
Let $M$ be an effect monoid. By ``homming into $M$'' one obtains an
adjunction:
$$\xymatrix{
\op{\big(\EMod_{M}\big)}\ar@/^1.5ex/[rr]^-{\Conv(-,M)} 
   & \top & 
   \Conv_{M} = \rlap{$\;\EM(\Dst_{M})$}
   \ar@/^1.5ex/[ll]^-{\EMod(-,M)}
}$$
\end{prop}

\begin{myproof}
Given a convex set $X\in\Conv_{M}$, the homset $\Conv(X,M)$ of affine
maps is an effect module, with $f\orthogonal g$ iff $\allin{x}{X}{f(x)
  \orthogonal g(x)}$ in $M$. In that case one defines $(f\ovee g)(x) =
f(x)\ovee g(x)$.  It is easy to see that this is again an affine
function. Similarly, the pointwise scalar product $(m\scalar f)(x) =
m\cdot f(x)$ yields an affine function. This mapping $X \mapsto
\Conv(X,M)$ gives a contravariant functor since for $h\colon
X\rightarrow X'$ in $\Conv_{M}$ pre-composition with $h$ yields a map
$(-) \after h \colon \Conv(X', M) \rightarrow \Conv(X, M)$ of
effect modules.

\auxproof{
For a formal convex sum $\ovee_{j}s_{j}x_{j}$ one has:
$$\begin{array}{rcl}
(f \ovee g)(\sum_{j}s_{j}x_{j})
& = &
f(\ovee_{j}s_{j}x_{j}) \ovee g(\ovee_{j}s_{j}x_{j}) \\
& = &
\big(\ovee_{j}s_{j}f(x_{j})\big) \ovee \big(\ovee_{j}s_{j}g(x_{j})\big) \\
& = &
\ovee_{j}s_{j}(f(x_{j}) \ovee g(x_{j})) \\
& = &
\ovee_{j}s_{j}(f\ovee g)(x_{j}) \\
(r\scalar f)(\ovee_{j}s_{j}x_{j})
& = &
r\cdot f(\ovee_{j}s_{j}x_{j}) \\
& = &
r\cdot (\ovee_{j}s_{j}f(x_{j})) \\
& = &
\ovee_{j}s_{j} r\cdot f(x_{j}) \\
& = &
\ovee_{j}s_{j}(r\scalar f)(x_{j}) \\
r \scalar (f \ovee g)
& = &
\lam{x}{r\cdot (f \ovee g)(x)} \\
& = &
\lam{x}{r\cdot (f(x) \ovee g(x))} \\
& = &
\lam{x}{r\cdot f(x) \ovee r \cdot g(x))} \\
& = &
\lam{x}{(r\scalar f)(x) \ovee (r\scalar g)(x)} \\
& = &
(r \scalar f) \ovee (r \scalar g) \\
(f \ovee g) \after h
& = &
\lam{x}{f(h(x)) \ovee g(h(x))} \\
& = &
\lam{x}{(f \after h)(x) \ovee (g \after h)(x)} \\
& = &
(f \after h) \ovee (f \after h) \\
(r \scalar f) \after h
& = &
\lam{x}{r \cdot f(h(x))} \\
& = &
\lam{x}{r \cdot (f \after h)(x)} \\
& = &
r \scalar (f \after h).
\end{array}$$
}

In the other direction, for an effect module $E\in\EMod_{M}$, the
homset $\EMod(E, M)$ of effect module maps yields a convex set: for a
formal convex sum $\sum_{j}m_{j}\ket{f_{j}}$, where $f_{j} \colon E
\rightarrow M$ in $\EMod_{M}$ and $m_{j}\in M$, we can define an
actual sum $f\colon E\rightarrow M$ by $f(y) = \bigovee_{j}\,
m_{j}\scalar f_{j}(y)$. This $f$ forms a map of effect modules. Again,
functoriality is obtained via pre-composition.

\auxproof{
$$\begin{array}{rcl}
f(y \ovee z)
& = &
\ovee_{j}\,m_{j}\cdot f_{j}(y\ovee z) \\
& = &
\ovee_{j}\,m_{j}\cdot (f_{j}(y) \ovee f_{j}(z)) \\
& = &
\ovee_{j}\,m_{j}\cdot f_{j}(y) \ovee m_{j}\cdot f_{j}(z)) \\
& = &
\ovee_{j}\,m_{j}\cdot f_{j}(y) \ovee \ovee_{j}\,m_{j}\cdot f_{j}(z)) \\
& = &
f(y) \ovee f(z) \\
f(1)
& = &
\ovee_{j}\,m_{j}\cdot f_{j}(1) \\
& = &
\ovee_{j}\,m_{j} \cdot 1 \\
& = &
\ovee_{j}\,m_{j} \\
& = &
1.
\end{array}$$

For a map $h\colon Y \rightarrow Y'$ of effect modules precomposition
with $h$ gives an affine map $(-) \after h \colon \EMod(Y', M)
\rightarrow \EMod(Y, M)$ since for a convex sum
$\sum_{j}m_{j}f_{j}$ of effect algebra maps $f_{j} \colon
Y'\rightarrow [0,1]$ we get:
$$\begin{array}{rcl}
(\ovee_{j}\,m_{j}f_{j}) \after h
& = &
\lam{y}{(\ovee_{j}\,m_{j}f_{j})(h(y))} \\
& = &
\lam{y}{\ovee_{j}\,m_{j}\cdot f_{j}(h(y))} \\
& = &
\lam{y}{\ovee_{j}\,m_{j}\cdot (f_{j} \after h)(y))} \\
& = &
\ovee_{j}\,m_{j}\cdot (f_{j} \after h).
\end{array}$$
}

The dual adjunction between $\EMod_M$ and $\Conv_M$ involves a
bijective correspondence that is obtained by swapping arguments. \QED

\auxproof{
For $X\in\Conv_{M}$ and $Y\in\EMod_{M}$, we have:
$$\begin{prooftree}
{\xymatrix{ 
X\ar[r]^-{f} & \EMod(Y,M) 
   \rlap{\hspace*{2em} in $\Conv$}}}
\Justifies
{\xymatrix{ Y\ar[r]_-{g} & \Conv(X, M)
   \rlap{\hspace*{2.3em} in $\EMod$}}}
\end{prooftree}$$

\noindent What needs to be checked is that for a map $f$ of convex
sets as indicated, the swapped version $\widehat{f} =
\lamin{y}{Y}{\lamin{x}{X}{f(x)(y)}} \colon Y \rightarrow
\Conv(X,M)$ is a map of effect modules---and similarly for $g$.

We show that $\widehat{f}(y) \colon X\rightarrow M$ is
a map of convex sets:
$$\begin{array}{rcl}
\widehat{f}(y)(\ovee_{j}\,m_{j}x_{j})
& = &
f(\ovee_{j}\,m_{j}x_{j})(y) \\
& = &
\big(\ovee_{j}\,m_{j}f(x_{j})\big)(y) \\
& = &
\ovee_{j}\,m_{j}f(x_{j})(y) \\
& = &
\ovee_{j}\,m_{j}\widehat{f}(y)(x_{j}).
\end{array}$$

\noindent Next we have to check that $\widehat{f}$ is a map
of effect modules.
$$\begin{array}{rcl}
\widehat{f}(y_{1} \ovee y_{2})
& = &
\lam{x}{f(x)(y_{1} \ovee y_{2})} \\
& = &
\lam{x}{f(x)(y_{1}) \ovee f(x)(y_{2})} \\
& = &
\lam{x}{\widehat{f}(y_{1})(x) \ovee \widehat{f}(y_{2})(x)} \\
& = &
\widehat{f}(y_{1}) \ovee \widehat{f}(y_{2}) \\
\widehat{f}(1)
& = &
\lam{x}{f(x)(1)} \\
& = &
\lam{x}{1} \\
& = &
1
\end{array}$$

Next, starting from $g$ we take $\widehat{g} =
\lamin{x}{X}{\lamin{y}{Y}{g(y)(x)}}$. First, $\widehat{g}(x)$ is a
map of effect modules.
$$\begin{array}{rcl}
\widehat{g}(x)(y_{1} \ovee y_{2})
& = &
g(y_{1} \ovee y_{2})(x) \\
& = &
(g(y_{1}) \ovee g(y_{2}))(x) \\
& = &
g(y_{1})(x) \ovee g(y_{2})(x) \\
& = &
\widehat{g}(x)(y_{1}) \ovee \widehat{g}(x)(y_{2}) \\
\widehat{g}(x)(1)
& = &
g(1)(x) \\
& = &
(\lam{z}{1})(x) \\
& = &
1.
\end{array}$$

\noindent Also, $\widehat{g}$ is affine:
$$\begin{array}{rcl}
\widehat{g}(\ovee_{j}\,m_{j}x_{j})
& = &
\lam{y}{g(y)(\ovee_{j}\,m_{j}x_{j})} \\
& = &
\lam{y}{\ovee_{j}\,m_{j}g(y)(x_{j})} \\
& = &
\lam{y}{\ovee_{j}\,m_{j}\widehat{g}(x_{j})(y)} \\
& = &
\ovee_{j}\,m_{j}\widehat{g}(x_{j}).
\end{array}$$
}
\end{myproof}

\section{Predicates and tests}\label{PredSec}

This section introduces predicates as maps of the form $X \rightarrow
1+1$ in a category, and more generally, tests as maps $X \rightarrow
n\cdot 1 = 1+\cdots+ 1$. At this stage our only aim is to illustrate
what such predicates/tests are in the categories of interest. The
algebraic structure of predicates will be investigated in the next
section.

In topos theory predicates are described as maps of the form $X
\rightarrow \Omega$, where $\Omega$ is a special object that carries
Heyting algebraic structure. The crucial topos property is that such
maps $X \rightarrow \Omega$ correspond to subobjects of $X$. Here we
use $1+1$ instead of $\Omega$. This looks extremely simple, but we
shall see that predicates $X \rightarrow 1+1$ include many interesting
examples, such as fuzzy predicates $X \rightarrow [0,1]$ taking values
in the unit interval $[0,1]$ of probabilities. In this section we
elaborate the overview of predicates given in Figure~\ref{ValidityFig}
in the introduction.

Let's be more explicit about notation: we write $0$ for the initial
object in a category, with unique map $!\colon 0 \rightarrow X$ to
each object $X$. The coproduct object $X+Y$ of two objects $X,Y$ comes
with two coprojections $\kappa_{1}\colon X \rightarrow X+Y$ and
$\kappa_{2}\colon Y \rightarrow X+Y$ which are universal: for each
pair of maps $f\colon X \rightarrow Z$ and $g\colon Y\rightarrow Z$
there is a unique cotuple $[f,g]\colon X+Y\rightarrow Z$ with
$[f,g]\after\kappa_{1} = f$ and $[f,g]\after \kappa_{2} = g$. For two
maps $h\colon X \rightarrow A$ and $k\colon Y\rightarrow B$ one writes
$h+k = [\kappa_{1} \after h, \kappa_{2} \after k]\colon X+Y
\rightarrow A+B$. The codiagonal map $X+X \rightarrow X$ is written as
$\nabla = [\idmap,\idmap]$. These coproducts, coprojections and
cotuples generalise to $n$-ary form: $X_{1}+\cdots+X_{n}$. If all
these objects are the same, we have a so-called copower, written as
$n\cdot X = X+\cdots+X$.

For the products $X\times Y$ we use standard notation, with
projections $X \stackrel{\pi_1}{\leftarrow} X\times Y
\stackrel{\pi_2}{\rightarrow} Y$ and tuples $\tuple{f,g} \colon Z
\rightarrow X\times Y$, for maps $f\colon Z \rightarrow X, g\colon Z
\rightarrow Y$. The empty product is a final object $1$, with unique
map $!\colon X \rightarrow 1$ for each object $X$. We shall frequently
use that products in a category $\cat{B}$ form coproducts in the
opposite category $\op{\cat{B}}$.

\auxproof{
We recall the following basic result (see
\textit{e.g.}~\cite{Jacobs13b}).

\begin{lem}
\label{CoproductComonadLem}
Let $\cat{B}$ be a category with coproducts $+$. For each natural
number $n > 0$, the $n$-fold copower functor $n\cdot (-) \colon
\cat{B} \rightarrow \cat{B}$ is a comonad, where
$$\begin{array}{rcl}
n\cdot X
& = &
\underbrace{X + \cdots + X}_{\textrm{$n$ times}}
\end{array}$$

\noindent The counit $\varepsilon \colon n\cdot X \rightarrow X$ and
comultiplication $\delta \colon n\cdot X \rightarrow n\cdot (n\cdot X)$ are
given by:
$$\begin{array}{rccclcrcccl}
\varepsilon
& = &
\nabla
& = &
[\idmap, \ldots, \idmap]
& \qquad &
\delta
& = &
\kappa_{1} + \cdots + \kappa_{n}
& = &
[\kappa_{i} \after \kappa_{i}]_{i\leq n}.
\end{array}\eqno{\qEd}$$
\end{lem}
}

\begin{defi}
\label{PredDef}
Let $\cat{B}$ be a category with coproducts and with a final object
$1\in\cat{B}$. Let $n\in\NNO$ be non-zero.
\begin{enumerate}
\item An \emph{$n$-test} on an object $X\in\cat{B}$ is map $p\colon X
  \rightarrow n\cdot 1 = 1+\cdots+1$.  


\item A 2-test on $X$ is also called a \emph{predicate} on $X$.  We
  write $\Pred(X) = \Hom(X,1+1)$ for the homset of predicates on
  $X$. There are true and false predicates $1_{X}, 0_{X} \in
  \Pred(X)$, defined as:
$$\xymatrix@C-.5pc{
1_{X} = \big(X\ar[r]^-{!_X} & 1\ar[r]^-{\kappa_{1}} & 1+1\big)
& &
0_{X} = \big(X\ar[r]^-{!_X} & 1\ar[r]^-{\kappa_{2}} & 1+1\big)
}$$

\noindent Each predicate $p\in\Pred(X)$ has an orthocomplement
$p^{\perp} \in \Pred(X)$ defined by swapping the outcomes:
$$\xymatrix@C+.5pc{
p^{\perp} = \big(X\ar[r]^-{p} & 1+1\ar[r]^-{[\kappa_{2},\kappa_{1}]}_-{\cong} 
   & 1+1\big)
}$$

\noindent The predicates $\Pred(1)$ on the final object $1\in\cat{B}$
play a special role and will be called \emph{scalars} (sometimes
\emph{probabilities}).

\item Each map $f\colon Y \rightarrow X$ in $\cat{B}$ yields a
  ``substitution'' or ``reindexing'' function $\Pred(f) = (-) \after f
  \colon \Pred(X) \rightarrow \Pred(Y)$. In this way we get a
  functor $\Pred \colon \op{\cat{B}} \rightarrow \Sets$.

Several different notations are in use for this map $\Pred(f)$, namely
$f^{*}$ in categorical logic, $\wp(f)$ for `weakest precondition' in
programming logic'. We shall use all of these interchangeably, so
$\Pred(f) = f^{*} = \wp(f) = (-) \after f$.
\end{enumerate}
\end{defi}

\noindent We often drop the subscripts for the predicates $0 = 0_{X}$ and $1 =
1_{X}$ when they are clear from the context. Obviously, $1^{\perp} =
0$ and $0^{\perp} = 1$. Also, $p^{\perp\perp} = p$. Thus we have
`double negation' built into our logic, in line with the right leg in
Figure~\ref{BAFig}. The substitution function $f^{*}$ clearly
preserves truth, false and orthocomplements.

\begin{exas}
\label{PredEx}
We describe some of the motivating examples.
\begin{enumerate}
\item \label{PredExSets} In the category $\Sets$ we identify the
  $n$-fold copower $n\cdot = 1 + \cdots + 1$ with the $n$-element set,
  commonly also written as $n$. An $n$-test on a set $X$ is thus a
  function $p\colon X \rightarrow n$ giving a partition of $X$ in $n$
  subsets $p^{-1}(i)\subseteq X$, for $i\in n$. A predicate, or 2-test,
  is a predicate described as characteristic function $X \rightarrow 2
  = \{0,1\}$. The probabilities $\Pred(1)$ are the Booleans $\{0,1\}$.

More generally, in a topos, the maps $X \rightarrow 1+1$ correspond to
the Boolean predicates. The set of scalars $\Pred(1)$ need not be as
trivial as $\{0,1\}$. For instance, in $\Sets^{2}$ one has four maps
$(1,1) \rightarrow (1,1)+(1,1) = (1+1, 1+1)$, namely $(\kappa_{i},
\kappa_{j})$ for $i,j\in\{1,2\}$.

\item \label{PredExKlD} Our next example involves the distribution
  monad $\Dst$ on $\Sets$ from Appendix~\ref{DiscProbSubsec}, which is
  used to model \emph{discrete} probabilistic computations. The
  associated Kleisli category $\Kl(\Dst)$ captures these computations,
  via stochastic matrices, as maps $X \rightarrow Y$ in $\Kl(\Dst)$,
  corresponding to functions $X \rightarrow \Dst(Y)$. An $n$-test $p$
  on a set $X\in\Kl(\Dst)$ is a function $p\colon X \rightarrow
  \Dst(n)$. It assigns to each element $x\in X$ an $n$-tuple $p(x)
  \colon n \rightarrow [0,1]$ of probabilities $p(x)(i) \in [0,1]$
  which add up to 1, that is, $\sum_{i} p(x)(i) = 1$. A predicate on
  $X$ can be identified with a fuzzy predicate $X \rightarrow
  [0,1]$. The scalars $\Pred(1)$ are in this case the usual
  probabilities in the unit interval $[0,1] \subseteq \R$. Thus we use
  the bijective correspondences:
$$\begin{prooftree}
\begin{prooftree}
\xymatrix{\mbox{Kleisli map } X\ar[r] & 2}
\Justifies
\xymatrix{\mbox{function } X\ar[r] & \Dst(2)
   \rlap{$\; \cong [0,1]$}}
\end{prooftree}
\Justifies
\mbox{fuzzy predicate in }[0,1]^{X}
\end{prooftree}$$

\item \label{PredExKlG} For \emph{continuous} probabilistic
  computation one uses the Giry monad $\Giry$ on the category $\Meas$
  of measurable spaces, see Appendix~\ref{ContProbSubsec}.  In the
  category $\Meas$, and also in the Kleisli category $\Kl(\Giry)$ of
  the Giry monad $\Giry$ on $\Meas$, the object $n\cdot 1$ is the
  $n$-element set with the discrete $\sigma$-algebra
  $\Pow(n)$. Therefore, $\Giry(n) \cong \Dst(n)$. Hence an $n$-test on
  $X\in\Meas$ is a measurable function $X \rightarrow \Giry(n)$. In
  particular, a predicate on $X\in\Meas$ is simply a measurable
  function $X \rightarrow [0,1]$, or a $[0,1]$-valued random variable,
  as used for instance in~\cite{Jacobs13a}.  As before we have the
  probabilities $[0,1]$ as scalars: $\Pred(1) \cong \Giry(2) \cong
  \Dst(2) \cong [0,1]$.

\item \label{PredExCstar} Let $\CstarPU$ be the category of (complex
  unital) $C^*$-algebras, with positive unital maps as morphisms
  between them, see Appendix~\ref{CstarSubsec}. The algebra $\C$ of
  complex numbers is initial in $\CstarPU$, since the only positive
  unital map $\C \rightarrow A$ to an arbitrary $C^*$-algebra $A$ is
  the map $z \mapsto z\cdot 1$, where $1\in A$ is the unit.  The
  category $\CstarPU$ has binary products $A_{1}\oplus A_{2}$, given
  by the cartesian product of the underlying spaces, with
  coordinate-wise operations, and supremum norm $\|(a_{1},a_{2})\| =
  \|a_{1}\| \vee \|a_{2}\|$. The singleton space $\{0\}$ is the final
  object in $\CstarPU$.

We shall be working in the opposite category $\op{(\CstarPU)}$. It
thus has coproducts $\oplus$ and a final object $\C$. The $n$-fold
coproduct $n\cdot 1$ is then the $n$-fold cartesian product $\C^{n}$.
An $n$-test on $A\in\CstarPU$ is a map $A \rightarrow \C^{n}$ in
$\op{(\CstarPU)}$, that is, a positive unital map $p\colon \C^{n}
\rightarrow A$. It can be identified with an $n$-tuple of effects
$e_{i} = p(\ket{i}) \in [0,1]_{A}$ with $e_{1} \ovee \cdots \ovee
e_{n} = 1$. Indeed, $\ket{i}\in\C^{n}$ is positive, and so $e_{i} =
p(\ket{i}) \geq 0$ in $A$. And $\bigovee_{i} e_{i} = \sum_{i}
p(\ket{i}) = p(\sum_{i}\ket{i}) = p(1) = 1$. A predicate on $A$ can be
identified with a pair of predicates $e_{1}, e_{2} \in [0,1]_{A}$ with
$e_{1} + e_{2} = 1$. Hence $e_{2} = 1 - e_{1}$, so we can identify the
predicate with a single effect $e\in [0,1]_{A}$, where $e^{\perp} = 1
- e$. Thus we have the fundamental correspondences:
$$\begin{prooftree}
\begin{prooftree}
\xymatrix{A\ar[r] & 2} \rlap{\hspace*{3.3em}in $\op{(\CstarPU)}$}
\Justifies
\xymatrix{\C^{2}\ar[r] & A} \rlap{\hspace*{3em}in $\CstarPU$}
\end{prooftree}
\Justifies
\mbox{effect in $[0,1]_{A} \subseteq A$}
\end{prooftree}\qquad\qquad$$

\noindent (This correspondence also works for \emph{completely}
positive unital maps, since a positive map $f\colon A \rightarrow B$
is automatically completely positive if either $A$ or $B$ is
commutative.)

The truth predicate in $[0,1]_{A}$ is obtained by applying the map $!
\, \after \pi_{1} \colon \C^{2} \rightarrow \C \rightarrow A$ to
$(1,0)\in \C^{2}$. It yields the top/unit element $1\in
   [0,1]_{A}$. Similarly, false is obtained by applying $!\,\after
   \pi_{2}$ to the pair $(1,0)$; this yields the bottom/zero element
   $0\in [0,1]_{A}$. The scalars for $C^*$-algebras are the predicates
   $\Pred(\C)$ on the complex numbers, and correspond to the usual
   probabilities $[0,1]$.

For a map $f\colon B \rightarrow A$ in $\op{(\CstarPU)}$ the
associated substitution function $f^{*} \colon \Pred(A) \rightarrow
\Pred(B)$ is defined by function application: an effect $e\in
     [0,1]_{A}$ is sent to the effect $f(e)\in [0,1]_{B}$.

We briefly describe what happens in the case we do not use positive
unital (PU) but the more common *-homomorphisms between
$C^*$-algebras, \textit{i.e.}~the multiplicative-involutive-unital
(MIU) maps. It is easy to see that predicates $A \rightarrow 2$ in
$\op{(\CstarMIU)}$ correspond to \emph{projections} in $A$, that is, to
idempotent effects $a\in [0,1]_{A}$. The resulting scalars are the
Booleans $\{0,1\}$. Hence in the MIU-case the logic is no longer
quantitative (with $[0,1]$ as scalars). But the projections themselves
do not form a Boolean algebra, in general.

\item \label{PredExRng} Let $\Rng$ be the category of rings (with
  unit), with ring homomorphisms between them. The set of integers
  $\Z$, with its standard ring structure, is initial in the category
  $\Rng$, and thus final in the opposite category $\op{\Rng}$. We
  shall write it simply as 1. The coproduct $1+1\in\op\Rng$ is then
  the cartesian product $\Z\times\Z$. We claim that there are
  bijective correspondences, for $R\in\Rng$,
$$\begin{prooftree}
\begin{prooftree}
\xymatrix{R\ar[r] & 1+1 \rlap{\hspace*{4em}in $\op{\Rng}$}}
\Justifies
\xymatrix{\Z^{2}\ar[r]|{f} & R \rlap{\hspace*{4.6em}in $\Rng$}}
\end{prooftree}
\Justifies
\mbox{idempotents $e\in R$}
\end{prooftree}$$

\noindent The correspondence arises as follows. Given a ring
homomorphism $f\colon \Z^{2} \rightarrow R$ we get an element $e_{f} =
f(1,0)\in R$ which is idempotent: $e_{f}^{2} = f(1,0)^{2} =
f((1,0)^{2}) = f(1,0) = e_{f}$. In the other direction, given an
idempotent $e\in R$ we define $f_{e} \colon \Z^{2} \rightarrow R$ by
$f_{e}(n,m) = n\cdot e + m\cdot (1-e)$. In order to prove that $f_e$
is a ring homomorphism, in particular that it preserves
multiplication, one uses that $e$ is an idempotent.

\auxproof{
$$\begin{array}{rcl}
f_{e}(0,0)
& = &
0 \cdot e + 0\cdot (1-e) \\
& = &
0 \\
f_{e}((n,m) + (n',m'))
& = &
f_{e}(n+n', m+m') \\
& = &
(n+n')\cdot e + (m+m')\cdot (1-e) \\
& = &
n\cdot e + m\cdot (1-e) + n'\cdot e + m'\cdot (1-e) \\
& = &
f_{e}(n,m) + f_{e}(n',m') \\
f_{e}(1,1)
& = &
1\cdot e + 1\cdot (1-e) \\
& = &
e + (1-e) \\
& = &
1 \\
f_{e}((n,m)\cdot (n',m'))
& = &
f_{e}(n\cdot n', m\cdot m') \\
& = &
(n\cdot n') \cdot e + (m\cdot m') \cdot (1-e) \\
& \smash{\stackrel{(*)}{=}} &
(n\cdot e)\cdot (n'\cdot e) + (n\cdot e) \cdot (m'\cdot (1-e)) \\
& & \qquad + \;
   (m\cdot (1-e)) \cdot (n'\cdot e) + (m\cdot (1-e))\cdot (m'\cdot (1-e)) \\
& = &
(n\cdot e + m\cdot (1-e))\cdot (n'\cdot e + m'\cdot (1-e)) \\
& = &
f_{e}(n,m) \cdot f_{e}(n',m')
\end{array}$$

\noindent The marked equation holds since:
$$\begin{array}{rcl}
(n\cdot e)\cdot (n'\cdot e)
& = &
(e + \cdots + e)\cdot (n'\cdot e) \\
& = &
e\cdot (n'\cdot e) + \cdots + e\cdot (n'\cdot e) \\
& = &
e\cdot (e + \cdots + e) + \cdots + e\cdot (e + \cdots + e) \\
& = &
(e^{2} + \cdots + e^{2}) + \cdots + (e^{2} + \cdots + e^{2}) \\
& = &
(e + \cdots + e) + \cdots + (e + \cdots + e) \\
& = &
(n\cdot n') \cdot e \\
(n\cdot e) \cdot (m'\cdot (1-e))
& = &
(e + \cdots + e) \cdot (m'\cdot (1-e)) \\
& = &
e\cdot (m'\cdot (1-e)) + \cdots + e\cdot (m'\cdot (1-e)) \\
& = &
e\cdot ((1-e) + \cdots + (1-e)) + \cdots + e\cdot ((1-e) + \cdots + (1-e)) \\
& = &
((e-e^{2}) + \cdots + (e-e^{2})) + \cdots + ((e-e^{2}) + \cdots + (e-e^{2})) \\
& = &
0 \\
(m\cdot (1-e))\cdot (m'\cdot (1-e)) 
& = &
((1-e) + \cdots + (1-e))\cdot (m'\cdot (1-e)) \\
& = &
(1-e)\cdot (m'\cdot (1-e)) + \cdots + (1-e)\cdot (m'\cdot (1-e)) \\
& = &
(1-e)\cdot ((1-e) + \cdots + (1-e)) \\
& & \qquad + \cdots + (1-e)\cdot ((1-e) + \cdots + (1-e)) \\
& = &
((1-e)^{2} + \cdots + (1-e)^{2}) + \cdots + ((1-e)^{2} + \cdots + (1-e)^{2}) \\
& = &
((1-e) + \cdots + (1-e)) + \cdots + ((1-e) + \cdots + (1-e)) \\
& & \qquad\mbox{since } (1-e)^{2} = 1 - 2e + e^{2} = 1 - e \\
& = &
(m\cdot m')\cdot (1-e).
\end{array}$$

We still have to check the bijective correspondence.
$$\begin{array}{rcl}
f_{e_f}(n, m)
& = &
n \cdot e_{f} + m\cdot (1-e_{f}) \\
& = &
n \cdot f(1,0) + m \cdot f(0,1) \\
& = &
f(1,0) + \cdots + f(1,0) + f(0,1) + \cdots + f(0,1) \\
& = &
f((1,0) + \cdots + (1,0)) + f((0,1) + \cdots + (0,1)) \\
& = &
f(n,0) + f(0,m) \\
& = &
f(n,m) \\
e_{f_e}
& = &
f_{e}(1,0) \\
& = &
1\cdot e + 0 \cdot (1-e) \\
& = &
e.
\end{array}$$
}

It is an elementary fact that the idempotent elements in a
\emph{commutative} ring form a Boolean algebra, see
\textit{e.g.}~\cite[I, 1.9]{Johnstone82}. We shall see that in the
non-commutative case these idempotents form an effect algebra, see
Proposition~\ref{RingPredProp}. The scalars are the idempotents of the
initial ring $\Z$, which are the elements $\{0,1\}$. An $n$-test in
$\op{\Rng}$, that is, a ring homomorphism $\Z^{n} \rightarrow R$ can
be identified with an $n$-tuple $e_{1}, \ldots, e_{n}$ of idempotent
elements of $R$ satisfying not only $e_{1} + \cdots + e_{n} = 1$ but
also $e_{i}\cdot e_{j} = 0$ for $i\neq j$. The latter multiplication
property follows from the fact that these maps $\Z^{n} \rightarrow R$
preserve multiplication.

\item \label{PredExDL} We write $\DL$ for the category of distributive
  lattices --- with finite joins $0, \vee$ and meets $1, \wedge$
  distributing over each other --- and lattices homomorphisms
  preserving these finite joins and meets. The 2-element lattice $2 =
  \{0,1\}$ is initial in $\DL$. In the opposite category $\op{\DL}$ a
  predicate $L\rightarrow 1+1$ is thus a lattice homomorphism $f\colon
  2\times 2 \rightarrow L$. It can be identified with an element $x\in
  L$ that has a complement: there is an $x'\in L$ with $x\wedge x' =
  0$ and $x\vee x' = 1$. It is easy to see that such a complement, if
  it exists, is necessarily unique. These `complementable' elements in
  $L$ form a Boolean sublattice of $L$. The scalars are the elements
  of the 2-element lattice $\{0,1\}$ --- which both have a
  complement. An $n$-test is an $n$-tuple of elements $x_{1}, \ldots,
  x_{n}$ with $x_{1} \vee \cdots \vee x_{n} = 1$ and $x_{i} \wedge
  x_{j} = 0$ for $i \neq j$.

\auxproof{
If both $x'$ and $x''$ are complements of $x$, then:
$$\begin{array}{rcl}
x'
& = &
x' \wedge 1 \\
& = &
x' \wedge (x \vee x'') \\
& = &
(x' \wedge x) \vee (x' \wedge x'') \\
& = &
0 \vee (x' \wedge x'') \\
& = &
x' \wedge x''
\end{array}$$

\noindent From this we get $x' \leq x''$. By symmetry we also get $x''
\leq x'$. Hence $x' = x''$.

The correspondence is given as follows.
\begin{itemize}
\item For $f\colon 2\times 2 \rightarrow L$ take $x = f(1,0)$ and $x' = 
f(0,1)$. Then:
$$\begin{array}{rcl}
x \wedge x'
& = &
f(1,0) \wedge f(0,1) \\
& = &
f((1,0) \wedge (0,1)) \\
& = &
f(1\wedge 0, 0 \wedge 1) \\
& = &
f(0,0) \\
& = &
0 \\
x\vee x'
& = &
f(1,0) \vee f(0,1) \\
& = &
f((1,0) \vee (0,1)) \\
& = &
f(1\vee 0, 0 \vee 1) \\
& = &
f(1,1) \\
& = &
1.
\end{array}$$

\item Given $x$ with complement $x'$ define $f\colon 2\times 2 \rightarrow L$
as $f(a,b) = (x\wedge a) \vee (x'\wedge b)$. This is a map of distributive
lattices:
$$\begin{array}{rcl}
f(0,0)
& = &
(x\wedge 0) \vee (x' \wedge 0) \\
& = &
0 \\
f((a_{1}, b_{1}) \vee (a_{2}, b_{2})) 
& = &
f(a_{1} \vee a_{2}, b_{1} \vee b_{2}) \\
& = &
(x \wedge (a_{1}\vee a_{2})) \vee (x' \wedge (b_{1} \vee b_{2})) \\
& = &
(x \wedge a_{1}) \vee (x\wedge a_{2}) \vee (x'\wedge b_{1}) \vee
  (x' \wedge b_{2}) \\
& = &
f(a_{1}, b_{1}) \vee f(a_{2}, b_{2}) \\
f(1,1)
& = &
(x\wedge 1) \vee (x' \wedge 1) \\
& = &
x \vee x' \\
& = &
1 \\
f((a_{1}, b_{1}) \wedge (a_{2}, b_{2})) 
& = &
f(a_{1} \wedge a_{2}, b_{1} \wedge b_{2}) \\
& = &
(x \wedge (a_{1}\wedge a_{2})) \vee (x' \wedge (b_{1} \wedge b_{2})) \\
& = &
((x\wedge a_{1}) \wedge (x\wedge a_{2})) \vee 
   ((x\wedge a_{1}) \wedge (x' \wedge b_{2})) \\
& & \qquad \vee
((x' \wedge b_{1}) \wedge (x\wedge a_{2})) \vee 
   ((x' \wedge b_{1}) \wedge (x' \wedge b_{2})) \\
& = &
((x\wedge a_{1}) \vee (x' \wedge b_{1})) \wedge
   ((x\wedge a_{2}) \vee (x' \wedge b_{2})) \\
& = &
f(a_{1}, b_{1}) \wedge f(a_{2}, b_{2})
\end{array}$$
\end{itemize}

Clearly, $0,1 \in L$ are each others complement. If $x$ has $x'$ as
complement, then we define $\neg x = x'$, which is again
complementable. Next assume $x,y$ have complements $x',y'$
respectively, then $x\vee y$ and $x\wedge y$ also have complements,
namely $x' \wedge y'$ and $x'\vee y'$ respectively, since:
$$\begin{array}{rcl}
(x\vee y) \wedge (x'\wedge y')
& = &
(x \wedge x' \wedge y') \vee (y \wedge x' \wedge y') \\
& = &
0 \vee 0 \\
& = &
0 \\
(x\vee y) \vee (x'\wedge y')
& = &
(x \vee y \vee x') \wedge (x \vee y \vee y') \\
& = &
1 \wedge 1 \\
& = &
1 \\
(x \wedge y) \wedge (x'\vee y') 
& = &
(x \wedge y \wedge x') \vee (x \wedge y \wedge y') \\
& = &
0 \vee 0 \\
& = & 
0 \\
(x \wedge y) \vee (x'\vee y') 
& = &
(x \vee x' \vee y') \wedge (y \vee x' \vee y') \\
& = &
1 \wedge 1 \\
& = &
1
\end{array}$$
}

\item \label{PredExBAEA} We conclude with three examples, as `sanity
  check'. In the category $\BA$ of Boolean algebras, the two element
  algebra $2 = \{0,1\}$ is initial, and maps of Boolean algebras
  $2\times 2 \rightarrow B$ correspond to elements of $B$. Hence all
  elements of $B$ occur as predicates $B\rightarrow 1+1$ in
  $\op{\BA}$. Similarly, in the category $\EMod$ of effect modules
  over $[0,1]$, this unit interval $[0,1]$ is initial, and effect
  module maps $[0,1]^{2} \rightarrow E$ correspond to elements of $E$.
  In the category $\Conv = \EM(\Dst)$ of convex sets, the coproduct
  $1+1$ is the unit interval $[0,1]$, since $1+1 \cong \Dst(1) +
  \Dst(1) \cong \Dst(1+1) \cong [0,1]$. Thus, the predicates $X
  \rightarrow 1+1$ in $\Conv$ are the affine maps $X \rightarrow
              [0,1]$, as used in the functor $\Hom(-,[0,1]) \colon
              \Conv \rightarrow \op{\EMod}$ in the adjunction $\Conv
              \rightleftarrows \op{\EMod}$ in
              Proposition~\ref{ConvEModAdjProp}.
\end{enumerate}
\end{exas}

\noindent One may wonder why several of the categories in the above list occur
in opposite form $\op{(-)}$. The appropriate way to see a morphism $L
\rightarrow K$ between two ``logical'' structures $L,K$ as a
computation is to interpret it as a predicate transformer computation
going in the other direction, from $K$ to $L$: it takes predicates on
the `post' state to predicates on the `pre' state of the computation.
For $C^*$-algebras this corresponds to Heisenberg's view on quantum
computation. In Dijkstra's weakest precondition
semantics~\cite{Dijkstra75} there is bijective correspondence between:
$$\begin{prooftree}
\xymatrix{\mbox{computations, as functions } X\ar[r]^-{f} & \Pow(Y)
   \mbox{ in the Kleisli category of the powerset monad}}
\Justifies
\xymatrix{\mbox{meet-preserving predicate transformers }
   \Pow(Y)\ar[r]_-{\wp(f)} & \Pow(X)}
\end{prooftree}$$

\noindent Notice that these predicate transformers operate in the
opposite direction. They can be described as a computations in the
opposite category $\op{(\CL_{\bigwedge})}$ of complete lattices and
meet-preserving functions. This is precisely like in the above
`logical' categories of the form $\op{(-)}$, see
also~\cite{Jacobs15b,Jacobs15a}.

In general the predicates $\Pred(1)$, \textit{i.e.}~the probabilities,
have a monoid structure. Later on, in Proposition~\ref{TwoMonoidProp},
it will be shown that this monoid is commutative in presence of a
(distributive) monoidal structure, following~\cite{KellyL80}. However,
this fact will not be used immediately.

\begin{lem}
\label{ScalarMonoidLem}
In a category $\cat{B}$ with coproduct $+$ and a final object $1$ the
collection $\Pred(1)$ of scalars, \textit{i.e.}~of predicates $1
\rightarrow 1+1$ on $1$, has a monoid structure, given by Kleisli
composition:
$$\xymatrix{
s\cdot r = s \klafter r = 
   \big(1\ar[r]^-{r} & 1+1\ar[r]^-{[s, \kappa_{2}]} & 1+1\big).
}$$


\noindent The truth predicate $1 = \kappa_{1}$ is unit for this
multiplication, and the falsum predicate $0 = \kappa_{2}$ is zero
element (meaning $0\cdot s = 0 = s\cdot 0$).
\end{lem}

\begin{myproof}
Easy, since this multiplication is Kleisli composition. \QED
\end{myproof}

On the Boolean scalars $\{0,1\}$ this multiplication is conjunction,
but on the unit interval $[0,1]$ of scalars in the other examples it
is multiplication.

\section{The effect algebra structure on predicates}\label{EAPredSec}

Having seen predicates and tests in the previous section we now give a
more precise description of the categorical properties that we
require. These requirements will be formulated below as our first out
of four assumptions. A category satisfying the first assumption is now
called an \emph{effectus}, since~\cite{JacobsWW15a}. We shall use
this name here too.

After the introduction of this first assumption, this section shows
two things:
\begin{itemize}
\item that these assumptions give us effect module structure
  on predicates;


\item that the example categories described in the previous section
  satisfy these assumptions.
\end{itemize}

\noindent Subsequent Sections~\ref{StatesSec} and~\ref{PredicateCoprodSec}
describe more consequences of the following assumptions.

In a category $\cat{B}$ with coproducts $+$ and a final object
$1\in\cat{B}$ the coproduct $X_{1}+X_{2}$ in $\cat{B}$ comes with
`partial projections' $\rhd_{i} \colon X_{1}+X_{2} \rightarrow X_{i}$,
described in $\cat{B}$ as:
\begin{equation}
\label{PartProjDiag}
\vcenter{\xymatrix{
X_{1}+1 & & X_{1}+X_{2}\ar[ll]_-{\rhd_{1} = \idmap+!}
   \ar[rr]^-{\rhd_{2} = [\kappa_{2}\after\,!, \kappa_{1}]} & & X_{2}+1
}}
\end{equation}

\noindent These partial projections form maps $\rhd_{i} \colon X_{1} +
X_{2} \rightarrow X_{i}$ in the Kleisli category $\KlL{\cat{B}}$ of
the lift monad $(-)+1$ on $\cat{B}$, where they can be described more
symmetrically as $\rhd_{1} = [\idmap, 0]$ and $\rhd_{2} = [0,
  \idmap]$. They are natural in this Kleisli category, since $\rhd_{i}
\klafter (f_{1}+f_{2}) = f_{i} \klafter \rhd_{i}$, where the $+$ on
the left-hand-side is understood in $\KlL{\cat{B}}$ --- which, in
$\cat{B}$, is given by $f_{1}+f_{2} = [(\kappa_{1}+\idmap) \after
  f_{1}, (\kappa_{2}+\idmap) \after f_{2}]$.  Obviously, there are
$n$-ary versions $\rhd_{i} \colon X_{1} + \cdots + X_{n} \rightarrow
X_{i}$ of these partial projections.


\auxproof{
$$\begin{array}{rcl}
\rhd_{1} \klafter (f_{1}+f_{2})
& = &
[\idmap+!, \kappa_{2}] \after [(\kappa_{1}+\idmap) \after f_{1}, 
   (\kappa_{2}+\idmap) \after f_{2}] \\
& = &
[[\kappa_{1}, \kappa_{2}] \after f_{1}, 
   [\kappa_{2}\after\,!, \kappa_{2}] \after f_{2}] \\
& = &
[f_{1}, \kappa_{2} \after [!,\idmap] \after f_{2}] \\
& = &
[f_{1}, \kappa_{2} \after \,!] \\
& = &
[f_{1}, \kappa_{2}] \after (\idmap+!) \\
& = &
f_{1} \klafter \rhd_{1}
\\
\rhd_{2} \klafter (f_{1}+f_{2})
& = &
[[\kappa_{2}\after\,!, \kappa_{1}],\kappa_{2}] \after 
   [(\kappa_{1}+\idmap) \after f_{1}, (\kappa_{2}+\idmap) \after f_{2}] \\
& = &
[[\kappa_{2}\after\,!, \kappa_{2}] \after f_{1}, 
   [\kappa_{1}, \kappa_{2}] \after f_{2}] \\
& = &
[\kappa_{2} \after \, !, f_{2}] \\
& = &
[f_{2}, \kappa_{2}] \after [\kappa_{2}\after\,!, \kappa_{1}] \\
& = &
f_{2} \klafter \rhd_{2}
\\
J(g_{1})+J(g_{2})
& = &
[(\kappa_{1}+\idmap) \after \kappa_{1} \after g_{1}, 
   (\kappa_{2}+\idmap) \after \kappa_{1} \after g_{2}] \\
& = &
[\kappa_{1} \after \kappa_{1} \after g_{1},
   \kappa_{1} \after \kappa_{2} \after g_{2}] \\
& = &
\kappa_{1} \after \kappa_{1} \after g_{1}, \kappa_{2} \after g_{2}] \\
& = &
\kappa_{1} \after (g_{1}+g_{2}) \\
& = &
J(g_{1}+g_{2}).
\end{array}$$
}

\begin{assumption}
\label{CoprodAss}
We use a category $\cat{B}$ that is an \emph{effectus}. This means
that it has a final object $1$ and finite coproducts $(+,0)$ such that
diagrams of the following form are pullbacks in $\cat{B}$.
\begin{equation}
\label{CoprodAssSquares}
\vcenter{\xymatrix{
A+ X\ar[r]^-{\idmap+ f}\ar[d]_{g+\idmap} & 
  A+ Y\ar[d]^{g+\idmap}
& &
A\ar@{=}[r]\ar[d]_{\kappa_1} & A\ar[d]^{\kappa_1}  \\
B+ X\ar[r]_-{\idmap+ f} & B+ Y
& &
A+ X\ar[r]_-{\idmap+ f} & A+ Y
}}
\end{equation}

\noindent Additionally, we require that the following two maps are
jointly monic in an effectus.
\begin{equation}
\label{CoprodOneAssJoint}
\vcenter{\xymatrix@C+2pc{
\llap{$(1$}+1)+1\ar@/^1ex/[rr]^-{\IV = [\rhd_{1},\kappa_{2}] = 
   [[\kappa_{1}, \kappa_{2}], \kappa_{2}]}
   \ar@/_1ex/[rr]_-{\XI = [\rhd_{2},\kappa_{2}] = 
   [[\kappa_{2}, \kappa_{1}], \kappa_{2}]} & & 1+1
}}
\end{equation}

\noindent The maps $\rhd_i$ are the partial projections
from~\eqref{PartProjDiag}. This joint monicity requirement means that
if $f,g$ satisfy $\IV \after f = \IV \after g$ and $\XI \after f = \XI
\after g$, then $f=g$. The symbols $\IV$ and $\XI$ should suggest what
these two maps do.
\end{assumption}

The pullback requirements are rather mild and hold in many situations.
The joint monicity of the maps~\eqref{CoprodOneAssJoint} however is
more restrictive, see also~\cite{Jacobs11c}. Later on in this section,
in Example~\ref{PowNonAssumpRem}, we show that it fails in the Kleisli
category of the powerset monad. Thus, non-deterministic computation,
in its standard form, does not fit in the current setting.
Lemma~\ref{JointMonicLem} elaborates on this joint monicity
requirement.  At this stage we note that it says that in the Kleisli
category of the lift monad the two partial projections $\rhd_{1},
\rhd_{2} \colon 1+1 \rightarrow 1$ are jointly monic.

At the end of Remark~\ref{LiftRem} we explain how the maps $\rhd_i$
may be seen as partial projections in the Kleisli category of the
lift monad, and how the joint monicity requirement translates
to this Kleisli category.

First we obtain some more pullbacks. We show that assumptions about
coproducts imply that they are \emph{disjoint}: the coprojections are
monic and have empty intersection (the pullback of different
coprojections is empty).  Hence our pullback requirements can be
understood as a form of `generalised' disjointness. In fact, the
situation is:
$$\xymatrix{
{\left(\begin{array}{c} \mbox{disjoint and stable} \\[-.4em]
   \mbox{coproducts} \end{array}\right)}\ar@{=>}[r] & 
{\left(\begin{array}{c} \mbox{pullback requirements} \\[-.4em]
   \mbox{in Assumption~\ref{CoprodAss}~\eqref{CoprodAssSquares}} 
   \end{array}\right)}\ar@{=>}[r] & 
{\left(\begin{array}{c} \mbox{disjoint} \\[-.4em]
   \mbox{coproducts} \end{array}\right)}
}$$

\noindent Categories with coproducts that are disjoint and stable
(under pullback) are also called \emph{extensive},
see~\cite{CarboniLW93}.

\begin{lem} 
\label{CoprodAssLem}
In an effectus $\cat{B}$, the initial object $0\in\cat{B}$ is strict,
the coprojections $\kappa_{i} \colon A_{i} \rightarrow A_{1}+A_{2}$
are monic and have empty intersection, and diagrams of the following
form are pullbacks.
\begin{equation}
\label{CoprodAssLemSquares}
\vcenter{\xymatrix@C-1pc{
A\ar[rr]^-{g}\ar[d]_{\kappa_1} & & B\ar[d]^{\kappa_1}  
& \qquad &
(A+X)+1\ar[rrr]^-{(g+\idmap)+\idmap}\ar[d]_{[\rhd_{1},\kappa_{2}]} & & &
   (B+X)+1\ar[d]^{[\rhd_{1},\kappa_{2}]}
\\
A+ X\ar[rr]_-{g+f} & & B+ Y
& &
A+1\ar[rrr]_-{g+\idmap} & & & B+1
}}
\end{equation}

\noindent (The diagram on the right is a bit simpler when formulated
in the Kleisli category $\KlL{\cat{B}}$ of the lift monad, but we
prefer to avoid confusion between $\cat{B}$ and $\KlL{\cat{B}}$.)
\end{lem}

\begin{myproof}
Strictness of the initial object $0\in\cat{B}$ means that each
map $f\colon X \rightarrow 0$ is an isomorphism. For such a map,
we have to prove that the composite $X \rightarrow 0 \rightarrow X$
is the identity. Consider the diagram:
$$\xymatrix@R-.5pc@C-1pc{
X\ar@/^2ex/[drrr]^{f}\ar@/_2ex/[ddr]_{\kappa_2}\ar@{-->}[dr] \\
& 0\pullback\ar@{=}[rr]\ar[d]_-{\kappa_{1}} & & 0\ar[d]^{\kappa_{1} = \kappa_{2}} \\
& 0+X\ar[rr]_-{\idmap+f} & & 0+0
}$$

\noindent The rectangle is a pullback, as
in~\eqref{CoprodAssSquares}. By initiality of $0$, we have
$\kappa_{1} = \kappa_{2} \colon 0 \rightarrow 0+0$, so that the outer
diagram commutes. Then we get the dashed map, as indicated, which must
be $f\colon X \rightarrow 0$. But then $!_{X} \after f = [!_{X},
  \idmap[X]] \after \kappa_{1} \after f = [!_{X}, \idmap[X]] \after
\kappa_{2} = \idmap[X]$.

We now prove a special case of the diagram on the left
in~\eqref{CoprodAssLemSquares} with $f=\idmap$:
\begin{equation}
\label{CoprodAssLemSquareId}
\vcenter{\xymatrix{
A\ar[r]^-{g}\ar[d]_{\kappa_1} & B\ar[d]^{\kappa_1} \\
A+ X\ar[r]_-{g+\idmap} & B+ X
}}
\end{equation}

\noindent This is a pullback since it can be obtained from the square
on the left in~\ref{CoprodAssSquares}:
$$\xymatrix@C+1pc{
A\pullback\ar[r]_-{\kappa_2}^-{\cong}\ar[d]_{\kappa_1}\ar@/^5ex/[rrr]^-{g} & 
   0+A\pullback\ar[r]_-{\idmap+g}\ar[d]_{!+\idmap} & 
   0+B\pullback\ar[r]_-{[!,\idmap]}^-{\cong}\ar[d]^{!+\idmap} & B\ar[d]^{\kappa_1}
\\
A+X\ar[r]^-{[\kappa_{2},\kappa_{1}]}_-{\cong}\ar@/_5ex/[rrr]_-{g+\idmap} &
   X+A\ar[r]^-{\idmap+g} & 
   X+B\ar[r]^-{[\kappa_{2},\kappa_{1}]}_-{\cong} &
   B+X
}$$

\noindent The square on the left in~\eqref{CoprodAssLemSquares} is a
pullback since it can be obtained from the square on the right
in~\eqref{CoprodAssSquares} and the previous square
in~\eqref{CoprodAssLemSquareId}:
$$\xymatrix@C+1pc{
A\pullback\ar[r]_-{g}\ar[d]_{\kappa_1}\ar@/^3ex/[rr]^-{g} & 
   B\pullback\ar@{=}[r]\ar[d]_{\kappa_{1}} & B\ar[d]^{\kappa_{1}}
\\
A+X\ar[r]^-{g+\idmap}\ar@/_3ex/[rr]_-{g+f} & B+X\ar[r]^-{\idmap+f} & B+Y
}$$

\noindent The first coprojection $\kappa_{1} \colon A \rightarrow A+X$ is
monic since, using the pullback on the right in~\eqref{CoprodAssSquares},
we get a diagram as on the left below.
$$\vcenter{\xymatrix@C+0pc@R-.5pc{
A\pullback\ar@{=}[r]\ar@{=}[d] & 
   A\pullback\ar@{=}[r]\ar[d]_{\kappa_1}^{\cong} & A\ar[d]^{\kappa_{1}}
\\
A\ar[r]^-{\kappa_1}_-{\cong}\ar@/_4ex/[rr]_-{\kappa_{1}} & 
   A+0\ar[r]^-{\idmap+\,!} & A+X
}}
\qquad\qquad
\vcenter{\xymatrix@C+0pc@R-.5pc{
0\ar[dr]^{\cong}\ar[dd]\ar[rr] & & 
   X\ar[d]^{\kappa_{2}}_{\cong}\ar@/^5ex/[dd]^{\kappa_2} \\
&  0+0\pullback\ar[r]\ar[d] & 0+X\ar[d]_{!+\idmap}
\\
A\ar[r]^-{\kappa_1}_-{\cong}\ar@/_4ex/[rr]_-{\kappa_{1}} & 
   A+0\ar[r]^-{\idmap+\,!} & A+X
}}$$

\noindent A similar argument applies to the second coprojection
$\kappa_2$.  The above diagram on the right shows that the
intersection (pullback) of $\kappa_{1}, \kappa_{2}$ is the initial
object $0$. 

We turn to the square on the right in~\eqref{CoprodAssLemSquares}.
Let $f\colon Y\rightarrow A+1$ and $h\colon Y \rightarrow (B+X)+1$
satisfy $(g+\idmap) \after f = [\rhd_{1}, \kappa_{2}] \after h$.
Consider the situation below, where in general, $\alpha$ is the
associativity isomorphism $U+(V+W)\rightarrow (U+V)+W$, given
explicitly by $\alpha = [\kappa_{1} \after \kappa_{1},
  \kappa_{2}+\idmap]$ and $\alpha^{-1} = [\idmap+\kappa_{1},
  \kappa_{2} \after \kappa_{2}]$.
$$\xymatrix@R-.5pc{
Y\ar@/_2ex/[ddr]_-{f}\ar@/^2ex/[drr]^{\alpha^{-1} \after h}\ar@{-->}[dr]^(0.6){k}
\\
& A+(X+1)\ar[r]^-{g+\idmap}\ar[d]_{\idmap+!}\pullback & 
   B+(X+1)\ar[d]^{\idmap+! = [\rhd_{1},\kappa_{2}] \after \alpha} 
\\
& A + 1\ar[r]_-{g+\idmap} & B+1
}$$
\auxproof{
\noindent The outer diagram commutes, since:
$$\begin{array}{rcl}
[\rhd_{1},\kappa_{2}] \after \alpha
& = &
[\idmap+!,\kappa_{2}] \after [\kappa_{1} \after \kappa_{1}, \kappa_{2}+\idmap] \\
& = &
[\kappa_{1} , [\kappa_{2}\after\,!!, \kappa_{2}]] \\
& = &
[\kappa_{1} , \kappa_{2} \after [!,\idmap]] \\
& = &
\idmap+!.
\end{array}$$
}

\noindent Then $k' = \alpha \after k \colon Z \rightarrow (X+A)+1$ is the
required mediating map. \QED

\auxproof{
Since:
$$\begin{array}{rcl}
[\rhd_{1}, \kappa_{2}] \after k'
& = &
[\rhd_{1}, \kappa_{2}] \after \alpha \after k \\
& = &
[\idmap+\,!, \kappa_{2}] \after 
   [\kappa_{1} \after \kappa_{1}, \kappa_{2}+\idmap] \after k\\
& = &
[\kappa_{1}, [\kappa_{2}\after\,!, \kappa_{2}]] \after k \\
& = &
[\kappa_{1}, \kappa_{2}\after [!, \kappa_{2}]] \after k \\
& = &
[\kappa_{1}, \kappa_{2}\after \, !] \after k \\
& = &
(\idmap+!) \after k \\
& = &
f 
\\
((g+\idmap)+\idmap) \after k'
& = &
((g+\idmap)+\idmap) \after \alpha \after k \\
& = &
\alpha \after (g+\idmap) \after k \\
& = &
\alpha \after \alpha^{-1} \after h \\
& = &
h.
\end{array}$$
}
\end{myproof}

The joint monicity requirement in~\eqref{CoprodOneAssJoint} says that
the two maps $\IV, \XI \colon (1+1)+1 \rightarrow 1+1$ are jointly
monic. This simple formulation can be generalised in two ways: using
arbitrary objects $X$ instead of $1$, and using arbitrary many $1$'s
instead of only 2. The equivalence of these generalisations to the
formulation in~\eqref{CoprodOneAssJoint} is due to Kenta Cho. (A
slightly stronger formulation is used in~\cite{JacobsWW15a}.)

\begin{lem}
\label{JointMonicLem}
Let $\cat{B}$ be a category with finite coproducts and a final object,
in which the squares in~\eqref{CoprodAssSquares} are pullbacks. Then
the following statements are equivalent.
\begin{enumerate}
\item \label{JointMonicLemTwoOnes} The two maps $\IV, \XI \colon
  (1+1)+1 \rightarrow 1+1$ in~\eqref{CoprodOneAssJoint} are jointly
  monic.

\item \label{JointMonicLemTwoX} The two maps $\IV, \XI \colon (X+X)+1
  \rightarrow X+1$ are jointly monic, where as before, $\IV =
              [\idmap+\;!, \kappa_{2}]$ and $\XI =
              [[\kappa_{2}\after\;!, \kappa_{1}], \kappa_{2}]$.

\item \label{JointMonicLemManyX} For each $n\geq 1$, the $n$ maps
  $[\rhd_{i}, \kappa_{2}] \colon n\cdot X + 1 \rightarrow X+1$ are
  jointly monic.

\auxproof{
For each non-zero $n\in\NNO$ the following families of maps are
disjoint monic.
\begin{equation}
\label{CoprodAssJoint}
\vcenter{\xymatrix@C+1pc{
n\cdot X + 1\ar[r]^-{[\rhd_{i}, \kappa_{2}]} & X+1,
}}
\end{equation}

\noindent where for $1\leq i\leq n$ the ``partial projection''
$\rhd_{i} \colon X + \cdots + X \rightarrow X+1$ is defined as:
$$\begin{array}{rcl}
\rhd_{i} \after \kappa_{j}
& = &
\left\{\begin{array}{ll}
\kappa_{1} \quad & \mbox{if } i=j \\
\kappa_{2} \after \, !_{X} \quad & \mbox{otherwise.}
\end{array}\right.
\end{array}$$

\noindent Thus, $\rhd_{i}$ sends the $i$-th coproduct component in
$n\cdot X = X + \cdots + X$ to $X$ in $X+1$ and all the rest to $1$.
}
\end{enumerate}
\end{lem}

\begin{myproof}
The implication $\eqref{JointMonicLemTwoX} \Rightarrow
\eqref{JointMonicLemManyX}$ works by induction. We concentrate on the
implication $\eqref{JointMonicLemTwoOnes} \Rightarrow
\eqref{JointMonicLemTwoX}$.  Let $f, g\colon Y \rightarrow (X+X)+1$
satisfy $[\rhd_{i}, \kappa_{2}] \after f = [\rhd_{i}, \kappa_{2}]
\after g \colon Y \rightarrow X+1$ for $i = 1,2$. Consider the
diagram:
$$\xymatrix@R-.5pc{
(X+X)+1\pullback\ar[rr]^-{(!+\idmap)+\idmap}\ar[d]_{(\idmap+\,!)+\idmap} & &
(1+X)+1\pullback\ar[rr]^-{[\rhd_{2},\kappa_{2}]}\ar[d]^{(\idmap+\,!)+\idmap}
& &
   X+1\ar[d]^{!+\idmap}
\\
(X+1)+1\pullback\ar[rr]^-{(!+\idmap)+\idmap}\ar[d]_{[\rhd_{1},\kappa_{2}]}
& &
   (1+1)+1\ar[rr]_-{[\rhd_{2},\kappa_{2}]}\ar[d]^{[\rhd_{1},\kappa_{2}]} & &
   1+1
\\
X+1\ar[rr]_-{!+\idmap} & & 1+1
}$$

\noindent The lower and right rectangles are pullbacks, via the
diagram on the right in~\eqref{CoprodAssLemSquares}. The upper left
rectangle is a pullback by the diagram on the left
in~\eqref{CoprodAssSquares}, using associativity of $+$. 

The equations $[\rhd_{i}, \kappa_{1}] \after f = [\rhd_{i}, \kappa_{1}]
\after g$ imply $[\rhd_{i}, \kappa_{1}] \after ((!+\,!)+\idmap)
\after f = [\rhd_{i}, \kappa_{1}] \after ((!+\,!)+\idmap) \after g$,
so that $((!+\,!)+\idmap) \after f = ((!+\,!)+\idmap) \after g$,
because the maps $[\rhd_{i}, \kappa_{1}] \colon (1+1)+1 \rightarrow
1+1$ are jointly monic, by assumption. But then:
$$\left\{\begin{array}{rcll}
((!+\idmap)+\idmap) \after f
& = &
((!+\idmap)+\idmap) \after g \qquad &
   \mbox{by the upper right pullback} \\
((\idmap+\,!)+\idmap) \after f
& = &
((\idmap+\,!)+\idmap) \after g &
   \mbox{by the lower left pullback.}
\end{array}\right.$$

\noindent Hence $f=g$ by the upper left pullback. \QED

\auxproof{
The implication $\eqref{JointMonicLemTwoX} \Rightarrow
\eqref{JointMonicLemManyX}$ is done by induction on $n$. The case $n =
1$, where $[\rhd_{1}, \kappa_{2}] = \idmap = [\rhd_{2}, \kappa_{2}]
\colon X+1 \rightarrow X+1$, holds, and the case $n=2$ holds by
induction. 

Now let $f, g \colon Y \rightarrow (n+1)\cdot X + 1 = (n\cdot X + X) +
1$ satisfy $[\rhd_{i}, \kappa_{2}] \after f = [\rhd_{i}, \kappa_{2}]
\after g$, for $i\leq n+1$. Take $f' = [\rhd_{1},\kappa_{2}] \after
f\colon Y \rightarrow (n+1)\cdot X +1 \rightarrow n\cdot X + 1$.  and
similarly, $g' = [\rhd_{1},\kappa_{2}] \after g$. Then for each $i\leq
n$, $[\rhd_{i},\kappa_{2}] \after f' = [\rhd_{i},\kappa_{2}] \after
g'$. Hence $f' = g'$. But since $[\rhd_{2},\kappa_{2}] \after f =
[\rhd_{n+1},\kappa_{2}] \after f = [\rhd_{n+1},\kappa_{2}] \after g =
[\rhd_{2},\kappa_{2}] \after g$, we get $f = g$
by~\eqref{JointMonicLemTwoX}.
}
\end{myproof}

The main reason for the assumptions about coproducts in the beginning
of this section lies in the following definition and result. The
specific definition of sums $\ovee$ via bounds that we use below goes
back to~\cite{Jacobs11c}. But such techniques to get partially
additive structure are much older, see~\cite{ArbibM86}.

\begin{defi}
\label{OveeDef}
Let $\cat{B}$ be a category satisfying Assumption~\ref{CoprodAss} and 
let $p,q\colon X \rightarrow 1+1$ be two predicates on $X\in\cat{B}$.
\begin{enumerate}
\item These two predicates are \emph{orthogonal}, written as
  $p\orthogonal q$, if there is a \emph{bound} $b \colon X \rightarrow
  (1+1)+1$ such that:
$$\begin{array}{rclcrcl}
\IV \after b
& = &
p
& \qquad\mbox{and}\qquad &
\XI \after b
& = &
q.
\end{array}$$

\noindent According to Assumption~\ref{CoprodAss} these two maps $\IV,
\XI \colon (1+1)+1 \rightarrow 1+1$ are jointly monic. Hence there is
at most one such a bound $b$.

\item If $p\orthogonal q$ via bound $b$ then we define a new ``sum''
  predicate $p \ovee q \colon X \rightarrow 1+1$ as:
$$\begin{array}{rcccl}
p\ovee q
& = &
(\nabla+\idmap) \after b
& = &
([\idmap,\idmap]+\idmap) \after b.
\end{array}$$
\end{enumerate}

\noindent Pictorially we have the following situation, with bound $b$
for predicates $p,q$.
$$\begin{array}{ccccc}
\vcenter{\xymatrix@R-.7pc@C-2.3pc{
& \qquad & & & & & & \llap{(}1\ar@{|->}[dr]  
   & \!\!+\!\! & 1\rlap{)}\ar@{|->}[dr] & +\!\! & 1\ar@{|->}[dl] \\
X\ar[urrrrr]^{b}\ar[rrrrrrrr]_{p} & & & & & & & & 1 & + & 1
}}
& \quad &
\vcenter{\xymatrix@R-.5pc@C-2.3pc{
& \qquad & & & & & & \llap{(}1\ar@{|->}[drrr]  
   & \!\!+\!\! & 1\rlap{)}\ar@{|->}[dl] & +\!\! & 1\ar@{|->}[dl] \\
X\ar[urrrrr]^{b}\ar[rrrrrrrr]_{q} & & & & & & & & 1 & + & 1
}}
& \quad &
\vcenter{\xymatrix@R-.5pc@C-2.3pc{
& \qquad & & & & & & \llap{(}1\ar@{|->}[dr]  
   & \!\!+\!\! & 1\rlap{)}\ar@{|->}[dl] & +\!\! & 1\ar@{|->}[dl] \\
X\ar[urrrrr]^{b}\ar[rrrrrrrr]_{p\ovee q} & & & & & & & & 1 & + & 1
}}
\\
p = \IV \after b
& &
1 = \XI \after b
& &
p\ovee q = (\nabla+\idmap) \after b.
\end{array}$$
\end{defi}

\noindent The next result shows how the pullback properties from
Assumption~\ref{CoprodAss} and Lemma~\ref{CoprodAssLem} yield effect
algebra structure on predicates.

\begin{prop}
\label{CoprodEAProp}
Let $\cat{B}$ be a category satisfying Assumption~\ref{CoprodAss}.
\begin{enumerate}
\item \label{CoprodEAPropObj} For each object $X\in\cat{B}$ the
  collection $\Pred(X)$ of predicates $X \rightarrow 1+1$ on $X$ forms
  an effect algebra, via $\ovee, 0, (-)^{\perp}$.

\item \label{CoprodEAPropMor} For each map $f\colon Y \rightarrow X$
  the substitution function $f^{*} = (-) \after f \colon \Pred(X)
  \rightarrow \Pred(Y)$ is a map of effect algebras.
\end{enumerate}

\noindent Thus we obtain a functor $\Pred\colon\cat{B} \rightarrow \op{\EA}$.
\end{prop}

\begin{myproof}
The operation $\ovee$ is commutative, since if $b$ is a bound for
$p,q$, then $b' = ([\kappa_{2}, \kappa_{1}] + \idmap) \after b$ is a
bound for $q,p$:
$$\begin{array}{rcl}
\IV \after b'
& = &
[\idmap,\kappa_{2}] \after ([\kappa_{2}, \kappa_{1}] + \idmap) \after b
\hspace*{\arraycolsep} = \hspace*{\arraycolsep}
[[\kappa_{2}, \kappa_{1}], \kappa_{2}] \after b 
\hspace*{\arraycolsep} = \hspace*{\arraycolsep}
q \\
\XI \after b'
& = &
[[\kappa_{2}, \kappa_{1}], \kappa_{2}] \after 
   ([\kappa_{2}, \kappa_{1}] + \idmap) \after b \\
& = &
[[\kappa_{2}, \kappa_{1}] \after [\kappa_{2}, \kappa_{1}], \kappa_{2}]
   \after b 
\hspace*{\arraycolsep} = \hspace*{\arraycolsep}
[\idmap, \kappa_{2}] \after b 
\hspace*{\arraycolsep} = \hspace*{\arraycolsep}
p.
\end{array}$$

\noindent This $b$ and $b'$ yield the same sum:
$$\begin{array}{rcl}
q\ovee p
\hspace*{\arraycolsep} = \hspace*{\arraycolsep}
(\nabla+\idmap) \after b' 
& = &
(\nabla+\idmap) \after ([\kappa_{2}, \kappa_{1}] + \idmap) \after b \\
& = &
((\nabla \after [\kappa_{2}, \kappa_{1}]) + \idmap) \after b 
\hspace*{\arraycolsep} = \hspace*{\arraycolsep}
(\nabla+\idmap) \after b 
\hspace*{\arraycolsep} = \hspace*{\arraycolsep}
p\ovee q.
\end{array}$$

\noindent The zero predicate $0 = \kappa_{2} \after\; !_{X} \colon X \rightarrow
1+1$ is neutral element for $\ovee$: the equation $0\ovee p = p$ holds
via the bound $b = (\kappa_{2}+\idmap) \after p \colon X \rightarrow
(1+1)+1$.

\auxproof{
$$\begin{array}{rcl}
[\idmap,\kappa_{2}] \after b
& = &
[\idmap,\kappa_{2}] \after (\kappa_{2}+\idmap) \after p \\
& = &
[\kappa_{2}, \kappa_{2}] \after p \\
& = &
\kappa_{2} \after \nabla \after p \\
& = &
\kappa_{2} \after \; !_{X} \\
& = &
0 \\
{[[\kappa_{2}, \kappa_{1}], \kappa_{2}]} \after b
& = &
{[[\kappa_{2}, \kappa_{1}], \kappa_{2}]} \after (\kappa_{2}+\idmap) 
   \after p \\
& = &
[\kappa_{1}, \kappa_{2}] \after p \\
& = &
p \\
0\ovee p
& = &
(\nabla+\idmap) \after b \\
& = &
(\nabla+\idmap) \after (\kappa_{2}+\idmap) \after p \\
& = &
p.
\end{array}$$
}

For associativity of $\ovee$ assume predicates $p,q,r\colon X \rightarrow
1+1$ where $p \orthogonal q$, say via bound $a$, and $(p \ovee q)
\orthogonal r$, via bound $b$. Thus we have equations:
$$\begin{array}{lcl}
\left\{\begin{array}{rcl}
\IV \after a & = & p \\
\XI \after a & = & q \\
(\nabla+\idmap) \after a & = & p\ovee q
\end{array}\right.
& \qquad\qquad &
\left\{\begin{array}{rcl}
\IV \after b & = & p\ovee q \\
\XI \after b & = & r \\
(\nabla+\idmap) \after b & = & (p\ovee q) \ovee r.
\end{array}\right.
\end{array}$$

\noindent We consider the following pullback, occurring
on the left in~\eqref{CoprodAssSquares},
$$\xymatrix@R-1.8pc{
X\ar[dd]_{b}\ar@/^2ex/[ddrrr]^-{a}\ar@{-->}[ddr]^(0.55){c} \\
\\
(1+1)+1\ar@/_2ex/[ddr]_-{\alpha^{-1}}^-{\cong} & 
   (1+1)+(1+1)\ar[dd]_{\nabla+\idmap}\ar[rr]^-{\idmap+\nabla}\pullback & &
   (1+1)+1\ar[dd]^{\nabla+\idmap} \\
& & \\
& 1+(1+1)\ar[rr]_-{\idmap+\nabla} & & 1+1
}$$

\noindent where $\alpha = [\kappa_{1} \after \kappa_{1}, \kappa_{2} +
  \idmap]$ is the obvious associativity isomorphism.

We now take $c' = [[\kappa_{2} , \kappa_{1} \after \kappa_{1}],
  \kappa_{2}+\idmap] \after c \colon X \rightarrow (1+1)+1$. This $c'$
is a bound for $q$ and $r$, giving $q\orthogonal r$. Next, the map
$c'' = [\kappa_{1}, \kappa_{2}+\idmap] \after c \colon X \rightarrow
(1+1)+1$ is a bound for $p$ and $q\ovee r = (\nabla+\idmap) \after
c'$. Finally we obtain the associativity of $\ovee$:
$$\begin{array}{rcl}
p\ovee (q\ovee r)
\hspace*{\arraycolsep} = \hspace*{\arraycolsep}
(\nabla+\idmap) \after c'' 
& = &
(\nabla+\idmap) \after [\kappa_{1}, \kappa_{2}+\idmap] \after c \\
& = &
[[\kappa_{1}, \kappa_{1}], [\kappa_{1}, \kappa_{2}]] \after c
\\
& = &
[[\kappa_{1},\kappa_{1}], \kappa_{2}] \after 
   [\kappa_{1} \after \kappa_{1} \after \nabla, \kappa_{2}+\idmap] 
   \after c \\
& = &
[[\kappa_{1},\kappa_{1}], \kappa_{2}] \after \alpha \after 
  (\nabla+\idmap) \after c \\
& = &
(\nabla+\idmap) \after b 
\hspace*{\arraycolsep} = \hspace*{\arraycolsep}
(p\ovee q) \ovee r.
\end{array}$$

\auxproof{ 
This $c'$ is a bound for $q$ and $r$.
$$\begin{array}{rcl}
[\idmap,\kappa_{2}] \after c'
& = &
[\idmap,\kappa_{2}] \after [[\kappa_{2} , \kappa_{1} \after
    \kappa_{1}], \kappa_{2}+\idmap] \after c \\
& = &
{[[\kappa_{2},\kappa_{1}],[\kappa_{2},\kappa_{2}]]} \after c \\
& = &
{[[\kappa_{2},\kappa_{1}],\kappa_{2}]} \after (\idmap+\nabla) \after c \\
& = &
{[[\kappa_{2},\kappa_{1}],\kappa_{2}]} \after a \\
& = &
q \\
{[[\kappa_{2}, \kappa_{1}], \kappa_{2}]} \after c' 
& = &
{[[\kappa_{2}, \kappa_{1}], \kappa_{2}]} \after 
   [[\kappa_{2} , \kappa_{1} \after
    \kappa_{1}], \kappa_{2}+\idmap] \after c \\
& = &
{[[\kappa_{2},\kappa_{2}], [\kappa_{1}, \kappa_{2}]]} \after c \\
& = &
{[[\kappa_{2}, \kappa_{1}], \kappa_{2}]} \after 
   [\kappa_{1} \after \kappa_{1} \after \nabla, \kappa_{2}+\idmap] 
   \after c \\
& = &
{[[\kappa_{2}, \kappa_{1}], \kappa_{2}]} \after \alpha \after
   (\nabla+\idmap) \after c \\
& = &
{[[\kappa_{2}, \kappa_{1}], \kappa_{2}]} \after b \\
& = &
r.
\end{array}$$

The map $c'' = [\kappa_{1}, \kappa_{2}+\idmap] \after c \colon X
\rightarrow (1+1)+1$ is then a bound for $p$ and $q\ovee r =
(\nabla+\idmap) \after c'$.
$$\begin{array}{rcl}
[\idmap,\kappa_{2}] \after c''
& = &
[\idmap,\kappa_{2}] \after [\kappa_{1}, \kappa_{2}+\idmap] \after c \\
& = &
[\idmap, [\kappa_{2},\kappa_{2}]] \after c \\
& = &
[\idmap, \kappa_{2}] \after (\idmap+\nabla) \after c \\
& = &
[\idmap, \kappa_{2}] \after a \\
& = &
p \\
{[[\kappa_{2}, \kappa_{1}], \kappa_{2}]} \after c''
& = &
{[[\kappa_{2}, \kappa_{1}], \kappa_{2}]} \after 
   [\kappa_{1}, \kappa_{2}+\idmap] \after c \\
& = &
{[[\kappa_{2}, \kappa_{1}], [\kappa_{1}, \kappa_{2}]]} \after c \\
& = &
[[\kappa_{1}, \kappa_{1}], \kappa_{2}] \after
    [[\kappa_{2}, \kappa_{1} \after \kappa_{1}], 
     [\kappa_{1} \after \kappa_{2}, \kappa_{2}]] \after c \\
& = &
(\nabla+\idmap) \after [[\kappa_{2} , \kappa_{1} \after
    \kappa_{1}], \kappa_{2}+\idmap] \after c \\
& = &
(\nabla+\idmap) \after c' \\
& = &
q\ovee r.
\end{array}$$
}

\noindent It is not hard to see that $b = \kappa_{1} \after p \colon X
\rightarrow (1+1)+1$ is a bound for $p$ and $p^{\perp} = p \after
            [\kappa_{2}, \kappa_{1}]$, yielding $p\ovee p^{\perp} =
            1$.

\auxproof{
$$\begin{array}{rcl}
[\idmap, \kappa_{2}] \after b
& = &
[\idmap, \kappa_{2}] \after \kappa_{1} \after p \\
& = &
p \\
{[[\kappa_{2}, \kappa_{1}], \kappa_{2}]} \after b
& = &
{[[\kappa_{2}, \kappa_{1}], \kappa_{2}]} \after \kappa_{1} \after p \\
& = &
[\kappa_{2}, \kappa_{1}] \after p \\
& = &
p^{\perp} \\
(\nabla+\idmap) \after b 
& = &
(\nabla+\idmap) \after \kappa_{1} \after p \\
& = &
\kappa_{1} \after \nabla \after p \\
& = &
\kappa_{1} \\
& = & 
1.
\end{array}$$
}

Next we show that $p^{\perp}$ is the only predicate with $p\ovee
p^{\perp} = 1$. So assume also for $q\colon X \rightarrow 1+1$ we have
$p\ovee q = 1$, say via bound $b$. Then $\IV \after b = p$, $\XI \after
b = q$ and $(\nabla+\idmap) \after b = 1 = \kappa_{1} \after\;
!_{X}$. We use the pullback~\eqref{CoprodAssLemSquares} in:
\begin{equation}
\label{TrueBoundDiag}
\vcenter{\xymatrix@R-1.8pc@C-1pc{
X\ar@/_2ex/[ddddrr]_{!_X}\ar@/^2ex/[ddrrrr]^-{b}\ar@{-->}[ddrr]^-{c} \\
\\
& & 1+1\ar[dd]_{\nabla}\ar[rr]^-{\kappa_1}\pullback & &
   (1+1)+1\ar[dd]^{\nabla+\idmap} \\
& & & \\
& & 1\ar[rr]_-{\kappa_1} & & 1+1
}}
\end{equation}

\noindent The fact that the bound $b$ is of the form $\kappa_{1}\after
c$ is enough to obtain $q = p^{\perp}$:
$$\begin{array}{rcl}
p^{\perp}
\hspace*{\arraycolsep} = \hspace*{\arraycolsep}
[\kappa_{2}, \kappa_{1}] \after p
\hspace*{\arraycolsep} = \hspace*{\arraycolsep}
[\kappa_{2}, \kappa_{1}] \after [\idmap,\kappa_{2}] \after b 
& = &
[[\kappa_{2}, \kappa_{1}], \kappa_{1}] \after \kappa_{1} \after c \\
& = &
[[\kappa_{2}, \kappa_{1}], \kappa_{2}] \after \kappa_{1} \after c \\
& = &
[[\kappa_{2}, \kappa_{1}], \kappa_{2}] \after b 
\hspace*{\arraycolsep} = \hspace*{\arraycolsep}
q.
\end{array}$$
Finally, suppose $1 \orthogonal p$; we must prove $p=0 \colon X
\rightarrow 1+1$. We may assume a bound $b$ with $\IV \after b = 1 =
\kappa_{1} \after\; !_{X}$ and $\XI \after b = p$.  Now we use the
pullback on the right in~\eqref{CoprodAssSquares} to obtain unique map
$X \rightarrow 1$ as mediating map in:
$$\xymatrix@R-1.8pc{
X\ar@/_2ex/[ddddr]_{!_X}\ar[rr]^-{b}\ar@{-->}[ddr]^(0.6){!_X} & &
   (1+1)+1\ar@/^2ex/[ddr]_-{\cong}^(0.75){\alpha^{-1}}
      \ar`r[rr]`[dddd]^{[\idmap,\kappa_{2}]}[rdddd] & & \\
\\
& 1\ar@{=}[dd]\ar[rr]^-{\kappa_1}\pullback & &
   1+(1+1)\ar[dd]^{\idmap+\nabla} & & \\
& & \\
& 1\ar[rr]_-{\kappa_1} & & 1+1
}$$

\noindent As a result, $b = \alpha \after \kappa_{1} \after\; !_{X} = 
\kappa_{1} \after \kappa_{1} \after\; !_{X}$. Then:
$$\begin{array}{rcccccccl}
p
& = &
\XI \after b
& = &
[[\kappa_{2}, \kappa_{1}], \kappa_{2}] \after \kappa_{1} \after \kappa_{1}
   \after\; !_{X}
& = &
\kappa_{2} \after\; !_{X}
& = &
0.
\end{array}$$\medskip

\noindent For the second point of the proposition we have to prove that
substitution preserves $\ovee$. So let $f\colon Y \rightarrow X$ be a
map in $\cat{B}$, and let $p,q\colon X \rightarrow 1+1$ be orthogonal
predicates, say via bound $b\colon X \rightarrow (1+1)+1$. Then $b
\after f$ is a bound for $p \after f = f^{*}(p)$ and $q \after f =
f^{*}(q)$. Hence:
$$\begin{array}{rcccccl}
f^{*}(p) \ovee f^{*}(q)
& = &
(\nabla+\idmap) \after b \after f
& = &
(p\ovee q) \after f
& = &
f^{*}(p\ovee q).
\end{array}\eqno{\qEd}$$
\end{myproof}\bigskip

\noindent In Lemma~\ref{ScalarMonoidLem} we have seen that the collection
$\Pred(1)$ of scalars (or probabilities) $1 \rightarrow 1+1$ carries a
monoid structure $s\cdot r = [s, \kappa_{2}] \after r$. By the
previous result $\Pred(1)$ also carries an effect algebra structure.
It turns out that these two structures interact appropriately to form
an effect monoid, that is, a monoid in the monoidal category of effect
algebras, see Remark~\ref{EffectMonoidalRem} and~\cite{JacobsM12a}.

\begin{lem}
\label{ScalarEffectMonoidLem}
In a category $\cat{B}$ satisfying Assumption~\ref{CoprodAss} the
scalars $\Pred(1) = \Hom(1, 1+1)$ form an effect monoid.
\end{lem}

\begin{myproof}
Let $s, r, r' \colon 1 \rightarrow 1+1$ be three predicates on $1$,
with $r \orthogonal r'$ via bound $b\colon 1 \rightarrow
(1+1)+1$. We define the following two new bounds.
$$\begin{array}{rclcrcl}
c & = & \left(\vcenter{\xymatrix@R-1pc{
   1\ar[d]^{s} \\ 
   1+1\ar[d]^{b+\idmap} \\
   ((1+1)+1)+1\ar[d]^{[\idmap,\kappa_{2}]} \\
   (1+1)+1 }}\right)
& \qquad &
d & = & \left(\vcenter{\xymatrix@R-1pc{
   1\ar[d]^-{b} \\
   (1+1)+1\ar[d]^{(s+s)+\idmap} \\
   ((1+1)+(1+1))+1
   \ar[d]^-{[[\kappa_{1}+\idmap, \kappa_{2}+\idmap], \kappa_{2}]} \\
   (1+1)+1 }}\right)
\end{array}$$

\noindent The bound $c$ shows $r\cdot s \orthogonal r'\cdot s$ and
$(r\cdot s) \ovee (r'\cdot s) = (r\ovee r') \cdot s$. Similarly, $d$
proves $s\cdot r \orthogonal s\cdot r'$ and $(s\cdot r) \ovee (s\cdot
r') = s\cdot (r\ovee r')$. \QED

\auxproof{
$$\begin{array}{rcl}
\IV \after c
& = &
[\idmap,\kappa_{2}] \after [\idmap,\kappa_{2}] \after (b+\idmap)
   \after p \\
& = &
[[\idmap,\kappa_{2}],\kappa_{2}] \after (b+\idmap) \after p \\
& = &
[[\idmap,\kappa_{2}] \after b,\kappa_{2}] \after p \\
& = &
[q,\kappa_{2}] \after p \\
& = &
q\cdot p 
\\
\XI \after c
& = &
[[\kappa_{2}, \kappa_{1}],\kappa_{2}] \after [\idmap,\kappa_{2}] \after 
   (b+\idmap) \after p \\
& = &
[[[\kappa_{2}, \kappa_{1}],\kappa_{2}],\kappa_{2}] \after 
   (b+\idmap) \after p \\
& = &
[[[\kappa_{2}, \kappa_{1}],\kappa_{2}] \after b,\kappa_{2}] \after p \\
& = &
[q',\kappa_{2}] \after p \\
& = &
q'\cdot p \\
(\nabla+\idmap) \after c 
& = &
(\nabla+\idmap) \after [\idmap,\kappa_{2}] \after (b+\idmap)
  x \after p \\
& = &
[(\nabla+\idmap),\kappa_{2}] \after (b+\idmap) \after p \\
& = &
[(\nabla+\idmap) \after b,\kappa_{2}] \after p \\
& = &
[q\ovee q',\kappa_{2}] \after p \\
& = &
(q\ovee q')\cdot p.
\end{array}$$

$$\begin{array}{rcl}
[\idmap, \kappa_{2}] \after d
& = &
[\idmap, \kappa_{2}] \after 
   [[\kappa_{1}+\idmap, \kappa_{2}+\idmap], \kappa_{2}] \after
   ((p+p)+\idmap) \after b \\
& = &
{[[[\kappa_{1},\kappa_{2}], [\kappa_{2},\kappa_{2}]], \kappa_{2}]} \after
   ((p+p)+\idmap) \after b \\
& = &
{[[\idmap, \kappa_{2}\after \nabla] \after (p+p), \kappa_{2}]} \after b \\
& = &
{[[p, \kappa_{2}\after \nabla \after p], \kappa_{2}]} \after b \\
& = &
{[[p, \kappa_{2}], \kappa_{2}]} \after b \\
& = &
{[p, \kappa_{2}]} \after [\idmap, \kappa_{2}] \after b \\
& = &
{[p, \kappa_{2}]} \after q \\
& = &
p\cdot q \\
{[[\kappa_{2}, \kappa_{1}], \kappa_{2}]} \after d 
& = &
{[[\kappa_{2}, \kappa_{1}], \kappa_{2}]} \after 
   [[\kappa_{1}+\idmap, \kappa_{2}+\idmap], \kappa_{2}] \after
   ((p+p)+\idmap) \after b \\
& = &
{[[[\kappa_{2},\kappa_{2}], [\kappa_{1},\kappa_{2}]], \kappa_{2}]} \after
   ((p+p)+\idmap) \after b \\
& = &
{[[\kappa_{2} \after \nabla, \idmap] \after (p+p), \kappa_{2}]}
   \after b \\
& = &
{[[\kappa_{2} \after \nabla \after p, p], \kappa_{2}]} \after b \\
& = &
{[[\kappa_{2}, p], \kappa_{2}]} \after b \\
& = &
{[p, \kappa_{2}]} \after [[\kappa_{2},\kappa_{1}], \kappa_{2}] \after b \\
& = &
{[p, \kappa_{2}]} \after q' \\
& = &
p\cdot q' \\
(p\cdot q) \ovee (p\cdot q')
& = &
(\nabla+\idmap) \after d \\
& = &
(\nabla+\idmap) \after
   [[\kappa_{1}+\idmap, \kappa_{2}+\idmap], \kappa_{2}] \after
   ((p+p)+\idmap) \after b \\
& = &
{[[\idmap+\idmap, \idmap+\idmap], \kappa_{2}]} \after
   ((p+p)+\idmap) \after b \\
& = &
{[\nabla, \kappa_{2}]} \after ((p+p)+\idmap) \after b \\
& = &
{[\nabla \after (p+p), \kappa_{2}]} \after b \\
& = &
{[p \after \nabla, \kappa_{2}]} \after b \\
& = &
{[p, \kappa_{2}]} \after (\nabla+\idmap) \after b \\
& = &
{[p, \kappa_{2}]} \after (q\ovee q') \\
& = &
p\cdot (q\ovee q').
\end{array}$$
}
\end{myproof}

There is more structure: the effect monoid $\Pred(1)$ of scalars acts
on every effect algebra $\Pred(X)$ of predicates in an appropriate
manner, turning $\Pred(X)$ into an effect module over $\Pred(1)$.
Recall that we write $\EMod_{M}$ for the category of effect modules
over an effect monoid $M$. This subscript $M$ will be dropped when
it is clear from the context.

For a scalar $s\colon 1 \rightarrow 1+1$ and a predicate $p\colon X
\rightarrow 1+1$ we define the predicate $s\scalar p \colon X
\rightarrow 1+1$ obtained by scalar multiplication as the following
(Kleisli) composite.
\begin{equation}
\label{ScalarEqn}
\xymatrix@C-1.5pc{
s\scalar p = \big(X\ar[rr]^-{p} & & 1+1\ar[rrr]^-{[s,\kappa_{2}]} & & &
   1+1\big).
}
\end{equation}

\begin{prop}
\label{CoprodEModProp}
If $\cat{B}$ is a category satisfying Assumption~\ref{CoprodAss}, then
each collection of predicates $\Pred(X)$ is an effect module over
$\Pred(1)$, via the scalar multiplication $\scalar $
from~\eqref{ScalarEqn}. For $X = 1$ this scalar multiplication
$\scalar$ is the same as multiplication from
Lemma~\ref{ScalarMonoidLem}.

The substitution functor $f^{*} = (-) \after f\colon \Pred(X)
\rightarrow \Pred(Y)$ induced by a map $f\colon Y\rightarrow X$ in
$\cat{B}$ preserves this effect module structure, in the sense that
$f^{*}(s\scalar p) = s\scalar f^{*}(p)$. Thus, the functor
$\Pred\colon \cat{B} \rightarrow \op{\EA}$ from
Proposition~\ref{CoprodEAProp} restricts to a functor $\Pred \colon
\cat{B} \rightarrow \op{(\EMod_{\Pred(1)})}$.
\end{prop}

\begin{myproof}
The scalar multiplication properties $1\scalar p = p$ and $r\scalar
(s\scalar p) = (r\cdot s) \scalar p$ follow from straightforward
calculations. The fact that $\scalar$ is a bihomomorphism of effect
algebras is shown as in the proof of Lemma~\ref{ScalarEffectMonoidLem}.
The rest is easy. \QED

\auxproof{
$$\begin{array}{rcl}
1\scalar p
& = &
[\idmap,\kappa_{2}] \after (\kappa_{1}+\idmap) \after p \\
& = &
[\kappa_{1},\kappa_{2}] \after p \\
& = &
p \\
r\scalar (s\scalar p)
& = &
[r,\kappa_{2}] \after [s,\kappa_{2}] \after p \\
& = &
[[r,\kappa_{2}] \after s,\kappa_{2}] \after p \\
& = &
[r\cdot s, \kappa_{2}] \after p \\
& = &
(r\cdot s)\scalar p \\
0\scalar p
& = &
[\kappa_{2}, \kappa_{2}] \after p \\
& = &
\kappa_{2} \after \nabla \after p \\
& = &
\kappa_{2} \\
& = &
0 
\\
s\scalar 0
& = &
[s,\kappa_{2}] \after \kappa_{2} \\
& = &
\kappa_{2} \\
& = &
0.
\end{array}$$

We have to prove that $\scalar$ is a bihomomorphism. So assume
$s\orthogonal r$, for $s,r\in \Pred(1)$, say via bound $b\colon 1
\rightarrow (1+1)+1$. We take $c = [b,\kappa_{2}] \after p \colon X
\rightarrow (1+1)+1$ and claim that $c$ is a bound for $s\scalar p$
and $r\scalar p$. Indeed,
$$\begin{array}{rcl}
\IV \after c
& = &
[\idmap,\kappa_{2}] \after [b,\kappa_{2}] \after p \\
& = &
[[\idmap,\kappa_{2}] \after b,\kappa_{2}] \after p \\
& = &
[s,\kappa_{2}] \after p \\
& = &
s\scalar p 
\\
\XI \after c
& = &
[[\kappa_{2},\kappa_{1}],\kappa_{2}] \after [b,\kappa_{2}] \after p \\
& = &
[[\kappa_{2},\kappa_{1}],\kappa_{2}] \after b,\kappa_{2}] \after p \\
& = &
[r,\kappa_{2}] \after p \\
& = &
r\scalar p 
\\
(s\scalar p) \ovee (r\scalar p)
& = &
(\nabla+\idmap) \after c \\
& = &
(\nabla+\idmap) \after [b,\kappa_{2}] \after p \\
& = &
[(\nabla+\idmap) \after b,\kappa_{2}] \after p \\
& = &
[s\ovee r,\kappa_{2}] \after p \\
& = &
(s\ovee r)\scalar p.
\end{array}$$

\noindent And if $b\colon X \rightarrow (1+1)+1$ is a bound for $p,q$,
then $c = [[\kappa_{1}+\idmap, \kappa_{2}+\idmap], \kappa_{2}] \after
((s+s)+\idmap) \after b \colon X \rightarrow (1+1)+1$ is a bound for
$s\scalar p$ and $s\scalar q$ since:
$$\begin{array}{rcl}
\IV \after c
& = &
[\idmap,\kappa_{2}] \after [[\kappa_{1}+\idmap, \kappa_{2}+\idmap], \kappa_{2}] 
   \after ((s+s)+\idmap) \after b \\
& = &
[[[\kappa_{1},\kappa_{2}], [\kappa_{2}, \kappa_{2}]], \kappa_{2}]
   \after ((s+s)+\idmap) \after b \\
& = &
[[\idmap, \kappa_{2} \after \nabla], \kappa_{2}]
   \after ((s+s)+\idmap) \after b \\
& = &
[[s, \kappa_{2} \after \nabla \after s], \kappa_{2}] \after b \\
& = &
[[s, \kappa_{2}], \kappa_{2}] \after b \\
& = &
[s, \kappa_{2}] \after [\idmap, \kappa_{2}] \after b \\
& = &
[s, \kappa_{2}] \after p \\
& = &
s\scalar p 
\\
\XI \after c
& = &
[[\kappa_{2},\kappa_{1}],\kappa_{2}] \after 
   [[\kappa_{1}+\idmap, \kappa_{2}+\idmap], \kappa_{2}] 
   \after ((s+s)+\idmap) \after b \\
& = &
[[[\kappa_{2},\kappa_{2}], [\kappa_{1}, \kappa_{2}]], \kappa_{2}]
   \after ((s+s)+\idmap) \after b \\
& = &
[[\kappa_{2} \after \nabla, \idmap], \kappa_{2}]
   \after ((s+s)+\idmap) \after b \\
& = &
[[\kappa_{2} \after \nabla \after s, s], \kappa_{2}] \after b \\
& = &
[[\kappa_{2}, s], \kappa_{2}] \after b \\
& = &
[s, \kappa_{2}] \after [[\kappa_{2},\kappa_{1}], \kappa_{2}] \after b \\
& = &
[s, \kappa_{2}] \after q \\
& = &
s\scalar q \\
(s\scalar p) \ovee (s\scalar q)
& = &
(\nabla+\idmap) \after c \\
& = &
(\nabla+\idmap) \after [[\kappa_{1}+\idmap, \kappa_{2}+\idmap], \kappa_{2}] 
   \after ((s+s)+\idmap) \after b \\
& = &
[[\idmap+\idmap, \idmap+\idmap], \kappa_{2}] 
   \after ((s+s)+\idmap) \after b \\
& = &
[\nabla, \kappa_{2}] 
   \after ((s+s)+\idmap) \after b \\
& = &
[\nabla \after (s+s), \kappa_{2}] \after b \\
& = &
[s \after \nabla, \kappa_{2}] \after b \\
& = &
[s, \kappa_{2}] \after (\nabla+\idmap) \after b \\
& = &
[s, \kappa_{2}] \after (p\ovee q) \\
& = &
s\scalar (p\ovee q).
\end{array}$$

For $X = 1$ we have, for $s,p\colon 1 \rightarrow 1+1$, satisfying:
$$\begin{array}{rcccl}
s\scalar p
& = &
[s, \kappa_{2}] \after p
& = &
s \cdot p.
\end{array}$$

And for $f\colon Y \rightarrow X$, 
$$\begin{array}{rcccccccl}
f^{*}(s\scalar p)
& = &
(s\scalar p) \after f 
& = &
[s, \kappa_{2}] \after p \after f
& = &
[s, \kappa_{2}] \after f^{*}(p)
& = &
s\scalar f^{*}(p).
\end{array}$$
}
\end{myproof}

At this stage we like to emphasise that the relatively weak structure
in Assumption~\ref{CoprodAss} already gives us quite some logical
structure, namely the indexed category of effect modules from the
previous result. Our next aim is to investigate our assumptions more
concretely in the examples listed in the previous section.

\begin{exa}
\label{CoprodEx}
The category of $\Sets$ is an effectus. In fact, every extensive
category, with disjoint and universal coproducts~\cite{CarboniLW93},
satisfies Assumption~\ref{CoprodAss}. Every topos is extensive, see
\textit{e.g}~\cite[Vol.3, \S\S 5.9]{Borceux94}).  The proofs involve
some elementary diagrammatic reasoning and will be skipped here. A
predicate of the form $X \rightarrow 1+1 = 2$ in $\Sets$ is a
characteristic function of a subset of $X$. For such predicates
$P,Q\subseteq X$ the effect algebra sum $P\ovee Q$ is defined if
$P\cap Q = \emptyset$, and in that case equals the union: $P\ovee Q =
P \cup Q$. The scalars in $\Sets$ are the Booleans $2 = \{0,1\}$, and
scalar multiplication is trivial: $0\scalar P = \emptyset$ and
$1\scalar P = P$.

\auxproof{ 
We recall that coproducts are called \emph{disjoint} if the
coprojections $\kappa_{1}, \kappa_{2}$ are monic and their
intersection is empty.  The latter means that the following diagram is
a pullback.
$$\xymatrix@R-.5pc@C-.5pc{
0\ar[r]^-{!}\ar[d]_{!}\pullback & Y\ar@{ >->}[d]^{\kappa_2} \\
X\ar@{ >->}[r]_-{\kappa_{1}} & X+Y
}$$

\noindent These coproducts are called \emph{universal} if for each map
$f\colon Z \rightarrow X+Y$ the two pullbacks below exist, and yield
an isomorphism $[k_{1},k_{2}] \colon Z_{1}+Z_{2} \conglongrightarrow Z$.
$$\xymatrix@R-.5pc{
Z_{1}\ar[d]_{f_1}\ar@{ >->}[rr]^-{k_{1}}\pullback & & 
   Z\ar[d]^{f} & & Z_{2}\ar[d]^{f_2}\ar@{ >->}[ll]_-{k_2}\pullback[dl] \\
X\ar@{ >->}[rr]_-{\kappa_1} & & X+Y & & Y\ar@{ >->}[ll]^-{\kappa_2}
}$$

\noindent A category with such disjoint and universal coproducts is
called \emph{extensive}.

We remark that the initial object $0$ of an extensive category is
\emph{strict}: any map $X\rightarrow 0$ is an isomorphism. In case an
extensive category has products $\times$, these products distribute
over coproducts, see~\cite{CarboniLW93}.

Given a map $f\colon X \rightarrow 0$, the two squares
below are obviously pullbacks:
$$\xymatrix@R-.5pc{
X\ar@{=}[r]\ar[d]_{f}\pullback & X\ar[d]^{f} & 
   X\ar@{=}[l]\ar[d]^{f}\pullback[dl] \\
0\ar@{=}[r] & 0 & 0\ar@{=}[l]
}$$

\noindent Hence $[\idmap,\idmap] \colon X+X \rightarrow X$ is
an isomorphism. Since $[\idmap,\idmap] \after \kappa_{i} = \idmap$,
the coprojections $\kappa_{i}$ are also isomorphisms. But then:
$$\begin{array}{rcccccl}
\idmap[X] 
& = &
[\idmap[X], !_{X} \after f] \after \kappa_{1}
& = &
[\idmap[X], !_{X} \after f] \after \kappa_{1}
& = &
!_{X} \after f.
\end{array}$$

\noindent The equation $f \after \;!_{X} = \idmap[0]$ trivially holds.

\bigskip

In an extensive category Assumption~\ref{CoprodAss} holds: all
listed squares are pullbacks, and the two maps $[\idmap,\kappa_{2}],
[[\kappa_{2},\kappa_{1}],\kappa_{2}]$ are jointly monic.

\bigskip

We show that diagram on the left in~\eqref{CoprodAssSquares} is a
pullback. Assume we have maps $h\colon Z \rightarrow B+X$ and $k\colon
Z \rightarrow A+Y$ with $(\idmap+g) \after h = (f+\idmap) \after
k$. We form the pullbacks:
$$\xymatrix@R-.5pc{
Z_{1}\ar[d]_{h_1}\ar@{ >->}[rr]^-{p_{1}}\pullback & & 
   Z\ar[d]^{h} & & Z_{2}\ar[d]^{h_2}\ar@{ >->}[ll]_-{p_2}\pullback[dl] \\
B\ar@{ >->}[rr]_-{\kappa_1} & & B+X & & X\ar@{ >->}[ll]^-{\kappa_2}
}$$

$$\xymatrix@R-.5pc{
W_{1}\ar[d]_{k_1}\ar@{ >->}[rr]^-{q_{1}}\pullback & & 
   Z\ar[d]^{k} & & W_{2}\ar[d]^{k_2}\ar@{ >->}[ll]_-{q_2}\pullback[dl] \\
A\ar@{ >->}[rr]_-{\kappa_1} & & A+Y & & Y\ar@{ >->}[ll]^-{\kappa_2}
}$$

\noindent Our aim is to show that there are isomorphisms
$\varphi_{i} \colon Z_{i} \conglongrightarrow W_{i}$. We get
$Z_{1} \cong W_{1}$ via the next two inclusions of monos.
$$\xymatrix@R-.5pc{
Z_{1}\ar@{-->}[dr]^-{\varphi_1}\ar@/_2ex/[dddr]_{h_1}
   \ar@{ >->}@/^2ex/[drr]^{p_1} & & 
& &
W_{1}\ar@{-->}[dr]\ar@{ >->}@/^2ex/[drr]^{q_1}\ar[dd]_{k_1} & & \\
& W_{1}\ar[d]_{k_1}\ar@{ >->}[r]^-{q_1}\pullback & Z\ar[d]^{k}
& &
& Z_{1}\ar[d]_{h_1}\ar@{ >->}[r]^-{p_1}\pullback & Z\ar[d]^{h} \\
& A\ar[d]_{f}\ar@{ >->}[r]^-{\kappa_1}\pullback & A+Y\ar[d]^{f+\idmap}
& &
A\ar@/_2ex/[dr]_{f} & B\ar@{=}[d]\ar@{ >->}[r]^-{\kappa_1}\pullback & 
   B+X\ar[d]^{\idmap+g} \\
& B\ar@{ >->}[r]^-{\kappa_1} & B+Y
& &
& B\ar@{ >->}[r]^-{\kappa_1} & B+Y
}$$

\noindent The outer diagrams commute since:
$$\begin{array}{ccc}
\begin{array}{rcl}
\lefteqn{(f+\idmap) \after k \after p_{1}} \\
& = &
(\idmap+g) \after h \after p_{1} \\
& = &
(\idmap+g) \after \kappa_{1} \after h_{1} \\
& = &
\kappa_{1} \after h_{1}
\end{array}
& \qquad &
\begin{array}{rcl}
\lefteqn{(\idmap+g) \after h \after q_{1}} \\
& = &
(f+\idmap) \after k \after q_{1} \\
& = &
(f+\idmap) \after \kappa_{1} \after k_{1} \\
& = &
\kappa_{1} \after f \after k_{1}
\end{array}
\end{array}$$

\noindent In a similar way we get $\varphi_{2} \colon Z_{2}
\conglongrightarrow W_{2}$:
$$\xymatrix@R-.5pc{
Z_{2}\ar@{-->}[dr]^-{\varphi_2}\ar[dd]_{k_2}
   \ar@{ >->}@/^2ex/[drr]^{p_2} & & 
& &
W_{2}\ar@{-->}[dr]\ar@{ >->}@/^2ex/[drr]^{q_2}\ar@/_2ex/[dddr]_{k_2} & & \\
& W_{2}\ar[d]_{k_2}\ar@{ >->}[r]^-{q_2}\pullback & Z\ar[d]^{k}
& &
& Z_{2}\ar[d]_{h_2}\ar@{ >->}[r]^-{p_2}\pullback & Z\ar[d]^{h} \\
X\ar@/_2ex/[dr]_{g} & 
   Y\ar@{=}[d]\ar@{ >->}[r]^-{\kappa_2}\pullback & A+Y\ar[d]^{f+\idmap}
& &
& X\ar[d]_{g}\ar@{ >->}[r]^-{\kappa_2}\pullback & 
   B+X\ar[d]^{\idmap+g} \\
& Y\ar@{ >->}[r]^-{\kappa_2} & B+Y
& &
& Y\ar@{ >->}[r]^-{\kappa_2} & B+Y
}$$

\noindent The outer diagrams commute since:
$$\begin{array}{ccc}
\begin{array}{rcl}
\lefteqn{(f+\idmap) \after k \after p_{2}} \\
& = &
(\idmap+g) \after h \after p_{2} \\
& = &
(\idmap+g) \after \kappa_{2} \after h_{2} \\
& = &
\kappa_{2} \after g \after h_{2}
\end{array}
& \qquad &
\begin{array}{rcl}
\lefteqn{(\idmap+g) \after h \after q_{2}} \\
& = &
(f+\idmap) \after k \after q_{2} \\
& = &
(f+\idmap) \after \kappa_{2} \after k_{2} \\
& = &
\kappa_{2} \after k_{2}
\end{array}
\end{array}$$

\noindent We notice that the following square commutes.
$$\xymatrix{
W_{1}+Z_{2}\ar[r]^-{\idmap+\varphi_{2}}_-{\cong}
   \ar[d]_{\varphi_{1}^{-1}+\idmap}^{\cong} &
   W_{1}+W_{2}\ar[d]^{[q_{1},q_{2}]}_{\cong} \\
Z_{1}+Z_{2}\ar[r]_-{[p_{1},p_{2}]}^-{\cong} & Z
}$$

We now take as mediating map:
$$\xymatrix@C+1pc{
\ell = \Big(Z\ar[r]^-{[p_{1},p_{2}]^{-1}}_-{\cong} & 
   Z_{1}+Z_{2}\ar[r]^-{\varphi_{1}+\idmap}_-{\cong} &
   W_{1}+Z_{2}\ar[r]^-{k_{1}+h_{2}} & A+X\Big)
}$$

\noindent Then:
$$\begin{array}{rcl}
(f+\idmap) \after \ell
& = &
(f+\idmap) \after (k_{1}+h_{2}) \after (\varphi_{1}+\idmap) \after
   [p_{1},p_{2}]^{-1} \\
& = &
((f\after k_{1} \after \varphi_{1})+h_{2}) \after [p_{1},p_{2}]^{-1} \\
& = &
(h_{1}+h_{2}) \after [p_{1},p_{2}]^{-1} \\
& = &
[\kappa_{1} \after h_{1}, \kappa_{2} \after h_{2}] \after
   [p_{1},p_{2}]^{-1} \\
& = &
[h \after p_{1}, h \after p_{2}] \after [p_{1},p_{2}]^{-1} \\
& = &
h \after [p_{1}, p_{2}] \after [p_{1},p_{2}]^{-1} \\
& = &
h \\
(\idmap+g) \after \ell
& = &
(\idmap+g) \after (k_{1}+h_{2}) \after (\varphi_{1}+\idmap) \after
   [p_{1},p_{2}]^{-1} \\
& = &
(\idmap+g) \after (k_{1}+h_{2}) \after (\idmap+\varphi_{2}^{-1}) \after
   [q_{1},q_{2}]^{-1} \\
& & \qquad \mbox{by the above square of isomorphisms} \\
& = &
(k_{1}+(g \after h_{2} \after \varphi_{2}^{-1})) \after
   [q_{1},q_{2}]^{-1} \\
& = &
(k_{1}+k_{2}) \after [q_{1},q_{2}]^{-1} \\
& = &
[\kappa_{1} \after k_{1}, \kappa_{2} \after k_{2}] \after
   [q_{1},q_{2}]^{-1} \\
& = &
[k \after q_{1}, k \after q_{2}] \after [q_{1},q_{2}]^{-1} \\
& = &
k \after [q_{1}, q_{2}] \after [q_{1},q_{2}]^{-1} \\
& = &
k.
\end{array}$$

Finally, assume $\ell'\colon Z\rightarrow A+X$ also satisfies
$(f+\idmap) \after \ell' = h$ and $(\idmap+g) \after \ell' = k$.
We form the pullbacks:
$$\xymatrix@R-.5pc{
Z'_{1}\ar[d]_{\ell_1}\ar@{ >->}[rr]^-{p'_{1}}\pullback & & 
   Z\ar[d]^{\ell'} & & 
   Z'_{2}\ar[d]^{\ell_2}\ar@{ >->}[ll]_-{p'_2}\pullback[dl] \\
A\ar@{ >->}[rr]_-{\kappa_1} & & A+X & & X\ar@{ >->}[ll]^-{\kappa_2}
}$$

\noindent Then we obtain maps $\psi_{i}\colon Z'_{i} \rightarrow
Z_{i}$ in:
$$\xymatrix@R-.5pc{
Z'_{1}\ar@{ >->}@/^2ex/[drrr]^{p'_1}\ar[d]_{\ell_1}\ar@{-->}[dr]^-{\psi_1} 
   & & & & & & 
   Z'_{2}\ar@{ >->}@/_2ex/[dlll]_{p'_2}\ar@{-->}[dl]_{\psi_2}
      \ar@/^2ex/[ddl]^{\ell_2} \\
A\ar@/_2ex/[dr]_{f} & Z_{1}\ar[d]_{h_1}\ar@{ >->}[rr]^-{p_{1}}\pullback & & 
   Z\ar[d]^{h} & & Z_{2}\ar[d]^{h_2}\ar@{ >->}[ll]_-{p_2}\pullback[dl] \\
& B\ar@{ >->}[rr]_-{\kappa_1} & & B+X & & X\ar@{ >->}[ll]^-{\kappa_2}
}$$

\noindent The two (outer) diagrams commute:
$$\begin{array}{ccc}
\begin{array}{rcl}
h \after p'_{1}
& = &
(f+\idmap) \after \ell' \after p'_{1} \\
& = &
(f+\idmap) \after \kappa_{1} \after \ell_{1} \\
& = &
\kappa_{1} \after f \after \ell_{1}
\end{array}
& \qquad &
\begin{array}{rcl}
h \after p'_{2}
& = &
(f+\idmap) \after \ell' \after p'_{2} \\
& = &
(f+\idmap) \after \kappa_{2} \after \ell_{2} \\
& = &
\kappa_{2} \after \ell_{2}.
\end{array}
\end{array}$$

\noindent We obtain maps $\psi_i$ satisfying:
$$\begin{array}{rclcrcl}
[p_{1}, p_{2}] \after (\psi_{1}+\psi_{2})
& = &
[p'_{1}, p'_{2}]
& \quad\mbox{and thus}\quad &
(\psi_{1}+\psi_{2}) \after [p'_{1}, p'_{2}]^{-1}
& = &
[p_{1}, p_{2}]^{-1}.
\end{array}$$

\noindent Now we obtain $\ell'=\ell$ via:
$$\begin{array}{rcl}
\ell
& = &
(k_{1}+h_{2}) \after (\varphi_{1}+\idmap) \after [p_{1},p_{2}]^{-1} \\
& = &
(k_{1}+h_{2}) \after (\varphi_{1}+\idmap) \after
   (\psi_{1}+\psi_{2}) \after [p'_{1}, p'_{2}]^{-1} \\
& = &
((k_{1} \after \varphi_{1} \after \psi_{1}) + (h_{2} \after \psi_{2}))
    \after [p'_{1}, p'_{2}]^{-1} \\
& \smash{\stackrel{(*)}{=}} &
(\ell_{1}+\ell_{2}) \after [p'_{1}, p'_{2}]^{-1} \\
& = &
[\kappa_{1} \after \ell_{1}, \kappa_{2} \after \ell_{2}] \after 
   [p'_{1}, p'_{2}]^{-1} \\
& = &
[\ell' \after p'_{1}, \ell' \after p'_{2}] \after [p'_{1}, p'_{2}]^{-1} \\
& = &
\ell' \after [p'_{1}, p'_{2}] \after [p'_{1}, p'_{2}]^{-1} \\
& = &
\ell'
\end{array}$$

\noindent The marked equation $\smash{\stackrel{(*)}{=}}$ holds since:
$$\begin{array}{rcl}
\kappa_{1} \after k_{1} \after \varphi_{1} \after \psi_{1}
& = &
k \after q_{1} \after \varphi_{1} \after \psi_{1} \\
& = &
k \after p_{1} \after \psi_{1} \\
& = &
k \after p'_{1} \\
& = &
(\idmap+g) \after \ell' \after p'_{1} \\
& = &
(\idmap+g) \after \kappa_{1} \after \ell_{1} \\
& = &
\kappa_{1} \after \ell_{1}.
\end{array}$$

\bigskip

We show that diagram on the left in~\eqref{CoprodAssSquares} Assume we
have maps $h\colon Z \rightarrow A+X$ and $k\colon Z\rightarrow A$
with $(\idmap+f) \after h = \kappa_{1} \after k$. We have to prove
that $\kappa_{1} \after k = h$. Form the pullbacks:
$$\xymatrix@R-.5pc{
Z_{1}\ar[d]_{h_1}\ar@{ >->}[rr]^-{p_{1}}\pullback & & 
   Z\ar[d]^{h} & & Z_{2}\ar[d]^{h_2}\ar@{ >->}[ll]_-{p_2}\pullback[dl] \\
A\ar@{ >->}[rr]_-{\kappa_1} & & A+X & & X\ar@{ >->}[ll]^-{\kappa_2}
}$$

\noindent in which $[p_{1}, p_{2}] \colon Z_{1} + Z_{2}
\conglongrightarrow Z$. Next we consider the situation:
$$\xymatrix@R-.5pc{
Z_{2}\ar@{-->}[dr]\ar[r]^-{p_2}\ar[d]_{h_2} & Z\ar@/^2ex/[dr]^-{k} \\
X\ar@/_2ex/[dr]_{f} & 0\ar[r]\ar[d]\pullback & A\ar@{ >->}[d]^{\kappa_1} \\
& Y\ar@{ >->}[r]_-{\kappa_2} & A+Y
}$$

\noindent The outer diagram commutes since:
$$\begin{array}{rcl}
\kappa_{2} \after f \after h_{2}
& = &
(\idmap+f) \after \kappa_{2} \after h_{2} \\
& = &
(\idmap+f) \after h \after p_{2} \\
& = &
\kappa_{1} \after k \after p_{2}.
\end{array}$$

\noindent Since $0$ is strict we get $Z_{2} \cong 0$ and thus $Z \cong
Z_{1}+Z_{2} \cong Z_{1}$. Thus $p_{1}$ is an isomorphism.  This gives
the required map $h_{1} \after p_{1}^{-1} \colon Z \rightarrow A$. We
have and equality of maps $Z \rightarrow A \rightarrow A+Y$:
$$\begin{array}{rcl}
\kappa_{1} \after h_{1} \after p_{1}^{-1}
& = &
(\idmap+f) \after \kappa_{1} \after h_{1} \after p_{1}^{-1} \\
& = &
(\idmap+f) \after h \\
& = &
\kappa_{1} \after k.
\end{array}$$

\noindent Since $\kappa_{1}$ is monic, this gives $h_{1} \after
p_{1}^{-1} = k$, and thus, as required,
$$\begin{array}{rcccccl}
\kappa_{1} \after k
& = &
\kappa_{1} \after h_{1} \after p_{1}^{-1}
& = &
h \after p_{1} \after p_{1}^{-1} 
& = &
h.
\end{array}$$

\bigskip

For the jointly monic requirement in Assumption~\ref{CoprodAss} we
first give a set-theoretic argument. Next we give a categorical
proof for $n=2$.

Assume functions $f,g\colon Y \rightarrow (n+1)\cdot X$ satisfy
$\nabla_{i} \after f = \nabla_{i} \after g$, for all $1 \leq i \leq
n$. Then, $f(y)$ is in the $i$-coproduct component iff $g(y)$ is in
the $i$-th component, and their values are the same in that case. But
this also means that if $f(y)$ is in the $n+1$-component, then it is
not in the $i$-th component, so $g(y)$ is not in the $i$-th component,
and so $g(y)$ must also be in the $n+1$-component. Moreover, in that
case the values $f(y) = \kappa_{n+1}x$ and $g(y) = \kappa_{n+1}x'$
must be equal, which simply follows from $\nabla_{1} \after f =
\nabla_{1} \after g$.

In a second step we show that the maps $[\idmap,\kappa_{2}],
[[\kappa_{2},\kappa_{1}],\kappa_{2}]$ are jointly monic. Let
$f,g\colon Y \rightarrow (X+X)+X$ be maps with $[\idmap,\kappa_{2}]
\after f = [\idmap,\kappa_{2}] \after g$ and
$[[\kappa_{2},\kappa_{1}], \kappa_{2}] \after f =
       [[\kappa_{2},\kappa_{1}], \kappa_{2}] \after g$. These maps
       $f,g$ can be split in three parts, corresponding to whether
       their result is in the first, second, or third component of
       $(X+X)+X$. In the first and second cases these splittings of
       $f,g$ are equal, since we have pullbacks:
$$\xymatrix@R-.5pc{
X\ar@{=}[d]\ar@{ >->}[r]^-{\kappa_{1}\after\kappa_{1}}\pullback & 
   (X+X)+X\ar[d]^{[\idmap,\kappa_{2}]} 
& &
X\ar@{=}[d]\ar@{ >->}[r]^-{\kappa_{1}\after\kappa_{2}}\pullback & 
   (X+X)+X\ar[d]^{[[\kappa_{2},\kappa_{1}], \kappa_{2}]} \\
X\ar@{ >->}[r]_-{\kappa_{1}} & X+X
& &
X\ar@{ >->}[r]_-{\kappa_{2}} & X+X
}$$

\noindent These diagrams are pullbacks since we can reorganise
them as a pullback on the left in~\eqref{CoprodAssLemSquares}, in:
$$\xymatrix@R-.5pc@C+1pc{
X\ar@{=}[d]\ar@{ >->}[r]_-{\kappa_{1}\after\kappa_{1}}
   \ar@/^4ex/[rr]^{\kappa_1}\pullback & 
   (X+X)+X\ar[d]^{[\idmap,\kappa_{2}]}\ar[r]^-{\alpha^{-1}}_-{\cong} &
   X + (X+X)\ar@/^2ex/[dl]^{\idmap+\nabla} \\
X\ar@{ >->}[r]_-{\kappa_{1}} & X+X
}$$

$$\xymatrix@R-.5pc@C+1pc{
X\ar@{=}[d]\ar@{ >->}[r]_-{\kappa_{1}\after\kappa_{2}}
   \ar@/^4ex/[rrr]^{\kappa_1}\pullback & 
   (X+X)+X\ar[d]^{[[\kappa_{2},\kappa_{1}], \kappa_{2}]}
      \ar[rr]^-{[[\kappa_{2}\after\kappa_{1},\kappa_{1}],
   \kappa_{2}\after\kappa_{2}]}_-{\cong} & &
   X + (X+X)\ar@/^2ex/[dll]^{\idmap+\nabla} \\
X\ar@{ >->}[r]_-{\kappa_1} & X+X
}$$

\noindent The difficulty lies in the third component. Here
we use that we have pullbacks.
$$\xymatrix@R-.5pc{
X+X\ar[d]_{\nabla}\ar@{ >->}[r]^-{\kappa_{2}+\idmap}\pullback & 
   (X+X)+X\ar[d]^{[\idmap,\kappa_{2}]} 
& &
X+X\ar[d]_{\nabla}\ar@{ >->}[r]^-{\kappa_{1}+\idmap}\pullback & 
   (X+X)+X\ar[d]^{[[\kappa_{2},\kappa_{1}], \kappa_{2}]} \\
X\ar@{ >->}[r]_-{\kappa_{2}} & X+X
& &
X\ar@{ >->}[r]_-{\kappa_{2}} & X+X
}$$

\noindent They are pullbacks, again by the pullback on the left
in~\eqref{CoprodAssLemSquares}:
$$\xymatrix@R-.5pc@C+1pc{
X+X\ar[d]_{\nabla}\ar@{ >->}[r]_-{\kappa_{2}+\idmap}\ar@/^4ex/[rr]^{\kappa_2}
   \pullback & 
   (X+X)+X\ar[d]^{[\idmap,\kappa_{2}]}\ar[r]^-{\alpha^{-1}}_-{\cong} &
   X + (X+X)\ar@/^2ex/[dl]^{\idmap+\nabla} \\
X\ar@{ >->}[r]_-{\kappa_{2}} & X+X
}$$

$$\xymatrix@R-.5pc@C+1pc{
X+X\ar[d]_{\nabla}\ar@{ >->}[r]_-{\kappa_{1}+\idmap}
   \ar@/^4ex/[rrr]^{\kappa_2}\pullback & 
   (X+X)+X\ar[d]^{[[\kappa_{2},\kappa_{1}], \kappa_{2}]}
      \ar[rr]^-{[[\kappa_{2}\after\kappa_{1},\kappa_{1}],
   \kappa_{2}\after\kappa_{2}]}_-{\cong} & &
   X + (X+X)\ar@/^2ex/[dll]^{\idmap+\nabla} \\
X\ar@{ >->}[r]_-{\kappa_2} & X+X
}$$

\noindent Thus the objects $U_{1}, U_{2}$ form pullbacks for both
maps $f,g$, in:
$$\xymatrix@C+1pc{
& V\ar@{ >->}[dl]_{p_1}\ar@{ >->}[dr]^{p_2}
   \ar[d]|{\ell = k_{1} \after p_{1} = k_{2} \after p_{2}}  & \\
U_{1}\ar[d]_{h_1}\ar@{ >->}[r]^-{k_1} & 
   Y\ar@<-.5ex>[d]_{f}\ar@<+.5ex>[d]^{g} & 
   U_{2}\ar@{ >->}[l]_-{k_2}\ar[d]^{h_2} \\
X+X\ar@{ >->}[r]_-{\kappa_{1}+\idmap} & (X+X)+X &
   X+X\ar@{ >->}[l]^-{\kappa_{2}+\idmap}
}$$

\noindent The object $V$ on top is also obtained as pullback.

In a next step we use the pullback square below, obtained again
via the pullback on the left in~\eqref{CoprodAssLemSquares}:
$$\xymatrix@R-.5pc@C+1pc{
X\ar@{ >->}[d]_{\kappa_2}\ar@{ >->}[r]^-{\kappa_2}\pullback &
   X+X\ar@{ >->}[d]^{\kappa_{2}+\idmap}\ar[dr]^{\kappa_2} \\
X+X\ar@{ >->}[r]^-{\kappa_{1}+\idmap}\ar@/_4ex/[rr]_{\idmap+\kappa_2} & 
   (X+X)+X\ar[r]^-{\alpha^{-1}}_-{\cong} &
   X + (X+X)
}$$

\noindent It yields:
$$\xymatrix{
V\ar@/_2ex/[ddr]_{h_{1}\after p_{1}}\ar@/^2ex/[drr]^{h_{2}\after p_{2}}
   \ar@{-->}[dr]^{a} & & \\
& X\ar@{ >->}[d]_{\kappa_2}\ar@{ >->}[r]^-{\kappa_2}\pullback &
   X+X\ar@{ >->}[d]^{\kappa_{2}+\idmap} \\
& X+X\ar@{ >->}[r]^-{\kappa_{1}+\idmap} & (X+X)+X
}$$

\noindent The outer rectangle commutes since:
$$\begin{array}{rcl}
(\kappa_{1}+\idmap) \after h_{1} \after p_{1}
& = &
f \after k_{1} \after p_{1} \\
& = &
f \after k_{2} \after p_{2} \\
& = &
(\kappa_{2}+\idmap) \after h_{2} \after p_{2}
\end{array}$$

\noindent This object $V$ gives the third component of both $f,g$,
since we have a pullback diagram:
$$\xymatrix{
V\ar[d]_{a}\ar@{ >->}[r]^-{\ell} & Y\ar[d]^{f} \\
X\ar@{ >->}[r]_-{\kappa_2} & (X+X)+X
}$$

\noindent It commutes since:
$$\begin{array}{rcl}
f \after \ell
& = &
f \after k_{1} \after p_{1} \\
& = &
(\kappa_{1}+\idmap) \after h_{1} \after p_{1} \\
& = &
(\kappa_{1}+\idmap) \after \kappa_{2} \after a \\
& = &
\kappa_{2} \after a.
\end{array}$$

\noindent It is a pullback: assume we have $u\colon A \rightarrow X$
and $v\colon A \rightarrow Y$ with $\kappa_{2} \after u = f \after v$.
The maps $\kappa_{2} \after u \colon A \rightarrow X+X$ and $v$ then
satisfy: $(\kappa_{1}+\idmap) \after \kappa_{2} \after u = \kappa_{2}
\after u = f \after v$. Similarly, $(\kappa_{2}+\idmap) \after
\kappa_{2} \after u = \kappa_{2} \after u = f \after v$. Thus we get
maps $u_{i} \colon A \rightarrow U_{i}$ with $h_{i} \after u_{i} =
\kappa_{2} \after u$ and $k_{i} \after u_{i} = v$. The latter yields
a map $w\colon A \rightarrow V$ with $p_{i} \after w = u_{i}$. 
This $w$ is the (unique) mediating map that we seek:
$$\begin{array}{rcl}
\ell \after w
& = &
k_{1} \after p_{1} \after w \\
& = &
k_{1} \after u_{i} \\
& = &
v \\
\kappa_{2} \after a \after w 
& = &
h_{1} \after p_{1} \after w \\
& = &
h_{1} \after u_{1} \\
& = &
\kappa_{2} \after u.
\end{array}$$
}

The Kleisli category $\Kl(\Dst)$ of the distribution monad is als an
effectus. We leave the verification of the pullback
properties~\eqref{CoprodAssSquares} to the interested reader and show
that the maps $\IV, \XI \colon (1+1) + 1 \rightarrow 1+1$
in~\eqref{CoprodOneAssJoint} are jointly monic in $\Kl(\Dst)$. So let
$\varphi, \psi\in\Dst(3)$ be distributions with equalities
$\Dst(\IV)(\varphi) = \Dst(\IV)(\psi)$ and $\Dst(\XI)(\varphi) =
\Dst(\XI)(\psi)$ in $\Dst(2)$ --- see~\eqref{DstMapEqn} for the
functor $\Dst$ on maps. In order to disambiguate matters we write $3 =
\{a, b, c\}$ and $2 = \{u, v\}$. The functor applications yield convex
sums:
$$\begin{array}{rclcrcl}
\Dst(\IV)(\varphi)
& = &
\varphi(a)\ket{u} + (\varphi(b)+\varphi(c))\ket{v}
& \quad\mbox{and}\quad &
\Dst(\XI)(\varphi)
& = &
\varphi(b)\ket{u} + (\varphi(a) + \varphi(c))\ket{v}.
\end{array}$$

\noindent Similarly for $\psi$. Hence the equations
$\Dst(\IV)(\varphi) = \Dst(\IV)(\psi)$ and $\Dst(\XI)(\varphi) =
\Dst(\XI)(\psi)$ immediately give us $\varphi(a) = \psi(a)$ and
$\varphi(b) = \psi(b)$. We still need to prove $\varphi(c) =
\psi(c)$. But $\varphi$ and $\psi$ are distributions, so their
probabilities add up to one:
$$\begin{array}{rcccccl}
\varphi(c)
& = &
1 - \varphi(a) - \varphi(b)
& = &
1 - \psi(a) - \psi(b)
& = &
\psi(c).
\end{array}$$

\auxproof{
$$\begin{array}{rcccccl}
\varphi(\kappa_{i}x)
& = &
\Dst([\rhd_{i},\kappa_{2}])(\varphi)(\kappa_{1}x)
& = &
\Dst([\rhd_{i},\kappa_{2}])(\psi)(\kappa_{1}x)
& = &
\psi(\kappa_{i}x).
\end{array}$$

\noindent In order to be able to conclude $\varphi=\psi$, we only need
the equation $\varphi(\kappa_{n+1}0) = \psi(\kappa_{n+1}0)$, for the
sole element $0$ of the singleton set $1 = \{0\}$. But this equations
follows from the above ones, since $\varphi$ and $\psi$ are
distributions, with probabilities adding up to 1:
$$\begin{array}{rcccccl}
\varphi(\kappa_{n+1}0)
& = &
1 - \sum_{u\in n\cdot X}\varphi(u)
& = &
1 - \sum_{u\in n\cdot X}\psi(u)
& = &
\psi(\kappa_{n+1}0).
\end{array}$$

We check that the first diagram in~\eqref{CoprodAssSquares} is a
pullback. We reason pointwise and assume $\varphi\in\Dst(B+X)$ and
$\psi\in\Dst(A+Y)$ satisfy $(\idmap+g)^{\#}(\varphi) =
(f+\idmap)^{\#}(\psi)$. Thus we have equal functions $B+Y\rightarrow
[0,1]$. This means for each $b\in B$ and $y\in Y$,
$$\begin{array}{rclcrcl}
\varphi(\kappa_{1}b)
& = &
\sum_{a\in A}\psi(\kappa_{1}a)\cdot f(a)(b)
& \qquad &
\sum_{x\in X}\varphi(\kappa_{2}x)\cdot g(x)(y)
& = &
\psi(\kappa_{2}y)
\end{array}$$

Just checking, $f^{\#} \colon \Dst(A) \rightarrow \Dst(B)$
is given by:
$$\begin{array}{rcl}
f^{\#}(\chi)(b)
& = &
\mu\big(\Dst(f)(\chi)\big)(b) \\
& = &
\sum_{\varphi}\Dst(f)(\chi)(\varphi)\cdot \varphi(b) \\
& = &
\sum_{\varphi}\big(\sum_{a, f(a)=\varphi}\chi(a)\big)\cdot \varphi(b) \\
& = &
\sum_{a}\chi(a)\cdot f(a)(b)
\end{array}$$

Hence we define $\chi\in\Dst(A+X)$ as $\chi(\kappa_{1}a) =
\psi(\kappa_{1}a)$ and $\chi(\kappa_{2}x) = \varphi(\kappa_{2}x)$. This
yields a probability distribution:
$$\begin{array}{rcl}
\sum_{u\in A+X}\chi(u)
& = &
\sum_{a} \chi(\kappa_{1}a) \;+\; \sum_{x}\chi(\kappa_{2}x) \\
& = &
\sum_{a} \psi(\kappa_{1}a) \;+\; \sum_{x}\varphi(\kappa_{2}x) \\
& = &
\sum_{a} \psi(\kappa_{1}a)\cdot (\sum_{b}f(a)(b)) \;+\; 
   \sum_{x}\varphi(\kappa_{2}x) \\
& = &
\sum_{a,b} \psi(\kappa_{1}a)\cdot f(a)(b) \;+\; 
   \sum_{x}\varphi(\kappa_{2}x) \\
& = &
\sum_{b}\varphi(\kappa_{1}b) \;+\; 
   \sum_{x}\varphi(\kappa_{2}x) \\
& = &
1.
\end{array}$$

\noindent Also, this $\chi$ is the appropriate mediating
distribution:
$$\begin{array}{rcl}
(f+\idmap)^{\#}(\chi)(\kappa_{1}b)
& = &
\sum_{a\in A}\chi(\kappa_{1}a)\cdot f(a)(b) \\
& = &
\sum_{a\in A}\psi(\kappa_{1}a)\cdot f(a)(b) 
\hspace*{\arraycolsep} = \hspace*{\arraycolsep}
\varphi(\kappa_{1}b) \\
(f+\idmap)^{\#}(\chi)(\kappa_{2}x)
& = &
\chi(\kappa_{2}x)
\hspace*{\arraycolsep} = \hspace*{\arraycolsep}
\varphi(\kappa_{2}x) \\
(\idmap+g)^{\#}(\chi)(\kappa_{1}a)
& = &
\chi(\kappa_{1}a) 
\hspace*{\arraycolsep} = \hspace*{\arraycolsep}
\psi(\kappa_{1}a) \\
(\idmap+g)^{\#}(\chi)(\kappa_{2}y)
& = &
\sum_{x} \chi(\kappa_{2}x)\cdot g(x)(y) \\
& = &
\sum_{x} \varphi(\kappa_{2}x)\cdot g(x)(y)
\hspace*{\arraycolsep} = \hspace*{\arraycolsep}
\psi(\kappa_{2}y).
\end{array}$$

\noindent Thus, $(f+\idmap)^{\#}(\chi) = \varphi$ and
$(\idmap+g)^{\#}(\chi) = \psi$. Uniqueness is left to the reader.

Next, for the diagram on the right in~\eqref{CoprodAssSquares}, assume
we have $\varphi\in\Dst(A+X)$ and $\psi\in\Dst(A)$ satisfying
$(\idmap+f)^{\#}(\varphi) = \Dst(\kappa_{1})(\psi)$. The aim is to
prove $\varphi = \Dst(\kappa_{1})(\psi)$. The assumptions mean for
each $y\in Y$ and $a\in B$,
$$\begin{array}{rclcrcl}
\varphi(\kappa_{1}a)
& = &
\psi(a)
& \qquad &
\sum_{x}\varphi(\kappa_{2}x)\cdot f(x)(y)
& = &
0.
\end{array}$$

\noindent The latter means $\varphi(\kappa_{2}x) = 0$ or
$f(x)(y)=0$. For each $x\in X$ we the support $\supp(f(x)) \subseteq
Y$ of $f(x)\in\Dst(Y)$ is non-empty. Hence we can choose an element
$y_{x}$ with $f(x)(y_{x}) \neq 0$. This means $\varphi(\kappa_{2}x) =
0$. Thus we have:
$$\begin{array}{rcl}
\Dst(\kappa_{1})(\psi)(\kappa_{1}a)
& = &
\psi(a) \\
& = &
\varphi(\kappa_{1}a) \\
\Dst(\kappa_{1})(\psi)(\kappa_{2}x)
& = &
0 \\
& = &
\varphi(\kappa_{2}x).
\end{array}$$
}

\noindent The sum $\ovee$ of fuzzy predicates $p,q\colon X \rightarrow 1+1$ is
defined if $p(x) + q(x) \leq 1$ for all $x\in X$. In that case one
takes the pointwise sum $(p\ovee q)(x) = p(x) + q(x)$. Multiplication
with a scalar $r\in [0,1]$ is also defined pointwise: $(r\scalar p)(x)
= r\cdot p(x)$. Substitution $f^{*}(p)$, for $f\colon Y \rightarrow X$
in $\Kl(\Dst)$, is given by Kleisli composition $\klafter$, and yields:
\begin{equation}
\label{KlDSubstEqn}
\begin{array}{rcccl}
f^{*}(p)(y)
& = &
(p \klafter f)(y)
& = &
\displaystyle\sum_{x} p(x) \cdot f(y)(x).
\end{array}
\end{equation}
Also the Kleisli category $\Kl(\Giry)$ of the Giry monad $\Giry$ is an
effectus. The verification of Assumption~\ref{CoprodAss} involve some
properties of integrals --- like splitting of an integral over a
coproducts into a sum of integrals, see~\cite[Lemma~12]{Jacobs13a} ---
that are skipped here.  Since $\Giry(2) \cong \Dst(2) \cong [0,1]$ we
obtain that predicates $X \rightarrow 1+1 = 2$ in $\Kl(\Giry)$
correspond to measurable functions $X \rightarrow [0,1]$. In
particular, the scalars are the probabilities $[0,1]$. These
measurable fuzzy predicates obviously carry the structure of an effect
module over $[0,1]$, via pointwise operations. For a Kleisli map $f
\colon Y \rightarrow \Giry(X)$, substitution $f^{*}(p)$ yields a
predicate $Y\rightarrow [0,1]$ on $Y$ defined via Kleisli composition
(written as $\klafter$) as integration (see
Section~\ref{ContProbSubsec}):
\begin{equation}
\label{KlGSubstEqn}
\begin{array}{rcccccl}
f^{*}(p)(y)
& = &
(p \klafter f)(y)
& = &
p_{*}(f(y))
& = &
\displaystyle\int p\intd f(y).
\end{array}
\end{equation}

\noindent This forms a `continuous' version of the `discrete'
formula~\eqref{KlDSubstEqn}.

\auxproof{
We start with the diagram on the left in~\eqref{CoprodAssSquares}, as in:
$$\xymatrix{
A+ X\ar[r]^-{\idmap+ f}\ar[d]_{g+\idmap} & 
  A+ Y\ar[d]^{g+\idmap}
& &
A\ar@{=}[r]\ar[d]_{\kappa_1} & A\ar[d]^{\kappa_1}  \\
B+ X\ar[r]_-{\idmap+ f} & B+ Y
& &
A+ X\ar[r]_-{\idmap+ f} & A+ Y
}$$

\noindent So let $f\colon X \rightarrow \Giry(Y)$ and $g\colon A
\rightarrow \Giry(B)$ be Kleisli maps, and $\phi\in\Giry(B+X), \psi\in
\Giry(A+Y)$ be probability measures satisfying $(\idmap+f)_{*}(\phi) =
(g+\idmap)_{*}(\psi)$.  This means for measurable subsets
$K\in\Sigma_{B+Y}$,
$$\begin{array}{rcl}
\int (\idmap+f)(-)(K)\intd\phi
& = &
\int (g+\idmap)(-)(K)\intd\psi
\end{array}$$

\noindent The function $\idmap+f \colon B+X \rightarrow \Giry(B+Y)$
is given by:
$$\begin{array}{rcl}
(\idmap+f)(\kappa_{1}b)(K)
& = &
\Giry(\kappa_{1})(\eta(b))(K) \\
& = &
\eta(b)(\kappa_{1}^{-1}(K)) \\
& = &
\left\{\begin{array}{ll}
\indic{M}(b) \quad & \mbox{if }K=\kappa_{1}M, \mbox{ for } M\in\Sigma_{B} \\
0 & \mbox{if }K=\kappa_{2}N \mbox{ for } N\in\Sigma_{Y}
\end{array}\right. \\
(\idmap+f)(\kappa_{2}x)(K)
& = &
\Giry(\kappa_{1})(f(x))(K) \\
& = &
f(x)(\kappa_{1}^{-1}(K)) \\
& = &
\left\{\begin{array}{ll}
f(x)(M) \quad & \mbox{if }K=\kappa_{1}M, \mbox{ for } M\in\Sigma_{B} \\
0 & \mbox{if }K=\kappa_{2}N \mbox{ for } N\in\Sigma_{Y}
\end{array}\right.
\end{array}$$

\noindent Similarly, $g+\idmap \colon A+Y \rightarrow \Giry(B+Y)$ is
given by:
$$\begin{array}{rcl}
(g+\idmap)(\kappa_{1}a)(K)
& = &
\left\{\begin{array}{ll}
g(a)(M) \quad & \mbox{if }K=\kappa_{1}M, \mbox{ for } M\in\Sigma_{B} \\
0 & \mbox{if }K=\kappa_{2}N \mbox{ for } N\in\Sigma_{Y}
\end{array}\right. \\
(g+\idmap)(\kappa_{2}y)(K)
& = &
\left\{\begin{array}{ll}
0 \quad & \mbox{if }K=\kappa_{1}M, \mbox{ for } M\in\Sigma_{B} \\
\indic{N}(y) & \mbox{if }K=\kappa_{2}N \mbox{ for } N\in\Sigma_{Y}
\end{array}\right.
\end{array}$$

\noindent The two integrals in the above equation can be split,
as in~\cite[Lemma~12]{Jacobs13a}:
$$\begin{array}{rcl}
\lefteqn{\int (\idmap+f)(-)(K)\intd\phi} \\
& = &
\phi(\kappa_{1}B) \cdot \int (\idmap+f)(\kappa_{1}-)(K)  \intd 
   \frac{\phi(\kappa_{1}-)}{\phi(\kappa_{1}B)} +
\phi(\kappa_{2}X) \cdot \int (\idmap+f)(\kappa_{2}-)(K) \intd 
   \frac{\phi(\kappa_{2}-)}{\phi(\kappa_{2}X)} \\
\end{array}$$

\noindent Hence we can separate the two cases $K = \kappa_{1}M, K =
\kappa_{2}N$ in:
$$\begin{array}{rcl}
\int (\idmap+f)(-)(\kappa_{1}M)\intd\phi
& = &
\phi(\kappa_{1}B) \cdot \int \eta(-)(M)  \intd 
   \frac{\phi(\kappa_{1}-)}{\phi(\kappa_{1}B)} \\
& = &
\phi(\kappa_{1}B) \cdot \int \indic{M}  \intd 
   \frac{\phi(\kappa_{1}-)}{\phi(\kappa_{1}B)} \\
& = &
\phi(\kappa_{1}B) \cdot \frac{\phi(\kappa_{1}M)}{\phi(\kappa_{1}B)} \\
& = &
\phi(\kappa_{1}M) \\
\int (\idmap+f)(-)(\kappa_{2}N)\intd\phi
& = &
\phi(\kappa_{2}X) \cdot \int f(-)(N) \intd 
   \frac{\phi(\kappa_{2}-)}{\phi(\kappa_{2}X)}.
\end{array}$$

\noindent Similarly we get:
$$\begin{array}{rcl}
\int (g+\idmap)(-)(\kappa_{1}M)\intd\psi
& = &
\psi(\kappa_{1}A)\cdot \int g(-)(M) \intd
   \frac{\psi(\kappa_{1}-)}{\psi(\kappa_{1}A)} \\
\int (g+\idmap)(-)(\kappa_{2}N)\intd\psi
& = &
\psi(\kappa_{2}N).
\end{array}$$

\noindent The original equation thus amounts to:
$$\begin{array}{rclcrcl}
\psi(\kappa_{1}A) \cdot 
   \int g(-)(M) \intd \frac{\psi(\kappa_{1}-)}{\psi(\kappa_{1}A)}
& = &
\phi(\kappa_{1}M)
& \qquad &
\phi(\kappa_{2}X) \cdot 
   \int f(-)(N) \intd \frac{\phi(\kappa_{2}-)}{\phi(\kappa_{2}X)}
& = &
\psi(\kappa_{2}N)
\end{array}$$

\noindent We need to define a probability measure $\chi \in \Giry(A+X)$.
For $M\in\Sigma_{A}$ and $N\in\Sigma_{X}$ we put:
$$\begin{array}{rclcrcl}
\chi(\kappa_{1}M)
& = &
\psi(\kappa_{1}M)
& \qquad &
\chi(\kappa_{2}N)
& = &
\phi(\kappa_{2}N).
\end{array}$$

\noindent It is easy to see that $\chi$ preserves sums of joins, and
also:
$$\begin{array}{rcl}
\chi(A+X)
& = &
\chi(\kappa_{1}A) + \chi(\kappa_{2}X) \\
& = &
\psi(\kappa_{1}A) + \phi(\kappa_{2}X) \\
& = &
1 - \psi(\kappa_{2}Y) + 1 - \phi(\kappa_{1}B) \\
& = &
2 - \phi(\kappa_{2}X) \cdot \int f(-)(Y) \intd 
   \frac{\phi(\kappa_{2}-)}{\phi(\kappa_{2}X)} - \phi(\kappa_{1}B) \\ 
& = &
2 - \phi(\kappa_{2}X) \cdot \int \indic{Y} \intd 
   \frac{\phi(\kappa_{2}-)}{\phi(\kappa_{2}X)} - \phi(\kappa_{1}B) \\ 
& = &
2 - \phi(\kappa_{2}X) \cdot \frac{\phi(\kappa_{2}Y)}{\phi(\kappa_{2}X)} - 
   \phi(\kappa_{1}B) \\ 
& = &
2 - \phi(\kappa_{2}Y) - \phi(\kappa_{1}B) \\ 
& = &
2 - \phi(B+Y) \\
& = & 
1.
\end{array}$$

\noindent This is the appropriate probability measure, since, for
$M\in\Sigma_{A}$ and $N\in\Sigma_{Y}$,
$$\begin{array}{rclcrcl}
(\idmap+f)_{*}(\chi)(\kappa_{1}M)
& = &
\int (\idmap+f)(-)(\kappa_{1}M) \intd \chi \\
& = &
\chi(\kappa_{1}A)\cdot \int \eta(-)(M) \int 
   \frac{\chi(\kappa_{1}-)}{\chi(\kappa_{1}A)} \\
& = &
\chi(\kappa_{1}A)\cdot \frac{\chi(\kappa_{1}M)}{\chi(\kappa_{1}A)} \\
& = &
\chi(\kappa_{1}M) \\
& = &
\psi(\kappa_{1}M) \\
(\idmap+f)_{*}(\chi)(\kappa_{2}N)
& = &
\int (\idmap+f)(-)(\kappa_{2}N) \intd \chi \\
& = &
\chi(\kappa_{2}X) \cdot \int f(-)(N) \intd 
   \frac{\chi(\kappa_{2}-)}{\chi(\kappa_{2}X)} \\
& = &
\phi(\kappa_{2}X) \cdot \int f(-)(N) \intd 
   \frac{\phi(\kappa_{2}-)}{\phi(\kappa_{2}X)} \\
& = &
\psi(\kappa_{2}N).
\end{array}$$

\noindent Similarly, for $M\in\Sigma_{B}$ and $N\in\Sigma_{X}$,
$$\begin{array}{rcl}
(g+\idmap)_{*}(\chi)(\kappa_{1}M)
& = &
\int (g+\idmap)(\kappa_{1}M) \intd \chi \\
& = &
\chi(\kappa_{1}A) \cdot \int g(-)(M) \intd
   \frac{\chi(\kappa_{1}-)}{\chi(\kappa_{1}A)} \\
& = &
\psi(\kappa_{1}A) \cdot \int g(-)(M) \intd
   \frac{\psi(\kappa_{1}-)}{\psi(\kappa_{1}A)} \\
& = &
\phi(\kappa_{1}M) \\
(g+\idmap)_{*}(\chi)(\kappa_{2}N)
& = &
\int (g+\idmap)(\kappa_{2}N) \intd \chi \\
& = &
\chi(\kappa_{2}X) \cdot \int \eta(-)(N) \intd
   \frac{\chi(\kappa_{2}-)}{\chi(\kappa_{2}X)} \\
& = &
\chi(\kappa_{2}X) \cdot \frac{\chi(\kappa_{2}N)}{\chi(\kappa_{2}X)} \\
& = &
\chi(\kappa_{2}N) \\
& = &
\phi(\kappa_{2}N)
\end{array}$$

\noindent Uniqueness is clear.

We turn to the second diagram in~\eqref{CoprodAssSquares}. So assume a
map $f\colon X \rightarrow \Giry(Y)$ and two probability measures
$\phi\in\Giry(A), \psi \in \Giry(A+X)$ with $(\kappa_{1})_{*}(\phi) =
(\idmap+f)_{*}(\psi)$. We need to prove $(\kappa_{1})_{*}(\phi) =
\psi$. The assumption can be split, for $M\in\Sigma_{A}$ and
$N\in\Sigma_{Y}$, and yields:
$$\begin{array}{rclcrcl}
(\kappa_{1})_{*}(\phi)(\kappa_{1}M)
& = &
(\idmap+f)_{*}(\psi)(\kappa_{1}M)
& \qquad &
(\kappa_{1})_{*}(\phi)(\kappa_{2}N)
& = &
(\idmap+f)_{*}(\psi)(\kappa_{2}N)
\end{array}$$

\noindent We elaborate these equations. We first note:
$$\begin{array}{rcccccl}
(\kappa_{1})_{*}(\phi)
& = &
(\eta \after \kappa_{1})_{*}(\phi)
& = &
\Giry(\kappa_{1})(\phi)
& = &
\phi \after \kappa_{1}^{-1}.
\end{array}$$

\noindent Hence:
$$\begin{array}{rccclcrcccccl}
(\kappa_{1})_{*}(\phi)(\kappa_{1}M)
& = &
\phi(\kappa_{1}^{-1}\kappa_{1}M)
& = &
\phi(M)
& \qquad &
(\kappa_{1})_{*}(\phi)(\kappa_{2}N)
& = &
\phi(\kappa_{1}^{-1}\kappa_{2}N)
& = &
\phi(\emptyset)
& = &
0.
\end{array}$$

\noindent Thus we have:
$$\begin{array}{rcl}
\phi(M)
& = &
(\idmap+f)_{*}(\psi)(\kappa_{1}M) \\
& = &
\int (\idmap+f)(-)(\kappa_{1}M) \intd \psi \\
& = &
\psi(\kappa_{1}A) \cdot \int \eta(-)(M) \intd 
   \frac{\psi(\kappa_{1}-)}{\psi(\kappa_{1}A)} \\
& = &
\psi(\kappa_{1}A) \cdot \frac{\psi(\kappa_{1}M)}{\psi(\kappa_{1}A)} \\
& = &
\psi(\kappa_{1}M) \\
0
& = &
(\idmap+f)_{*}(\psi)(\kappa_{2}N) \\
& = &
\int (\idmap+f)(-)(\kappa_{2}N) \intd \psi \\
& = &
\psi(\kappa_{2}X) \cdot \int f(-)(N) \intd 
   \frac{\psi(\kappa_{2}-)}{\psi(\kappa_{2}X)}.
\end{array}$$

\noindent We need to prove $\psi(\kappa_{2}U) = 0$ for each
$U\in\Sigma_{X}$. For each $x\in U$ we have $f(x)(Y) = 1$, so
$\indic{U} \leq f(-)(Y)$. Thus
$\frac{\psi(\kappa_{1}U)}{\psi(\kappa_{2}X)} \leq \int f(-)(Y) \intd
\frac{\psi(\kappa_{2}-)}{\psi(\kappa_{2}X)} = 0$. This implies
$\psi(\kappa_{2}U) = 0$.

We turn to the jointly monic requirement in~\eqref{CoprodAssJoint}.
So assume $\phi,\psi\in\Giry(n\cdot X+1)$ with
$[\rhd_{i},\kappa_{2}]_{*}(\phi) = [\rhd_{i},\kappa_{2}]_{*}(\psi)$
for each $i\leq n$.  Here, $[\rhd_{i},\kappa_{2}] \colon n\cdot X +1
\rightarrow X+1$ in $\Kl(\Giry)$ is a cotuple of maps $\eta \after
\kappa_{i}$, and thus of the form $\eta \after [\rhd_{i},\kappa_{2}]$
for $[\rhd_{i},\kappa_{2}]$ in $\Sets$. Thus
$([\rhd_{i},\kappa_{2}])_{*}(\phi) =
\Giry([\rhd_{i},\kappa_{2}])(\phi) = \phi \after
     [\rhd_{i},\kappa_{2}]^{-1} \colon \Sigma_{X+1} \rightarrow
     [0,1]$. Now we note that for $M\in\Sigma_{X}$ and $N\in\Sigma_{1}
     = \{\emptyset,1\}$,
$$\begin{array}{rcl}
[\rhd_{i},\kappa_{2}]^{-1}(\kappa_{1}M)
& = &
\kappa_{i}M \\
{[\rhd_{i},\kappa_{2}]}^{-1}(\kappa_{2}\emptyset)
& = &
\emptyset \\
{[\rhd_{i},\kappa_{2}]}^{-1}(\kappa_{2}1)
& = &
\kappa_{1}X \cup \cdots \cup \kappa_{i-1}X \cup \kappa_{i+1}X \cup \cdots \cup
   \kappa_{n}X \cup \kappa_{n+1}1.
\end{array}$$

\noindent Thus for $i\leq n$,
$$\begin{array}{rcl}
\phi(\kappa_{i}M)
& = &
\phi\big([\rhd_{i},\kappa_{2}]^{-1}(\kappa_{1}M)\big) \\
& = &
[\rhd_{i},\kappa_{2}]_{*}(\phi)(\kappa_{1}M) \\
& = &
[\rhd_{i},\kappa_{2}]_{*}(\psi)(\kappa_{1}M) \\
& = &
\psi\big([\rhd_{i},\kappa_{2}]^{-1}(\kappa_{1}M)\big) \\
& = &
\psi(\kappa_{i}M).
\end{array}$$

\noindent Further,
$$\begin{array}{rcl}
\phi(\kappa_{n+1}1)
& = &
\phi\big([\rhd_{1},\kappa_{2}]^{-1}(\kappa_{2}1)\big) - \phi(\kappa_{2}X) - \cdots
   - \phi(\kappa_{n}X) \\
& = &
[\rhd_{i},\kappa_{2}]_{*}(\phi)(\kappa_{2}1) - \phi(\kappa_{2}X) - \cdots
   - \phi(\kappa_{n}X) \\
& = &
[\rhd_{i},\kappa_{2}]_{*}(\psi)(\kappa_{2}1) - \psi(\kappa_{2}X) - \cdots
   - \psi(\kappa_{n}X) \\
& = &
\psi(\kappa_{n+1}1).
\end{array}$$



}

Next we look at opposite categories $\op{(\CstarPU)}$, $\op{\Rng}$ and
$\op{\DL}$ of $C^*$-algebras, rings, and distributive lattices. We
shall start reasoning in $\Rng$, since the arguments carry over to the
other categories in a straightforward manner.

We wish to prove that the squares of~\eqref{CoprodAssLemSquares} are
pullbacks in the category $\op{\Rng}$. Thus we have to show that they
are pushouts in $\Rng$, using products $\times$ instead of coproducts
$+$:
$$\xymatrix{
A\times X\ar[r]^-{\idmap\times f}\ar[d]_{g\times\idmap} & 
  A\times Y\ar[d]^{g\times\idmap}
& &
A\times X\ar[r]^-{\idmap\times f}\ar[d]_{\pi_1} & A\times Y\ar[d]^{\pi_1} \\
B\times X\ar[r]_-{\idmap\times f} & B\times Y
& &
A\ar@{=}[r] & A
}$$

\noindent We start with the one on the left. Suppose we have ring
homomorphisms $\alpha\colon A\times Y\rightarrow C$ and $\beta \colon
B\times X \rightarrow C$ with $\alpha \after (\idmap\times f) = \beta
\after (g\times\idmap)$. Explicitly, this means that:
$$\begin{array}{rcl}
\alpha(a, f(x))
& = &
\beta(g(a), x), \quad\mbox{for all } a\in A, x\in X.
\end{array}$$

\noindent Then we define $\gamma \colon B\times Y \rightarrow C$ by:
$$\begin{array}{rcl}
\gamma(b,y)
& = &
\alpha(0,y) + \beta(b,0).
\end{array}$$

\noindent It is easy to see that $\gamma$ preserves addition. It
preserves the unit, since:
$$\begin{array}{rcl}
\gamma(1,1)
\hspace*{\arraycolsep} = \hspace*{\arraycolsep}
\alpha(0,1) + \beta(1,0) 
& = &
\alpha(0, f(1)) + \beta(1,0) \\
& = &
\beta(g(0), 1) + \beta(1,0) 
\hspace*{\arraycolsep} = \hspace*{\arraycolsep}
\beta(0,1) + \beta(1,0)
\hspace*{\arraycolsep} = \hspace*{\arraycolsep}
\beta(1,1) 
\hspace*{\arraycolsep} = \hspace*{\arraycolsep}
1.
\end{array}$$

\noindent Via a similar trick one shows that $\gamma$ preserves
multiplication: the `cross terms' disappear since:
$$\begin{array}{rcccccl}
\alpha(0, y)\cdot \beta(b, 0)
& = &
\alpha(0, y) \cdot \alpha(0, f(1))\cdot \beta(b, 0)
& = &
\alpha(0, y) \cdot \beta(g(0), 1) \cdot \beta(b, 0)
& = &
0.
\end{array}$$

\noindent Clearly, this $\gamma$ is the unique mediating map.

\auxproof{
This $\gamma$ is a linear map:
$$\begin{array}{rcl}
\gamma((b_{1},y_{1}) + (b_{2},y_{2}))
& = &
\gamma(b_{1}+b_{2}, y_{1}+y_{2}) \\
& = &
\alpha(0, y_{1}+y_{2}) + \beta(b_{1}+b_{2}, 0) \\
& = &
\alpha(0, y_{1}) + \alpha(0, y_{2}) + \beta(b_{1}, 0) + \beta(b_{2}, 0) \\
& = &
\alpha(0, y_{1}) + \beta(b_{1}, 0) + \alpha(0, y_{2}) + \beta(b_{2}, 0) \\
& = &
\gamma(b_{1},y_{1}) + \gamma(b_{2},y_{2}) \\
\gamma(z\cdot (b,y))
& = &
\gamma(z\cdot b, z\cdot y) \\
& = &
z\cdot \alpha(0, y) + z\cdot \beta(b, 0) \\
& = &
z\cdot (\alpha(0, y) + \beta(b, 0)) \\
& = &
z\cdot \gamma(b,y).
\end{array}$$

In the ring-theoretic case it also preserves multiplication since:
$$\begin{array}{rcl}
\gamma(b_{1},y_{1}) \cdot \gamma(b_{2},y_{2})
& = &
(\alpha(0,y_{1}) + \beta(b_{1},0)) \cdot (\alpha(0,y_{2}) + \beta(b_{2},0)) \\
& = &
\alpha(0,y_{1})\cdot \alpha(0,y_{2}) + \alpha(0,y_{1})\cdot \beta(b_{2},0)) \\
& & \qquad +
   \beta(b_{1},0) \cdot \alpha(0,y_{2}) + \beta(b_{1},0) \cdot \beta(b_{2},0)) \\
& = &
\alpha(0,y_{1})\cdot \alpha(0,y_{2}) + 
   \alpha(0,y_{1}\cdot 1)\cdot \beta(b_{2},0)) \\
& & \qquad +
   \beta(b_{1},0) \cdot \alpha(0,1\cdot y_{2}) + 
   \beta(b_{1},0) \cdot \beta(b_{2},0)) \\
& = &
\alpha(0, y_{1}\cdot y_{2}) +
   \alpha(0,y_{1})\cdot \alpha(0,1)\cdot \beta(b_{2},0)) \\
& & \qquad +
   \beta(b_{1},0) \cdot \alpha(0,1) \cdot \alpha(0,y_{2}) + 
   \beta(b_{1}\cdot b_{2}, 0) \\
& = &
\alpha(0, y_{1}\cdot y_{2}) +
   \alpha(0,y_{1})\cdot \alpha(0,f(1))\cdot \beta(b_{2},0)) \\
& & \qquad +
   \beta(b_{1},0) \cdot \alpha(0,f(1)) \cdot \alpha(0,y_{2}) + 
   \beta(b_{1}\cdot b_{2}, 0) \\
& = &
\alpha(0, y_{1}\cdot y_{2}) +
   \alpha(0,y_{1})\cdot \beta(g(0),1)\cdot \beta(b_{2},0)) \\
& & \qquad +
   \beta(b_{1},0) \cdot \beta(g(0),1) \cdot \alpha(0,y_{2}) + 
   \beta(b_{1}\cdot b_{2}, 0) \\
& = &
\alpha(0, y_{1}\cdot y_{2}) +
   \alpha(0,y_{1})\cdot \beta(0\cdot b_{2},1 \cdot 0)) \\
& & \qquad +
   \beta(b_{1}\cdot 0,0\cdot 1) \cdot \alpha(0,y_{2}) + 
   \beta(b_{1}\cdot b_{2}, 0) \\
& = &
\alpha(0, y_{1}\cdot y_{2}) + \alpha(0,y_{1})\cdot 0 +
   0 \cdot \alpha(0,y_{2}) + \beta(b_{1}\cdot b_{2}, 0) \\
& = &
\alpha(0, y_{1}\cdot y_{2}) + \beta(b_{1}\cdot b_{2}, 0) \\
& = &
\gamma(b_{1} \cdot b_{2}, y_{1} \cdot y_{2}) \\
& = &
\gamma((b_{1},y_{1}) \cdot (b_{2},y_{2}))
\end{array}$$

\noindent This map $\gamma$ is positive, since if $(b,y) \geq 0$,
then $b\geq 0$ and $y\geq 0$, so $(0,y) \geq 0$ and $(b,0) \geq 0$,
which yields $\alpha(0,y) \geq 0$ and $\beta(b,0) \geq 0$, and
thus $\gamma(b,y) = \alpha(0,y) + \beta(b,0) \geq 0$.

This map $\gamma$ satisfies:
$$\begin{array}{rcl}
(\gamma \after (g\times\idmap))(a,y)
& = &
\gamma(g(a), y) \\
& = &
\alpha(0, y) + \beta(g(a), 0) \\
& = &
\alpha(0, y) + \alpha(a, f(0)) \\
& = &
\alpha(a,y) \\
(\gamma \after (\idmap\times f))(b,x)
& = &
\gamma(b, f(x)) \\
& = &
\alpha(0, f(x)) + \beta(b, 0) \\
& = &
\beta(g(0), x) + \beta(b, 0) \\
& = &
\beta(b,x).
\end{array}$$

Finally, if $\gamma' \colon B\times Y \rightarrow C$ also satisfies
$\gamma' \after (g\times\idmap) = \alpha$ and $\gamma' \after
(\idmap\times f) = \beta$, then:
$$\begin{array}{rcl}
\gamma(b,y)
& = &
\alpha(0,y) + \beta(b,0) \\
& = &
\gamma'(g(0), y) + \gamma'(b, f(0)) \\
& = &
\gamma'(b,y).
\end{array}$$
}

We turn to the above square on the right: given $\alpha \colon A
\rightarrow B$ and $\beta \colon A\times Y \rightarrow B$ with $\alpha
\after \pi_{1} = \beta \after (\idmap\times f)$. The only possible
mediating map is $\alpha \colon A \rightarrow B$. Hence we need to
prove $\beta = \alpha \after \pi_{1} \colon A\times Y \rightarrow B$,
that is, $\beta(a, y) = \alpha(a)$, for all $a\in A, y\in Y$.
Consider the function $\beta' \colon Y \rightarrow B$ defined by
$\beta'(y) = \beta(0,y)$. This $\beta'$ preserves finite sums $(+,0)$
and multiplication $\cdot$, and satisfies:
$$\begin{array}{rcccccccl}
\beta'(1)
& = &
\beta(0, 1)
& = &
\beta(0, f(1))
& = &
\alpha(0)
& = &
0.
\end{array}$$

\auxproof{
$$\begin{array}{rcl}
\beta'(0)
& = &
\beta(0,0) \\
& = &
0 \\
\beta'(y+y')
& = &
\beta(0, y+y') \\
& = &
\beta(0,y) + \beta(0,y') \\
& = &
\beta'(y) + \beta'(y') \\
\beta'(y\cdot y') 
& = &
\beta(0, y\cdot y') \\
& = &
\beta(0,y) \cdot \beta(0,y') \\
& = &
\beta'(y) \cdot \beta'(y') \\
\beta'(z\cdot y)
& = &
\beta(0, z\cdot y) \\
& = &
z\cdot \beta(0,y) \\
& = &
z\cdot \beta'(y).
\end{array}$$
}

\noindent Then $\beta' = 0$, since $\beta'(y) = \beta'(1\cdot y) =
\beta'(1)\cdot \beta'(y) = 0\cdot \beta'(y) = 0$. But then we are done
since:
$$\begin{array}{rcccccccl}
\beta(a,y)
& = &
\beta(a, 0) + \beta(0, y)
& = &
\beta(a, f(0)) + \beta'(y)
& = &
\alpha(a) + 0 
& = &
\alpha(a).
\end{array}$$

\noindent Next we investigate the jointly monic requirement
from~\eqref{CoprodOneAssJoint} in the category $\op{\Rng}$. When we
translate this back to $\Rng$ we have to prove the corresponding maps
$\IV, \XI \colon \Z\times \Z \rightarrow \Z\times\Z\times\Z$ are
jointly epic, where now $\IV = \tuple{\idmap, \pi_{2}}$ and $\XI =
\tuple{\tuple{\pi_{2},\pi_{1}}, \pi_{2}}$. Explicitly:
$$\begin{array}{rclcrcl}
\IV(k,m)
& = &
(k, m, m)
& \qquad\mbox{and}\qquad &
\XI(k,m)
& = &
(m, k, m).
\end{array}$$

\noindent Let $f, g\colon \Z\times\Z\times\Z \rightarrow A$ satisfy $f
\after \IV = g \after \IV$ and $f \after \XI = g \after \XI$. We
obtain $f=g$ from:
$$\begin{array}{rcl}
f(k,m,\ell)
& = &
f(k, \ell, \ell) + f(0, m - \ell, 0) \\
& = &
\big(f \after \IV\big)(k, \ell) + \big(f \after \XI\big)(m - \ell, 0) \\
& = &
\big(g \after \IV\big)(k, \ell) + \big(g \after \XI\big)(m - \ell, 0) \\
& = &
g(k, m, \ell).
\end{array}$$

\auxproof{
These maps are given by:
$$\begin{array}{rcl}
[\rhd_{i},\pi_{2}](a, k)
& = &
(k\cdot 1, \ldots, k\cdot 1, a, k\cdot 1, \ldots, k\cdot 1, k),
\end{array}$$

\noindent where the element $a$ is at the $i$-th position. Now let
$f,g\colon A^{n}\times \Z \rightarrow B$ be ring homomorphisms with
$f \after [\rhd_{i},\pi_{2}] = g \after [\rhd_{i},\pi_{2}]$ for each $i$.
We prove $f=g$ in the following manner.
$$\begin{array}{rcl}
\lefteqn{f(a_{1}, \ldots, a_{n}, k)} \\
& = &
f(a_{1}, k\cdot 1, \ldots, k\cdot 1, k) + f(0, a_{2} - k\cdot 1, 0, \ldots, 0)
  + \cdots + f(0, \ldots, 0, a_{n} - k\cdot 1, 0) \\
& = &
g(a_{1}, k\cdot 1, \ldots, k\cdot 1, k) + g(0, a_{2} - k\cdot 1, 0, \ldots, 0)
  + \cdots + g(0, \ldots, 0, a_{n} - k\cdot 1, 0) \\
& = &
g(a_{1}, \ldots, a_{n}, k).
\end{array}$$
}



\noindent The structure on predicates in a ring will be described
separately in Proposition~\ref{RingPredProp} below.

This line of reasoning showing that the category $\op{\Rng}$ is an
effectus can be transferred directly to the opposite $\op{\DL}$ of the
category of distributive lattices, and also to the opposite
$\op{(\CstarPU)}$ of the category of $C^*$-algebras with positive
unital maps. The only point, in the latter case, where some care is
needed is in the demonstration that $\beta'(1) = 0$ implies
$\beta'=0$, where $\beta'$ is a positive unital map as defined
above. It is an elementary result in the theory of $C^*$-algebras that
a linear positive map $f$ with $f(1) = 0$ must be the zero map. (For a
proof, use that the operator norm $\|f\|$ equals $\|f(1)\|$, so
$\|f(a)\| \leq \|f(1)\|\cdot \|a\| = 0$, for each element $a$.)

\auxproof{
In $\op{\DL}$ we have $[\rhd_{i},\pi_{2}] \colon A\times 2 \rightarrow
A^{n}\times 2$ given by, for $k\in 2 = \{0,1\}$,
$$\begin{array}{rcl}
[\rhd_{i},\pi_{2}](a, k)
& = &
(k, \ldots, k, a, k, \ldots, k, k).
\end{array}$$

\noindent Now let $f,g\colon A^{n}\times 2 \rightarrow B$ be
homomorphisms of distributive lattices with $f \after
[\rhd_{i},\pi_{2}] = g \after [\rhd_{i},\pi_{2}]$ for each $i$.  We
prove $f=g$ in the following manner.
$$\begin{array}{rcl}
f(a_{1}, \ldots, a_{n}, 0) 
& = &
f(a_{1}, 0, \ldots, 0, 0) \vee f(0, a_{2}, 0, \ldots, 0, 0)
  \vee \cdots \vee f(0, \ldots, 0, a_{n}, 0) \\
& = &
g(a_{1}, 0, \ldots, 0, 0) \vee g(0, a_{2}, 0, \ldots, 0, 0)
  \vee \cdots \vee g(0, \ldots, 0, a_{n}, 0) \\
& = &
g(a_{1}, \ldots, a_{n}, 0) \\
f(a_{1}, \ldots, a_{n}, 1) 
& = &
f(a_{1}, 1, \ldots, 1, 1) \wedge f(1, a_{2}, 1, \ldots, 1, 1)
  \wedge \cdots \wedge f(1, \ldots, 1, a_{n}, 1) \\
& = &
g(a_{1}, 1, \ldots, 1, 1) \wedge g(1, a_{2}, 1, \ldots, 1, 1)
  \wedge \cdots \wedge g(1, \ldots, 1, a_{n}, 1) \\
& = &
g(a_{1}, \ldots, a_{n}, 1).
\end{array}$$
}

The effect algebra structure from Definition~\ref{OveeDef},
specialised to the subset of effects $[0,1]_{A} \subseteq A$ of a
$C^*$-algebra $A$, is the standard one, with $e\ovee d$ defined if
$e+d \leq 1$, and in that case equal to $e+d$. The resulting order is
the usual one for $C^*$-algebras.


\end{exa}

We now formulate the structure of predicates on rings. Recall from
Example~\ref{PredEx}~\eqref{PredExRng} that predicates $R \rightarrow
1+1$ in the opposite $\op{\Rng}$ of the category of rings are in
one-to-one correspondence with idempotents of the ring $R$.


\begin{prop}
\label{RingPredProp}
The effect algebra of idempotents $\Pred(R) = \setin{e}{R}{e^{2} = e}$
in a ring $R$ can be described as follows. Two idempotents $e,d\in
\Pred(R)$ are orthogonal (\textit{i.e.}~$e \orthogonal d$) iff $e\cdot
d = 0 = d\cdot e$. In that case their effect algebra sum $e \ovee d
\in \Pred(R)$ equals their sum $e\ovee d = e+d$ in $R$. One has
$e^{\perp} = 1 - e$. The induced order relation is: $e \leq d$ iff
$e\cdot d = e = d\cdot e$.

\auxproof{
For abstract reasons the set of idempotents $\Pred(R)$ of a ring $R$ is an
effect algebra. We check the requirements concretely. 
\begin{itemize}
\item Clearly, if $e \orthogonal d$, then $e\ovee d = e+d$ is a
  idempotent again, since $(e\ovee d)^{2} = e^{2} + ed + de + d^{2} =
  e + d = e\ovee d$.

\item If $e \orthogonal d$, then $ed = 0 = de$, so that also $d
  \orthogonal e$ and $e\ovee d = e+d = d+e = d\ovee e$.

\item Also, $e\orthogonal 0$, since $e0 = 0 = 0e$, and $e\ovee 0 = e+0 = e$.

\item Assume $e \orthogonal d$ and $(e\ovee d) \orthogonal f$. We need
  to prove $d\orthogonal f$ and $e \orthogonal (d \ovee f)$. So, by
  assumption, $ed = 0 = de = 0$ and $(e+d)f = 0 = f(e+d) = 0$.
\begin{itemize}
\item For $d\orthogonal f$ we need to prove $df = 0 = fd$. 
$$\begin{array}{rcccccccccl}
df
& = &
(0 + d^{2})f
& = &
(de + d^{2})f 
& = &
d(e+d)f
& = &
d0
& = &
0.
\end{array}$$

\noindent And similarly,
$$\begin{array}{rcccccccccl}
fd
& = &
f(0 + d^{2})f
& = &
f(ed + d^{2})
& = &
f(e+d)d
& = &
0d
& = &
0.
\end{array}$$

By the same reasoning we get $ef = 0 = fe$, namely:
$$\begin{array}{rcccccccccl}
ef
& = &
(e^{2}+0)f
& = &
(e^{2}+ed)f 
& = &
e(e+d)f
& = &
e0
& = &
0.
\end{array}$$

\noindent And:
$$\begin{array}{rcccccccccl}
fe
& = &
f(e^{2}+0)
& = &
f(e^{2} + de)
& = &
f(e+d)e
& = &
0e
& = &
0.
\end{array}$$

\item For $e \orthogonal (d \ovee f)$ we need to prove $e(d+f) = 0 =
  (d+f)e$. But this is easy by the previous point, since:
$$\begin{array}{rcccccccl}
e(d+f)
& = &
ed + ef
& = & 
0+0
& = &
de + fe
& = &
(d+f)e.
\end{array}$$
\end{itemize}

\item Clearly, $e \orthogonal (1-e)$ since $e(1-e) = e - e^{2} = 0 = e - e^{2}
= (1-e)e$, and $e \ovee (1-e) = e + (1-e) = 1$. Also, $1-e$ is unique
with this property, since if $de = 0 = ed$ and $e+d=1$, then $d = 1-e$.

\item If $e\orthogonal 1$, then $1e = 0 = e1$, so that $e=0$.
\end{itemize}

\noindent Just curiosity: how difficult is it to prove directly that
$e \leq d$ iff $ed = e = de$, for $e,d$ idempotents is a partial
order?
\begin{itemize}
\item $e \leq e$ holds since $e^{2} = e$

\item If $e \leq d \leq f$, then $ed = e = de$ and $df = d = fd$. We
have to prove $ef = 0 = fe$. This is easy:
$$\begin{array}{rcccccccccccl}
ef
& = &
edf
& = &
e0
& = &
0
& = &
0e
& = &
fde
& = &
fe.
\end{array}$$

\item If $e \leq d$ and $d \leq e$, then $ed = e = de$ and
$de = d = ed$ from which we immediately get $e=d$.
\end{itemize}
}

In the special case that the ring $R$ is commutative, the effect
algebra $\Pred(R)$ of idempotents is actually a Boolean algebra with
$e \wedge d = ed$ and $e\vee d = e + d - ed$.

\auxproof{
We first prove that $ed=de$ is the meet.
\begin{itemize}
\item $ed \leq e$ since $e(ed) = e^{2}d = ed$

\item $ed \leq d$ since $d(ed) = ed^{2} = ed$

\item if $x \leq e,d$ for an idempotent $x\in R$ then $x \leq ed$
since $x(ed) = xed = xd = x$.
\end{itemize}

\noindent The resulting join can now be described as:
$$\begin{array}{rcccccccl}
e \vee d
& = &
(e^{\perp} \wedge d^{\perp})^{\perp}
& = &
1 - (1-e)\cdot (1-d)
& = &
1 - (1 - d - e + ed)
& = &
e + d - ed.
\end{array}$$

\noindent We have that $\wedge$ distributes over $\vee$ since:
$$\begin{array}{rcl}
e \wedge (d \vee f)
& = &
e\cdot (d + f - df) \\
& = &
ed + ef - edf \\
& = &
ed + ef - (ed)(ef) \\
& = &
ed \vee ef \\
& = &
(e\wedge d) \vee (e\wedge f)
\end{array}$$

\noindent Hence we get a Boolean algebra.
}

These mappings can be made functorial, and yield a diagram:
\begin{equation}
\label{RingPredDiag}
\vcenter{\xymatrix@R-.5pc{
\Rng\ar[d]_{\Pred} & & \,\CRng\ar@{_(->}[ll]\ar[d]^{\Pred} \\
\EA & & \,\BA\ar@{_(->}[ll]
}}
\end{equation}

\noindent where $\CRng$ is the category of commutative rings.

\auxproof{
Clearly, a ring homomorphism $f\colon R \rightarrow S$ restricts to
an effect algebra homomorphism between subsets of idempotents, since:
\begin{enumerate}
\item if $e\in R$ is an idempotent, so is $f(e) \in S$, since
$f(e)\cdot f(e) = f(e^{2}) = f(e)$

\item if $e \orthogonal d$, then $f(e) \orthogonal f(d)$, since
$e \orthogonal d$ means $ed+de = 0$, so that $f(e)f(d) + f(d)f(e) =
f(ed+de) = f(0) = 0$.

\item in that case $f(e\ove d) = f(e+d) = f(e) + f(d) = f(e) \ovee
  f(d)$.

\item obviously, $f(1) = 1$.
\end{enumerate}
}
\end{prop}

The second part of this result --- that the idempotents of a
commutative ring form a Boolean algebra --- is well-known, see for
instance~\cite[V, Lemma~2.3]{Johnstone82}. The first part however ---
that in the general case the idempotents form an effect algebra --- is
less familiar, but the essence was already observed
in~\cite{Katrnoska92}.


\begin{myproof}
For idempotent elements $e,d\in R$, corresponding to functions $f_{e},
f_{d} \colon \Z^{2} \rightarrow R$ as in
Example~\ref{PredEx}~\eqref{PredExRng}, we have $e\orthogonal d$ iff
there is a bound $b\colon \Z^{3} \rightarrow R$ with $b \after
\tuple{\idmap,\pi_{2}} = f_{e}$ and $b \after
\tuple{\tuple{\pi_{2},\pi_{1}}, \pi_{2}} = f_{d}$. This means that
$b(n,m,m) = f_{e}(n,m)$ and $b(m,n,m) = f_{d}(n,m)$. Moreover, one
obtains an idempotent $x = b(0,0,1) \in R$ satisfying:
$$\begin{array}{rcccccccl}
e + d + x
& = &
f_{e}(1,0) + f_{d}(1,0) + x
& = &
b(1,0,0) + b(0,1,0) + b(0,0,1)
& = &
b(1,1,1)
& = &
1.
\end{array}$$

\noindent But then $x = 1 - e - d$. Using that $b$ is a ring
homomorphism we get:
$$\begin{array}{rcccccccl}
e\cdot d
& = &
b(1,0,0)\cdot b(0,1,0)
& = &
b\big((1,0,0)\cdot(0,1,0)\big)
& = &
b(0,0,0)
& = &
0.
\end{array}$$

\noindent In the same way we can get $d\cdot e = 0$.

In the reverse direction, if $e\cdot d = 0 = d\cdot e$, then we define
a map $b\colon \Z^{3} \rightarrow R$ by:
$$\begin{array}{rcl}
b(n,m,k)
& = &
n\cdot e + m \cdot d + k\cdot (1-e-d).
\end{array}$$

\noindent One immediately gets $b(n,m,m) = f_{e}(n,m)$ and $b(m,n,m) =
f_{d}(n,m)$. With a bit more effort, using $e\cdot d = 0 = d\cdot e$,
one proves that $b$ is a ring homomorphism. The resulting effect
algebra sum $e\ovee d\in R$ is given by their sum $e\ovee d = e+d$ in
$R$, corresponding to the function $f_{e\ovee d} \colon \Z^{2}
\rightarrow R$ described by $f_{e\ovee d}(n,m) = b(n,n, m)$, as in
Definition~\ref{OveeDef}.

\auxproof{
$$\begin{array}{rcl}
b(n,m,m)
& = &
n\cdot e + m \cdot d + m\cdot (1-e-d) \\
& = &
n\cdot e - m\cdot (1-e) \\
& = &
f_{e}(n,m) \\
b(m,n,m)
& = &
m\cdot e + n \cdot d + m\cdot (1-e-d) \\
& = &
n\cdot d + m\cdot (1-d) \\
& = &
f_{d}(n,m).
\end{array}$$

$$\begin{array}{rcl}
b(0,0,0)
& = &
0\cdot e + 0 \cdot d + 0\cdot (1-e-d) \\
& = &
0 \\
b((n, m, k) + (n',m',k'))
& = &
b(n+n', m+m', k+k') \\
& = &
(n+n')\cdot e + (m+m') \cdot d + (k+k')\cdot (1-e-d) \\
& = &
n\cdot e + m\cdot d + k\cdot (1-e-d) 
   + n'\cdot e + m' \cdot d + k'\cdot (1-e-d) \\
& = &
b(n,m,k) + b(n',m',k') \\
b(1,1,1)
& = &
1\cdot e + 1 \cdot d + 1\cdot (1-e-d) \\
& = &
e + d + 1 - e - d \\
& = &
1.
\end{array}$$

\noindent We do the most difficult case separately:
$$\begin{array}{rcl}
\lefteqn{b(n,m,k) \cdot b(n',m',k')} \\
& = &
\big(n\cdot e + m\cdot d + k\cdot (1-e-d)\big) \cdot
   \big(n'\cdot e + m' \cdot d + k'\cdot (1-e-d)\big) \\
& = &
(n\cdot e)\cdot(n'\cdot e) + (n\cdot e)\cdot(m'\cdot d) + 
   (n\cdot e)\cdot(k'\cdot (1-e-d)) \\
& & \quad +
(m\cdot d)\cdot(n'\cdot e) + (m\cdot d)\cdot(m'\cdot d) + 
   (m\cdot d)\cdot(k'\cdot (1-e-d)) \\
& & \quad +
(k\cdot (1-e-d))\cdot(n'\cdot e) + (k\cdot (1-e-d))\cdot(m'\cdot d) + 
   (k\cdot (1-e-d))\cdot(k'\cdot (1-e-d)) \\
& = &
(n\cdot n')\cdot (e\cdot e) + (n\cdot m') \cdot (e\cdot d) +
   (n\cdot k')\cdot (e - e^{2} - ed) \\
& & \quad
   (m\cdot n')\cdot (d\cdot e) + (m\cdot m')\cdot (d\cdot d) +
   (m \cdot k') \cdot (d - de - d^{2}) \\
& & \quad
   (k\cdot n')\cdot (e - e^{2} - de) + (k\cdot m')\cdot (d - ed - d^{2}) +
   (k\cdot k')\cdot (1 - e - d - e + e^{2} + ed - d + de + d^{2}) \\
& = &
(n\cdot n')\cdot e + (m\cdot m') \cdot d + (k\cdot k')\cdot (1-e-d) \\
& = &
b(n\cdot n', m\cdot m', k\cdot k') \\
& = &
b((n, m, k) \cdot (n',m',k')).
\end{array}$$

The sum $e+d$ of orthogonal idempotents $e,d$ is idempotent again since:
$$\begin{array}{rcccccl}
(e+d)^{2}
& = &
e^{2} + ed + de + d^{2}
& = &
e + 0 + d
& = &
e + d.
\end{array}$$
}

Next we prove the characterisation of the order $\leq$ on
predicates/idempotents. Assume $e\cdot d = e = d\cdot e$. We need to
find an idempotent $x\in R$ with $e\cdot x = 0 = x\cdot e$ and $e + x
= d$. Clearly, the only possible choice is $x = d - e$. This $x$ is
idempotent and satisfies $e \orthogonal x$ since
$$\begin{array}{rcccccccccccccccl}
e\cdot x
& = &
e\cdot (d-e)
& = &
e\cdot d - e^{2}
& = &
e - e
& = &
0 
& = &
e - e
& = &
d\cdot e - e^{2}
& = &
(d-e)\cdot e
& = &
x\cdot e.
\end{array}$$

\auxproof{
This $x$ is idempotent since $x^{2} = (d-e)(d- e) = d^{2} - de - ed + e^{2} 
= d - e - e + e = d-e$.
}

\noindent In the other direction, assume $e \leq d$. Then there is an
idempotent $x\in R$ with $e \orthogonal x$ and $e + x = d$. The latter
equation yields $x = d - e$. Orthogonality $e\orthogonal x$ gives
$e\cdot (e-d) = 0 = (e-d)\cdot e$, so that $e\cdot d = e = d\cdot
e$. \QED

\auxproof{ 
$$\begin{array}{rcccl}
1 - (1-e)\cdot(1-d)
& = &
1 - (1 - d - e + ed) 
& = &
e + d - ed.
\end{array}$$

\noindent Then $e,d \leq e\vee d$ since:
$$\begin{array}{rcl}
e(e\vee d)
& = &
e(d + e - ed) \\
& = &
ed + e^{2} - e^{2}d \\
& = &
ed + e - ed \\
& = &
e \\
d(e\vee d)
& = &
d(d + e - ed) \\
& = &
d^{2} + de - ded \\
& = &
d + de - d^{2}e \\
& = &
d + de - de \\
& = &
d.
\end{array}$$

\noindent And if $e \leq x$ and $d \leq x$, then $e\vee d \leq x$,
since $ex = e$ and $dx = d$ implies:
$$\begin{array}{rcl}
(e\vee d)x
& = &
(d + e - ed)x \\
& = &
dx + ex - edx \\
& = &
d + e - ed \\
& = &
e\vee d.
\end{array}$$

Similarly, $e\cdot d$ is the greatest lower bound of $e,d$, since
$e\cdot d \leq e,d$ is shown in:
$$\begin{array}{rcl}
(e \wedge d)e
& = &
ede \\
& = &
e^{2}d \\
& = &
ed \\
& = &
e\wedge d \\
(e \wedge d)d
& = &
ed^{2} \\
& = &
ed \\
& = &
e\wedge d.
\end{array}$$

\noindent And if $x \leq e$ and $x \leq d$, then $x\leq e \wedge d$
since:
$$\begin{array}{rcccccl}
x(e \wedge d)
& = &
xed 
& = &
xd 
& = &
x.
\end{array}$$

We also have $e \wedge e^{\perp} = 0$ and $e\vee e^{\perp} = 1$
since:
$$\begin{array}{rcl}
e \vee e^{\perp}
& = &
e + (1-e) - e(1-e) \\
& = &
1 - e + e^{2} \\
& = &
1 \\
e \wedge e^{\perp}
& = &
e(1-e) \\
& = &
e - e^{2} \\
& = &
0.
\end{array}$$

And, distributivity holds:
$$\begin{array}{rcl}
e \wedge (d \vee d')
& = &
e \cdot (d + d' - dd') \\
& = &
ed + ed' - edd' \\
& = &
ed + ed' - eedd' \\
& = &
ed + ed' - eded' \\
& = &
(e\wedge d) \vee (e\wedge d')
\end{array}$$
}
\end{myproof}

Recall from Example~\ref{PredEx}~\eqref{PredExRng} that for a ring $R$
an $n$-test $R \rightarrow n\cdot 1$ in $\op{\Rng}$ consists of a ring
homomorphism $f\colon \Z^{n} \rightarrow R$, corresponding to $n$
idempotents $e_{i} = f(\ket{i}) \in R$ with $e_{1} + \cdots + e_{n} =
1$ and $e_{i}\cdot e_{j} = 0$ for $i\neq j$. Such an $n$-test is an
essential ingredient of the Peirce decomposition $R
\conglongrightarrow \bigoplus_{i,j} e_{i}Re_{j}$ of the ring $R$.

Although the requirements in Assumption~\ref{CoprodAss} look rather
mild, they do not hold in all categories. Specifically, in the Kleisli
category of the powerset monad they fail. This is significant, because
this powerset monad is used for the semantics of non-deterministic
computation. It is thus out of scope.

\begin{rem}
\label{PowNonAssumpRem}
The Kleisli category $\Kl(\Pow)$ of the powerset monad $\Pow$ can be
identified with the category of sets and relations between them.  It
does not satisfy the jointly monic requirement in
Assumption~\ref{CoprodAss}. We construct a counterexample for $n=2$.

For a singleton set $1$ we have $1+1+1 \cong \{a,b,c\}$ and $1+1 \cong
\{u,v\}$. We can describe the two maps $\IV,
\XI\colon 3 \rightrightarrows 2$
in~\eqref{CoprodOneAssJoint} as:
$$\begin{array}{rclcrccclcrccclcrcl}
\IV(a)
& = &
u
& \quad &
\IV(b)
& = &
\IV(c)
& = &
v 
& \qquad\qquad &
\XI(a) 
& = &
\XI(c) 
& = &
v
& \quad &
\XI(b)
& = & 
u.
\end{array}$$

\noindent The two subset $S = \{a, b\}$ and $T = \{a,b,c\}$ of $3$
are clearly different but they satisfy:
$$\begin{array}{rccclcrcccl}
\Pow(\IV)(S) & = & \{u,v\} & = & \Pow(\IV)(T)
& \quad\mbox{and}\quad &
\Pow(\XI)(S) & = & \{u,v\} & = & \Pow(\XI)(T).
\end{array}$$

\noindent Hence the two maps $\IV$ and $\XI$ are \emph{not} jointly
monic in the category $\Kl(\Pow)$.

On a different level, notice that the empty set $0$ is a zero object
in $\Kl(\Pow)$: it is both initial and final. Hence $0+0 \cong 0$ is
also a zero object, so that there is only one predicate $X \rightarrow
0+0$. The category $\Kl(\Pow)$ is not of interest in the current
setting.
\end{rem}

We finish this section with some observations about multiple sums
$e_{1} \ovee \cdots \ovee e_{n}$.

\begin{rem}
\label{ManyOveeRem}
Let $\cat{B}$ be an effectus \textit{i.e.}~a category for which
Assumption~\ref{CoprodAss} holds. Assume we have multiple predicates
$p_{1}, \ldots, p_{n} \colon X \rightarrow 1+1$. We would like to
express what it means that these $p_{i}$ are jointly orthogonal,
\textit{i.e.}~that their sum $p_{1} \ovee \cdots \ovee p_{n} \colon X
\rightarrow 1+1$ exists.

This can be done by using the $n$-ary maps $[\rhd_{i},\kappa_{2}]
\colon n\cdot 1 + 1 \rightarrow 1+1$ from Assumption~\ref{CoprodAss}.
We say that $p_{1}, \ldots, p_{n} \colon X \rightarrow 1+1$ are
pairwise orthogonal if there is a single bound $b\colon 1 \rightarrow
(n+1)\cdot 1$ with $[\rhd_{i},\kappa_{2}] \after b = p_{i}$, for each
$i\in \{1,\ldots, n\}$. The sum is then defined as:
$$\xymatrix{
p_{1} \ovee \cdots \ovee p_{n} = \big(X \ar[r]^-{b} &
   (n+1)\cdot 1 = n\cdot 1 + 1\ar[r]^-{\nabla+\idmap} & 1+1\big).
}$$
\end{rem}\bigskip

\noindent With these $n$-ary sums we can give an alternative description of
$n$-tests, in line with how they are commonly understood, namely as
predicates that add up to 1. This generalises the bijective
correspondences that we have in the context of $C^*$-algebras between:
$$\begin{prooftree}
\begin{prooftree}
\begin{prooftree}
\mbox{$n$-tests } A \longrightarrow n\cdot 1 \mbox{ in }\op{(\CstarPU)}
\Justifies
\mbox{maps } \C^{n} \longrightarrow A \mbox{ in }\CstarPU 
\end{prooftree}
\Justifies
\mbox{maps of effect modules } [0,1]^{n} \longrightarrow [0,1]_{A} 
\end{prooftree}
\Justifies
\mbox{effects } e_{1}, \ldots, e_{n} \in [0,1]_{A} \mbox{ with }
   e_{1} \ovee \cdots \ovee e_{n} = 1
\end{prooftree}$$

\noindent Here is the general formulation.

\begin{lem}
\label{TestLem}
In a category satisfying Assumption~\ref{CoprodAss} there is a
bijective correspondence between:
\begin{equation}
\label{TestCor}
\begin{prooftree}
\mbox{$n$-tests } q \colon X \longrightarrow n\cdot 1
\Justifies
\mbox{predicates } p_{1}, \ldots, p_{n} \colon X \longrightarrow 1+1
   \mbox{ with } p_{1} \ovee \cdots \ovee p_{n} = 1
\end{prooftree}
\end{equation}
\end{lem}

\begin{myproof}
First suppose we have such a collection of predicates $p_{i} \colon X
\rightarrow 1+1$ with $p_{1} \ovee \cdots \ovee p_{n} = 1$. If $b
\colon X \rightarrow (n+1)\cdot 1 = n\cdot 1 + 1$ with
$[\rhd_{i},\kappa_{2}] \after b = p_{i}$ is the bound, like in
Remark~\ref{ManyOveeRem}, then $(\nabla+\idmap) \after b = \kappa_{1}
\after\; !_{X}$.  We proceed like in diagram~\eqref{TrueBoundDiag} in
the proof of Proposition~\ref{CoprodEAProp}:
$$\xymatrix@R-1.8pc@C-1pc{
X\ar@/_2ex/[ddddrr]_{!_X}\ar@/^2ex/[ddrrrr]^-{b}\ar@{-->}[ddrr]^(0.6){q} \\
\\
& & n\cdot 1\ar[dd]_{!}\ar[rr]^-{\kappa_1}\pullback & &
   (n+1)\cdot 1\rlap{$\;=n\cdot 1 + 1$}\ar[dd]^{\nabla+\idmap =\, !+\idmap} \\
& & & \\
& & 1\ar[rr]_-{\kappa_1} & & 1+1
}$$

\noindent Thus we obtain an $n$-test $q\colon X \rightarrow n\cdot 1$.
It satisfies:
$$\begin{array}{rcccl}
[\rhd_{i},\kappa_{2}] \after \kappa_{1} \after q
& = &
[\rhd_{i},\kappa_{2}] \after b
& = &
p_{i}.
\end{array}$$

This equation suggests how to proceed in the other direction: given an
$n$-test $q\colon X \rightarrow n\cdot 1$, define $p_{i} = \rhd_{i}
\after q = [\rhd_{i},\kappa_{2}] \after \kappa_{1} \after q \colon X
\rightarrow 1+1$. Then, by construction, $\kappa_{1} \after q \colon X
\rightarrow (n+1)\cdot 1$ is a bound for these $p_i$, with:
$$\begin{array}{rcccccccl}
p_{1} \ovee \cdots \ovee p_{n}
& = &
(\nabla+\idmap) \after \kappa_{1} \after q
& = &
\kappa_{1} \after \nabla \after q
& = &
\kappa_{1} \after\; !_{X} 
& = &
1.
\end{array}\eqno{\qEd}$$

\end{myproof}

\section{States}\label{StatesSec}

Having seen predicates as maps of the form $X\rightarrow 1+1$, we turn
to states, and describe them as maps of the form $1\rightarrow X$.
They are sometimes called points. We define validity $\models$ for
states and predicates, via an abstract Born rule, and show how these
predicates and states give rise to a state-and-effect triangle.

\begin{defi}
\label{StateDef}
In a category $\cat{B}$ with final object $1\in\cat{B}$ we define a
\emph{state} to be a map of the form $1 \rightarrow X$. We write $\Stat(X)
= \Hom(1, X)$ for the set of states.

For a map $f\colon X \rightarrow Y$ in $\cat{B}$ we get a function
$\Stat(f) = f_{*} = f \after (-) \colon \Stat(X) \rightarrow
\Stat(Y)$. This yields a functor $\Stat \colon \cat{B} \rightarrow
\Sets$.
\end{defi}

In the category $\Sets$, states of an object (set) $X$ are just
elements of $X$. But in the Kleisli category $\Kl(\Dst)$, states $1
\rightarrow X$ correspond to functions $1 \rightarrow \Dst(X)$, and
thus to distributions $\varphi\in\Dst(X)$, see Example~\ref{StateEx}
below. Similarly, states in the Kleisli category $\Kl(\Giry)$ of the
Giry monad correspond to probability distributions
$\phi\in\Giry(X)$. And states $1 \rightarrow A$ in the category
$\op{(\CstarPU)}$ are positive unital maps $A \rightarrow \C$, and
thus states as they are commonly used in the theory of $C^*$-algebras.

Typically, states are closed under convex combinations. This can also
be shown in the current context, where the scalars are of the form $1
\rightarrow 1+1$. It requires the generalised notion of convex set
with respect to an effect monoid, as described towards the end of
section~\ref{PrelimSec}.

\begin{lem}
\label{StatConvexLem}
Let $\cat{B}$ be an effectus. For each object $X\in\cat{B}$ the set of
states $\Stat(X) = \Hom(1, X)$ is a convex set over the effect monoid
$\Pred(1) = \Hom(1, 1+1)$ of scalars in $\cat{B}$.

For each map $f\colon X \rightarrow Y$ in $\cat{B}$ the associated
function $\Stat(f) = f_{*} = f \after (-) \colon \Stat(X) \rightarrow
\Stat(Y)$ is affine. Thus, taking states yields a functor $\Stat
\colon \cat{B} \rightarrow \Conv_{\Pred(1)}$.
\end{lem}

\begin{myproof}
Assume we have $n$ scalars $r_{1}, \ldots, r_{n} \colon 1 \rightarrow
1+1$ with $r_{1} \ovee \cdots \ovee r_{n} = 1$ together with $n$
states $\omega_{1}, \ldots, \omega_{n} \colon 1 \rightarrow X$. The
scalars correspond by~\eqref{TestCor} to an $n$-test $q \colon 1
\rightarrow n\cdot 1$. Then we define a convex sum of these states as:
$$\xymatrix{
\bigovee_{i} r_{i}\omega_{i} = \big(1\ar[r]^-{q} & 
   n\cdot 1\ar[rr]^-{[\omega_{1}, \ldots, \omega_{n}]} & & X\big).
}$$

\noindent More formally, we describe this definition as a function
$\alpha \colon \Dst_{M}(\Stat(X)) \rightarrow \Stat(X)$, given
by $\alpha(\sum_{i}r_{i}\ket{\omega_{i}}) = \bigovee_{i} r_{i}\omega_{i}$.
It forms an an Eilenberg-Moore algebra since:
$$\begin{array}{rcl}
\big(\alpha \after \eta\big)(\omega)
& = &
\alpha(1\ket{\omega}) \\
& = &
[\omega] \after \idmap \\
& = &
\omega \\
\big(\alpha \after \mu\big)\big(
   \sum_{i}r_{i}\ket{\sum_{j\in J_{i}}s_{ij}\ket{\omega_{ij}}}\big)
& = &
\alpha\big(\sum_{ij} r_{i}\cdot s_{ij}\ket{\omega_{ij}}\big) \\
& = &
\bigovee_{ij} (r_{i}\cdot s_{ij}) \omega_{ij} \\
& = &
\bigovee_{i}r_{i}(\bigovee_{j\in J_{i}}s_{ij}\omega_{ij}) \\
& = &
\alpha\big(\sum_{i}r_{i}\ket{
   \alpha(\sum_{j\in J_{i}}s_{ij}\ket{\omega_{ij}})}\big) \\
& = &
\big(\alpha \after \Dst_{M}(\alpha)\big)\big(
   \sum_{i}r_{i}\ket{\sum_{j\in J_{i}}s_{ij}\ket{\omega_{ij}}}\big).
\end{array}$$

\noindent These convex sums make the set $\Stat(X)$ into a convex
set. The sums are preserved by the functions $\Stat(f)$ since:
$$\begin{array}[b]{rcl}
\Stat(f)\big(\bigovee_{i} r_{i}\ket{\omega_{i}}\big)
& = &
f \after [\omega_{1}, \ldots, \omega_{n}] \after q \\
& = &
[f \after \omega_{1}, \ldots, f \after \omega_{n}] \after q
\hspace*{\arraycolsep}=\hspace*{\arraycolsep}
\bigovee_{i} r_{i}\ket{\Stat(f)(\omega_{i})}.
\end{array}\eqno{\qEd}$$
\end{myproof}\medskip

\noindent Next we define validity via composition. This definition is extremely
simple but turns out to be very powerful since it can be interpreted
in many categories and can mean many different things, depending on
the meaning of composition, see the overview in
Figure~\ref{ValidityFig}.

\begin{defi}
\label{ValidityDef}
Given a state $\omega\colon 1 \rightarrow X$ on an object $X$ in an
effectus, and a predicate $p\colon X \rightarrow 1+1$ on that same
object, we can define the \emph{validity probability} $\omega\models
p$ in $\Pred(1)$ as the scalar:
\begin{equation}
\label{ModelsEqn}
\begin{array}{rcl}
\omega\models p
& \;=\; &
p \after \omega \;\colon\; 1 \longrightarrow 1+1.
\end{array}
\end{equation}
\end{defi}\bigskip

\noindent When the homset of scalars $\Hom(1,1+1)$ is $\{0,1\}$, this expression
$\omega\models p$ is a truth value. But $\omega\models p$ may also be
a probability, when $[0,1]$ is set of scalars, as illustrated
below. The definition~\eqref{ModelsEqn} may be seen as \emph{Born's
  rule}, in most elementary shape. It takes various forms, depending
on the category in which $\models$ is interpreted. This will be
elaborated next.

\begin{exas}
\label{StateEx}
We shall see what the validity probability $\omega\models p$ amounts
to in our running examples. We encounter the formulations occurring
in~\cite[Table~1]{dHondtP06a}, see also Figure~\ref{ValidityFig}. The
idea of using integration as logical validity goes back
to~\cite{Kozen81,Kozen85}.

We first consider the Kleisli category $\Kl(\Dst)$ of the distribution
monad $\Dst$ on $\Sets$. A state $\omega$ on a set $X\in\Kl(\Dst)$ is
a function $1 \rightarrow \Dst(X)$, which corresponds to a discrete
probability distribution $\omega\in\Dst(X)$.  Thus $\Stat(X) =
\Dst(X)$ and $f_{*} = \Stat(f)$ is the Kleisli extension map $\Dst(X)
\rightarrow \Dst(Y)$, given by $f_{*}(\omega)(y) =
\sum_{x}\omega(x)\cdot f(x)(y)$, for $f\colon X \rightarrow
\Dst(Y)$. A predicate $p\colon X \rightarrow 1+1$ corresponds to
function $p\colon X \rightarrow [0,1]$. The validity probability is
obtained by Kleisli composition and looks as follows.
$$\begin{array}{rcl}
\omega\models p
& \;=\; &
\sum_{x}\omega(x)\cdot p(x) \;\in\; [0,1].
\end{array}$$

In the Kleisli category $\Kl(\Giry)$ of the Giry monad this formula
becomes integration: for a measurable space $X$, a state $\omega
\colon 1 \rightarrow X$ in $\Kl(\Giry)$ is a probability measure
$\omega\in\Giry(X)$.  A predicate $p$ on $X$ is a measurable map
$p\colon X \rightarrow [0,1]$.  Then:
\begin{equation}
\label{KlGValidityEqn}
\begin{array}{rcl}
\omega\models p
& \;=\; &
\displaystyle\int p\intd \omega \;\in\; [0,1].
\end{array}
\end{equation}
In the category $\op{\Rng}$ a state of a ring $R$ is a ring
homomorphism $\omega \colon R \rightarrow \Z$. Such a map is sometimes
called a $\Z$-point in algebraic geometry. A predicate is an
idempotent $e\in R$, see Proposition~\ref{RingPredProp}. Validity
becomes function application:
$$\begin{array}{rcl}
\omega\models e
& \;=\; &
\omega(e) \;\in\; \{0,1\}.
\end{array}$$

\noindent The outcome $\omega(e) \in \Z$ is in $\{0,1\}$ since $e$ is
an idempotent, and $\omega$ is a ring homomorphism, so that
$\omega(e)^{2} = \omega(e)$.

In the category $\op{\DL}$ of distributive lattices a state of $L$ is
a lattice homomorphism $\omega \colon L \rightarrow 2$. As is
well-known, it may be identified with a \emph{prime filter} $U =
\omega^{-1}(1) \subseteq L$, that is, with a non-empty upset $U$ that
is closed under meets and satisfies: $0\not\in U$ and $x\vee y \in U$
implies either $x\in U$ or $y\in U$. We can formulate validity
$\models$ either in terms of states $\omega\colon L \rightarrow 2$ or
in terms of prime filters $U\subseteq L$, as:
$$\begin{array}{rclcrcl}
\omega\models e
& \;=\; &
\omega(e) \;\in\{0,1\}
& \qquad\mbox{or as}\qquad &
U\models e
& \;=\; &
(e\in U).
\end{array}$$

\noindent where the predicate $e\in L$ has a complement, see
Example~\ref{PredEx}~\eqref{PredExDL}. Of course, instead of prime
filters one may equivalently use prime ideals, corresponding to
kernels $\omega^{-1}(0)$ of states $\omega$. 

The same approach works in the opposite of the category of Boolean
algebras, where a state of a Boolean algebra $B$ is a map of Boolean
algebra $B\rightarrow 2$, which corresponds to an ultrafilter of $B$.

Next, we consider the opposite $\op{(\CstarPU)}$ of the category of
$C^*$-algebras with positive unital maps. A state on a $C^*$-algebra
$A$ is, as usual, a positive unital map $\omega\colon A \rightarrow
\C$. It is a map $1 \rightarrow A$ in this opposite category
$\op{(\CstarPU)}$. For an effect/predicate $a\in [0,1]_{A} =
\setin{x}{A}{0 \leq x \leq 1}$ we obtain the validity probability
simply by function application:
$$\begin{array}{rcl}
\omega\models a
& \;=\; &
\omega(a) \;\in\; [0,1].
\end{array}$$

\noindent The outcome $\omega(a)$ is in the unit interval $[0,1] \subseteq
\R$ because: $\omega(a) \geq 0$ since $a\geq 0$ and $\omega$ is positive,
and $\omega(a) \leq 1$ since $a\leq 1$ and $\omega(1) = 1$.

We briefly look at the special case where our $C^*$-algebra $A$ is the
algebra $\B(\H)$ of bounded linear maps on a finite-dimensional
Hilbert space $\H$. We recall that the subset $\DM(\H) \hookrightarrow
\B(\H)$ of density matrices contains the positive operators with trace
equal to one. Each such density matrix $\rho\in\DM(\H)$ gives rise to
a state $\omega_{\rho} \colon \B(H) \rightarrow \C$, namely
$\omega_{\rho}(f) = \tr(f\rho)$ where $\tr$ is the trace
operation. Then, for an effect $E\in\Ef(\H) = [0,1]_{\B(\H)} =
\set{A\colon \H \rightarrow \H}{0 \leq A \leq \idmap}$ on $\H$ we have
the usual probability formula (see~\textit{e.g}~\cite{dHondtP06a}):
$$\begin{array}{rcccl}
(\omega_{\rho} \models E)
& = &
\omega_{\rho}(E)
& = &
\tr(E\rho).
\end{array}$$

\noindent If we take the spectral decomposition of the density matrix
$\rho$, as $\rho = \sum_{i} \lambda_{i}\ket{v_i}\bra{v_i}$, where the
$\lambda_{i} \in [0,1]$ are the eigenvalues and the $\ket{v_i}\in\H$
are the (orthonormal) eigenvectors, then we get the standard Born
rule:
$$\begin{array}{rcccl}
(\omega_{\rho} \models E)
& = &
\tr(E\rho)
& = &
\sum_{i} \lambda_{i} \bra{v_i}E\ket{v_i}.
\end{array}$$

\noindent What we use is the isomorphism between states and density
matrices (for finite-dimensional $\H$):
\begin{equation}
\label{HilbStatEqn}
\begin{array}{rcl}
\Stat(\B(\H)) 
& \cong &
\DM(\H).
\end{array}
\end{equation}
\end{exas}\bigskip

\noindent At this stage we emphasise how much structure we get from the
relatively mild effectus requirements in Assumption~\ref{CoprodAss}:
we do not only have state and effect functors, but, if we involve the
adjunction $\op{(\EMod_{M})} \leftrightarrows \Conv_{M}$ between
effect modules and convex sets from Proposition~\ref{ConvEModAdjProp},
we also get a ``state-and-effect'' triangle of the form:
\begin{equation}
\label{GeneralTriangleDiag}
\vcenter{\xymatrix{
\op{(\EMod_{M})}\ar@/^1.5ex/[rr]^-{\Hom(-,M)} & \top & 
   \Conv_{M}\ar@/^1.5ex/[ll]^-{\Hom(-,M)} \\
& \cat{B}\ar[ul]^{\Hom(-,1+1)=\Pred\quad}\ar[ur]_{\;\Stat=\Hom(1,-)} &
}} 
\qquad\mbox{where } M = \Pred(1) = \Stat(1+1).
\end{equation}

\noindent Such triangles are quite common in program semantics
see~\cite{Jacobs15a,Jacobs15b} for more information. They capture the
correspondence between programs, in the base category, state
transformers in the upper right category (in Schr\"odinger style), and
predicate transformers in the upper left category (in Heisenberg
style).

In general, the two triangles in~\eqref{GeneralTriangleDiag} do not
commute, but we do have canonical natural transformations between
them, given by validity $\models$.

\begin{prop}
For an effectus $\cat{B}$ there are, in the setting of
diagram~\eqref{GeneralTriangleDiag}, natural transformations:
\begin{equation}
\label{GeneralTriangleNatrosDiag}
\vcenter{\xymatrix{
\op{(\EMod_{M})} & & \Conv_{M} \\
& \cat{B}\urtwocell^{\Stat\quad}_{\hspace*{5em}\Hom(\Pred(-), M)}{\beta}
  \ultwocell_{\quad\Pred}^{\Hom(\Stat(-),M)\hspace*{5em}}{\alpha}
}}
\end{equation}

\noindent given as follows. For a predicate $p\colon X \rightarrow
1+1$ and state $\omega\colon 1 \rightarrow X$, both $\alpha$ and $\beta$
are defined via validity $\models$, as in:
$$\begin{array}{rclcrcl}
\alpha_{X}(p)(\omega)
& = &
\omega \models p
& = &
p \after \omega
& = &
\beta_{X}(\omega)(p).
\end{array}$$
\end{prop}\medskip

\noindent The maps $\alpha_{X}$ and $\beta_{X}$ are each other's transposes
in the adjunction from Proposition~\ref{ConvEModAdjProp}:
$$\begin{prooftree}
\xymatrix{\Pred(X)\ar[r]^-{\alpha_X} & \Hom(\Stat(X), M)}
\Justifies
\xymatrix{\Stat(X)\ar[r]_-{\beta_X} & \Hom(\Pred(X), M)}
\end{prooftree}$$

\noindent Note the reversal of direction for $\alpha$, due to the use
of the \emph{opposite} of the category of effect modules
in~\eqref{GeneralTriangleNatrosDiag}.  These $\alpha$ and $\beta$ may
be called the \emph{Born} natural transformations because they are
given by the Born rule~\eqref{ModelsEqn}.

\begin{myproof}
It is easy to see that $\alpha$ and $\beta$ are natural
transformations. Some basic properties have to be checked.  For
instance, each $\alpha_{X}(p) = p \after (-) = p_{*} \colon \Stat(X)
\rightarrow \Stat(1+1) = M$ is an affine map, by
Lemma~\ref{StatConvexLem}. And each $\alpha_{X} \colon \Pred(X)
\rightarrow \Hom(\Stat(X), M)$ is a map of effect modules, where the
homset inherits the effect module structure from $M = \Pred(1)$ in a
pointwise manner, see Propositions~\ref{ConvEModAdjProp}
and~\ref{CoprodEModProp}. Similarly for $\beta$. \QED

\auxproof{
For naturality of $\alpha$ and $\beta$, let $f\colon X \rightarrow Y$
be a map in $\cat{B}$. We must show that the following squares commute.
$$\xymatrix{
\Pred(X)\ar[r]^-{\alpha_X} & \Hom(\Stat(X), M)
&
\Stat(X)\ar[d]_{\Stat(f)}\ar[r]^-{\beta_X} & 
   \Hom(\Pred(X), M)\ar[d]^{(-) \after \Pred(f)}
\\
\Pred(Y)\ar[u]^{\Pred(f) = f^{*}}\ar[r]_-{\alpha_{Y}} & 
   \Hom(\Stat(Y), M)\ar[u]_{(-) \after \Stat(f)}
&
\Stat(Y)\ar[r]_-{\beta_Y} & \Hom(\Pred(Y), M)
}$$

\noindent Indeed, for $q\in\Pred(Y)$ and $\omega\in\Stat(X)$,
$$\begin{array}{rcl}
\big(((-) \after \Stat(f)) \after \alpha_{Y}\big)(q)(\omega)
& = &
\big(\alpha_{Y}(q) \after \Stat(f)\big)(\omega) \\
& = &
\alpha_{Y}(q)(f \after \omega) \\
& = &
q \after f \after \omega \\
& = &
f^{*}(q) \after \omega \\
& = &
\alpha_{X}(f^{*}(q))(\omega) \\
& = &
\big(\alpha_{X} \after f^{*}\big)(q)(\omega) \\
\big(((-) \after \Pred(f)) \after \beta_{X}\big)(\omega)(q)
& = &
\big(\beta_{X}(\omega) \after f^{*}\big)(q) \\
& = &
\beta_{X}(\omega)(f^{*}(q)) \\
& = &
f^{*}(q) \after \omega \\
& = &
q \after f \after \omega \\
& = &
q \after \Stat(f)(\omega) \\
& = &
\beta_{Y}(\Stat(f)(\omega))(q) \\
& = &
\big(\beta_{Y} \after \Stat(f)\big)(\omega)(q).
\end{array}$$

For $X\in\cat{B}$, $\alpha_{X} = \lam{p}{p \after (-)} \colon \Pred(X)
\rightarrow \Hom(\Stat(X), M)$ is a map of effect modules. For instance:
$$\begin{array}{rcl}
\alpha_{X}(p_{1} \ovee p_{2})
& = &
\lam{\omega}{(p_{1} \ovee p_{2}) \after \omega} \\
& = &
\lam{\omega}{\omega^{*}(p_{1} \ovee p_{2})} \\
& = &
\lam{\omega}{\omega^{*}(p_{1}) \ovee \omega^{*}(p_{2})} \\
& = &
\lam{\omega}{\omega^{*}(p_{1})} \ovee \lam{\omega}{\omega^{*}(p_{2})} \\
& = &
\alpha_{X}(p_{1}) \ovee \alpha_{X}(p_{2})
\end{array}$$
}
\end{myproof}

The triangle~\eqref{GeneralTriangleDiag} commutes up-to-isomorphism
when these maps $\alpha, \beta$ from~\eqref{GeneralTriangleNatrosDiag}
are both isomorphisms. This requires further assumptions, which we
briefly illustrate (but do not impose).  In our two main examples
(with $M = [0,1]$) we do \emph{not} have such commutation:
$$\xymatrix{
\op{\EMod}\ar@/^1.5ex/[rr]^-{\Hom(-,[0,1])} & \top & 
   \Conv\ar@/^1.5ex/[ll]^-{\Hom(-,[0,1])}
& &
\op{\EMod}\ar@/^1.5ex/[rr]^-{\Hom(-,[0,1])} & \top & 
   \Conv\ar@/^1.5ex/[ll]^-{\Hom(-,[0,1])} \\
& \Kl(\Dst)\ar[ul]^{\Pred}\ar[ur]_{\Stat} &
& &
& \op{(\CstarPU)}\ar[ul]^{\Pred}\ar[ur]_{\Stat} &
}$$

\noindent In the triangle on the left we do have:
$$\begin{array}{rcccccl}
\Hom(\Stat(X), [0,1])
& = &
\Conv(\Dst(X), [0,1])
& \cong &
[0,1]^{X}
& = &
\Pred(X).
\end{array}$$

\noindent But in the other direction we have $\Hom(\Pred(X), [0,1])
\cong \Dst(X)$ when $X$ is finite. For details, see~\cite{JacobsM12b},
where the mapping $X \mapsto \Hom(\Pred(X), [0,1]) = \EMod([0,1]^{X},
[0,1])$ is called the expectation monad $\Exp$. This expectation monad
is a monad on $\Sets$ that extends the distribution monad, in the
sense that $\Exp(X) \cong \Dst(X)$ for finite $X$.

For the $C^*$-algebra example on the right we do have commutation in
the other direction:
$$\begin{array}{rcccccl}
\Hom(\Pred(A), [0,1])
& = &
\EMod([0,1]_{A}, [0,1]_{\C})
& \cong &
\CstarPU(A, \C)
& = &
\Stat(A),
\end{array}$$

\noindent using that the effect functor $[0,1]_{(-)} \colon \CstarPU
\rightarrow \EMod$ is full and faithful, see~\cite{FurberJ13a}. In the
other direction we need to involve the (compact Hausdorff) topology on
states of a $C^*$-algebra. Via the validity relation $\models$ one can
put a weak-*-like topology on the homset of states. The restriction
to convex \emph{compact Hausdorff} spaces gives an equivalence with
suitably complete effect modules, which is sometimes called Kadison
duality after~\cite{Kadison51}, see also~\cite{JacobsM12b}.

A natural question is to which extent effects and states determine
each other.

\begin{defi}
\label{SeparationDef}
In an effectus $\cat{B}$ we say:
\begin{enumerate}
\item \emph{predicates can be separated} if for each pair of predicates
$p_{1}, p_{2} \colon X \rightarrow 1+1$,
$$\begin{array}{rcl}
\Big(\all{\omega\colon 1 \rightarrow X}
   {(\omega\models p_{1}) = (\omega\models p_{2})}\Big)
& \Longrightarrow &
p_{1} = p_{2}.
\end{array}$$

\noindent Equivalently, if $p_{1}\neq p_{2}$, then there is a state
$\omega$ with $(\omega\models p_{1}) \neq (\omega\models p_{2})$.
This expresses that the Born natural transformation $\alpha$ 
in~\eqref{GeneralTriangleNatrosDiag} is injective.

\item \emph{states can be separated} if for each pair of states
$\omega_{1}, \omega_{2} \colon 1 \rightarrow X$,
$$\begin{array}{rcl}
\Big(\all{p\colon X \rightarrow 1+1}
   {(\omega_{1}\models p) = (\omega_{2}\models p)}\Big)
& \Longrightarrow &
\omega_{1} = \omega_{2}.
\end{array}$$

\noindent This says that $\beta$ in~\eqref{GeneralTriangleNatrosDiag} is
injective.
\end{enumerate}
\end{defi}

\noindent Separation of predicates is the more interesting property. For
instance for distributive lattices it says that for (complementable)
elements $a,b$ in a distributive lattice with $a\neq b$, there is a
prime filter (or ideal) that contains one, but not the other.  This
separation property holds, see~\cite[I, Prop.~2.5]{Johnstone82} using
the Prime Ideal Theorem.  

In $C^*$-algebras predicates can also be separated: let $a, b \in
[0,1]_{A}$ be effects in a $C^{*}$-algebra $A$ with $a \neq b$. Then
$a-b$ is a non-zero self-adjoint element. Hence
by~\cite[Thm.~7.1.2]{Arveson81} there is a state $\omega \colon A
\rightarrow \C$ with $|\omega(a-b)| = \|a - b\|$. From $a - b \neq 0$
we get $\|a - b\| \neq 0$, and thus $\omega(a) \neq \omega(b)$.


Separation of states is an easier property. For instance, if we have
two states $\omega_{1}, \omega_{2}\in \Giry(X)$ in the Kleisli
category $\Kl(\Giry)$ of the Giry monad with $(\omega_{1}\models p) =
(\omega_{2}\models p)$ for each predicate $p$, then we take for $p$
indicator functions $\indic{M} \colon X \rightarrow [0,1]$ for
$M\in\Sigma_{X}$. We get by~\eqref{KlGValidityEqn}:
$$\begin{array}{rcccccccccl}
\omega_{1}(M)
& = &
\int \indic{M}\intd \omega_{1}
& = &
(\omega_{1}\models\indic{M})
& = &
(\omega_{2}\models\indic{M})
& = &
\int \indic{M}\intd \omega_{2}
& = &
\omega_{2}(M).
\end{array}$$

\noindent Hence $\omega_{1}=\omega_{2}$. Also, states can be separated
in $C^*$-algebras: let $\omega_{1}, \omega_{2} \colon A \rightarrow
\C$ be two states of a $\C$-algebra $A$ with $\omega_{1} \neq
\omega_{2}$. Then there is an element $a\in A$ with $\omega_{1}(a)
\neq \omega_{2}(a)$. Since each element in a $C^*$-algebra can be
written as a linear combination of (four) positive elements --- see
also the proof of Lemma~\ref{EffectCommutationLem} later on --- we may
assume $a$ is positive (and $a\neq 0$). We take $e =
\frac{1}{\|a\|}\cdot a$, which is an effect since $0 \leq e \leq
\|e\|\cdot 1 = 1$. Moreover,
$$\begin{array}{rcccccccccl}
(\omega_{1}\models e)
& = &
\omega_{1}(e)
& = &
\frac{1}{\|a\|}\cdot\omega_{1}(a)
& \neq &
\frac{1}{\|a\|}\cdot\omega_{2}(a)
& = &
\omega_{2}(e)
& = &
(\omega_{2}\models e).
\end{array}$$

\noindent Separation of predicates and/or states is thus a sensible
additional requirement to add in axiomatisation.

\subsection{Programs, states, and predicates}\label{StatProgrPredSubsec}

Let's take a closer look at the state-and-effect
triangle~\eqref{GeneralTriangleDiag} of an effectus $\cat{B}$. We use
the following interpretations of morphisms in $\cat{B}$.
\begin{equation}
\label{StatProgrPredEqn}
\left\{{\renewcommand\arraystretch{1}\begin{array}{ll}
\mbox{states} & \omega \colon 1 \longrightarrow X \\
\mbox{programs} & f \colon X \longrightarrow Y \\
\mbox{predicates} \quad & q \colon Y \longrightarrow 1+1
\end{array}}\right.
\end{equation}

\noindent Such a program, also called a computation, $f\colon X
\rightarrow Y$ yields two `transformer' maps:
$$\left\{{\renewcommand\arraystretch{1}\begin{array}{ll}
\mbox{state transformer} & f_{*} = f \after (-) \colon \Stat(X) 
   \longrightarrow \Stat(Y) \\
\mbox{predicate transformer} \quad & f^{*} = (-) \after f \colon
   \Pred(Y) \longrightarrow \Pred(X) 
\end{array}}\right.$$

\noindent In a quantum setting the state transformer $f_*$ captures
Schr\"odinger's style of forward computation, whereas the predicate
transformer $f^*$, going backwards, corresponds to Heisenberg's style.
It turns a predicate on the codomain of $f$ into a predicate on the
domain of $f$, acting like a weakest precondition operation $\wp(f)$.


Notice that the scalars $\omega\models f^{*}(q)$ and $f_{*}(\omega)
\models q$ are the same, giving a Galois/adjoint style correspondence:
\begin{equation}
\label{PredStatProbEqn}
\xymatrix@C-.5pc{
\Big(\omega\models f^{*}(q)\Big) = \Big(1\ar[r]^-{\omega} &
   X\ar[r]^-{f} & Y\ar[r]^-{q} & 1+1\Big) = \Big(f_{*}(\omega) \models q\Big).
}
\end{equation}

\noindent These scalars $(f_{*}(\omega) \models q) = (\omega \models
f^{*}(q))$ are interpreted as the expected probability that predicate
$q$ holds after running program $f$ in state $\omega$. This Galois
correspondence~\eqref{PredStatProbEqn} is described more concretely
for $\Kl(\Giry)$ in~\cite[Prop.~6]{Jacobs13a}, where it relates
Kleisli composition and substitution for $\Giry$. This correspondence
is reminiscent of the `satisfaction condition' of institutions,
see~\cite{GoguenB92}.

When the base category $\cat{B}$ in the
triangle~\eqref{GeneralTriangleDiag} is $\op{(\CstarPU)}$ the
interpretation of the direction of computations is more subtle,
because the category $\op{(\CstarPU)}$ already captures Heisenberg's
picture of quantum computation --- also known as matrix mechanics. Its
restriction to observables may be understood as a quantum version of
weakest precondition calculation, see~\cite{dHondtP06a}. State
transformer semantics for Hilbert spaces, considered as $C^*$-algebras
via their bounded maps $\B(-)$, can be described as transformations of
density matrices, via~\eqref{HilbStatEqn}.

In the remainder of this paper we shall describe several constructions
on programs and on predicates. Ideally, these constructions should
come with appropriate calculation rules that break down the
constructions to calculations with scalars/probabilities. For instance
we already know that substitution maps $f^{*}$ preserve the effect
module structure.  In particular, $(\omega \models (-)) = \omega^{*}
\colon \Pred(X) \rightarrow \Pred(1)$ is a map of effect
modules. Hence we can reduce validity probabilities $\models$
involving $\ovee$ and $(-)^{\perp}$ on predicates to $\ovee$ and
$(-)^{\perp}$ on scalars, as in:
$$\begin{array}{rcl}
\big(\omega\models f^{*}(q_{1}\ovee q_{2})\big)
& = &
\big(\omega\models f^{*}(q_{1})\big) \ovee
   \big(\omega\models f^{*}(q_{2})\big) \\
\big(\omega \models f^{*}(q^{\perp})\big)
& = &
\big(\omega \models f^{*}(q)\big)^{\perp}.
\end{array}$$

\section{Predicates and coproducts}\label{PredicateCoprodSec}

This section investigates the interaction in an effectus between
coproducts, as in Assumption~\ref{CoprodAss}, and predicates. The
elements of an effect module may be seen as formulas in a logic.
There are truth and falsum formulas $1,0$, and formulas are closed
under sum $\ovee$ and orthocomplement $(-)^{\perp}$. In the current
categorical context where we have effect modules $\Pred(X)$ of
predicates $X\rightarrow 1+1$ on an object $X$ in a category, we can
form additional formulas. Notably, if we have two predicates $p\colon
X \rightarrow 1+1$ and $q\colon Y \rightarrow 1+1$ we can form a
cotuple predicate $[p,q] \colon X+Y \rightarrow 1+1$. This cotuple
looks like a new predicate constructor. But the
isomorphism~\eqref{CotupleIsoPredDiag} in the next result shows that
we get nothing new.

\begin{lem}
\label{CotupleIsoPredLem}
Let $\cat{B}$ be an effectus, with its predicate functor $\Pred \colon
\cat{B} \rightarrow \op{(\EMod_{M})}$ from
Proposition~\ref{CoprodEModProp}, where $M = \Pred(1) = \Hom(1, 1+1)$
is the effect monoid of scalars. This functor $\Pred$ preserves finite
coproducts and the final object $1$, \textit{i.e.}~it sends:
\begin{itemize}
\item coproducts $(+,0)$ in $\cat{B}$ to products $(\times, 1)$ in
  $\EMod_{M}$;

\item the final object $1\in\cat{B}$ to the initial object $M = \Pred(1)$
in $\EMod_M$.
\end{itemize}

\noindent In particular, for objects $X,Y\in\cat{B}$ there is an
isomorphism of effect modules:
\begin{equation}
\label{CotupleIsoPredDiag}
\xymatrix{
\Pred(X+Y)\ar@/^2ex/[rr]^-{\tuple{\kappa_{1}^{*}, \kappa_{2}^{*}}} & \cong &
   \Pred(X)\times\rlap{$\Pred(Y)$}\ar@/^2ex/[ll]^{[-,-]}
}
\end{equation}

\noindent On the right-hand-side we use the obvious componentwise
effect module structure. This isomorphism is natural in $X,Y$.
\end{lem}

\begin{myproof}
Clearly, $\Pred(0) = \Hom(0,1+1)$ is a singleton set, which,
understood as trivial effect module, is final in the category
$\EMod_M$. The two functions $\tuple{\kappa_{1}^{*}, \kappa_{2}^{*}}$
and $[-,-]$ in~\eqref{CotupleIsoPredDiag} are each other's
inverses. We have already seen in Proposition~\ref{CoprodEModProp}
that substitution functors --- in particular these $\kappa_{i}^*$ ---
preserve the effect module structure. So we only have to prove that
the cotuple map $[-,-] \colon \Pred(X)\times\Pred(Y) \rightarrow
\Pred(X+Y)$ preserves the effect module structure.

\auxproof{
$$\begin{array}{rcl}
\big(\tuple{\kappa_{1}^{*}, \kappa_{2}^{*}} \after [-,-]\big)(p,q)
& = &
\tuple{\kappa_{1}^{*}, \kappa_{2}^{*}}([p,q]) \\
& = &
(\kappa_{1}^{*}([p,q]), \kappa_{2}^{*}([p,q])) \\
& = &
([p,q] \after \kappa_{1}, [p,q] \after \kappa_{2}) \\
& = &
(p,q) \\
\big([-,-] \after \tuple{\kappa_{1}^{*}, \kappa_{2}^{*}}\big)(r)
& = &
[\kappa_{1}^{*}(r), \kappa_{2}^{*}(r)] \\
& = &
[r\after \kappa_{1}, r \after \kappa_{2}] \\
& = &
r\after [\kappa_{1}, \kappa_{2}] \\
& = &
r.
\end{array}$$
}

Obviously, $[-,-]$ preserves truth:
$$\begin{array}{rcccccccl}
[1,1]
& = &
[\kappa_{1} \after\; !_{X}, \kappa_{1} \after\; !_{Y}]
& = &
\kappa_{1} \after [\, !_{X}, \, !_{Y}] 
& = &
\kappa_{1} \after\; !_{X+Y} 
& = &
1.
\end{array}$$

\noindent If we have $p_{1} \orthogonal p_{2}$ for $p_{i}\in\Pred(X)$
via bound $b$, and $q_{1} \orthogonal q_{2}$ for $q_{i}\in\Pred(Y)$
via bound $c$, then one easily checks that $[b,c] \colon X+Y
\rightarrow (1+1)+1$ proves $[p_{1}, q_{1}] \orthogonal [p_{2},
  q_{2}]$ and also $[p_{1}, q_{1}] \ovee [p_{2}, q_{2}] = [p_{1} \ovee
  p_{2}, q_{1} \ovee q_{2}]$. Preservation of scalar multiplication
from~\eqref{ScalarEqn} is obvious:
$$\begin{array}{rcccccl}
s \scalar [p,q]
& = &
[s, \kappa_{2}] \after [p,q]
& = &
[[s, \kappa_{2}] \after p, [s, \kappa_{2}] \after q]
& = &
[s\scalar p, s\scalar q].
\end{array}$$

\auxproof{
$$\begin{array}{rcl}
[\idmap,\kappa_{2}] \after [b,c]
& = &
[[\idmap,\kappa_{2}] \after b, [\idmap,\kappa_{2}] \after c] \\
& = &
[p_{1}, q_{1}] \\
{[[\kappa_{2},\kappa_{1}],\kappa_{2}]} \after [b,c]
& = &
[[[\kappa_{2},\kappa_{1}],\kappa_{2}] \after b, [\idmap,\kappa_{2}] \after c] \\
& = &
[p_{2}, q_{2}] \\
{[p_{1}, q_{1}] \ovee [p_{2}, q_{2}]}
& = &
(\nabla+\idmap) \after [b,c] \\
& = &
[(\nabla+\idmap) \after b, (\nabla+\idmap) \after c] \\
& = &
[p_{1} \ovee p_{2}, q_{1} \ovee q_{2}].
\end{array}$$
}

\noindent Naturality of the isomorphism in~\eqref{CotupleIsoPredDiag}
is easy and left to the reader. Finally, $\Pred(1) = M$ is clearly an
initial object in $\EMod_M$, because for each $E\in\EMod_M$ there is a
unique map $M \rightarrow E$ of effect modules, namely $r \mapsto r
\scalar 1$. \QED

\auxproof{
For $f\colon X' \rightarrow X$ and $g\colon Y'\rightarrow Y$,
$$\begin{array}{rcl}
\tuple{\kappa_{1}^{*}, \kappa_{2}^{*}} \after (f+g)^{*}
& = &
\tuple{\kappa_{1}^{*} \after (f+g)^{*}, \kappa_{2}^{*} \after (f+g)^{*}} \\
& = &
\tuple{((f+g) \after \kappa_{1})^{*}, ((f+g) \after \kappa_{2})^{*}} \\
& = &
\tuple{(\kappa_{1} \after f)^{*}, (\kappa_{2} \after g)^{*}} \\
& = &
\tuple{f)^{*} \after \kappa_{1}^{*}, g^{*} \after \kappa_{2}^{*}} \\
& = &
(f^{*}\times g^{*}) \after \tuple{\kappa_{1}^{*}, \kappa_{2}^{*}}.
\end{array}$$
}
\end{myproof}

The isomorphism~\eqref{CotupleIsoPredDiag} is known from many
instances.
\begin{itemize}
\item In $\Sets$ it amounts to $\{0,1\}^{X+Y} \cong \{0,1\}^{X} \times
\{0,1\}^{Y}$, which, in terms of powersets instead of characteristic
functions becomes: $\Pow(X+Y) \cong \Pow(X)\times\Pow(Y)$.

\item In $\Kl(\Dst)$ it means $[0,1]^{X+Y} \cong [0,1]^{X} \times [0,1]^{Y}$.

\item And for $C^*$-algebras $A,B$ it says $[0,1]_{A\oplus B} \cong
  [0,1]_{A} \times [0,1]_{B}$.
\end{itemize}

\noindent Because $\Pred \colon \cat{B} \rightarrow \op{(\EMod_{M})}$
preserves finite coproducts we have isomorphisms of effect modules:
$$\begin{array}{rcl}
\Pred(X_{1}+\cdots+X_{n}) 
& \cong &
\Pred(X_{1}) \times \cdots \times \Pred(X_{n}).
\end{array}$$

\noindent In particular, the right-to-left cotuple is a map of effect
modules. This yields a useful equation:
\begin{equation}
\label{NaryCotupleOveeEqn}
\begin{array}{rcl}
[p_{1}, \ldots, p_{n}] \ovee [q_{1}, \ldots, q_{n}]
& = &
[p_{1}\ovee q_{1}, \ldots, p_{n}\ovee q_{n}].
\end{array}
\end{equation}

\auxproof{
Let $b_{i} \colon X_{i} \rightarrow (1+1)+1$ be a bound for $p_{i}
\orthogonal q_{i}$. We take $b = [b_{1}, \ldots, b_{n}] \colon
X_{1} + \cdots + X_{n} \rightarrow (1+1)+1$. It satisfies:
$$\begin{array}{rcl}
[\idmap,\kappa_{2}] \after b
& = &
[[\idmap,\kappa_{2}] \after b_{1}, \ldots, [\idmap,\kappa_{2}] \after b_{n}] \\
& = &
[p_{1}, \ldots, p_{n}] \\
{[[\kappa_{2},\kappa_{1}],\kappa_{2}]} \after b
& = &
[[[\kappa_{2},\kappa_{1}],\kappa_{2}] \after b_{1}, \ldots, 
    [\idmap,\kappa_{2}] \after b_{n}] \\
& = &
[q_{1}, \ldots, q_{n}] \\
{[p_{1}, \ldots, p_{n}] \ovee [q_{1}, \ldots, q_{n}]}
& = &
(\nabla+\idmap) \after b \\
& = &
[(\nabla+\idmap) \after b_{1}, \ldots, (\nabla+\idmap) \after b_{n}] \\
& = &
[p_{1} \ovee q_{1}, \ldots, p_{n} \ovee q_{n}].
\end{array}$$
}

\noindent As immediate consequence we get $[p,q] = [p,0] \ovee [0,q]$.
The predicates $[p,0]$ and $[0,q]$ occurring in this equation play a
special role. Informally, the predicate $[p,0]$ says that an element
$z\colon X+Y$ must be in the first component $X$ of the coproduct
$X+Y$ and in that case $p$ holds for $z$. The predicate $[p,1]$
involves an implication: \emph{if} $z$ is in the first component $X$,
\emph{then} $p$ holds for $z$. We define special operators:
$$\begin{array}{rclcrcl}
\FstAnd(p)
& = &
[p,0]
& \qquad\qquad &
\FstThen(p)
& = &
[p,1].
\end{array}$$

\noindent Similarly we have:
$$\begin{array}{rclcrcl}
\SndAnd(q)
& = &
[0,q]
& \qquad\qquad &
\SndThen(p)
& = &
[1,q].
\end{array}$$

\noindent Lemma~\ref{CotupleIsoPredLem} implies $[p,q] = \FstAnd(p)
\ovee \SndAnd(q)$. It is easy to see that the ``and'' and ``then''
versions satisfy the De Morgan equalities:
\begin{equation}
\label{FstSndDeMorganEqn}
\begin{array}{rclcrcl}
\FstAnd(p^{\perp})^{\perp}
& = &
\FstThen(p)
& \qquad\mbox{and}\qquad &
\SndAnd(p^{\perp})^{\perp}
& = &
\SndThen(p).
\end{array}
\end{equation}

\auxproof{
$$\begin{array}{rcl}
\FstAnd(p^{\perp})^{\perp}
& = &
[\kappa_{2},\kappa_{1}] \after [[\kappa_{2},\kappa_{1}] \after p, 0] \\
& = &
[p, [\kappa_{2},\kappa_{1}] \after 0] \\
& = &
[p, 1] \\
& = &
\FstThen(p)
\end{array}$$
}

\noindent The relevance of these operators lies in the following
result. It states that substitution functors $\kappa_{i}^*$ for
coprojections automatically have adjoints. They will be useful later
on in Section~\ref{TestSec}.

\begin{lem}
\label{CasesAdjLem}
Let category $\cat{B}$ satisfy Assumption~\ref{CoprodAss}. Then there
are (order) adjunctions:
\begin{equation}
\label{CaseAdjDiag}
\begin{array}{ccc}
\FstAnd \dashv \kappa_{1}^{*} \dashv \FstThen
& \qquad\mbox{and}\qquad &
\SndAnd \dashv \kappa_{2}^{*} \dashv \SndThen
\end{array}
\end{equation}

\noindent in a situation:
$$\xymatrix{
\Pred(X)\ar@/_3ex/[rr]_-{\FstThen = [-,1]}^-{\bot}
   \ar@/^3ex/[rr]^-{\FstAnd = [-,0]}_-{\bot} & &
\Pred(X+Y)\ar[ll]|(0.5){\;\kappa_{1}^*\,}\ar[rr]|(0.5){\;\kappa_{2}^*\,} & & 
   \Pred(Y)\ar@/_3ex/[ll]_-{\SndAnd = [0,-]}^-{\bot}
   \ar@/^3ex/[ll]^-{\SndThen = [1,-]}_-{\bot}
}$$

\noindent These adjoints commute with substitution, in the sense that
for each map $f\colon Z \rightarrow X$ in $\cat{B}$ the following
diagrams commute.
$$\xymatrix@C+.5pc{
\Pred(X)\ar[r]^-{\FstAnd}\ar[d]_{f^*} & \Pred(X+Y)\ar[d]^{(f+\idmap)^*}
& &
\Pred(X)\ar[r]^-{\FstThen}\ar[d]_{f^*} & \Pred(X+Y)\ar[d]^{(f+\idmap)^*}
\\
\Pred(Z)\ar[r]_-{\FstAnd} & \Pred(Z+Y)
& &
\Pred(Z)\ar[r]_-{\FstThen} & \Pred(Z+Y)
}$$

\noindent (And similarly for Snd.) In the language of fibrations ---
see~\cite{Jacobs99a} --- this means that these adjoints form a
\emph{fibred} adjunction.
\end{lem}

The adjoints in~\eqref{CaseAdjDiag} are \emph{order} adjoints (or
Galois connections). In particular, the adjoints to $\kappa_{i}^*$ are
not themselves maps of effect modules, but only monotone maps.  Their
properties will be described in greater detail in
Lemma~\ref{CasesPropertiesLem} below. In the categorical logic left
and right adjoints to substitution are often written as $\coprod$ and
$\prod$ respectively, see~\cite{Jacobs99a}. Hence one may also write
$\FstAnd = \coprod_{\kappa_{1}}$ and $\FstThen = \prod_{\kappa_{1}}$,
and similarly for the second coprojection $\kappa_2$.

\begin{myproof}
We prove that there is a bijective correspondence:
$$\begin{prooftree}
\textstyle\FstAnd(p) \;\leq\; q
\Justifies
p \;\leq \kappa_{1}^{*}(q)
\end{prooftree}$$

\noindent where $q$ is a predicate on $X+Y$.

First, assume $\FstAnd(p) \leq q$, say via $\FstAnd(p) \ovee s = q$,
where $s\in\Pred(X+Y)$. Let $b\colon X+Y \rightarrow (1+1)+1$ be a
bound for the sum $\FstAnd(p) \ovee s$, so that $\IV \after b =
\FstAnd(p)$, $\XI \after b = s$, and $(\nabla+\idmap) \after b =
\FstAnd(p) \ovee s = q$. We take $c = b \after \kappa_{1} \colon X
\rightarrow (1+1)+1$. Then one easily checks that $c$ is a bound
showing $p \orthogonal \kappa_{1}^{*}(s)$ and $p \ovee
\kappa_{1}^{*}(s) = \kappa_{1}^{*}(q)$. Hence $p \leq
\kappa_{1}^{*}(q)$.

\auxproof{
$$\begin{array}{rcl}
[\idmap,\kappa_{2}] \after c
& = &
[\idmap,\kappa_{2}] \after b \after \kappa_{1} \\
& = &
\FstAnd(p) \after \kappa_{1} \\
& = &
[p, \kappa_{2} \after\; !_{Y}] \after \kappa_{1} \\
& = &
p \\
{[[\kappa_{2},\kappa_{1}],\kappa_{2}]} \after c
& = &
[[\kappa_{2},\kappa_{1}],\kappa_{2}] \after b \after \kappa_{1} \\
& = &
s \after \kappa_{1} \\
& = &
\kappa_{1}^{*}(s) \\
(\nabla+\idmap) \after c
& = &
(\nabla+\idmap) \after b \after \kappa_{1} \\
& = &
q \after \kappa_{1} \\
& = &
\kappa_{1}^{*}(q).
\end{array}$$
}

In the other direction, assume $p \leq \kappa_{1}^{*}(q)$, via $p
\ovee s = \kappa_{1}^{*}(q)$ with bound $b \colon X \rightarrow
(1+1)+1$. We now take $c = [b, (\kappa_{2}+\idmap) \after q \after
  \kappa_{2}] \colon X+Y \rightarrow (1+1)+1$. Then $c$ proves
$\FstAnd(p) \orthogonal t$, where $t = [s, q \after \kappa_{2}] \colon
X+Y \rightarrow 1+1$, and $\FstAnd(p) \ovee t = q$. Hence $\FstAnd(p)
\leq q$.

\auxproof{
$$\begin{array}{rcl}
[\idmap, \kappa_{2}] \after c
& = &
[\idmap, \kappa_{2}] \after [b, (\kappa_{2}+\idmap) \after q \after \kappa_{2}] \\
& = &
[[\idmap, \kappa_{2}] \after b, 
   [\kappa_{2}, \kappa_{2}] \after q \after \kappa_{2}] \\
& = &
[p, \kappa_{2} \after \nabla \after q \after \kappa_{2}] \\
& = &
[p, \kappa_{2} \after\; !_{Y}] \\
& = &
\FstAnd(p) 
\\
{[[\kappa_{2},\kappa_{1}], \kappa_{2}]} \after c
& = &
[[\kappa_{2},\kappa_{1}], \kappa_{2}] \after 
   [b, (\kappa_{2}+\idmap) \after q \after \kappa_{2}] \\
& = &
[[[\kappa_{2},\kappa_{1}], \kappa_{2}] \after b,
   [\kappa_{1}, \kappa_{2}] \after q \after \kappa_{2}] \\
& = &
[s, q \after \kappa_{2}] \\
& = &
t 
\\
(\nabla+\idmap) \after c
& = &
(\nabla+\idmap) \after [b, (\kappa_{2}+\idmap) \after q \after \kappa_{2}] \\
& = &
[(\nabla+\idmap) \after b, (\idmap+\idmap) \after q \after \kappa_{2}] \\
& = &
[q \after \kappa_{1}, q \after \kappa_{2}] \\
& = &
q.
\end{array}$$
}

Having established the adjunction $\FstAnd \dashv \kappa_{1}^{*}$ we use
the De Morgan equation to prove $\kappa_{1}^{*} \dashv \FstThen$ in:
$$\begin{prooftree}
\begin{prooftree}
q \;\leq\; \textstyle\FstThen(p) \rlap{$\; = \FstAnd(p^{\perp})^{\perp}$} 
\Justifies
\textstyle\FstAnd(p^{\perp}) \;\leq\; q^{\perp}
\end{prooftree}
\Justifies
\begin{prooftree}
p^{\perp} \;\leq\; \kappa_{1}^{*}(q^{\perp}) \rlap{$\;=\kappa_{1}^{*}(q)^{\perp}$}
\Justifies
\kappa_{1}^{*}(q) \;\leq\; p
\end{prooftree}
\end{prooftree}$$

\noindent Finally, for $f\colon Z \rightarrow X$ one has:
$$\begin{array}{rcccccccl}
(f+\idmap)^{*}(\FstAnd(p))
& = &
[p, 0] \after (f+\idmap) 
& = &
[p \after f, 0]
& = &
[f^{*}(p), 0]
& = &
\FstAnd(f^{*}(p)).
\end{array}\eqno{\qEd}$$\smallskip
\end{myproof}

\noindent There is a bit more to say about these adjoints $\FstAnd$ and
$\FstThen$.  We restrict ourselves to properties of $\FstAnd$, because
the corresponding properties for $\FstThen$ are easily obtained since
$\FstThen$ can be expressed as the De Morgan dual of $\FstAnd$,
see~\eqref{FstSndDeMorganEqn}.

\newpage
\begin{lem}
\label{CasesPropertiesLem}
The adjoint $\FstAnd \dashv \kappa_{1}^{*}$ from Lemma~\ref{CasesAdjLem}
satisfies:
\begin{enumerate}
\item $\kappa_{1}^{*} \after \FstAnd = \idmap$;

\item $\FstAnd$ preserves and reflects the order: $p\leq p'$ iff
  $\FstAnd(p) \leq \FstAnd(p')$;

\item \label{CasesPropertiesLemZero} $\FstAnd(0) = 0$;

\item \label{CasesPropertiesLemOvee} $\FstAnd(p_{1}\ovee p_{2}) =
\FstAnd(p_{1}) \ovee \FstAnd(p_{2})$, when $p_{1} \orthogonal p_{2}$;

\item \label{CasesPropertiesLemScalar} $\FstAnd(s\scalar p) = s\scalar
  \FstAnd(p)$ where $\scalar$ is scalar multiplication
  from~\eqref{ScalarEqn}.
\end{enumerate}
\end{lem}

\begin{myproof}
The results follow from the isomorphism $\Pred(X+Y) \cong
\Pred(X)\times\Pred(Y)$ in Lemma~\ref{CotupleIsoPredLem}.
\begin{enumerate}
\item Directly from the isomorphism~\ref{CotupleIsoPredDiag}.

\auxproof{
For a predicate $p\colon X \rightarrow 1+1$ we have:
$$\begin{array}{rcccccl}
\kappa_{1}^{*}(\FstAnd(p))
& = &
\FstAnd(p) \after \kappa_{1}
& = &
[p, \kappa_{2} \after\; !_{X}] \after \kappa_{1}
& = &
p.
\end{array}$$
}

\item Idem, since the isomorphism preserves and reflects the order, so:
$p \leq p'$ iff $[p, 0] \leq [p', 0]$.

\auxproof{
Let $p \leq p'$; the unit of the adjunction gives $p \leq p' \leq
\kappa_{1}^{*}\FstAnd(p')$, and thus $\FstAnd(p) \leq \FstAnd(p')$.
For the other direction, use the previous point.
}

\item Because~\eqref{CotupleIsoPredDiag} is an isomorphism of effect
modules. Explicitly: $\FstAnd(0_{X}) = [0_{X}, 0_{Y}] = 0_{X+Y}$.

\auxproof{
We calculate in $\Pred(X+Y)$,
$$\begin{array}{rcccccccccl}
\FstAnd(0)
& = &
[0, 0]
& = &
[\kappa_{2} \after\; !_{X}, \kappa_{2} \after\; !_{Y}]
& = &
\kappa_{2} \after [!_{X}, \, !_{Y}] 
& = &
\kappa_{2} \after\; !_{X+Y}
& = &
0.
\end{array}$$
}

\item Similarly: $\FstAnd(p_{1} \ovee p_{2}) = [p_{1} \ovee p_{2}, 0] =
[p_{1} \ovee p_{2}, 0 \ovee 0] = [p_{1}, 0] \ovee [p_{2}, 0] = \FstAnd(p_{1})
\ovee \FstAnd(p_{2})$.

\auxproof{
Assume $p_{1},p_{2}\in\Pred(X)$ are orthogonal via bound
  $b\colon X \rightarrow (1+1)+1$. Then $\FstAnd(p_{1})$ and
  $\FstAnd(p_{2})$ in $\Pred(X+Y)$ are orthogonal via bound $c = [b,
  \kappa_{2} \after\; !_{Y}] \colon X+Y \rightarrow (1+1)+1$, giving
  $\FstAnd(p_{1}) \ovee \FstAnd(p_{2}) = \FstAnd(p_{1}\ovee p_{2})$.
$$\begin{array}{rcl}
[\idmap,\kappa_{2}] \after c
& = &
[\idmap,\kappa_{2}] \after [b, \kappa_{2} \after\; !_{Y}] \\
& = &
[[\idmap,\kappa_{2}] \after b, \kappa_{2} \after\; !_{Y}] \\
& = &
[p_{1}, \kappa_{2} \after\; !_{Y}] \\
& = &
\FstAnd(p_{1}) \\
{[[\kappa_{2},\kappa_{1}],\kappa_{2}]} \after c
& = &
[[\kappa_{2},\kappa_{1}],\kappa_{2}] \after [b, \kappa_{2} \after\; !_{Y}] \\
& = &
[[[\kappa_{2},\kappa_{1}],\kappa_{2}] \after b, \kappa_{2} \after\; !_{Y}] \\
& = &
[p_{2}, \kappa_{2} \after\; !_{Y}] \\
& = &
\FstAnd(p_{2}) \\
\FstAnd(p_{1}) \ovee \FstAnd(p_{2})
& = &
(\nabla+\idmap) \after c \\
& = &
(\nabla+\idmap) \after [b, \kappa_{2} \after\; !_{Y}] \\
& = &
[p_{1}\ovee p_{2}, \kappa_{2} \after\; !_{Y}] \\
& = &
\FstAnd(p_{1}\ovee p_{2}).
\end{array}$$
}

\item Idem, using that $s \scalar 0 = 0$. \QED

\auxproof{
For a scalar $s\colon 1\rightarrow 1+1$ and a predicate $p\colon
  X \rightarrow 1+1$ we have:
$$\begin{array}[b]{rcl}
\FstAnd(s\scalar p)
\hspace*{\arraycolsep} = \hspace*{\arraycolsep}
[s\scalar p, \kappa_{2} \after\; !_{X}] 
& = &
[[s,\kappa_{2}] \after p, \kappa_{2} \after\; !_{X}] \\
& = &
[s,\kappa_{2}] \after [p, \kappa_{2} \after\; !_{X}] \\
& = &
[s,\kappa_{2}] \after \FstAnd(p)
\hspace*{\arraycolsep} = \hspace*{\arraycolsep}
s\scalar \FstAnd(p).
\end{array}\eqno{\qEd}$$
}
\end{enumerate}
\end{myproof}

\begin{rem}
\label{LiftRem}
As already mentioned, the coproducts $+$ and the final object $1$ in
an effectus $\cat{B}$ give rise to the `lift' or `maybe' monad $\Lift
= (-)+1$ on $\cat{B}$. Its unit maps $X \rightarrow X+1$ are the first
coprojections $\kappa_{1}$, and its multiplication maps $(X+1)+1
\rightarrow X+1$ are of the form $[\idmap,\kappa_{2}]$. The associated
category of Eilenberg-Moore algebras $\EM(\Lift)$ is isomorphic to the
category with states $\omega \colon 1 \rightarrow X$ as object,
corresponding to an algebra $[\idmap,\omega] \colon X+1 \rightarrow
X$. Maps correspond to commuting triangles between states.

The Kleisli category of the lift monad, written as $\KlL{\cat{B}}$,
typically captures \emph{partial} computations that may fail to
produce an output, for instance via non-terminating behaviour. The
Kleisli category of the lift monad $(-)+1$ on $\Sets$ has maps of the
form $X \rightarrow Y+1$, corresponding to partial functions $X
\rightarrow Y$. With the `flat' order, these maps form a directed
complete partial order (dcpo).

We can also consider the lift monad $\Lift = (-)+1$ on the
Kleisli category of the distribution monad $\Dst$. A map $X
\rightarrow Y$ in the associated Kleisli category (of lift) is a map
$X \rightarrow Y+1$ in the Kleisli category $\Kl(\Dst)$ of the
distribution monad. Hence it is a function $X \rightarrow
\Dst(Y+1)$. A crucial observations is that a distribution $\varphi
\in\Dst(Y+1)$ on $Y+1$ is a \emph{subdistribution} on $Y$: a formal
convex sum $\sum_{i}r_{i}\ket{y_{i}}$ of elements $y_{i}\in Y$ and
probabilities $r_{i}\in [0,1]$ whose sum is \emph{below} 1, as in:
$\sum_{i}r_{i} \leq 1$. These subdistributions can be ordered
pointwise, and have the everywhere-zero subdistribution as least
element. In case infinite supports are allowed, they are directed
complete.

A similar phenomenon exists for $C^*$-algebras: the lift monad
$(-)\oplus\C$ on the opposite category $\op{(\CstarPU)}$ yields a
correspondence between maps in $\CstarPU$.
$$\begin{prooftree}
\xymatrix{\mbox{positive unital } A\oplus \C\ar[r]^-{f} & B}
\Justifies
\xymatrix{\mbox{positive subunital }A\ar[r]_-{g} & B}
\end{prooftree}$$

\noindent Subunitality of $g\colon A \rightarrow B$ means $0 \leq g(1)
\leq 1$. This correspondence is given as follows. For $f \colon
A\oplus \C \rightarrow B$ take $\overline{f}(a) = f(a, 0)$, and for
$g\colon A \rightarrow B$, take $\overline{g}(a,z) = g(a) + z\cdot
(1-g(1))$. When $A,B$ are $W^*$-algebras, the normal completely
positive subunital maps $A\rightarrow B$ carry a dcpo structure,
see~\cite{Cho14a,Rennela14a}, which can be used to interpret loops,
like while or recursion, see also~\cite{RennelaS15a} for more
information.


\auxproof{
We first check that these functions are well-defined.
\begin{itemize}
\item $\overline{f}$ is positive, since if $a \geq 0$, then $(0, a)
  \geq 0$, so that $\overline{f}(a) = f(0,a) \geq 0$. It is subunital,
since $\overline{f}(1) = f(0,1) \leq 1$, since $f(1,0) \geq 0$ and
$f(0,1) + f(1,0) = f(1,1) = 1$.

\item $\overline{g}$ is positive, since if $(z,a) \geq 0$, then $z\geq
  0$ and $a \geq 0$. Hence $g(a) \geq 0$ and since $0 \leq g(1) \leq
  1$ we have $1 - g(1) \geq 0$ and also $z\cdot (1-g(1)) \geq
  0$. Hence $\overline{g}(z,a) = g(a) + z\cdot (1-g(1)) \geq 0$. Also,
  by construction, $\overline{g}(1,1) = g(1) + 1\cdot (1-g(1)) = g(1)
  + 1 - g(1) = 1$.
\end{itemize}

$$\begin{array}{rcl}
\overline{\overline{f}}(z,a)
& = &
\overline{f}(a) + z \cdot (1 - \overline{f}(1)) \\
& = &
f(0, a) + z\cdot (1 - f(0, 1)) \\
& = &
f(0, a) + z\cdot (f(1,1) - f(0, 1)) \\
& = &
f(0, a) + z\cdot f(1,0) \\
& = &
f(0, a) + f(z,0) \\
& = &
f(z,a) \\
\overline{\overline{g}}(a)
& = &
\overline{g}(0,a) \\
& = &
g(a) + 0\cdot (1- g(1)) \\
& = &
g(a).
\end{array}$$
}

A systematic investigation of partiality in effectuses can be found
in~\cite{Cho15a}. There, the main result shows that an effectus can
equivalently be described as a finitely partially additive category
(as in~\cite{ArbibM86}), with suitable effect algebra structure. This
shows that partiality is an intrinsic feature of the theory of
effectuses.


\end{rem}

\begin{rem}
\label{StatesCopropPresRem}
The starting point of this section is Lemma~\ref{CotupleIsoPredLem},
showing that the predicate functor $\Pred \colon \cat{B} \rightarrow
\op{\EMod}$ of an effectus $\cat{B}$ automatically preserves
coproducts. A natural question to ask is: does the states functor
$\Stat \colon \cat{B} \rightarrow \Conv$ in a state-and-effect
triangle~\eqref{GeneralTriangleDiag} also preserve coproducts?

This matter is addressed in~\cite{JacobsWW15a}. There it is shown that
preservation of $+$ by the states functor $\Stat$ does not come for
free, but corresponds to an important property, namely
\emph{normalisation} of states. This is relevant in conditional
probability, as we briefly illustrate in the Kleisli category
$\Kl(\Dst)$ of the distribution monad $\Dst$ on $\Sets$, where this
preservation/normalisation property holds. If we have a distribution
$\varphi \in \Dst(X+Y)$ on a coproduct set $X+Y$, then we can write
$\varphi$ as convex combination $r\cdot \Dst(\kappa_{1})(\varphi|X) +
(r-1)\cdot \Dst(\kappa_{2})(\varphi|Y)$, for certain `normalised'
distributions $\varphi|X\in\Dst(X), \varphi|Y\in \Dst(Y)$ and $r =
\sum_{x\in X} \varphi(\kappa_{1}x) \in [0,1]$.  Such convex
combinations form the coproduct in the category $\Conv$,
see~\cite{JacobsWW15a} for details.

\auxproof{
Take $r = \sum_{x\in X}\varphi(\kappa_{1}x) \in [0,1]$. Then $1 - r =
\sum_{y\in Y} \varphi(\kappa_{2}y)$. Distinguish:
\begin{itemize}
\item If $r = 0$, then take $\varphi|Y(y) = \varphi(\kappa_{2}y)$.

\item If $r = 1$, take $\varphi|X(x) = \varphi(\kappa_{1}x)$.

\item If $s\neq 0, s\neq 1$, then take $\varphi|X(x) =
  \frac{1}{r}\cdot \varphi(\kappa_{1}x)$ and $\varphi|Y =
  \frac{1}{1-r}\cdot \varphi(\kappa_{2}y)$. These $\varphi|X$ and
  $\varphi|Y$ are both distributions, and:
$$\begin{array}{rcl}
\varphi
& = &
\sum_{x\in X}\varphi(\kappa_{1}x)\ket{\kappa_{1}x} + 
   \sum_{y\in Y} \varphi(\kappa_{2}y)\ket{\kappa_{2}y} \\
& = &
r\cdot \sum_{x\in X}\varphi|X(x)\ket{\kappa_{1}x} + 
   (1-r)\cdot \sum_{y\in Y} \varphi|Y(y)\ket{\kappa_{2}y} \\
& = &
r\cdot \Dst(\kappa_{1})(\varphi|X) + (r-1)\cdot \Dst(\kappa_{2})(\varphi|Y).
\end{array}$$
\end{itemize}
}
\end{rem}

\section{Measurement instruments}\label{MeasurementSec}

Our second assumption introduces (discrete) measurement instruments as
certain maps in our base category. Such instruments have been
introduced in~\cite{DaviesL70}, see also~\cite{Ozawa84}
and~\cite{HeinosaariZ12}.  Characteristic aspects are:
\begin{enumerate}
\item different measurement outcomes can be distinguished
\item state changes caused by measurement are taken into account.
\end{enumerate}

\noindent Assumption~\ref{MeasurementAss} below describes our
categorical formalisation of discrete instruments for finitary
measurements, as maps of the form $X \rightarrow X+\cdots+X$, induced
by a test. The original references~\cite{DaviesL70,Ozawa84} deal with
the non-finitary, continuous case, where effects are indexed via a
measurable space, see
Remark~\ref{MeasurementRem}~\eqref{MeasurementRemInstr} below. The
different measurement outcomes in the first point correspond to the
different coproduct options in the codomain $X+\cdots+X$. The state
changes are captured by the side-effect associated with such an
instrument, which we define as the composite with the codiagonal
$\nabla = [\idmap,\cdots,\idmap]$ in:
$$\xymatrix{
X\ar[r] & X+\cdots+X\ar[r]^-{\nabla} & X
}$$

\noindent There is no side-effect if this map $X \rightarrow X$ is the
identity. We use the term `side-effect-free' for this situation. These
side-effects will be studied systematically in the next section.  In
this section we concentrate on the definition of instruments and on
examples, and show how they can be used for guarded test programs.

We should emphasise that the instrument assumption below requires the
presence of certain maps, satifying certain properties. It is an open
question, under which conditions, these properties fully determine the
instrument maps. Hence these instruments form \emph{structure} in the
category, and their presence is not a \emph{property} of the category
--- as long as we cannot show that they are uniquely determined. This
is not uncommon: our third assumption is about the presence of tensors
$\otimes$, which are also non-uniquely determined structure.

\begin{assumption}
\label{MeasurementAss}
Let $\cat{B}$ be an effectus, \textit{i.e.}~a category satisfying
Assumption~\ref{CoprodAss}. For each $n$-test $p\colon X \rightarrow
n\cdot 1$ there is a \emph{measurement instrument} map $\instr_{p}
\colon X \rightarrow n\cdot X$ in $\cat{B}$ making the following
diagram commute.
$$\xymatrix@C+1pc{
X\ar[r]^-{p}\ar@{..>}[dr]_{\instr_p} & n\cdot 1 \\
& n\cdot X\ar[u]_{n\cdot !}
}$$

\noindent Notice that the `$1$' in the codomain of a test $X
\rightarrow n\cdot 1$ is replaced by `$X$' in an instrument $X
\rightarrow n\cdot X$. Thus, where a test only captures the different
output options (and no outcome, in $1$), an instrument has both
options and an outcome (in $X$). The above diagram says that
annihilating these outcomes in $X$ via the unique map $!\colon
X\rightarrow 1$ returns the original test.

We do not require uniqueness of these instruments, but we
do require the following properties.
\begin{enumerate}
\item \label{MeasurementAssInjection} (Injection) For an injective
  function $\phi \colon \{1, \ldots, n\} \rightarrowtail \{1,\ldots,
  m\}$ and an object $A$, write $[\phi] \colon n\cdot A \rightarrow
  m\cdot A$ for the obvious map $[\phi] = [\kappa_{\phi(1)}, \ldots,
    \kappa_{\phi(n)}]$. Then for each $n$-test $p\colon X \rightarrow
  n\cdot 1$ the following diagram commutes.
$$\xymatrix@R-.5pc{
X\ar[rr]^-{\instr_{p}}\ar[drr]_{\instr_{[\phi] \after p}} & & 
   n\cdot X\ar[d]^{[\phi]} \\
& & m\cdot X
}$$

\noindent This property says that reordering elements of a test and
adding zero (falsum) predicates is reflected in the corresponding
instrument. It implies for instance: $\instr_{p^{\perp}} =
[\kappa_{2}, \kappa_{1}] \after \instr_{p}$, for a predicate $p\colon
X \rightarrow 1+1$, where $p^{\perp} = [\kappa_{2}, \kappa_{1}] \after
p$ is the orthocomplement.



\item \label{MeasurementAssCotuple} (Cotuple) For two $n$-tests
  $p\colon X \rightarrow n\cdot 1$ and $q\colon Y \rightarrow n\cdot
  1$, the following diagram commutes:
$$\xymatrix@R-.5pc{
X+Y\ar[rr]^-{\instr_{[p,q]}}\ar[drr]_{\instr_{p}+\instr_{q}\quad} & & 
   n\cdot (X+Y)\ar[d]_{\cong}^{[\kappa_{1}+\kappa_{1}, \ldots,
    \kappa_{n}+\kappa_{n}]} \\
& & n\cdot X + n\cdot Y
}$$

\item \label{MeasurementAssCollapse} (Collapse) For each $n+m$-test
  $p\colon X \rightarrow n\cdot 1 + m\cdot 1$ one can ``collapse'' the
  $m$-options in two different ways, but this gives the same outcome,
  as expressed by the following commuting diagram.
$$\xymatrix@R-1pc@C+1pc{
& n\cdot X + m\cdot X\ar@/^1ex/[dr]^(0.6){\idmap+\,!} & \\
X\ar@/^1ex/[ur]^(0.3){\instr_p}\ar@/_1ex/[dr]_(0.4){\instr_{(\idmap+\,!) \after p}\quad} 
   & & n\cdot X\rlap{$\, +\; 1$} \\
& n\cdot X+X\ar@/_1ex/[ur]_(0.55){\idmap+\,!}
}$$

\medskip

\item \label{MeasurementAssSef} (Side-effect-free) Given a map
  $q\colon X \rightarrow n\cdot X$ in $\cat{B}$ which is
  side-effect-free, that is $\nabla \after q = \idmap$, then $q =
  \instr_{p}$ for the $n$-test $p = (n\cdot !) \after q \colon X
  \rightarrow n\cdot 1$. In diagrams:
$$\mbox{if}\quad
\vcenter{\xymatrix@R-.5pc{
X\ar[r]^-{q}\ar@{=}[dr] & n\cdot X\ar[d]^{\nabla} \\
 & X
}}
\quad\mbox{then}\quad
\vcenter{\xymatrix@R-.5pc{
X\ar[r]^-{q}\ar[dr]_{p} & n\cdot X\ar[d]^{n\cdot !} \\
 & n\cdot 1
}}
\quad\mbox{satisfies}\quad
\vcenter{\xymatrix@C+.5pc{
X\ar[r]^-{q = \instr_{p}} & n\cdot X
}}
$$

\end{enumerate}
\end{assumption}

\noindent When $\nabla \after q = \idmap$ as in the last property, then $q
\colon X \rightarrow X+\cdots + X$ is merely selecting paths among the
various coproduct options, without having any computational
effect. Therefor we call such a $q$ side-effect-free, see the next
section. The above property~\eqref{MeasurementAssSef} says that such
path-selector maps $q$ arise as instruments from the corresponding
path-selector test $(n\cdot !) \after q \colon X \rightarrow
1+\cdots+1$. This path property is non-trivial, especially for
$C^*$-algebras, see Corollary~\ref{MeasurementCstarSefCor}. We
explicitly mention two consequences.

\begin{lem}
\label{MeasurementAssSefLem}
In a category satisfying Assumption~\ref{MeasurementAss} we have:
\begin{enumerate}
\item \label{MeasurementAssSefLemOne} $\instr_{q} = q$ for each
  $n$-test $q\colon 1 \rightarrow n\cdot 1$ on the final object $1$;

\item \label{MeasurementAssSefLemCoprojection} $\instr_{\kappa_{i}
  \after\, !} = \kappa_{i} \colon X \rightarrow n\cdot X$, for each $1
  \leq i \leq n$.
\end{enumerate}

\noindent In particular, $\instr_{q} = q$ when $q$ is a scalar, and
$\instr_{1} = \kappa_{1}, \instr_{0} = \kappa_{2}$, for the true and
false predicates $1,0 \colon X \rightarrow 1+1$.
\end{lem}

\begin{myproof}
We use property~\eqref{MeasurementAssSef} from
Assumption~\ref{MeasurementAss} in both cases. For $q\colon 1
\rightarrow n\cdot 1$ we clearly have $\nabla \after q = \idmap$ by
finality of $1$. But also $p = (n\cdot !) \after q = q$, so that
$\instr_{q} = q$. For the second point, notice that the coprojection
$\kappa_{i} \colon X \rightarrow n\cdot X$ satisfies $\nabla \after
\kappa_{i} = \idmap$. The $n$-test $(n\cdot !_{X}) \after \kappa_{i} =
\kappa_{i} \after \, !_{X} \colon X \rightarrow n\cdot 1$ gives
$\instr_{\kappa_{i} \after\, !_{X}} = \kappa_{i}$. \QED
\end{myproof}

We illustrate these instruments in our leading examples.

\begin{exas}
\label{MeasurementEx}
In the category $\Sets$, an $n$-test $p\colon X \rightarrow n\cdot 1 =
\{1,2, \ldots, n\}$ corresponds to a partition of $X$, that is, to a
cover of $X$ given by disjoint subsets. Hence we can define an
instrument $\instr_{p} \colon X \rightarrow n\cdot X$ that picks the
$i$-th coproduct option if the input $x$ is in the $i$-th subset:
\begin{equation}
\label{MeasurementSetsEqn}
\begin{array}{rclcrcl}
\instr_{p}(x)
& = &
\kappa_{i}x 
& \qquad\mbox{iff}\qquad &
p(x)
& = &
i.
\end{array}
\end{equation}

\noindent Clearly, $n\cdot ! \after \instr_{p} = p$.

For an $n$-test $p\colon X \rightarrow \Dst(n)$ in the Kleisli
category $\Kl(\Dst)$ the associated instrument $\instr_{p} \colon
X \rightarrow \Dst(n\cdot X)$ works as before, but now gives a distribution
over the various coproduct options, via the following convex sum:
\begin{equation}
\label{MeasurementDstEqn}
\begin{array}{rcl}
\instr_{p}(x)
& = &
p(x)(1)\ket{\kappa_{1}x} + \cdots + p(x)(n)\ket{\kappa_{n}x}
\end{array}
\end{equation}

\auxproof{
We check the condition:
$$\begin{array}{rcccccccl}
\Dst(n\cdot !)\big(\instr_{p}(x)\big)
& = &
\sum_{i} p(x)(i)\ket{(n\cdot !)(\kappa_{i}x)}
& = &
\sum_{i} p(x)(i)\ket{\kappa_{i}*}
& = &
\sum_{i} p(x)(i)\ket{i}
& = &
p(x).
\end{array}$$

We also check
Assumption~\ref{MeasurementAss}~\eqref{MeasurementAssCollapse} in
this case. So assume we have an $n+m$-test $p\colon X \rightarrow
\Dst(n+m)$. Then $p' = (\idmap+\,!) \klafter p \colon X \rightarrow
\Dst(n+1)$ is given by $p'(x)(i) = p(x)(i)$ if $i\leq n$ and
$p'(x)(n+1) = \sum_{j}p(x)(n+j)$. Then we get an equality of maps
$X \rightarrow \Dst(n\cdot X+1)$ in:
$$\begin{array}{rcl}
\big((\idmap+\,!) \klafter \instr_{p})(x)
& = &
p(x)(1)\ket{\kappa_{1}x} + \cdots + p(x)(n)\ket{\kappa_{n}x} +
   (\sum_{j}p(x)(n+j))\ket{\kappa_{n+1}*} \\
& = &
\big((\idmap+\,!) \klafter \instr_{p'})(x).
\end{array}$$

\noindent Notice that in this case we already have an equality of maps
$X \rightarrow \Dst(n\cdot X + X)$, so collapsing the last +-component
$X$ to $1$ is not needed in this probabilistic case.

We check the coherence conditions.
$$\begin{array}{rcl}
\lefteqn{\big([\kappa_{1}+\kappa_{1}, \cdots, \kappa_{n}+\kappa_{n}] \klafter
   \instr_{[p,q]}\big)(\kappa_{1}x)} \\
& = &
\Dst([\kappa_{1}+\kappa_{1}, \cdots, \kappa_{n}+\kappa_{n}])
   (\sum_{i}p(x)(i)\ket{\kappa_{i}\kappa_{1}x}) \\
& = &
\sum_{i}p(x)(i)\ket{\kappa_{1}\kappa_{i}x} \\
& = &
\Dst(\kappa_{1})(\instr_{p}(x)) \\
& = &
(\instr_{p}+\instr_{q})(\kappa_{1}x).
\end{array}$$

\noindent And for an injection $\phi \colon \{1, \ldots, n\} 
\rightarrowtail \{1,\ldots, m\}$,
$$\begin{array}{rcl}
\big([\phi] \klafter \instr_{p}\big)(x)
& = &
\Dst([\phi])(\sum_{i} p(x)(i)\ket{\kappa_{i}x}) \\
& = &
\sum_{i} p(x)(i)\ket{\kappa_{\phi(i)}x} \\
& = &
\instr_{[\phi] \klafter p}(x),
\end{array}$$

\noindent where we use that for $p\colon X \rightarrow n\cdot X$, written
as $p(x) = \sum_{i}p(x)(i)\ket{\kappa_{i}*}$ we have:
$$\begin{array}{rcccl}
\big([\phi] \klafter p\big)(x)
& = &
\Dst([\phi])(\sum_{i}p(x)(i)\ket{\kappa_{i}*}) 
& = &
\sum_{i}p(x)(i)\ket{\kappa_{\phi(i)}*}
\end{array}$$

If $q\colon X \rightarrow \Dst(n\cdot X)$ satisfies $\nabla \klafter q
= \idmap$, then $\sum_{i}q(x)(\kappa_{i}x') = 1$ if $x=x'$ and $0$ if
$x\neq x'$. The predicate $p = (n\cdot !) \klafter q \colon X
\rightarrow \Dst(n)$ is given by $p(x)(i) = q(x)(\kappa_{i}x)$. Then 
indeed:
$$\begin{array}{rcl}
\instr_{p}(x)
& = &
\sum_{i} p(x)(i)\ket{\kappa_{i}x} \\
& = &
\sum_{i} q(x)(\kappa_{i}x)\ket{\kappa_{i}x} \\
& = &
\sum_{i,x'} q(x)(\kappa_{i}x')\ket{\kappa_{i}x'} \\
& = &
q(x)
\end{array}$$
}

In the case of continuous probability, given by the Kleisli category
$\Kl(\Giry)$ of the Giry monad $\Giry$, an $n$-test is a measurable
map $p\colon X \rightarrow \Giry(n) \cong \Dst(n)$. Hence we can use
basically the same construction as before: the instrument $\instr_{p}
\colon X \rightarrow \Giry(n\cdot X)$ sends $x\in X$ to the
probability measure $\instr_{p}(x) \colon \Sigma_{n\cdot X}
\rightarrow [0,1]$ defined on basic measurable subsets
$\kappa_{i}M\in\Sigma_{n\cdot X}$ as:
\begin{equation}
\label{MeasurementGiryEqn}
\begin{array}{rcccl}
\instr_{p}(x)(\kappa_{i}M)
& = &
p(x)(i)\cdot\indic{M}(x)
& = &
\left\{\begin{array}{ll}
p(x)(i) \quad & \mbox{if }x\in M \\
0 & \mbox{otherwise.}
\end{array}\right.
\end{array}
\end{equation}

\noindent We sketch the proof of requirement~\eqref{MeasurementAssSef}
because it is slightly subtle. Let $q\colon X \rightarrow \Giry(n\cdot
X)$ be a measurable map satisfying $\nabla \klafter q = \idmap$, where
$\klafter$ is Kleisli composition. One easily checks that this amounts
to, for $x\in X$ and $M\in\Sigma_{X}$,
$$\begin{array}{rcccccl}
\indic{M}(x)
& = &
\eta(x)(M) 
& = &
\idmap(x)(M)
& = &
q(x)(\kappa_{1}M) + \cdots + q(x)(\kappa_{n}M)
\end{array}$$

\noindent This equality gives some useful information. First, if
$x\not\in M$, then $\indic{M}(x) = 0$, so $q(x)(\kappa_{i}M) = 0$ for
each $i$.  Second, if $M\subseteq N$ then $q(x)(\kappa_{i}M) \leq
q(x)(\kappa_{i}N)$ by monotonicity of the probability measure
$q(x)$. But if $x\in M$ this $\leq$ must be an equality, since also $1
= \indic{N}(x) = q(x)(\kappa_{1}N) + \cdots +
q(x)(\kappa_{n}N)$. Third, $q(x)(\kappa_{i}M) = q(x)(\kappa_{i}X)
\cdot \indic{M}(x)$.  Indeed, if $x\not\in M$ then both sides are $0$
by the first observation. And if $x\in M$, then both sides are equal
by the second one. Next we note that the test $(n\cdot !) \klafter q
\colon X \rightarrow \Giry(n)$ can be described as:
$$\begin{array}{rcccccl}
((n\cdot !) \klafter q)(x)(i)
& = &
\int (n\cdot !)(-)(i) \intd q(x) 
& = &
\int \indic{\kappa_{i}X} \intd q(x) 
& = &
q(x)(\kappa_{i}X).
\end{array}$$

\noindent Finally we can prove $\instr_{(n\cdot !) \klafter q} = q$ via:
$$\begin{array}{rcccccl}
\instr_{(n\cdot !) \klafter q}(x)(\kappa_{i}M)
& = &
((n\cdot !) \klafter q)(x)(i) \cdot \indic{M}(x) 
& = &
q(x)(\kappa_{i}X) \cdot \indic{M}(x) 
& = &
q(x)(\kappa_{i}M).
\end{array}$$

\auxproof{
We first have to check that $\instr_{p}(x)$ is a probability measure.
$$\begin{array}{rcl}
\instr_{p}(x)(n\cdot X)
& = &
\instr_{p}(x)(\bigcup_{i} \kappa_{i}X) \\
& = &
\sum_{i} \instr_{p}(x)(\kappa_{i}X) \qquad
   \mbox{via extension of~\eqref{MeasurementGiryEqn}} \\
& = &
\sum_{i} p(x)(i)\cdot \indic{X}(x) \\
& = &
\sum_{i} p(x)(i) \\
& = &
1.
\end{array}$$

\noindent If we have a collection $N_{i\in I} \in \Sigma_{n\cdot X}$
of pairwise disjoint subsets, then we can split it in $n$-collections
forming a disjoint union $I = I_{1} \ovee \cdots \ovee I_{n}$, with
measurable subsets of the form $\kappa_{j}M_{i}$, for $i\in I_{j}$.  Then:
$$\begin{array}{rcl}
\instr_{p}(x)(\bigovee_{i\in I}N_{i}) 
& = &
\instr_{p}(x)(\bigovee_{j\leq n}\bigovee_{i\in I_j} \kappa_{i}M_{i}) \\
& = &
\sum_{j} \instr_{p}(x)(\bigovee_{i\in I_j} \kappa_{j}M_{i}) \\
& = &
\sum_{j} p(x)(j) \cdot \indic{\bigovee_{i\in I_j} M_{i}} \\
& = &
\sum_{j} p(x)(j) \cdot (\sum_{i\in I_j}\indic{M_{i}}) \\
& = &
\sum_{i,j} p(x)(j) \cdot \indic{M_{i}} \\
& = &
\sum_{i,j} \instr_{p}(x)(\kappa_{j}M_{i}) \\
& = &
\sum_{i\in I} \instr_{p}(x)(N_{i}).
\end{array}$$

Next, we show that the annihilation property holds. This requires that
we have an equality of maps $X \rightarrow \Giry(n\cdot 1)$. An
element of $\Giry(n\cdot 1)$ is a probability measure $\Sigma_{n\cdot
  1} = \Sigma_{n} = \Pow(n) \rightarrow [0,1]$. Such measures
correspond to functions $n \rightarrow [0,1]$, by application to
$\{i\}$. In the calculation below we identify $i\in n$ with the
singleton $\{i\}\in\Pow(n)$.
$$\begin{array}{rcl}
\big((n\cdot !) \klafter \instr_{p}\big)(x)(\{i\})
& = &
(n\cdot !)_{*}(\instr_{p}(x))(\{i\}) \\
& = &
\int (n\cdot !)(-)(\{i\}) \intd \instr_{p}(x) \\
& = &
\int \indic{\kappa_{i}X} \intd \instr_{p}(x) \\
& = &
\instr_{p}(x)(\kappa_{i}X) \\
& = &
p(x)(i) \cdot \indic{X}(x) \\
& = &
p(x)(i).
\end{array}$$

\noindent We continue with the four properties from
Assumption~\ref{MeasurementAss}.
\begin{enumerate}
\item For an injection $\phi \colon \{1,
\ldots, n\} \rightarrowtail \{1, \ldots, m\}$ we have $[\phi] \colon
n\cdot A \rightarrow \Giry(m\cdot A)$ given by:
$$\begin{array}{rcl}
[\phi](\kappa_{i}a)(\kappa_{j}M)
& = &
\left\{\begin{array}{ll}
\indic{\kappa_{i}M}(\kappa_{i}a) \quad & \mbox{if }\phi(i)=j \\
0 & \mbox{otherwise} 
\end{array}\right.
\end{array}$$

\noindent Thus for an $n$-test $p\colon X \rightarrow \Giry(n)\cong \Dst(n)$,
and $\{j\} \in \Pow(m) = \Sigma_{m}$,
$$\begin{array}{rcl}
\big([\phi] \klafter p\big)(x)(\{j\})
& = &
[\phi]_{*}(p(x))(\{j\}) \\
& = &
\int [\phi](-)(\{j\}) \intd p(x) \\
& = &
\left\{\begin{array}{ll}
\int \indic{\{i\}} \intd p(x) \quad & \mbox{if }\phi(i) = j \\
0 & \mbox{otherwise}
\end{array}\right. \\
& = &
\left\{\begin{array}{ll}
p(x)(i) \quad & \mbox{if }\phi(i) = j \\
0 & \mbox{otherwise}
\end{array}\right. \\
\end{array}$$

\noindent and so:
$$\begin{array}{rcl}
\big([\phi] \klafter \instr_{p}\big)(x)(\kappa_{j}M)
& = &
[\phi]_{*}(\instr_{p}(x))(\kappa_{j}M) \\
& = &
\int [\phi](-)(\kappa_{j}M) \intd \instr_{p}(x) \\
& = &
\left\{\begin{array}{ll}
\int \indic{\kappa_{i}M} \intd \instr_{p}(x) \quad
   & \mbox{if }\phi(i) = j \\
\int \indic{\emptyset} \intd \instr_{p}(x) & \mbox{otherwise}
\end{array}\right. \\
& = &
\left\{\begin{array}{ll}
\instr_{p}(x)(\kappa_{i}M) \quad  & \mbox{if }\phi(i) = j \\
0 & \mbox{otherwise}
\end{array}\right. \\
& = &
\left\{\begin{array}{ll}
p(x)(i)\cdot \indic{M}(x) \quad  & \mbox{if }\phi(i) = j \\
0 \cdot \indic{M}(x) & \mbox{otherwise}
\end{array}\right. \\
& = &
([\phi]\klafter p)(x)(j) \cdot \indic{M}(x) \\
& = &
\instr_{[\phi] \klafter p}(x)(\kappa_{j}M)
\end{array}$$

\item For measurable predicates $p,q \colon X \rightarrow [0,1]$ we have
for $M\in\Sigma_{X}$,
$$\begin{array}{rcl}
\lefteqn{\big([\kappa_{1}+\kappa_{1}, \ldots, \kappa_{n}+\kappa_{n}] \klafter 
   \instr_{[p,q]}\big)(\kappa_{j}x)(\kappa_{j'}\kappa_{i}M)} \\
& = &
[\kappa_{1}+\kappa_{1}, \ldots, \kappa_{n}+\kappa_{n}]_{*}
   (\instr_{[p,q]}(\kappa_{j}x))(\kappa_{j'}\kappa_{i}M) \\
& = &
\int [\kappa_{1}+\kappa_{1}, \ldots, \kappa_{n}+\kappa_{n}](-)
   (\kappa_{j'}\kappa_{i}M) \intd \instr_{[p,q]}(\kappa_{j}x) \\
& = &
\int \big(\lamin{u}{n\cdot (X+Y)}
   {\eta([\kappa_{1}+\kappa_{1}, \ldots, \kappa_{n}+\kappa_{n}](u))
    (\kappa_{j'}\kappa_{j}M)}\big) \intd \instr_{[p,q]}(\kappa_{j}x) \\
& = &
\int \big(\lamin{u}{n\cdot (X+Y)}{\indic{\kappa_{j'}\kappa_{i}M}(
   [\kappa_{1}+\kappa_{1}, \ldots, \kappa_{n}+\kappa_{n}](u))}\big) 
   \intd \instr_{[p,q]}(\kappa_{j}x) \\
& = &
\int \big(\lamin{u}{n\cdot (X+Y)}{\indic{\kappa_{i}\kappa_{j'}M}(u)}\big)
   \intd \instr_{[p,q]}(\kappa_{j}x) \\
& = &
\int \indic{\kappa_{i}\kappa_{j'}M} \intd \instr_{[p,q]}(\kappa_{j}x) \\
& = &
\instr_{[p,q]}(\kappa_{j}x)(\kappa_{i}\kappa_{j'}M) \\
& = &
[p,q](\kappa_{j}x)(i) \cdot \indic{\kappa_{j'}M}(\kappa_{j}x) \\
& = &
\left\{\begin{array}{ll}
p(x)(i) \cdot \indic{M}(x) \quad & \mbox{if } j = j' = 1 \\
q(x)(i) \cdot \indic{M}(x) \quad & \mbox{if } j = j' = 2 \\
0 & \mbox{otherwise}
\end{array}\right. \\
& = &
\left\{\begin{array}{ll}
\instr_{p}(x)(\kappa_{i}M) \quad & \mbox{if } j = j' = 1 \\
\instr_{q}(x)(\kappa_{i}M) \quad & \mbox{if } j = j' = 2 \\
0 & \mbox{otherwise}
\end{array}\right. \\
& = &
\big(\instr_{p}+\instr_{q}\big)(\kappa_{j}x)(\kappa_{j'}\kappa_{i}M).
\end{array}$$

\item Let $p\colon X \rightarrow \Giry(n+m)$ be given. Then,
for $1 \leq i \leq n+1$,
$$\begin{array}{rcl}
\big((\idmap+!) \klafter p\big)(x)(\{i\})
& = &
(\idmap+!)_{*}(p(x))(\{i\}) \\
& = &
\int (\idmap+!)(-)(\{i\}) \intd p(x) \\
& = &
\left\{\begin{array}{ll}
\int \indic{\{i\}} \intd p(x) \quad & \mbox{if } i \leq n \\
\int \indic{\{n+1, \ldots, n+m\}} \intd p(x) \quad & \mbox{if } i = n+1 \\
\end{array}\right. \\
& = &
\left\{\begin{array}{ll}
p(x)(i) \quad & \mbox{if } i \leq n \\
p(\{n+1, \ldots, n+m\}) \quad & \mbox{if } i = n+1 \\
\end{array}\right. \\
& = &
\left\{\begin{array}{ll}
p(x)(i) \quad & \mbox{if } i \leq n \\
p(n+1) + \cdots + p(n+m) \quad & \mbox{if } i = n+1 \\
\end{array}\right. 
\end{array}$$

\noindent Thus:
$$\begin{array}{rcl}
\big((\idmap+!) \klafter \instr_{p})(x)(\kappa_{i}M)
& = &
(\idmap+!)_{*}(\instr_{p}(x))(\kappa_{i}M) \\
& = &
\int (\idmap+!)(-)(\kappa_{i}M) \intd \instr_{p}(x) \\
& = &
\left\{\begin{array}{ll}
\int \indic{\kappa_{i}M} \intd \instr_{p}(x) \quad 
   & \mbox{if } i \leq n \\
\int \indic{\kappa_{n+1}X \cup \cdots \cup \kappa_{n+m}X} \intd \instr_{p}(x) \quad
   & \mbox{if } i = n+1 \\
\end{array}\right. \\
& = &
\left\{\begin{array}{ll}
\instr_{p}(x)(\kappa_{i}M)  & \mbox{if } i \leq n \\
\instr_{p}(x)(\kappa_{n+1}X \cup \cdots \cup \kappa_{n+m}X) \quad
   & \mbox{if } i = n+1 \\
\end{array}\right. \\
& = &
\left\{\begin{array}{ll}
p(x)(i) \cdot \indic{\kappa_{i}M}(x)  & \mbox{if } i \leq n \\
p(x)(n+1) + \cdots + p(n+m) \quad  & \mbox{if } i = n+1 \\
\end{array}\right. \\
\end{array}$$

\noindent At the same time, 
$$\begin{array}{rcl}
\big((\idmap+!) \klafter \instr_{(\idmap+!)\klafter p})(x)(\kappa_{i}M)
& = &
(\idmap+!)_{*}(\instr_{(\idmap+!)\klafter p}(x))(\kappa_{i}M) \\
& = &
\int (\idmap+!)_{*}(-)(\kappa_{i}M) \intd \instr_{(\idmap+!)\klafter p}(x) \\
& = &
\left\{\begin{array}{ll}
\int \indic{\kappa_{i}M} \intd \instr_{p}(x) \quad 
   & \mbox{if } i \leq n \\
\int \indic{\kappa_{n+1}X} \intd \instr_{(\idmap+!)\klafter p}(x) \quad
   & \mbox{if } i = n+1 \\
\end{array}\right. \\
& = &
\left\{\begin{array}{ll}
\instr_{p}(x)(\kappa_{i}M) \quad & \mbox{if } i \leq n \\
\instr_{(\idmap+!)\klafter p}(x)(\kappa_{n+1}X) \quad  & \mbox{if } i = n+1 \\
\end{array}\right. \\
& = &
\left\{\begin{array}{ll}
p(x)(i) \cdot \indic{\kappa_{i}M}(x)  & \mbox{if } i \leq n \\
((\idmap+!)\klafter p)(x)(n+1) \cdot \indic{X}(x) \quad  & \mbox{if } i = n+1 \\
\end{array}\right. \\
& = &
\left\{\begin{array}{ll}
p(x)(i) \cdot \indic{\kappa_{i}M}(x)  & \mbox{if } i \leq n \\
p(x)(n+1) + \cdots + p(n+m) \quad  & \mbox{if } i = n+1 \\
\end{array}\right. \\
\end{array}$$

\item Next, assume that $q\colon X \rightarrow \Giry(n\cdot X)$ satisfies
$\nabla \klafter q = \idmap$. This means that:
$$\begin{array}{rcl}
\indic{M}(x)
& = &
\eta(x)(M) \\
& = &
(\nabla \klafter q)(x)(M) \\
& = &
\nabla_{*}(q(x))(M) \\
& = &
\int \nabla(-)(M) \intd q(x) \\
& = &
\int \indic{\kappa_{1}M \cup \cdots \cup \kappa_{n}M} \intd q(x) \\
& = &
q(x)(\kappa_{1}M \cup \cdots \cup \kappa_{n}M) \\
& = &
q(x)(\kappa_{1}M) + \cdots + q(x)(\kappa_{n}M)
\end{array}$$

\noindent From this we conclude:
\begin{itemize}
\item If $x\not\in M$, then $\indic{M}(x) = 0$, so $q(x)(\kappa_{i}M)
  = 0$ for each $i$;

\item If $M\subseteq N$ and $x\in M$, then $q(x)(\kappa_{i}M) =
  q(x)(\kappa_{i}N)$. By monotonicity we have $q(x)(\kappa_{i}M) \leq
  q(x)(\kappa_{i}N)$. Now suppose $q(x)(\kappa_{i}M) + \varepsilon_{i} = 
  q(x)(\kappa_{i}N)$ for $\varepsilon_{i} \geq 0$. But then:
$$\begin{array}{rcl}
1
& = &
\indic{N}(x) \\
& = &
q(x)(\kappa_{1}N) + \cdots + q(x)(\kappa_{n}N) \\
& = &
q(x)(\kappa_{1}M) + \varepsilon_{1} + \cdots + 
   q(x)(\kappa_{n}M) + \varepsilon_{n} \\
& = &
1 + \sum_{i}\varepsilon_{i}
\end{array}$$

\noindent Hence $\varepsilon_{i} = 0$.

\item In particular, $q(x)(\kappa_{i}M) = q(x)(\kappa_{i}X) \cdot
  \indic{M}(x)$.  Indeed, if $x\not\in M$ then both sides are $0$ by
  the first bullet. And if $x\in M$, then both sides are equal by the
  second bullet.
\end{itemize}

\noindent We intend to show that $\instr_{(n\cdot !) \klafter q} = q$. But
we first compute:
$$\begin{array}{rcl}
((n\cdot !) \klafter q)(x)(i)
& = &
\int (n\cdot !)(-)(i) \intd q(x) \\
& = &
\int \indic{\kappa_{i}X} \intd q(x) \\
& = &
q(x)(\kappa_{i}X).
\end{array}$$

\noindent Finally:
$$\begin{array}{rcl}
\instr_{(n\cdot !) \klafter q}(x)(\kappa_{i}M)
& = &
((n\cdot !) \klafter q)(x)(i) \cdot \indic{M}(x) \\
& = &
q(x)(\kappa_{i}X) \cdot \indic{M}(x) \\
& = &
q(x)(\kappa_{i}M) \qquad \mbox{by the third bullet above.}
\end{array}$$
\end{enumerate}
}\bigskip

\noindent We continue with rings and show that instruments exist for
\emph{commutative} rings. For a ring $R$ an $n$-test $R \rightarrow
n\cdot 1$ consists of $n$ idempotents $\vec{e} = e_{1}, \ldots, e_{n}
\in R$ with $e_{1} + \cdots + e_{n} = 1$ and $e_{i}\cdot e_{j} = 0$
for $i\neq j$. If $R$ is commutative we get an instrument
$\instr_{\vec{e}} \colon R \rightarrow n\cdot R$ in the opposite
category $\op{\CRng}$. In the category $\CRng$ it takes the form of a
function $\instr_{\vec{e}} \colon R^{n} \rightarrow R$ defined by:
\begin{equation}
\label{MeasurementRingEqn}
\begin{array}{rcl}
\instr_{\vec{e}}(x_{1}, \ldots, x_{n})
& = &
e_{1}\cdot x_{1} + \cdots + e_{n}\cdot x_{n}.
\end{array}
\end{equation}

\noindent It clearly preserves sums $(+,0)$. It also preserves
multiplications $(\cdot, 1)$ by the properties of an $n$-test and
commutativity. Of the required properties in
Assumption~\ref{MeasurementAss} we only do
number~\eqref{MeasurementAssSef}: if $q\colon R^{n} \rightarrow R$ is
a `side-effect-free' ring homomorphism with $q \after \Delta =
\idmap$, then $q(x,\ldots,x) = x$. Moreover, the $n$-test $p = q
\after \, !^{n} \colon \Z^{n} \rightarrow R$ is given by $p(m_{1},
\ldots, m_{n}) = q(m_{1}\cdot 1, \ldots, m_{n}\cdot 1)$, and
corresponds to idempotents $e_{i} = q(\ket{i}) \in R$, satisfying:
$$\begin{array}{rcl}
\instr_{\vec{e}}(x_{1}, \ldots, x_{n})
& = &
e_{1}\cdot x_{1} + \cdots + e_{n}\cdot x_{n} \\
& = &
q(\ket{1})\cdot q(x_{1}, \ldots, x_{1}) + \cdots + 
   q(\ket{n})\cdot q(x_{n},\ldots,x_{n}) \\
& = &
q(\ket{1}\cdot (x_{1}, \ldots, x_{1})) + \cdots + 
   q(\ket{n}\cdot (x_{n},\ldots,x_{n})) \\
& = &
q(x_{1}, 0, \ldots, 0) + \cdots + q(0,\ldots,0,x_{n}) \\
& = &
q(x_{1}, \ldots, x_{n}).
\end{array}$$

\auxproof{
$$\begin{array}{rcl}
\instr_{\vec{e}}(\vec{x})\cdot \instr_{\vec{e}}(\vec{y})
& = &
(\sum_{i}e_{i}x_{i})\cdot (\sum_{j}e_{j}y_{j}) \\
& = &
\sum_{i,j} e_{i}e_{j}x_{i}x_{j} \\
& = &
\sum_{i} e_{i}x_{i}y_{i} \\
& = &
\instr_{\vec{e}}(x_{1}y_{1}, \ldots, x_{n}y_{n}) \\
& = &
\instr_{\vec{e}}(\vec{x}\vec{y}).
\end{array}$$

$$\begin{array}{rcl}
\big(\instr_{\vec{e}} \after \, !^{n}\big)(m_{1}, \ldots, m_{n})
& = &
\instr_{\vec{e}}(m_{1}\cdot 1, \ldots, m_{n}\cdot 1) \\
& = &
e_{1}\cdot (m_{1}\cdot 1) + \cdots + e_{n}\cdot (m_{n}\cdot 1) \\
& = &
(e_{1} + \cdots + e_{1}) + \cdots + (e_{n} + \cdots + e_{n}) \\
& = &
m_{1}\cdot e_{1} + \cdots + m_{n}\cdot e_{n} \\
& = &
m_{1}\cdot f_{\vec{e}}(\ket{1}) + \cdots + m_{n}\cdot f_{\vec{e}}(\ket{1}) \\
& = &
f_{\vec{e}}(m_{1}\ket{1}) + \cdots + f_{\vec{e}}(m_{n}\ket{1}) \\
& = &
f_{\vec{e}}(m_{1}, \ldots, m_{n})
\end{array}$$

For two $n$-tests $e_{i}\in R$ and $d_{i}\in S$ we get an equality
$R^{n} \times S^{n} \rightarrow R\times S$, namely:
$$\begin{array}{rcl}
\big(\instr_{\vec{e},\vec{d}} \after \tuple{\pi_{1}\times\pi_{1}, \ldots,
   \pi_{n}\times\pi_{n}}\big)(\vec{x}, \vec{y}) 
& = &
\instr_{\vec{e},\vec{d}}((x_{1}, y_{1}), \ldots, (x_{n},y_{n})) \\
& = &
(e_{1},d_{1})\cdot (x_{1},y_{1}) + \cdots + (e_{n},d_{n})\cdot (x_{n},y_{n}) \\
& = &
(e_{1}x_{1} + \cdots e_{n}x_{n}, d_{1}y_{1} + \cdots + d_{n}y_{n}) \\
& = &
\tuple{\instr_{\vec{e}}, \instr_{\vec{d}}}(\vec{x}, \vec{y}).
\end{array}$$

For $\phi \colon \{1,2\} \rightarrowtail \{1, 2, 3, 4, 5\}$ say with
$\phi(1) = 3, \phi(2) = 5$ we get a map $\tuple{\phi} =
\tuple{\pi_{3}, \pi_{5}} \colon R^{5} \rightarrow R^{2}$. Thus, for
idempotents $e_{1},e_{2} \in R$ with $e_{1}+e_{2} = 1$ and $e_{1}\cdot
e_{2} = 0$ we have a $5$-test $\tuple{\phi}(\vec{e}) = (0,0,e_{1}, 0,
e_{2})$, and:
$$\begin{array}{rcl}
\big(\instr_{\vec{e}} \after \tuple{\phi}\big)(x_{1}, x_{2}, x_{3}, x_{4}, x_{5})
& = &
\instr_{\vec{e}}(x_{3}, x_{5}) \\
& = &
e_{1}x_{3} + e_{2}x_{5} \\
& = &
0x_{1} + 0x_{2} + e_{1}x_{3} + 0x_{4} + e_{2}x_{5} \\
& = &
\instr_{\tuple{\phi}(\vec{e})}(x_{1}, x_{2}, x_{3}, x_{4}, x_{5}).
\end{array}$$
}\bigskip

\noindent Along the same lines one shows that instruments exist in the opposite
$\op{\DL}$ of the category of distributive lattices: for such a
lattice $L$, with an $n$-test, consisting of elements $e_{1}, \ldots,
e_{n} \in L$ satisfying $e_{1} \vee \cdots \vee e_{n} = 1$, and $e_{i}
\wedge e_{j} = 0$ for $i\neq j$, one defines an instrument
$\instr_{\vec{e}} \colon L^{n} \rightarrow L$ by:
$$\begin{array}{rcl}
\instr_{\vec{e}}(x_{1}, \ldots, x_{n})
& = &
(e_{1} \wedge x_{1}) \vee \cdots \vee (e_{n} \wedge x_{n}).
\end{array}$$

\auxproof{
We check that $\instr_{\vec{e}}$ is a map of distributive lattices:
$$\begin{array}{rcl}
\instr_{\vec{e}}(0, \ldots, 0)
& = &
(e_{1} \wedge 0) \vee \cdots \vee (e_{n} \wedge 0) \\
& = &
0 \vee \cdots \vee 0 \\
& = &
0 \\
\instr_{\vec{e}}(x_{1}, \ldots, x_{n}) \vee \instr_{\vec{e}}(y_{1}, \ldots, y_{n})
& = &
(e_{1} \wedge x_{1}) \vee \cdots \vee (e_{n} \wedge x_{n}) \\
& & \qquad \vee (e_{1} \wedge y_{1}) \vee \cdots \vee (e_{n} \wedge y_{n}) \\
& = &
(e_{1} \wedge (x_{1}\vee y_{1})) \vee \cdots \vee 
   (e_{n} \wedge (x_{n}\vee y_{n})) \\
& = &
\instr_{\vec{e}}(x_{1}\vee y_{1}, \ldots, x_{n}\vee y_{n}) \\
\instr_{\vec{e}}(1, \ldots, 1)
& = &
(e_{1} \wedge 1) \vee \cdots \vee (e_{n} \wedge 1) \\
& = &
e_{1} \vee \cdots \vee e_{n} \\
& = &
1 \\
\instr_{\vec{e}}(x_{1}, \ldots, x_{n}) \wedge \instr_{\vec{e}}(y_{1}, \ldots, y_{n})
& = &
\big((e_{1} \wedge x_{1}) \vee \cdots \vee (e_{n} \wedge x_{n})\big) \\
& & \qquad \wedge 
   \big((e_{1} \wedge y_{1}) \vee \cdots \vee (e_{n} \wedge y_{n})\big) \\
& = &
\bigvee_{i,j} (e_{i}\wedge x_{i}) \wedge (e_{j} \wedge y_{j}) \\
& = &
\bigvee_{i} (e_{i} \wedge x_{i} \wedge y_{i}) \\
& = &
\instr_{\vec{e}}(x_{1}\wedge y_{1}, \ldots, x_{n}\wedge y_{n}) \\
\end{array}$$
}

We turn to $C^*$-algebras. In a $C^*$-algebra $A$, with an $n$-test
$\vec{e} = e_{1}, \ldots, e_{n} \in [0,1]_{A}$ where $e_{1} + \cdots +
e_{n} = 1$, we have an instrument $\instr_{\vec{e}} \colon A^{n}
\rightarrow A$ given by:
\begin{equation}
\label{MeasurementCstarEqn}
\begin{array}{rcl}
\instr_{\vec{e}}(a_{1}, \ldots, a_{n})
& = &
\sqrt{e_1}\cdot a_{1}\cdot \sqrt{e_1} + \cdots + 
   \sqrt{e_n}\cdot a_{n}\cdot \sqrt{e_n}.
\end{array}
\end{equation}

\noindent This map is not only positive, but completely positive.  We
notice that if the $C^*$-algebra $A$ is commutative, then
$\instr_{\vec{e}}(a_{1}, \ldots, a_{n}) = e_{1}\cdot a_{1} + \cdots +
e_{n}\cdot a_{n}$, like for rings in~\eqref{MeasurementRingEqn}.  This
formula~\eqref{MeasurementCstarEqn} is sometimes called the
(generalised) L{\"u}ders rule, for unsharp predicates (effects). It
can be found, for instance, in~\cite[Eq.(1.3)]{BuschS98}.

Of the different properties
\eqref{MeasurementAssInjection}--\eqref{MeasurementAssSef} in
Assumption~\ref{MeasurementAss} for instruments we shall
demonstrate~\eqref{MeasurementAssCollapse}. So assume we have an
$n+m$-test, of the form $p\colon \C^{n}\times\C^{m} \rightarrow A$.
It corresponds to $n+m$ effects $e_{1}, \ldots, e_{n}, d_{1}, \ldots,
d_{m} \in [0,1]_{A}$, where $e_{i} = p(\ket{i}, \vec{0})$ and $d_{j} =
p(\vec{0}, \ket{j})$, with $\sum_{i}e_{i} + \sum_{j}d_{j} = 1$. The
test $p' = p \after (\idmap\times\,!) \colon \C^{n}\times\C
\rightarrow A$ is given by $p'(z_{1}, \ldots, z_{n}, w) = p(z_{1},
\ldots, z_{n}, w, \ldots, w)$. Hence $p'$ corresponds to the
$n+1$-tuple of effects $e_{1}, \ldots, e_{n}, \sum_{j}d_{j}$. We can
now compute the outcomes of the two paths in the diagram in
Assumption~\ref{MeasurementAss}~\eqref{MeasurementAssCollapse} in the
category $\CstarPU$. We show that the two associated maps
$A^{n}\times\C \rightarrow A$ are equal. The upper composite is:
$$\begin{array}{rcl}
\big(\instr_{p} \after (\idmap\times\,!)\big)(a_{1}, \ldots, a_{n}, w)
& = &
\instr_{\vec{e}, \vec{d}}(a_{1}, \ldots, a_{n}, w\cdot 1, \ldots, w\cdot 1) \\
& = &
\big(\sum_{i}\sqrt{e_i}\cdot a_{i}\cdot \sqrt{e_i}\big) + 
   \big(\sum_{j} \sqrt{d_j} \cdot (w\cdot 1) \cdot \sqrt{d_j}\big) \\
& = &
\big(\sum_{i}\sqrt{e_i}\cdot a_{i}\cdot \sqrt{e_i}\big) + w\cdot (\sum_{j}d_{j}).
\end{array}$$

\noindent We get the same outcome if we follow the lower path:
$$\begin{array}{rcl}
\big(\instr_{p'} \after (\idmap\times\,!)\big)(a_{1}, \ldots, a_{n}, w)
& = &
\instr_{\vec{e}, \sum_{j}d_{j}}(a_{1}, \ldots, a_{n}, w\cdot 1) \\
& = &
\big(\sum_{i}\sqrt{e_i}\cdot a_{i}\cdot \sqrt{e_i}\big) + 
   \big(\sqrt{\sum_{j} d_j} \cdot (w\cdot 1) \cdot \sqrt{\sum_{j} d_j}\big) \\
& = &
\big(\sum_{i}\sqrt{e_i}\cdot a_{i}\cdot \sqrt{e_i}\big) + w\cdot (\sum_{j}d_{j}).
\end{array}$$

\noindent The verification of property~\eqref{MeasurementAssSef} for
$C^*$-algebras will be postponed to
Corollary~\ref{MeasurementCstarSefCor} in the next section.

\auxproof{
We have to prove that the following diagram in $\CstarPU$
commutes.
$$\xymatrix@C+1pc{
A & \C^{n}\ar[l]_-{\vec{e}}\ar[d]^{!^{n}} \\
& A^{n}\ar@{..>}[ul]^{\instr_{\vec{e}}}
}$$

\noindent This works, since for $z_{i}\in\C$,
$$\begin{array}{rcl}
\big(\instr_{\vec{e}} \after \; !^{n}\big)(z_{1}, \ldots, z_{n})
& = &
\instr_{\vec{e}}(z_{1}\cdot 1, \ldots, z_{n}\cdot 1) \\
& = &
\sqrt{e_1}\cdot (z_{1}\cdot 1)\cdot \sqrt{e_1} + \cdots + 
   \sqrt{e_n}\cdot (z_{n}\cdot 1) \cdot \sqrt{e_n} \\
& = &
z_{1}\cdot e_{1} + \cdots + z_{n}\cdot e_{n} \\
& = &
z_{1}\cdot f_{\vec{e}}(\ket{1}) + \cdots + z_{n}\cdot f_{\vec{e}}(\ket{n}) \\
& = &
f_{\vec{e}}(z_{1}, \ldots, z_{n}).
\end{array}$$

We check the ``coherence'' equations. First, for $n$-tests $\vec{e}$ in $A$
and $\vec{d}$ in $B$ we have a commuting diagram:
$$\xymatrix@R-.5pc{
A\times B & & (A\times B)^{n}\ar[ll]_-{\instr_{\tuple{\vec{e},\vec{d}}}} \\
& & A^{n}\times B^{n}\ar[u]^{\cong}_{\tuple{\pi_{1}\times\pi_{1}, \ldots,
    \pi_{n}\times\pi_{n}}}\ar[ull]^{\instr_{\vec{e}}\times\instr_{\vec{d}}\quad}
}$$

\noindent Since:
$$\begin{array}{rcl}
\lefteqn{\big(\instr_{\tuple{\vec{e},\vec{d}}} \after 
    \tuple{\pi_{1}\times\pi_{1}, \ldots, \pi_{n}\times\pi_{n}}\big)
      ((a_{1}, \ldots, a_{n}), (b_{1}, \ldots b_{n}))} \\
& = &
\instr_{\tuple{\vec{e},\vec{d}}}((a_{1},b_{1}), \ldots, (a_{n},b_{n})) \\
& = &
\sum_{i} \sqrt{(e_{i},d_{i})}\cdot (a_{i},b_{i}) \cdot \sqrt{(e_{i},d_{i})} \\
& = &
\sum_{i} \tuple{\sqrt{e_{i}}\cdot a_{i} \cdot \sqrt{e_i},
   \sqrt{d_{i}}\cdot b_{i} \cdot \sqrt{d_{i}}} \\
& = &
\tuple{\sum_{i} \sqrt{e_{i}}\cdot a_{i} \cdot \sqrt{e_i},
   \sum_{i} \sqrt{d_{i}}\cdot b_{i} \cdot \sqrt{d_{i}}} \\
& = &
\tuple{\instr_{\vec{e}}(a_{1}, \ldots, a_{n}), \instr_{\vec{d}}(b_{1}, \ldots, b_{n})} \\
& = &
(\instr_{\vec{e}}\times\instr_{\vec{d}})((a_{1}, \ldots, a_{n}), (b_{1}, \ldots b_{n})).
\end{array}$$

Also, for an injection $\phi \colon \{1, \ldots, n\}
\conglongrightarrow \{1,\ldots, m\}$, we have $\tuple{\phi} =
\tuple{\pi_{\phi(i)}}_{1\leq i\leq n} \colon A^{m} \rightarrow A^{n}$
and an $m$-test $\vec{d} = \vec{e}\tuple{\phi}$ determined by
$f_{\vec{e}} \after \tuple{\phi} \colon \C^{m} \rightarrow A$ with
$$\begin{array}{rcl}
d_{j}
& = &
\left\{\begin{array}{ll}
e_{i} \quad & \mbox{if }\phi(i) = j \\
0 & \mbox{otherwise.}
\end{array}\right.
\end{array}$$

\noindent The following diagram
$$\xymatrix@R-.5pc{
A & & A^{n}\ar[ll]_-{\instr_{\vec{e}}} \\
& & A^{m}\ar[u]^{\cong}_{\tuple{\phi}}\ar[ull]^{\instr_{\vec{e}\tuple{\phi}}\quad}
}$$

\noindent commutes, since:
$$\begin{array}{rcl}
\big(\instr_{\vec{e}} \after \tuple{\phi}\big)(a_{1},\ldots,a_{m})
& = &
\sum_{i}\sqrt{e_i}\cdot a_{\phi(i)}\cdot \sqrt{e_i} \\
& = &
\sum_{j}\sqrt{d_j}\cdot a_{j}\cdot \sqrt{d_j} \\
& = &
\instr_{\vec{d}}(a_{1},\ldots,a_{m}) \\
& = &
\instr_{\vec{e}\tuple{\phi}}(a_{1},\ldots,a_{m}).
\end{array}$$
}

We briefly consider the special cases where the $C^*$-algebra is
$\B(\H)$, for a Hilbert space $\H$. The effects $[0,1]_{\B(\H)}$ of
$\B(\H)$ are the positive maps $\H\rightarrow\H$ below the
identity. Suppose our $n$-test of effects $p_{1}, \ldots, p_{n}\in
[0,1]_A$ consists of projections: $p_{i}^{2} = p_i$. Then $\sqrt{p_i}
= p_{i}$. Hence the associated instrument $\instr_{\vec{e}} \colon
\B(\H)^{n} \rightarrow \B(\H)$ takes the familiar form:
$$\begin{array}{rcl}
\instr_{\vec{p}}(A_{1}, \ldots, A_{n})
& = &
p_{1}A_{1}p_{1} + \cdots p_{n}A_{n}p_{n}.
\end{array}$$


\end{exas}

\begin{rem}
\label{MeasurementRem}
The above list of examples gives rise to the following observations.
\begin{enumerate}
\item \label{MeasurementRemInstr} For an $n$-test $\vec{e}$ in a
  $C^*$-algebra $A$ we have an instrument $\instr_{\vec{e}}(a_{1},
  \ldots, a_{n}) = \sum_{i}\sqrt{e_i}\cdot a_{1} \cdot \sqrt{e_i}$ as
  defined in~\eqref{MeasurementCstarEqn}. We briefly show how it forms
  a (discrete) instrument as defined in~\cite[\S4]{Ozawa84}, see
  also~\cite[\S5.1.2]{HeinosaariZ12}. The map $\instr_{\vec{e}}$ can
  be seen as a combination of $n$ completely positive subunital maps
  $A \rightarrow A$, namely $f_{i}(a) = \sqrt{e_i}\cdot a_{1}\cdot
  \sqrt{e_i}$. These maps $f_i$ are subunital since $0 \leq f_{i}(1) =
  e_{i} \leq 1$. In fact, let $\Sigma_{n} = \Pow(n)$ be the set of
  measurable subsets of the discrete measurable space $n$. Then we can
  define for each (measurable) subset $U\subseteq n$ a completely
  positive subunital map $f_{U} \colon A \rightarrow A$ via $f_{U}(a)
  = \sum_{i\in U}\sqrt{e_i}\cdot a_{1}\cdot \sqrt{e_i}$. As required
  for a (discrete) instrument, this mapping $U \mapsto f_{U}$ sends
  disjoint unions $U\ovee V$ to sums $f_{U} + f_{V}$. Moreover, the
  function $f_{n}$ associated with the largest (measurable) subset
  $n\subseteq n$ is unital. 

Thus we have constructed a function from a $\sigma$-algebra $\Sigma$
to the completely positive subunital maps $A \rightarrow A$ that is
additive, and sends the top element to a unital function. This is (the
$C^*$-algebraic version of) an instrument as defined
in~\cite{DaviesL70,Ozawa84,HeinosaariZ12}.



\item \label{MeasurementRemRing} We have defined instruments
  in~\eqref{MeasurementRingEqn} only for commutative rings. For a
  non-commutative ring $R$ with an $n$-test $\vec{e}$ one can define a
  map $f\colon R^{n} \rightarrow R$ as $f(x_{1}, \ldots, x_{n}) =
  \sum_{i} e_{i}\cdot x_{i} \cdot e_{i}$. However, this map does not
  preserve multiplication. Hence one may try to adapt the notion of
  `ring homomorphism' to something weaker, just like PU is a weakening
  of MIU for $C^*$-algebras. This requires a careful balancing act,
  because one would still like to have a reasonable notion of
  predicate as this new kind of map $R \rightarrow 1+1$. There is no
  obvious way to do that. This balancing act does work for
  $C^*$-algebras, by not using MIU but PU maps.

\item \label{MeasurementRemDatabase} We have used the computer science
  terminology `side-effect' for what is more commonly called `observer
  effect' in physics. It is an interesting question to what extent
  instruments, as defined above are useful to capture side-effects in
  computer science. Here is a first thought. One can imagine a
  database with state space $X$ for which one may observe certain
  content, say via a function $X \rightarrow n$, where what is
  observed is described here simply as a finite set $n$ of natural
  numbers. Associated with such an observer map $X \rightarrow n$ one
  can construct an `instrument' $X \rightarrow n\cdot X$ which has a
  side-effect on the state space $X$ of the database, for instance by
  logging the observation itself. This can be done for accountability
  reasons, like in a medical database. The instrument used for
  observation thus changes the state of the database by adding a log
  entry.

\item \label{MeasurementRemRequirements}
  Assumption~\ref{MeasurementAss} lists four coherence requirements
  for instruments. They suffice in the present setting. However, we do
  not claim that these four requirements are `complete'. Possibly,
  other formulations are more adequate, and guarantee that instruments
  are unique. This issue is still an area of active research.  For
  instance, in~\cite{ChoJWW15} instruments are (equivalently)
  described via partial `assert' maps $X \rightarrow X+1$, like in
  point~\eqref{MeasurementRemInstr}.  It is illustrated that in
  several example categories these assert maps can be obtained from
  quotient and comprehension maps, but the precise defining property
  remains elusive. Recent, unpublished research shows that in
  `Boolean' and in `classical' effectuses these assert/instrument maps
  are uniquely determined; but the `quantum' case remains open.
\end{enumerate}
\end{rem}

\noindent We conclude this section by showing how instruments $\instr$ as in
Assumption~\ref{MeasurementAss} give rise to a natural programming
language construct, namely a general \emph{test map} in the base
category of computations: for an $n$-test $p\colon X \rightarrow
n\cdot 1$ and $n$ parallel maps $f_{1}, \ldots, f_{n} \colon X
\rightarrow Y$ we can define a test map $X \rightarrow Y$ as:
\begin{equation}
\label{TestMapEqn}
\xymatrix@C-1pc{
p?[f_{1}, \ldots, f_{n}] = \big(X\ar[rr]^-{\instr_p} & &
    n\cdot X = X + \cdots + X\ar[rrr]^-{[f_{1}, \ldots, f_{n}]} & & & Y\big).
}
\end{equation}

\noindent The idea is that this test map $p?[f_{1}, \ldots, f_{n}]$
first performs a $p$-test, and then, depending on the outcome, does a
suitable combination of the maps $f_{1}, \ldots, f_{n}$. This
combination depends on the coproduct options in $n\cdot X$ that $p$
produces and acts on the state after $p$, incorporating the
side-effect of $p$, if any. 


Via the bijective correspondence~\eqref{TestCor} between an $n$-test
$p\colon X \rightarrow n\cdot 1$ and $n$ separate predicates $p_{1},
\ldots, p_{n} \colon X \rightarrow 1+1$ with $p_{1} \ovee \cdots \ovee
p_{n} = 1$, we can understand this test map $p?[f_{1}, \ldots, f_{n}]$
in~\eqref{TestMapEqn} also as a guarded command in the style
of~\cite{Dijkstra75}, written alternatively as:
\begin{equation}
\label{GuardedIfEqn}
\vcenter{{\renewcommand\arraystretch{.9}\begin{array}[t]{lcl}
\textsf{if} \\
\quad | \; p_{1} & \rightarrow & f_{1} \\[-.3em]
& \vdots & \\
\quad | \; p_{n} & \rightarrow & f_{n}  \\
\textsf{fi}
\end{array}}}
\qquad\qquad
\begin{array}{c}
\mbox{or, as we shall write} \\[-.3em]
\mbox{later, in Figure~\ref{SuperDenseFig},}
\end{array}
\qquad\qquad
\vcenter{{\renewcommand\arraystretch{.9}\begin{array}[t]{lcl}
\rlap{$\textsf{begin test } x$} \\
\quad | \; p_{1} & \rightarrow & f_{1} \\[-.3em]
& \vdots & \\
\quad | \; p_{n} & \rightarrow & f_{n}  \\
\rlap{\textsf{end test}}
\end{array}}}
\end{equation}

\noindent In the set-theoretic case these $p_i$ describe a partition
of $X$, so that for each $x\in X$ precisely one $f_i$ is selected. In
the probabilistic case a suitable convex combination of the $f_i$ is
executed. And in the quantum case this test map involves a more
complicated combination using square roots of the effects, as
in~\eqref{MeasurementCstarEqn}.

Appropriate logical ``test'' operators to reason about these
guarded commands will be described in Section~\ref{TestSec}. At this
stage we observe the following obvious property:
\begin{equation}
\label{TestMapPostEqn}
\begin{array}{rcl}
g \after \big(p?[f_{1}, \ldots, f_{n}]\big)
& = &
p?[g \after f_{1}, \ldots, g \after f_{n}].
\end{array}
\end{equation}

\noindent A corresponding pre-composition property holds for
so-called pure maps, see~\eqref{TestMapPreEqn}.

\auxproof{
$$\begin{array}{rcl}
g \after p?[f_{1}, \ldots, f_{n}]
& = &
g \after [f_{1}, \ldots, f_{n}] \after \instr_{p} \\
& = &
[g \after f_{1}, \ldots, g \after f_{n}] \after \instr_{p} \\
& = &
p?[g \after f_{1}, \ldots, g \after f_{n}].
\end{array}$$
}

The convex sum $\bigovee_{i}r_{i}\omega_{i}$ of states $\omega_i$
defined in the proof of Lemma~\ref{StatConvexLem} is an instance of
this test programming construct~\eqref{TestMapEqn}:
$$\begin{array}{rcccl}
\bigovee_{i}r_{i}\omega_{i}
& = &
[\omega_{1}, \ldots, \omega_{n}] \after p
& = &
{\renewcommand\arraystretch{.9}\begin{array}[t]{lcl}
\textsf{if} \\
\quad | \; r_{1} & \rightarrow & \omega_{1} \\[-.3em]
& \vdots & \\
\quad | \; r_{n} & \rightarrow & \omega_{n}  \\
\textsf{fi}
\end{array}}
\end{array}$$

\noindent This follows from the equation $\instr_{p} = p$ for $p\colon
1 \rightarrow n\cdot 1$ in
Lemma~\ref{MeasurementAssSefLem}~\eqref{MeasurementAssSefLemOne}.

\section{Side-effect free predicates and pure maps}\label{SefSec}

One of the advantages of the current, general approach is that it
allows us to express the side-effect --- sometimes called the observer
effect --- of an instrument. We already encountered it briefly in the
previous section. Here we add some further results about
side-effect-freeness, notably in relation to commutativity in
$C^*$-algebras. Also, we introduce an abstract notion of purity for
maps.


\begin{defi}
\label{SideEffectFreeDef}
Given a choice of mesurement instruments, the \emph{side-effect} of an
$n$-test $p\colon X \rightarrow n\cdot 1$ is the endomap $X
\rightarrow X$ obtained in:
$$\xymatrix@C-.5pc{
X\ar[rr]^-{\instr_p} & & 
   X+\cdots+X\ar[rrr]^-{\nabla = [\idmap, \ldots, \idmap]} & & & X
}$$

\noindent This test $p$ is called \emph{side-effect-free} if its
side-effect map $\nabla \after \instr_p$ is the identity.
\end{defi}

Tests of the form $1 \rightarrow n\cdot 1$, including scalars for
$n=2$, and the true and false predicates $X \rightarrow 1+1$, are
side-effect-free, simply because the composite $1 \rightarrow n\cdot 1
\rightarrow 1$ is necessarily the identity. Also, $p^{\perp}$ has the
same side-effect as $p$, by the (Injection) property in
Assumption~\ref{MeasurementAss}. The side-effect associated with an
effect $e\in [0,1]_{A}$ in a $C^{*}$-algebra $A$ is the map $A
\rightarrow A$ given by $a \mapsto \sqrt{e}\cdot a\cdot \sqrt{e} +
\sqrt{1-e}\cdot a\cdot \sqrt{1-e}$.  If $A$ is commutative, this map
is the identity. One of the main contributions of this section is
establishing a connection between side-effect-freeness and
commutativity in $C^*$-algebras. We recall that commutative
$C^*$-algebra are not `quantum' but `probabilistic' models of
computation. Formulated differently, having side-effects of
predicates/tests is a characteristic feature of the quantum
world. Indeed, an observation at the quantum level may change
(disturb) the state of the system that is being observed.

But first we show that there are no side-effects in classical and
probabilistic logic.

\begin{prop}
\label{ClassProbSideEffectFreeProp}
In the categories $\Sets$, $\Kl(\Dst)$ and $\Kl(\Giry)$ all tests are
side-effect-free. In those cases, there is a bijective correspondence:
\begin{equation}
\label{SideEffecFreeCor}
\begin{prooftree}
\xymatrix{ \mbox{$n$-tests } X\ar[r]^-{p} & n\cdot 1 }
\Justifies
\xymatrix{ \mbox{maps }X\ar[r]_-{f} & n\cdot X \mbox{ with } 
   \nabla \after f = \idmap }
\end{prooftree}
\end{equation}

\noindent Downwards, the map $f$ corresponding to test $p$ is the
instrument $f = \instr_p$ described in Example~\ref{MeasurementEx}.
Upwards, the test $p$ is obtained from $f$ as in
Assumption~\ref{MeasurementAss}~\eqref{MeasurementAssSef}, namely as
$p = (n\cdot !) \after f$.
\end{prop}

\begin{myproof}
The properties in Assumption~\ref{MeasurementAss} imply the bijective
correspondence, once we have proved that tests are side-effects-free.
In the set-theoretic case it is easy to see that the
instrument~\eqref{MeasurementSetsEqn} satisfies $\nabla \after
\instr_{p} = \idmap$. We shall skip $\Kl(\Dst)$ and do the calculation
in $\Kl(\Giry)$. We must show that the instrument $\instr_{p} \colon X
\rightarrow \Giry(n\cdot X)$ from~\eqref{MeasurementGiryEqn}, for a
test $p\colon X \rightarrow \Giry(n)$, satisfies $\nabla \klafter
\instr_{p} = \idmap$, with Kleisli composition in $\Kl(\Giry)$ written
as $\klafter$. Well, for $x\in X$ and $M\in\Sigma_{X}$,
$$\hspace*{-1em}\begin{array}[b]{rcl}
\big(\nabla \klafter \instr_{p}\big)(x)(M)
& = &
\nabla_{*}\big(\instr_{p}(x)\big)(M) 
    \hspace*{3.8em} \mbox{where $(-)_*$ is Kleisli extension} \\
& = &
\int \nabla(-)(M) \intd\, \instr_{p}(x) 
   \hspace*{2em}\mbox{see~\eqref{GiryKleisliExtEqn}} \\
& = &
\int \indic{\kappa_{1}M \cup \cdots \cup \kappa_{n}M} \intd\, \instr_{p}(x) \\
& = &
\instr_{p}(x)(\kappa_{1}M \cup \cdots \cup \kappa_{n}M) \\
& = &
\instr_{p}(x)(\kappa_{1}M) + \cdots + \instr_{p}(x)(\kappa_{n}M) \\
& = &
p(x)(1) \cdot \indic{M}(x) + \cdots + p(x)(n) \cdot \indic{M}(x) 
   \qquad \mbox{see~\eqref{MeasurementGiryEqn}} \\
& = &
\big(p(x)(1) + \cdots + p(x)(n)\big) \cdot \indic{M}(x) \\
& = &
\indic{M}(x) \\
& = &
\eta(x)(M) \\
& = &
\idmap(x)(M).
\end{array}\eqno{\qEd}$$

\auxproof{
In the set-theoretic case we have, for $p\colon X \rightarrow n$,
$$\begin{array}{rcl}
\big(\nabla \after \instr_{p}\big)(x)
& = &
\nabla(\kappa_{i}x) \quad \mbox{if } p(x) = i \\
& = &
x.
\end{array}$$

\noindent The bijective correspondence follows, since the operations
$\overline{p} = \instr_{p}$ and $\overline{f} = (n\cdot !) \after f$ are
each other's inverse:
$$\begin{array}{rcl}
\overline{\overline{p}}
& = &
(n\cdot !) \after \instr_{p} \\
& = &
p \\
\overline{\overline{f}}
& = &
\instr_{(n\cdot !) \after f} \\
& = &
f \quad \mbox{since $\nabla \after f = \idmap$, by (side-effects-free)}.
\end{array}$$

In the discrete probabilistic case the instrument $\instr_{p}
\colon X \rightarrow \Dst(n\cdot X)$ from~\eqref{MeasurementDstEqn},
for a test $p\colon X \rightarrow \Dst(n)$, satisfies $\nabla \klafter
\instr_{p} = \idmap$, with Kleisli composition in $\Kl(\Dst)$ written
as $\klafter$, since:
$$\begin{array}[b]{rcl}
\big(\nabla \klafter \instr_{p}\big)(x)
& = &
\Dst(\nabla)\big(p(x)(1)\ket{\kappa_{1}x} + \cdots + 
   p(x)(n)\ket{\kappa_{n}x}\big) \\
& = &
p(x)(1)\ket{\nabla(\kappa_{1}x)} + \cdots + p(x)(n)\ket{\nabla(\kappa_{n}x)} \\
& = &
p(x)(1)\ket{x} + \cdots + p(x)(n)\ket{x} \\
& = &
(p(x)(1) + \cdots + p(x)(n))\ket{x} \\
& = &
1\ket{x} \qquad \mbox{since } p(x)\in \Dst(n) \\
& = &
\idmap(x).
\end{array}$$
}\medskip
\end{myproof}

\noindent The following classical result in the theory of $C^*$-algebras will be
used without proof.

\begin{theorem}[Tomiyama~\cite{Tomiyama57a}]
\label{TomiyamaThm}
Let $A,B$ be $C^*$-algebras with two maps $f,g$ in:
$$\vcenter{\xymatrix@R-.5pc{
A\ar@{=}[drr]\ar[rr]^-{\text{$g$, MIU}} & & B\ar[d]^-{\text{$f$, PU}} \\
& & A
}}
\qquad\mbox{\textit{i.e.}}\qquad
\left\{{\renewcommand\arraystretch{1}\begin{array}{l}
\mbox{$g$ preserves multiplication, involution, unit (MIU)} \\
\mbox{$f$ preserves unit, positive elements (PU)} \\
f \after g = \idmap
\end{array}}\right.$$

\noindent Then $f$ is a map of modules, in the sense that, for all
$a\in A, b\in B$,
\begin{equation}
\label{TomiyamaEqn}
\begin{array}{rclcrcl}
a \cdot f(b)
& = &
f(g(a) \cdot b)
\qquad\mbox{and}\qquad
f(b)\cdot a
& = &
f(b \cdot g(a)).
\end{array}
\end{equation}
\end{theorem}

\noindent This result is also used in~\cite{FurberJ13b}. We shall
apply it here with the diagonal $\Delta = \tuple{\idmap, \ldots,
  \idmap} \colon A \rightarrow A^{n}$ for $g$. It clearly preserves
multiplication. This map $f$ is sometimes called a conditional
expectation.

%

The following is a key result in linking side-effect-freeness in
$\op{(\CstarPU)}$ to commutativity. Recall that the center
$\mathcal{Z}(A)\subseteq A$ of a $C^*$-algebra $A$ is the
sub-$C^*$-algebra of elements that commute with everything:
$\mathcal{Z}(A) = \setin{a}{A}{\allin{x}{A}{a\cdot x = x\cdot a}}$.

\begin{lem}
\label{CenterPULem}
Let $A$ be a $C^*$-algebra. There is a bijective correspondence between:
$$\begin{prooftree}
\mbox{$n$-tests }e_{1}, \ldots, e_{n} \in [0,1]_{A} \cap \mathcal{Z}(A)
   \mbox{ of effects $e_i$ in the center}
\Justifies
\mbox{positive unital maps $f \colon A^{n} \rightarrow A$
   with $f \after \Delta = \idmap$}
\end{prooftree}$$
\end{lem}


\begin{myproof}
For effects $\vec{e} = e_{1}, \ldots, e_{n}\in [0,1]_{A}$ in the
center $\mathcal{Z}(A)$ with $e_{1} + \cdots + e_{n} = 1$ one can
define a linear map $f_{\vec{e}} \colon A^{n} \rightarrow A$ by:
\begin{equation}
\label{CenterPUEqn}
\begin{array}{rcl}
f_{\vec{e}}(a_{1}, \ldots, a_{n})
& = &
\sum_{i} e_{i}\cdot a_{i}.
\end{array}
\end{equation}

\noindent We use that the effects $e_i$ are in the center in order to
show that this function $f_{\vec{e}}$ is positive: if $a_{i} =
x_{i}^{*}x_{i}$, then $f_{\vec{e}}(a_{1},\ldots,a_{n})$ is positive, being
the sum of positive elements:
$$\begin{array}{rcll}
f_{\vec{e}}(a_{1},\ldots,a_{n})
\hspace*{\arraycolsep} = \hspace*{\arraycolsep}
\sum_{i} e_{i}\cdot a_{i} 
& = &
\sum_{i} e_{i}\cdot x_{i}^{*} \cdot x_{i} \\
& = &
\sum_{i} x_{i}^{*} \cdot e_{i}\cdot x_{i} 
   & \mbox{since $e_i$ is in the center} \\
& = &
\sum_{i} x_{i}^{*} \cdot (\sqrt{e_i})^{*} \cdot \sqrt{e_i} \cdot x_{i}
   & \mbox{since $e_i$ is positive} \\
& = &
\sum_{i} (\sqrt{e_i} \cdot x_{i})^{*}\cdot (\sqrt{e_i} \cdot x_{i}).
\end{array}$$

\noindent Clearly, 
$$\begin{array}{rcccccccccl}
\big(f_{\vec{e}} \after \Delta\big)(a)
& = &
f_{\vec{e}}(a,\ldots,a)
& = &
\sum_{i} e_{i}\cdot a
& = &
(\sum_{i} e_{i}) \cdot a
& = &
1\cdot a
& = &
a.
\end{array}$$

\noindent In particular, $f_{\vec{e}}$ is unital. The effects $e_{i}$
can be retrieved from $f_{\vec{e}}$ as $e_{i} = f_{\vec{e}}(\ket{i})$,
where $\ket{i}\in A^{n}$ is the sequence with $1\in A$ at position $i$
and $0\in A$ elsewhere.

Conversely, let $f\colon A^{n} \rightarrow A$ be a positive unital map
with $f \after \Delta = \idmap$. We obtain effects $e_{i} = f(\ket{i})
\in [0,1]_{A}$ in the standard way.  And we get $f = f_{\vec{e}}$ via
Tomiyama's Theorem~\ref{TomiyamaThm}:
$$\begin{array}{rcl}
f(a_{1},\ldots, a_{n})
\hspace*{\arraycolsep} = \hspace*{\arraycolsep}
f\big(\sum_{i}\ket{i}\cdot\Delta(a_{i})\big) 
\hspace*{\arraycolsep} = \hspace*{\arraycolsep}
\sum_{i}f\big(\ket{i}\cdot\Delta(a_{i})\big) 
& = &
\sum_{i}f(\ket{i})\cdot a_{i} 
   \qquad \mbox{see~\eqref{TomiyamaEqn}} \\
& = &
\sum_{i}e_{i}\cdot a_{i} \\
& = &
f_{\vec{e}}(a_{1},\ldots, a_{n}).
\end{array}$$

\noindent We still need to prove that the effects $e_{i} = f(\ket{i})$
are in the center. Again, this follows from Tomiyama's Theorem:
$$\begin{array}[b]{rcll}
e_{i}\cdot x
\hspace*{\arraycolsep} = \hspace*{\arraycolsep}
f(\ket{i})\cdot x 
& = &
f\big(\ket{i}\cdot\Delta(x)\big) \qquad & \mbox{see~\eqref{TomiyamaEqn}} \\
& = &
f\big(\Delta(x)\cdot \ket{i}\big) \\
& = &
x\cdot f(\ket{i}) & \mbox{see~\eqref{TomiyamaEqn} again} \\
& = &
x\cdot e_{i}.
\end{array}\eqno{\qEd}$$
\end{myproof}

\begin{cor}
\label{CenterSefCor}
Let $\vec{e} = e_{1}, \ldots, e_{n}\in [0,1]_{A}$ be an $n$-test in a
$C^*$-algebra $A$. Then:
$$\begin{array}{rcl}
\vec{e} \mbox{ is side-effect-free}
& \Longleftrightarrow &
\mbox{each effect $e_i$ is in the center }\mathcal{Z}(A).
\end{array}$$
\end{cor}

\begin{myproof}
First we observe that if $e\in [0,1]_{A} \cap \mathcal{Z}(A)$, then
also $\sqrt{e}\in [0,1]_{A} \cap \mathcal{Z}(A)$. That the square root
$\sqrt{e}$ is an effect is obvious. It is also in the center
$\mathcal{Z}(A) \hookrightarrow A$ since this center is
sub-$C^*$-algebra, so the construction of $\sqrt{e}$ for
$e\in\mathcal{Z}(A)$ via the functional calculus takes place entirely
within $\mathcal{Z}(A)$. Secondly we recall that the instrument
$\instr_{\vec{e}} \colon A^{n} \rightarrow A$ for the $n$-test
$\vec{e}$ is given by $\instr_{\vec{e}}(\vec{a}) = \sum_{i}\sqrt{e_i}
\cdot a_{i} \cdot \sqrt{e_i}$, see~\eqref{MeasurementCstarEqn}.

For the implication $(\Longleftarrow)$ assume that the effects $e_i$
are in the center $\mathcal{Z}(A)$. Then so are the $\sqrt{e_i}$, and thus:
$$\begin{array}{rcccccccl}
\instr_{\vec{e}}(\vec{a})
& = &
\sum_{i}\sqrt{e_i} \cdot a_{i} \cdot \sqrt{e_i}
& = &
\sum_{i}\sqrt{e_i} \cdot \sqrt{e_i} \cdot a_{i} 
& = &
\sum_{i}e_{i} \cdot a_{i} 
& = &
f_{\vec{e}}(\vec{a}),
\end{array}$$

\noindent where the function $f_{\vec{e}}$ is as defined
in~\eqref{CenterPUEqn}. This yields $\instr_{\vec{e}} \after \Delta =
f_{\vec{e}} \after \Delta = \idmap$, as shown in the proof of
Lemma~\ref{CenterPULem}, so that the $n$-test $\vec{e}$ is
side-effect-free.

For $(\Longrightarrow)$, if $\instr_{\vec{e}} \after \Delta = \idmap$,
then by Lemma~\ref{CenterPULem} $\instr_{\vec{e}} = f_{\vec{d}}$ where
$d_{i} \in [0,1]_{A} \cap \mathcal{Z}(A)$ is defined as: $d_{i} =
\instr_{\vec{e}}(\ket{i}) = \sqrt{e_i}\cdot 1 \cdot \sqrt{e_i} = e_i$. \QED
\end{myproof}

There is another conclusion we can draw at this stage: the instruments
$\instr_{\vec{e}}$ for $C^*$-algebras defined
in~\eqref{MeasurementCstarEqn} satisfy the `side-effect' property from
Assumption~\ref{MeasurementAss}~\eqref{MeasurementAssSef}. This was
postponed at the end of Example~\ref{MeasurementEx}.

\begin{cor}
\label{MeasurementCstarSefCor}
The category $\op{(\CstarPU)}$ satisfies the measurement
Assumption~\ref{MeasurementAss}; in particular, the Side-effect-free
property holds.
\end{cor}

\begin{myproof}
Given a positive unital map $f\colon A^{n} \rightarrow A$ with $f
\after \Delta = \idmap$, we obtain an $n$-test $q = f \after \, !^{n}
\colon \C^{n} \rightarrow A$ corresponding to effects $e_{i} =
q(\ket{i}) = f(\ket{i})$ in the center. But then $f = \instr_{\vec{e}}$
as shown in the previous proof, and thus $\instr_{q} = \instr_{\vec{e}}
= f$. \QED
\end{myproof}

For our next result it is useful to make the following observation
explicit.

\begin{lem}
\label{EffectCommutationLem}
Let $a\in A$ be an element of a $C^*$-algebra $A$. Then:
$$\begin{array}{rcl}
\allin{x}{A}{a\cdot x = x\cdot a}
& \Longleftrightarrow &
\allin{e}{[0,1]_A}{a\cdot e = e\cdot a}.
\end{array}$$

\noindent In words: $a$ commutes with all elements $x\in A$ if and
only if it commutes with all effects $e\in[0,1]_{A}$.
\end{lem}

\begin{myproof}
The direction $(\Rightarrow)$ is obvious, so we concentrate on
$(\Leftarrow)$.  A standard result in the theory of $C^*$-algebras
says that we can write each $x\in A$ as linear combination $x = x_{1}
- x_{2} +ix_{3} - ix_{4}$ of positive elements $x_{i} \geq 0$. By
suitable scaling we can write $x = \frac{1}{r}(rx_{1} - rx_{2}
+irx_{3} - irx_{4})$, where $rx_{i}\leq 1$ is an effect, for $r >
0$. By assumption, each $rx_{i}\in [0,1]_{A}$ commutes with $a$. But
then $x$ itself also commutes with $a$ since:
$$\begin{array}[b]{rcl}
a\cdot x
& = &
a\cdot \frac{1}{r}(rx_{1} - rx_{2} +irx_{3} - irx_{4}) \\
& = &
\frac{1}{r}(a\cdot rx_{1} - a\cdot rx_{2} + ia\cdot rx_{3}) - 
   ia\cdot rx_{4}) \\
& = &
\frac{1}{r}(rx_{1}\cdot a - rx_{2}\cdot a +irx_{3}\cdot a - irx_{4}\cdot a) \\
& = &
\frac{1}{r}(rx_{1} - rx_{2} +irx_{3} - irx_{4})\cdot a \\
& = &
x \cdot a.
\end{array}\eqno{\qEd}$$
\end{myproof}

\begin{cor}
\label{CCstarSefCor}
A $C^*$-algebra is commutative if and only if all its tests are
side-effect-free.
\end{cor}

\begin{myproof}
If a $C^*$-algebra $A$ is commutative, then $\mathcal{Z}(A) = A$, so
in particular each effect $e\in [0,1]_{A}$ is in the center, so that
each test $\vec{e}$ is side-effect-free by
Corollary~\ref{CenterSefCor}. Conversely, if each test is
side-effect-free, then each effect $e\in [0,1]_{A}$, considered as
predicate (\textit{i.e.}~2-test), is in the center.  An arbitrary
element $a\in A$ thus commutes with each effect $e$. But this $a$ then
commutes with all other elements, by
Lemma~\ref{EffectCommutationLem}. This makes $A$ commutative. \QED
\end{myproof}

In~\cite{CoeckeP08} so-called von Neumann projective measurements on
Hilbert spaces are described as (Eilenberg-Moore) coalgebras of a
comonad. This result can be adapted to the current setting, but
requires side-effect-freeness.

\begin{defi}
\label{vNeumannTestDef}
A \emph{von Neumann} $n$-test is an $n$-test $e_{1}, \ldots, e_{n}$ of
effects in a $C^{*}$-algebra satisfying, for each $i,j$,
\begin{equation}
\label{vNeumannTestEqn}
\begin{array}{rcl}
e_{i} \cdot e_{j}
& = &
\left\{\begin{array}{ll}
e_{i} \quad & \mbox{if }i=j \\
0 & \mbox{otherwise.}
\end{array}\right.
\end{array}
\end{equation}

\noindent Such test will be called a \emph{central von Neumann} test
if each $e_i$ is in the center.
\end{defi}

Notice that each effect $e_i$ in a von Neumann test is a projection:
$e_{i}^{2} = e_{i}$. As a result, $\sqrt{e_i} = e_{i}$.

In~\cite{Jacobs13b} it is observed that, in general, for a category
with coproducts $+$, for each $n>0$ taking $n$-fold copowers $X
\mapsto n\cdot X = X+\cdots+X$ is a comonad, with counit $\varepsilon$
and comultiplication $\delta$:
$$\xymatrix@C+.5pc{
n\cdot X\ar[rr]^-{\varepsilon = \nabla = [\idmap,\ldots,\idmap]} & & X
&
n\cdot X\ar[rr]^-{\delta = \kappa_{1}+\cdots+\kappa_{n}} & & n\cdot (n\cdot X)
}$$

\noindent We shall be using this comonad on the opposite category
$\op{(\CstarPU)}$. It forms a monad $(-)^{n}$ on $\CstarPU$, with unit
$\Delta$ and multiplication $\pi_{1}+\cdots+\pi_{n}$. We then have the
following result, describing the instruments of central von Neumann
tests as coalgebras of a comonad.

\begin{theorem}
\label{vNeumannTestComonadProp}
For a $C^*$-algebra $A$ there is a bijective correspondence between:
$$\begin{prooftree}
\mbox{central von Neumann $n$-tests }e_{1}, \ldots, e_{n} \in [0,1]_{A}
\Justifies
\mbox{coalgebras $f \colon A \rightarrow n\cdot A$ in $\op{(\CstarPU)}$}
\end{prooftree}$$

\noindent This is a restriction of the bijective correspondence
in Lemma~\ref{CenterPULem}.
\end{theorem}

\begin{myproof}
Given a central von Neumann test $\vec{e} = e_{1}, \ldots, e_{n}$ we
already know that the corresponding function $f_{\vec{e}}$
from~\eqref{CenterPUEqn} satisfies $f_{\vec{e}} \after \Delta =
\idmap$. We show that the orthogonal projection
requirement~\eqref{vNeumannTestEqn} is equivalent to the (co)algebra
equation:
\begin{equation}
\label{vNeumannTestMuEqn}
\begin{array}{rcl}
f_{\vec{e}} \after (f_{\vec{e}})^{n}
& = &
f_{\vec{e}} \after (\pi_{1}\oplus\cdots\oplus\pi_{n})
  \;\colon\; \big(A^{n}\big)^{n} \longrightarrow A.
\end{array}
\end{equation}

\noindent This equivalence is the essence
of~\cite[Thm.~16.6]{CoeckeP08}.  Here it occurs in $C^*$-algebraic
terms.

So assume $\vec{e}$ is a (central) von Neumann test. Then for $n$
sequences $\vec{a_i} = (a_{i1}, \ldots, a_{in})\in A^{n}$ we can
prove~\eqref{vNeumannTestMuEqn}:
$$\begin{array}{rcll}
\big(f_{\vec{e}} \after (f_{\vec{e}})^{n}\big)(\vec{a_1}, \ldots, \vec{a_n})
& = &
f_{\vec{e}}\big(f_{\vec{e}}(\vec{a_1}), \ldots, f_{\vec{e}}(\vec{a_n})\big) \\
& = &
\sum_{i} e_{i} \cdot (\sum_{j} e_{j} \cdot a_{ij}) \\
& = &
\sum_{i,j} e_{i} \cdot e_{j} \cdot a_{ij} \\
& = &
\sum_{i} e_{i} \cdot a_{ii} & \mbox{by~\eqref{vNeumannTestEqn}} \\
& = &
f_{\vec{e}}(a_{11}, \ldots, a_{nn}) \\
& = &
f_{\vec{e}}\big(\pi_{1}(\vec{a_1}), \ldots, \pi_{n}(\vec{a_n})\big) \\
& = &
\big(f_{\vec{e}} \after (\pi_{1}\oplus\cdots\oplus\pi_{n})\big)
   (\vec{a_1}, \ldots, \vec{a_n}).
\end{array}$$

\noindent In the reverse direction, assuming~\eqref{vNeumannTestMuEqn} we
prove~\eqref{vNeumannTestEqn}. Take the sequence of sequences:
$$\begin{array}{rcl}
s
& = &
(\vec{0}, \ldots, \vec{0}, \ket{j}, \vec{0}, \ldots, \vec{0})
   \;\in\; \big(A^{n}\big)^{n}
\end{array}$$

\noindent with $\ket{j}\in A^{n}$ at the $i$-th position. Then:
$$\begin{array}[b]{rcll}
e_{i} \cdot e_{j}
& = &
e_{i} \cdot f_{\vec{e}}(\ket{j}) \\
& = &
f_{\vec{e}}\big(0, \ldots, 0, f_{\vec{e}}(\ket{j}), 0, \ldots, 0\big) \\
& = &
f_{\vec{e}}\big(f_{\vec{e}}(\vec{0}), \ldots, f_{\vec{e}}(\vec{0}), 
   f_{\vec{e}}(\ket{j}), f_{\vec{e}}(\vec{0}), \ldots, f_{\vec{e}}(\vec{0})\big) \\
& = &
\big(f_{\vec{e}} \after (f_{\vec{e}})^{n}\big)(s) \\
& = &
\big(f_{\vec{e}} \after (\pi_{1}\oplus\cdots\oplus\pi_{n}))(s) 
   & \mbox{by~\eqref{vNeumannTestMuEqn}} \\
& = &
\left\{\begin{array}{ll}
f_{\vec{e}}(0, \ldots, 0, 1, 0 \ldots, 0) \quad & \mbox{if }i=j \\
f_{\vec{e}}(0, \ldots, 0, 0, 0 \ldots, 0) & \mbox{otherwise}
\end{array}\right. \\
& = &
\left\{\begin{array}{ll}
e_{i} \quad & \mbox{if }i=j \\
0 & \mbox{otherwise.}
\end{array}\right.
\end{array}\eqno{\qEd}$$
\end{myproof}

\noindent We conclude this section by introducing `pure' maps, which commute
appropriately with instruments. 


\begin{defi}
\label{PureDef}
In an effectus with chosen instruments, as in
Assumption~\ref{MeasurementAss}, we call a map $f\colon X \rightarrow
Y$ \emph{pure} if for each $n$-test $q\colon Y \rightarrow n\cdot 1$
the following diagram commutes.
$$\xymatrix@R-.5pc{
X\ar[d]_{f}\ar[rr]^-{\instr_{f^{*}(q)}} & & n\cdot X\ar[d]^{n\cdot f} \\
Y\ar[rr]_{\instr_{q}} & & n\cdot Y
}$$

\noindent This makes $f$ a map of coalgebras, for the functor $n\cdot
(-)$.
\end{defi}

We should point out that `purity' depends on the instruments that
are used. The next result shows that by restricting ourselves to
pure maps we still have a category with coproducts.

\begin{lem}
\label{PureCatLem}
The pure maps form a subcategory: the identity map is pure, and if $f,
g$ are composable pure maps, then $g \after f$ is pure too. Also, the
coprojections are pure, and if $f_{i} \colon X_{i} \rightarrow Y$ is
pure for $i=1,2$, then so is the cotuple $[f_{1},f_{2}]$.
\end{lem}

\begin{myproof}
Closure under composition is obvious. We shall prove that
coprojections are pure. Closure under cotuples is done similarly,
using
Assumption~\ref{MeasurementAss}~\eqref{MeasurementAssCotuple}. Consider
a coprojection $\kappa_{1} \colon X \rightarrow X+Y$ and an $n$-test
$q = [q_{1}, q_{2}] \colon X+Y \rightarrow n\cdot 1$ on its codomain.
Then:
$$\begin{array}[b]{rcl}
\instr_{q} \after \kappa_{1}
& = &
[\kappa_{1}+\kappa_{1}, \ldots, \kappa_{n}+\kappa_{n}]^{-1} \after
   (\instr_{q_1} + \instr_{q_2}) \after \kappa_{1} 
   \qquad\mbox{by 
   Assumption~\ref{MeasurementAss}~\eqref{MeasurementAssCotuple}} \\
& = &
[n\cdot \kappa_{1}, n\cdot \kappa_{2}] \after \kappa_{1} \after
   \instr_{q_1} \\
& = &
(n\cdot \kappa_{1}) \after \instr_{\kappa_{1}^{*}(q)}.\rlap{\hbox to
      269 pt{\hfill\qEd}}
\end{array}$$
\auxproof{
We first prove closure under composition.  Suppose $f\colon X
\rightarrow Y$ and $g\colon Y \rightarrow Z$ are both pure. We show
that $g \after f$ is also pure. So assume a test $q \colon Z
\rightarrow n\cdot 1$. Then $g^{*}(q) = q \after g$ is also a test,
and:
$$\begin{array}{rcl}
(n\cdot (g\after f)) \after \instr_{(g \after f)^{*}(q)} 
& = &
(n \cdot g) \after (n\cdot f) \after \instr_{f^{*}g^{*}(q)} \\
& = &
(n\cdot g) \after \instr_{g^{*}(q)} \after f \\
& = &
\instr_{q} \after g \after f.
\end{array}$$
Next, for closure under cotuples, assume pure maps $f_{i} \colon X_{i}
\rightarrow Y$ with an $n$-test $q\colon Y \rightarrow n\cdot 1$.
$$\begin{array}{rcl}
\instr_{q} \after [f,g]
& = &
[\instr_{q} \after f, \instr_{q} \after g] \\
& = &
[(n\cdot f) \after \instr_{f^{*}(q)}, (n\cdot g) \after \instr_{g^{*}(q)}] 
   \qquad \mbox{because $f,g$ are pure} \\
& = &
[n\cdot f, n\cdot g] \after (\instr_{q\after f} + \instr_{q\after g}) \\
& = &
(n\cdot [f,g]) \after [n\cdot\kappa_{1}, n\cdot \kappa_{2}] 
   \after (\instr_{q\after f} + \instr_{q\after g}) \\
& = &
(n\cdot [f,g]) \after [\kappa_{1}+\kappa_{1}, \ldots, \kappa_{n}+\kappa_{n}]^{-1} 
   \after (\instr_{q\after f} + \instr_{q\after g}) \\
& = &
(n\cdot [f,g]) \after \instr_{[q \after f, q \after g]}
   \qquad\mbox{by 
   Assumption~\ref{MeasurementAss}~\eqref{MeasurementAssCotuple}} \\
& = &
(n\cdot [f,g]) \after \instr_{q \after [f, g]} \\
& = &
(n\cdot [f,g]) \after \instr_{[f, g]^{*}(q)}.
\end{array}$$}
\end{myproof}

\noindent The next result describes pure maps in several examples. Note that we
only give sufficient conditions, and no characterisations of pure maps
in these cases.

\begin{prop}
\label{PureExProp}
In our standard examples we have the following situation.
\begin{enumerate}
\item \label{PureExPropSets} In $\Sets$ all maps are pure.

\item \label{PureExPropKlD} In $\Kl(\Dst)$ all maps in the image of
  the inclusion $J\colon \Sets \rightarrow \Kl(\Dst)$ are pure.

\item \label{PureExPropKlG} Similarly, in $\Kl(\Giry)$ all maps in the
  image of the inclusion $J\colon \Meas \rightarrow \Kl(\Giry)$ are
  pure.

\item \label{PureExPropCstar} In $\CstarPU$ all MIU-maps are pure,
  \textit{i.e.}~all maps in the image of $\CstarMIU \rightarrow
  \CstarPU$.
\end{enumerate}
\end{prop}

\begin{myproof}
For $f\colon X \rightarrow Y$ in $\Sets$ and $q\colon Y \rightarrow
n$ we have:
$$\begin{array}{rcl}
\big(n\cdot f \after \instr_{q \after f}\big)(x) = \kappa_{i}y
& \Longleftrightarrow &
\instr_{q \after f}(x) = \kappa_{i}x \mbox{ and } 
   (n\cdot f)(\kappa_{i}x) = \kappa_{i}y \\
& \Longleftrightarrow &
q(f(x)) = i \mbox{ and } f(x) = y \\
& \Longleftrightarrow &
\instr_{q}(f(x)) = \kappa_{i}f(x) \mbox{ and } f(x) = y \\
& \Longleftrightarrow &
\big(\instr_{q} \after f\big)(x) = \kappa_{i}y.
\end{array}$$\medskip

\noindent In $\Kl(\Dst)$ a map $f\colon X \rightarrow Y$ of the form $f = \eta
\after g$ is pure, for a function $g\colon X \rightarrow Y$ in
$\Sets$, where $\eta$ is the unit $Y \rightarrow \Dst(Y)$ given by
$\eta(y) = 1\ket{y}$. Then $n\cdot f \colon n\cdot X \rightarrow
\Dst(n\cdot Y)$ is given by $(n\cdot f)(\kappa_{i}x) =
1\ket{\kappa_{i}g(x)}$. And for a predicate $q\colon Y \rightarrow
\Dst(n)$ we have $q \klafter f = q \after g$. Hence:
$$\begin{array}{rcll}
\big(n\cdot f \klafter \instr_{q \klafter f}\big)(x)
& = &
(n\cdot f)_{*}(\sum_{i} q(g(x))(i)\ket{\kappa_{i}x}) \quad
   & \mbox{with $(-)_*$ from~\eqref{DstKleisliExtEqn}} \\
& = &
\sum_{i} q(g(x))(i)\ket{\kappa_{i}g(x)} \\
& = &
\instr_{q}(g(x)) \\
& = &
\big(\instr_{q} \klafter f\big)(x).
\end{array}$$

\auxproof{
For a map $f\colon X \rightarrow Y$ in $\Meas$ and $q\colon Y
\rightarrow \Giry(n\cdot Y)$ we have $\instr_{q}(y)(\kappa_{i}N) =
q(y)(i) \cdot \indic{N}(y)$. Further,
$$\begin{array}{rcccl}
(\eta\after f)^{*}(q)(x)
& = &
(q \klafter (\eta \after f))(x)
& = &
q(f(x)).
\end{array}$$

\noindent Thus:
$$\begin{array}{rcl}
\big((n\cdot f) \klafter \instr_{q \after f}\big)(x)(\kappa_{i}N)
& = &
\Giry(n\cdot f)(\instr_{q \after f}(x))(\kappa_{i}N) \\
& = &
\instr_{q \after f}(x)((n\cdot f)^{-1}(\kappa_{i}N)) \\
& = &
\instr_{q \after f}(x)(\kappa_{i}f^{-1}(N)) \\
& = &
q(f(x))(i)\cdot \indic{f^{-1}(N)}(x) \\
& = &
q(f(x))(i)\cdot \indic{N}(f(x)) \\
& = &
\instr_{q}(f(x))(\kappa_{i}N) \\
& = &
\big(\instr_{q} \after f\big)(x)(\kappa_{i}N) \\
& = &
\big(\instr_{q} \klafter (\eta \after f)\big)(x)(\kappa_{i}N)
\end{array}$$
}

\noindent The Giry monad $\Giry$ is handled similarly. Finally, in the
case of $C^*$-algebras, a MIU-map $f\colon B \rightarrow A$ preserves
not only multiplication but also square roots. Hence, for an $n$-test
$\C^{n} \rightarrow B$, corresponding to effects $e_{i} = q(\ket{i})
\in B$ with $\bigovee_{i}e_{i} = 1$, we have:
$$\begin{array}[b]{rcl}
\big(\instr_{f \after \vec{e}} \after f^{n}\big)(b_{1}, \ldots, b_{n})
& = &
\instr_{\vec{f(e_{i})}}(f(b_{1}), \ldots, f(b_{n})) \\
& = &
\sum_{i} \sqrt{f(e_{i})} \cdot f(b_{i}) \cdot \sqrt{f(e_{i})} \\
& = &
\sum_{i} f(\sqrt{e_{i}}) \cdot f(b_{i}) \cdot f(\sqrt{e_{i}}) \\
& = &
\sum_{i} f(\sqrt{e_{i}} \cdot b_{i} \cdot \sqrt{e_{i}}) \\
& = &
f(\sum_{i} \sqrt{e_{i}} \cdot b_{i} \cdot \sqrt{e_{i}}) \\
& = &
\big(f \after \instr_{\vec{e}}\big)(b_{1}, \ldots, b_{n}).
\end{array}\eqno{\qEd}$$
\end{myproof}

\noindent Below we give an illustration of a map that is not of the form
$\eta\after f$ in $\Kl(\Dst)$ and which is not pure. Thus the
condition in point~\eqref{PureExPropKlD} is necessary. It is not clear
if the MIU condition in point~\ref{PureExPropCstar} is also
necessary. For states, it is known that a MIU-map $\omega \colon A
\rightarrow \C$ is pure in the sense that each linear map $f\colon A
\rightarrow \C$ for which both $f$ and $\omega-f$ are positive is of
the form $f = r\cdot \omega$ for some $r\in [0,1]$. This latter
formulation is equivalent to $\omega$ being an extreme point in the
convex set of PU maps $A\rightarrow\C$.

\begin{exa}
\label{KlDNonPureStateEx}
Let state $\omega\colon 1 \rightarrow \Dst(\{a,b\})$ in $\Kl(\Dst)$ be
given by: $\omega = \frac{1}{4}\ket{a} + \frac{3}{4}\ket{b}$ and
predicate $q\colon \{a,b\} \rightarrow [0,1]$ by $q(a) = \frac{2}{5},
q(b) = \frac{1}{3}$. The instrument $\instr_{q} \colon \{a,b\}
\rightarrow \Dst(\{a,b\} + \{a,b\})$ is:
$$\begin{array}{rclcrcl}
\instr_{q}(a)
& = &
\frac{2}{5}\ket{\kappa_{1}a} + \frac{3}{5}\ket{\kappa_{2}a}
& \qquad &
\instr_{q}(b)
& = &
\frac{1}{3}\ket{\kappa_{1}b} + \frac{2}{3}\ket{\kappa_{2}b}.
\end{array}$$

\noindent Thus $\instr_{q} \klafter \omega \colon 1 \rightarrow
\Dst(\{a,b\} + \{a,b\})$ is:
$$\begin{array}{rcl}
\instr_{q} \klafter \omega
& = &
\frac{1}{10}\ket{\kappa_{1}a} + \frac{3}{20}\ket{\kappa_{2}a} + 
   \frac{1}{4}\ket{\kappa_{1}b} + \frac{1}{2}\ket{\kappa_{2}b}.
\end{array}$$

\noindent On the other hand, $\omega^{*}(q) \colon 1 \rightarrow
          [0,1]$ is given by $\omega^{*}(q) =
          \frac{1}{4}\cdot\frac{2}{5} + \frac{3}{4}\cdot\frac{1}{3} =
          \frac{1}{10} + \frac{1}{4} = \frac{7}{20}$. The associated
          instrument $\instr_{\omega^{*}(q)} \colon 1 \rightarrow
          \Dst(1+1)$ is $\instr_{\omega^{*}(q)} = \frac{7}{20}\ket{0} +
            \frac{13}{20}\ket{1}$.  Then $(\omega+\omega) \klafter
            \instr_{\omega^{*}(q)} \colon 1 \rightarrow \Dst(\{a,b\} +
            \{a,b\})$ is the distribution
$$\begin{array}{rcl}
(\omega+\omega) \klafter \instr_{\omega^{*}(q)}
& = &
\frac{7}{80}\ket{\kappa_{1}a} + \frac{21}{80}\ket{\kappa_{1}b}
    + \frac{13}{80}\ket{\kappa_{2}a} + \frac{39}{80}\ket{\kappa_{2}b}.
\end{array}$$

\noindent Thus $\instr_{q} \klafter \omega \neq (\omega+\omega)
\klafter \instr_{\omega^{*}(q)}$, so that $\omega$ is not pure.
\end{exa}

Finally, we observe that for a \emph{pure} map $g$ one has the
following pre-composition property for test maps, as defined
in~\eqref{TestMapEqn}
\begin{equation}
\label{TestMapPreEqn}
\begin{array}{rcl}
\big(p?[f_{1}, \ldots, f_{n}]\big) \after g
& = &
(g^{*}(p))?[f_{1} \after g, \ldots, f_{n} \after g].
\end{array}
\end{equation}

\auxproof{
$$\begin{array}{rcl}
\big(p?[f_{1}, \ldots, f_{n}]\big) \after g
& = &
[f_{1}, \ldots, f_{n}] \after \instr_{p} \after g \\
& = &
[f_{1}, \ldots, f_{n}] \after n\cdot g \after \instr_{g^{*}(p)} \\
& = &
[f_{1} \after g, \ldots, f_{n} \after g] \after \instr_{g^{*}(p)} \\
& = &
(g^{*}(p))?[f_{1} \after g, \ldots, f_{n} \after g].
\end{array}$$
}

\section{Test operators in dynamic logic}\label{TestSec}

At this stage both our assumptions will be used to define test
operators on predicates. It will provide the basis for a dynamic logic
whose semantic basis will be investigated in this section. We do not
provide a proper logic with deduction rules, but only describe what
holds in categorical models.

\begin{defi}
\label{TestOperatorDef}
Let $\cat{B}$ be an effectus with a choice of instruments $\instr_{p}$
as in Assumptions~\ref{CoprodAss} and~\ref{MeasurementAss}. For two
predicates $p,q\colon X \rightarrow 1+1$ on the same object
$X\in\cat{B}$ we define two new ``test'' predicates on $X$, namely:
$$\begin{array}{rcccll}
\andthen{p}{q}
& = &
(\instr_{p})^{*}(\FstAnd(q))
& = &
[q, 0] \after \instr_{p} 
   & \colon X \longrightarrow X+X \longrightarrow 1+1 \\
\implies{p}{q}
& = &
(\instr_{p})^{*}(\FstThen(q))
& = &
[q, 1] \after \instr_{p}
   & \colon X \longrightarrow X+X \longrightarrow 1+1.
\end{array}$$

\noindent We call $\cat{B}$ a \emph{commutative} model if
$\andthen{p}{q} = \andthen{q}{p}$ for all predicates $p,q$. Also, we
call a predicate $p$ a projection if $\andthen{p}{p} = p$. Such
projections are sometimes called `sharp' predicates, and, in contrast,
effects are then the `unsharp' ones.
\end{defi}

We pronounce $\andthen{p}{q}$ as ``test $p$, and then $q$'', and
$\implies{p}{q}$ as ``test $p$, then $q$''. These $\andthen{p}{q}$ and
$\implies{p}{q}$ are test operators, like in dynamic logic,
see~\cite{HarelKT00}. Using the test map notation
from~\eqref{GuardedIfEqn} we may also describe these test operators
as:
$$\begin{array}{rclcrcl}
\andthen{p}{q}
& = &
{\renewcommand\arraystretch{.9}\begin{array}[t]{lcl}
\textsf{if} \\
\quad | \; p & \rightarrow & q \\
\quad | \; p^{\perp} & \rightarrow & 0  \\
\textsf{fi}
\end{array}}
& \qquad\qquad &
\implies{p}{q}
& = &
{\renewcommand\arraystretch{.9}\begin{array}[t]{lcl}
\textsf{if} \\
\quad | \; p & \rightarrow & q \\
\quad | \; p^{\perp} & \rightarrow & 1  \\
\textsf{fi}
\end{array}}
\end{array}$$

\noindent These notations emphasise the difference between $p$ and $q$
in these test operators, since $p$ is used for an action, namely
measurement, via the associated instrument $\instr_p$, after which $q$
is evaluated. In dynamic logic the test $p?$ has no side-effect, that
is, it does not change any state. But the current setting such a test
$p?$ may have a side-effect, see explicitly in
Corollary~\ref{TwoTestSefCor} below. Therefore, commutativity does not
hold in general, as the $C^*$-algebra example below demonstrates. We
now review the meaning of $\andthen{p}{q}$ and $\implies{p}{q}$ in our
running examples.  We will see that familiar logical operations emerge
from the above general definitions.

\begin{exas}
\label{TestOperatorEx}
In the category $\Sets$, for predicates $p,q\colon X \rightarrow 2$ we
have:
$$\begin{array}{rcccccl}
\andthen{p}{q}(x)
& = &
\La{q}(\instr_{p}(x)) 
& = &
\left\{\begin{array}{ll}
q(x) \;\; & \mbox{if }p(x) = 1 \\
0 & \mbox{otherwise}
\end{array}\right\}
& = &
p(x)\cdot q(x) \\
\implies{p}{q}(x)
& = &
\Ra{q}(\instr_{p}(x)) 
& = &
\left\{\begin{array}{ll}
q(x) \;\; & \mbox{if }p(x) = 1 \\
1 & \mbox{otherwise}
\end{array}\right\}
& = &
p(x)\cdot q(x) + (1-p(x)) \\
\end{array}$$

\noindent These formulas describe ordinary conjunction and implication
in terms of characteristic functions. Since $p(x) \in \{0,1\}$ we have
$p(x)\cdot p(x) = p(x)$, so that each predicate is a projection.  This
set-theoretic model is clearly commutative: $\andthen{p}{q} =
\andthen{p}{q}$ holds since multiplication on $\{0,1\}$ is
commutative. 

In the Kleisli category $\Kl(\Dst)$ of the distribution monad we get
the same formulas, but this time interpreted not in the set of
Booleans $\{0,1\}$ but in the unit interval $[0,1]$. For $p,q \colon X
\rightarrow [0,1]$ we elaborate $\andthen{p}{q}$ as function $X
\rightarrow \Dst(2)$. We first write the cotuple $[q, 0] \colon X+X
\rightarrow \Dst(2)$ as $\La{q}(\kappa_{1}x) = q(x)\ket{1} +
(1-q(x))\ket{0}$ and $\La{q}(\kappa_{2}x) = 1\ket{0}$. Then:
$$\begin{array}{rcl}
\andthen{p}{q}(x)
& = &
(\La{q} \klafter \instr_{p})(x) \\
& = &
\La{q}_{*}\big(p(x)\ket{\kappa_{1}x} +
   (1-p(x))\ket{\kappa_{2}x}\big) 
  \qquad \mbox{with $(-)_{*}$ from~\eqref{DstKleisliExtEqn}} \\
& = &
p(x)\cdot q(x)\ket{1} + p(x)\cdot (1-q(x))\ket{0} 
   + (1-p(x))\cdot 1\ket{0} \\
& = &
p(x)\cdot q(x)\ket{1} + (1-p(x)\cdot q(x))\ket{0}.
\end{array}$$

\noindent Thus, as fuzzy predicate $\andthen{p}{q} \colon X
\rightarrow [0,1]$ we can write $\andthen{p}{q}(x) = p(x)\cdot q(x)$
like in the set-theoretic case. Hence, $\Kl(\Dst)$ is a commutative
model.

The projections $p$ are the predicates with $p(x)^{2} = p(x)$, so that
$p(x) \in \{0,1\}$, and are thus the Boolean predicates. One can check
that the then-test is given by $\implies{p}{q}(x) = p(x)\cdot q(x) +
(1-p(x))$. In probability theory the latter formula for the then-test
is called the Reichenbach implication, see~\cite{Reichenbach49}.

Precisely the same formulas hold for predicates in the Kleisli
category $\Kl(\Giry)$ of the Giry monad, for \emph{measurable} maps $X
\rightarrow \Giry(2) \cong [0,1]$ where $X$ is a measurable space,
see also~\cite[Prop.~13]{Jacobs13a}.

\auxproof{
We check for $p,q \colon X \rightarrow \Giry(2)$, $x\in X$ and $i\in 2$,
$$\begin{array}{rcl}
\andthen{p}{q}(x)(i)
& = &
\big([q,0] \klafter \instr_{p}\big)(x)(i) \\
& = &
\int [q,0](-)(i) \intd \instr_{p}(x) \\
& = &
\instr_{p}(x)(\kappa_{1}X) \cdot \int [q,0](\kappa_{1}-)(i) \intd 
   \frac{\instr_{p}(x)(\kappa_{1}-)}{\instr_{p}(x)(\kappa_{1}X)} \\
& & \qquad +
\instr_{p}(x)(\kappa_{2}X) \cdot \int [q,0](\kappa_{2}-)(i) \intd 
   \frac{\instr_{p}(x)(\kappa_{2}-)}{\instr_{p}(x)(\kappa_{2}X)} \\
& = &
p(x)(1)\cdot\indic{X}(x) \cdot \int q(-)(i) \intd 
   \frac{p(x)(1)\cdot\indic{(-)}(x)}{p(x)(1)\cdot\indic{X}(x)} \\
& & \qquad +
p(x)(2)\cdot\indic{X}(x) \cdot \int 0(-)(i) \intd 
   \frac{p(x)(2)\cdot\indic{(-)}(x)}{p(x)(2)\cdot\indic{X}(x)} \\
& = &
p(x)(1)\cdot \int q(-)(i) \intd \indic{(-)}(x) + 
   p(x)(2)\cdot \int 0(-)(i) \intd \indic{(-)}(x) \\
& = &
p(x)(1)\cdot \int q(-)(i) \intd \eta(x) + 
   p(x)(2)\cdot \int 0(-)(i) \intd \eta(x) \\
& = &
p(x)(1)\cdot q(x)(i) + p(x)(2)\cdot 0(i).
\end{array}$$

\noindent Thus:
$$\begin{array}{rcl}
\andthen{p}{q}(x)(1)
& = &
p(x)(1)\cdot q(x)(1) + p(x)(2)\cdot 0(1) \\
& = &
p(x)(1)\cdot q(x)(1) \\
\andthen{p}{q}(x)(2)
& = &
p(x)(1)\cdot q(x)(2) + p(x)(2)\cdot 0(2) \\
& = &
p(x)(1)\cdot q(x)(2) + p(x)(2) \\
& = &
p(x)(1)\cdot (1-q(x)(1)) + (1-p(x)(1)) \\
& = &
1 - p(x)(1)\cdot q(x)(1) \\
& = &
1 - \andthen{p}{q}(x)(1).
\end{array}$$

\noindent In the same way we compute:
$$\begin{array}{rcl}
\implies{p}{q}(x)(i)
& = &
p(x)(1)\cdot \int q(-)(i) \intd \eta(x) + 
   p(x)(2)\cdot \int 1(-)(i) \intd \eta(x) \\
& = &
p(x)(1)\cdot q(x)(i) + p(x)(2)\cdot 0(i).
\end{array}$$

\noindent Thus:
$$\begin{array}{rcl}
\andthen{p}{q}(x)(1)
& = &
p(x)(1)\cdot q(x)(1) + p(x)(2)\cdot 1(1) \\
& = &
p(x)(1)\cdot q(x)(1) + (1 - p(x)(1)) \\
\andthen{p}{q}(x)(2)
& = &
p(x)(1)\cdot q(x)(2) \\
& = &
p(x)(1)\cdot (1-q(x)(1)) \\
& = &
p(x)(1) - p(x)(1)\cdot q(x)(1) \\
& = &
1 - \Big(p(x)(1)\cdot q(x)(1) + (1 - p(x)(1))\Big) \\
& = &
1 - \implies{p}{q}(x)(1).
\end{array}$$
}

In the (opposite of the) category of commutative rings, assume we have
two predicates (idempotents) $r,s\in R$, considered as ring
homomorphisms $f_{r}, f_{s} \colon \Z^{2} \rightarrow R$, like in
Example~\ref{PredEx}~\eqref{PredExRng}. The instrument $\instr_{r}
\colon R^{2} \rightarrow R$ associated with $r$ is given by
$\instr_{r}(x,y) = r\cdot x + (1-r)\cdot y$. The element $\andthen{r}{s}
\in R$ is defined as $\andthen{r}{s} = f(1,0)$ for the function
$f \colon \Z^{2} \rightarrow R$ given by:
$$\begin{array}{rcl}
f(n,m)
\hspace*{\arraycolsep} = \hspace*{\arraycolsep}
\big(\instr_{r} \after \tuple{f_{s}, f_{0}}\big)(n,m)
& = &
\instr_{r}(n\cdot s + m \cdot (1-s), m\cdot 1) \\
& = &
r\cdot (n\cdot s + m \cdot (1-s)) + (1-r) \cdot (m\cdot 1) \\
& = &
n \cdot (r\cdot s) + m \cdot (1 - r\cdot s).
\end{array}$$

\noindent Hence $\andthen{r}{s} = r\cdot s$. Similarly,
$\implies{r}{s} = r\cdot s + 1 - r$. In a same manner one obtains, for
complementable elements $x,y$ in a distributive lattice,
$\andthen{x}{y} = x\wedge y$ and $\implies{x}{y} = (x\wedge y) \vee
x^{\perp} = y \vee x^{\perp}$.

\auxproof{
We have $\implies{r}{s} = f(1,0)$ for the function $f$ given by:
$$\begin{array}{rcl}
f(n,m)
\hspace*{\arraycolsep} = \hspace*{\arraycolsep}
\big(\instr_{r} \after \tuple{f_{s}, f_{1}}\big)(n,m)
& = &
\instr_{r}(n\cdot s + m \cdot (1-s), n\cdot 1) \\
& = &
r\cdot (n\cdot s + m \cdot (1-s)) + (1-r) \cdot (n\cdot 1) \\
& = &
n \cdot (r\cdot s + 1 - r) + m \cdot (r - r\cdot s).
\end{array}$$
}

In a $C^*$-algebra $A$, for effects $e,d\in [0,1]_{A}$, we can compute
$\andthen{e}{d}$ via the map $f\colon \C^{2} \rightarrow A$ given by:
$$\begin{array}{rcl}
f(z,w)
& = &
\instr_{e}(\tuple{d, 0}(z,w)) \\
& = &
\instr_{e}(z\cdot d + w\cdot (1-d), w\cdot 1) \\
& = &
\sqrt{e}\cdot (z\cdot d + w\cdot (1-d)) \cdot \sqrt{e} +
   \sqrt{1-e}\cdot (w\cdot 1) \cdot \sqrt{1-e} \\
& = &
z\cdot \sqrt{e}\cdot d \cdot \sqrt{e} +
  w\cdot \sqrt{e}\cdot (1-d) \cdot \sqrt{e} +
  w \cdot (1-e).
\end{array}$$

\noindent We get the corresponding effect by taking $z=1,w=0$. Then:
\begin{equation}
\label{CstarTestOperatorEqn}
\begin{array}{rccclcrcl}
\andthen{e}{d}
& = &
f(1,0)
& = &
\sqrt{e}\cdot d \cdot \sqrt{e}
& \quad\mbox{and similarly}\quad &
\implies{e}{d}
& = &
\sqrt{e}\cdot d \cdot \sqrt{e} + 1 - e.
\end{array}
\end{equation}

\noindent This formula for $\andthen{e}{d}$ is precisely the one that
is used in~\cite{GudderG02,GudderN02}, for effects on a Hilbert space,
as instance of the notion of `sequential product' on effect algebras.

By definition, in the current context, a projection is an
effect $e\in [0,1]_A$ satisfying $e = \andthen{e}{e}$. This means $e =
\sqrt{e}\cdot e \cdot \sqrt{e} =
\sqrt{e}\cdot\sqrt{e}\cdot\sqrt{e}\cdot\sqrt{e} = e \cdot e = e^{2}$.
Hence projections are idempotent effects. Equivalently, following the
more common description, an idempotent is an element $a\in A$
satisfying $a^{2} = a = a^{*}$. Notice that for such a projection $a$
one has $\sqrt{a} = a$ and thus $\andthen{a}{\andthen{a}{e}} =
\andthen{a}{e}$ for each effect $e$, since:
$$\begin{array}{rcccccccl}
\andthen{a}{\andthen{a}{e}}
& = &
\sqrt{a}\cdot (\sqrt{a}\cdot e \cdot \sqrt{a}) \cdot \sqrt{a}
& = &
a \cdot e \cdot a 
& = &
\sqrt{a}\cdot e \cdot \sqrt{a}
& = &
\andthen{a}{e}.
\end{array}$$

\noindent Hence a double andthen-test $\andthen{a}{-}$ with a
projection $a$ is the same as a single test. This holds not only for
andthen-tests, but also for then-tests $\implies{a}{-}$, by the De
Morgan property from
Lemma~\ref{TestOperatorLem}~\eqref{TestOperatorLemMorgan} below.

\auxproof{
If $a^{2} = a = a^{*}$, then $a \geq 0$ since $a = a\cdot a = a^{*}
\cdot a$. Also, $a\leq 1$ since $1-a \geq 0$ because $1-a$ is also a
projection: $(1-a)\cdot (1-a) = 1 -a - a + a^{2} = 1 - a - a + a = 1
-a$ and $(1-a)^{*} = 1^{*} - a^{*} = 1 - a$.

$$\begin{array}{rcl}
\implies{a}{e}
& = &
\big(\andthen{a}{e^\perp}\big)^{\perp} \\
& = &
\big(\andthen{a}{\andthen{a}{e^\perp}}\big)^{\perp} \\
& = &
\implies{a}{\big(\andthen{a}{e^\perp}\big)^{\perp}} \\
& = &
\implies{a}{\implies{a}{e}}.
\end{array}$$
}
\end{exas}

Next we collect some basic result about our test operators.

\begin{lem}
\label{TestOperatorLem}
The andthen-test and then-test operators $\andthen{-}{-}$ and
$\implies{-}{-}$ in a category satisfying
Assumption~\ref{MeasurementAss} satisfy:
\begin{enumerate}[label=\enspace(\arabic*)]
\item \label{TestOperatorLemTrue} $\andthen{1}{p} = p = \andthen{p}{1}$;

\item \label{TestOperatorLemFalse} $\andthen{0}{p} = 0 =
  \andthen{p}{0}$;

\item \label{TestOperatorLemOvee} $\andthen{p}{q_{1} \ovee q_{2}} =
  \andthen{p}{q_{1}} \ovee \andthen{p}{q_{2}}$;

\item \label{TestOperatorLemScalar} $\andthen{p}{s\scalar q} =
  s\scalar \andthen{p}{q}$, for a scalar $s$;

\item \label{TestOperatorLemOrder} $\andthen{p}{q} \leq p$;

\item \label{TestOperatorLemMonotone} $q_{1} \leq q_{2}$ implies
$\andthen{p}{q_1} \leq \andthen{p}{q_2}$;

\item \label{TestOperatorLemMorgan} $\implies{p}{q} =
  \andthen{p}{q^{\perp}}^{\perp}$;

\item $\andthen{p}{q} \ovee p^{\perp} = \implies{p}{q}$, and so
  $\andthen{p}{q} \leq \implies{p}{q}$;

\item $\andthen{[p_{1}, p_{2}]}{q} = [\andthen{p_{1}}{q\after \kappa_{1}},
   \andthen{p_{1}}{q \after \kappa_{2}}]$;

\item $f^{*}(\andthen{p}{q}) = \andthen{f^{*}(p)}{f^{*}(q)}$ if $f$ is
  a pure map.
\end{enumerate}
\end{lem}

Via point~\eqref{TestOperatorLemMorgan} (and the other points) one
obtains various properties of $\implies{p}{-}$. For instance,
$$\begin{array}{rccclcrcccccl}
\implies{p}{0}
& \smash{\stackrel{\eqref{TestOperatorLemMorgan}}{=}} &
\andthen{p}{1}^{\perp}
& \smash{\stackrel{\eqref{TestOperatorLemTrue}}{=}} &
p^{\perp}
& \qquad &
\implies{p}{1}
& \smash{\stackrel{\eqref{TestOperatorLemMorgan}}{=}} &
\andthen{p}{0}^{\perp}
& \smash{\stackrel{\eqref{TestOperatorLemFalse}}{=}} &
0^{\perp}
& = &
1.
\end{array}$$

\begin{myproof}
We describe some relevant steps, using the properties of the $\FstAnd$
operator from Lemma~\ref{CasesPropertiesLem}.
\begin{enumerate}[label=\enspace(\arabic*)]
\item Obviously, $\andthen{1}{p} = (\instr_{1})^{*}(\FstAnd(p)) =
  \FstAnd(p) \after \kappa_{1} = [p,0] \after \kappa_{1} = p$. And:
$$\begin{array}{rcccccccl}
\andthen{p}{1}
& = &
\FstAnd(1) \after \instr_{p}
& = &
[\kappa_{1} \after\; !, \kappa_{2} \after\; !] \after \instr_{p}
& = &
(!+!) \after \instr_{p}
& = &
p.
\end{array}$$

\item Like in the previous point, $\andthen{0}{p} =
  (\instr_{0})^{*}([p,0]) = [p, 0] \after \kappa_{2} = 0$. And for the
  equation $\andthen{p}{0} = 0$, we use
  Lemma~\ref{CasesPropertiesLem}~\eqref{CasesPropertiesLemZero} in:
$$\begin{array}{rcccccl}
\andthen{p}{0}
& = &
(\instr_{p})^{*}(\FstAnd(0))
& = &
(\instr_{p})^{*}(0)
& = &
0.
\end{array}$$

\item Lemma~\ref{CasesPropertiesLem}~\eqref{CasesPropertiesLemOvee}
  says that $\FstAnd$ preserves $\ovee$, and thus:
$$\begin{array}{rcl}
\andthen{p}{q_{1}\ovee q_{2}}
& = &
(\instr_{p})^{*}(\FstAnd(q_{1}\ovee q_{2})) \\
& = &
(\instr_{p})^{*}(\FstAnd(q_{1})\ovee \FstAnd(q_{2})) \\
& = &
(\instr_{p})^{*}(\FstAnd(q_{1})) \ovee (\instr_{p})^{*}(\FstAnd(q_{2})) \\
& = &
\andthen{p}{q_{1}} \ovee \andthen{p}{q_{2}}.
\end{array}$$

\item Similarly, using
  Lemma~\ref{CasesPropertiesLem}~\eqref{CasesPropertiesLemScalar}.

\item Using points~\eqref{TestOperatorLemOvee}
  and~\eqref{TestOperatorLemTrue} we have:
$$\begin{array}{rcccccl}
\andthen{p}{q} \ovee \andthen{p}{q^{\perp}}
& = &
\andthen{p}{q \ovee q^{\perp}}
& = &
\andthen{p}{1}
& = &
p.
\end{array}$$

\noindent Hence $\andthen{p}{q} \leq p$.

\item Assume $q_{1} \leq q_{2}$, say via $q_{1} \ovee r =
  q_{2}$. Then, again using point~\eqref{TestOperatorLemOvee} we get:
$$\begin{array}{rcccl}
\andthen{p}{q_1} \ovee \andthen{p}{r}
& = &
\andthen{p}{q_{1} \ovee r}
& = &
\andthen{p}{q_2}.
\end{array}$$

\noindent Hence $\andthen{p}{q_1} \leq \andthen{p}{q_2}$.

\item Using the De Morgan duality between $\FstAnd$ and $\FstThen$ we
  get:
$$\begin{array}{rcl}
\implies{p}{q}
\hspace*{\arraycolsep} = \hspace*{\arraycolsep}
(\instr_{p})^{*}(\FstThen(q))
& = &
(\instr_{p})^{*}(\FstAnd(q^{\perp})^{\perp}) \\
& = &
\big((\instr_{p})^{*}(\FstAnd(q^{\perp}))\big)^{\perp}
\hspace*{\arraycolsep} = \hspace*{\arraycolsep}
\andthen{p}{q^{\perp}}^{\perp}.
\end{array}$$

\item We define the bound $b\colon X \rightarrow (1+1)+1$ by $b =
  [\kappa_{1}+\idmap, \kappa_{1} \after \kappa_{2}] \after (q+!)
  \after \instr_{p}$. It proves $\andthen{p}{q} \orthogonal p^{\perp}$
  and $\andthen{p}{q} \ovee p^{\perp} = \implies{p}{q}$.

\auxproof{
$$\begin{array}{rcl}
[\idmap, \kappa_{2}] \after b
& = &
[\idmap, \kappa_{2}] \after 
   [\kappa_{1}+\idmap, \kappa_{1} \after \kappa_{2}] \after 
   (q+!) \after \instr_{p} \\
& = &
[[\kappa_{1}, \kappa_{2}], \kappa_{2}] \after
   (q+!) \after \instr_{p} \\
& = &
[q, \kappa_{2} \after\; !_{X}] \after \instr_{p} \\
& = &
\andthen{p}{q} 
\\
{[[\kappa_{2},\kappa_{1}], \kappa_{2}]} \after b
& = &
[[\kappa_{2},\kappa_{1}], \kappa_{2}] \after 
   [\kappa_{1}+\idmap, \kappa_{1} \after \kappa_{2}] \after 
   (q+!) \after \instr_{p} \\
& = &
[[\kappa_{2}, \kappa_{2}], \kappa_{1}] \after
   (q+!) \after \instr_{p} \\
& = &
[\kappa_{2} \after \nabla \after q, \kappa_{1} \after\; !] 
  \after \instr_{p} \\
& = &
[\kappa_{2} \after\; !, \kappa_{1} \after\; !] \after \instr_{p} \\
& = &
[\kappa_{2}, \kappa_{1}] \after 
  [\kappa_{1} \after\; !, \kappa_{2} \after\; !] \after \instr_{p} \\
& = &
[\kappa_{2}, \kappa_{1}] \after p \\
& = &
p^{\perp} 
\\
(\nabla+\idmap) \after b
& = &
(\nabla+\idmap) \after 
   [\kappa_{1}+\idmap, \kappa_{1} \after \kappa_{2}] \after 
   (q+!) \after \instr_{p} \\
& = &
[\idmap+\idmap, \kappa_{1} \after \nabla \after \kappa_{2}] \after 
   (q+!) \after \instr_{p} \\
& = &
[q, \kappa_{1} \after\; !] \after \instr_{p} \\
& = &
\implies{p}{q}
\end{array}$$
}

\item We use
  Assumption~\ref{MeasurementAss}~\eqref{MeasurementAssCotuple} in:
$$\begin{array}[b]{rcl}
\andthen{[p_{1}, p_{2}]}{q}
& = &
[q, 0] \after \instr_{[p_{1}, p_{2}]} \\
& = &
[q, 0] \after [\kappa_{1}+\kappa_{1}, \kappa_{2}+\kappa_{2}] \after
   (\instr_{p_{1}} + \instr_{p_{2}}) \\
& = &
[[q \after \kappa_{1}, 0 \after \kappa_{1}] \after \instr_{p_1}, 
   [q \after \kappa_{2}, 0 \after \kappa_{2}] \after \instr_{p_2}] \\
& = &
[[q \after \kappa_{1}, 0] \after \instr_{p_1}, 
   [q \after \kappa_{2}, 0] \after \instr_{p_2}] \\
& = &
[\andthen{p_{1}}{q\after \kappa_{1}},
   \andthen{p_{2}}{q \after \kappa_{2}}]
\end{array}$$

\item If $f\colon Y \rightarrow X$ is pure, then for $p,q\in\Pred(X)$,
$$\begin{array}[b]{rcl}
f^{*}(\andthen{p}{q})
\hspace*{\arraycolsep} = \hspace*{\arraycolsep}
[q,0] \after \instr_{p} \after f 
& = &
[q,0] \after (f+f) \after \instr_{p \after f} 
   \qquad \mbox{see Definition~\ref{PureDef}} \\
& = &
[q \after f, 0] \after \instr_{f^{*}(p)} \\
& = &
\andthen{f^{*}(p)}{f^{*}(q)}.
\end{array}\eqno{\qEd}$$
\end{enumerate}
\end{myproof}

\noindent There are some further properties of test operators that we do not
list in Lemma~\ref{TestOperatorLem} but describe separately because
they deserve more attention. We first consider the substitution
operation $f^{*} = \Pred(f) = \wp(f) = (-) \after f$ from
Proposition~\ref{CoprodEAProp}~\eqref{CoprodEAPropMor} that makes the
mapping $X \mapsto \Pred(X)$ functorial. We show that for the special
case where the map $f$ is an instrument $\instr_p$, substitution can
be defined in logical terms. This result will be useful immediately,
in a subsequent lemma, but may also be useful for Lindenbaum models
where the logic has test operators.

The result will be formulated first for predicates (2-tests), and then
more generally for $n$-tests. We separately mention version for $n=2$
because it has an easier proof and because it is useful on its own.
It describes a form of ``if-then-else''.

\begin{prop}
\label{MeasSubstProp}
Substitution for instruments, in a category satisfying
Assumption~\ref{MeasurementAss}, can be described as follows.
\begin{enumerate}
\item \label{MeasSubstPropTwo} For predicates
$p\colon X \rightarrow 1+1$ and $q\colon X+X \rightarrow 1+1$,
$$\begin{array}{rcl}
(\instr_{p})^{*}(q)
& = &
\andthen{p}{q\after\kappa_{1}} \ovee
   \andthen{p^{\perp}}{q \after \kappa_{2}}
\end{array}$$

\item \label{MeasSubstPropMany} For an $n$-test $p\colon X \rightarrow
  n\cdot 1$, where $n\geq 2$, and a predicate $q\colon n\cdot X
  \rightarrow 1+1$ one has:
$$\begin{array}{rcl}
(\instr_{p})^{*}(q)
& = &
\andthen{p_{1}}{q_{1}} \ovee \cdots \ovee \andthen{p_{n}}{q_{n}},
\end{array}$$

\noindent for predicates $p_{i}$ corresponding to the $n$-test $p$ as
in Lemma~\ref{TestLem} --- via $p_{i} = \rhd_{i} \after p \colon X
\rightarrow 1+1$, with $n$-ary partial projections $\rhd_i$ like
in~\eqref{PartProjDiag} --- and predicates $q_{i} = q \after
\kappa_{i} \colon X \rightarrow 1+1$ on $X$.
\end{enumerate}
\end{prop}

\begin{myproof}
For the first point we calculate:
$$\begin{array}{rcll}
(\instr_{p})^{*}(q)
& = &
(\instr_{p})^{*}([q \after \kappa_{1}, q \after \kappa_{2}]) \\
& = &
(\instr_{p})^{*}([q \after \kappa_{1}, 0] \ovee [0, q \after \kappa_{2}])
   & \mbox{by~\eqref{NaryCotupleOveeEqn}} \\
& = &
(\instr_{p})^{*}([q \after \kappa_{1}, 0]) \ovee 
   (\instr_{p})^{*}([0, q \after \kappa_{2}]) \\
& & \qquad \mbox{since substitution is a map of effect algebras} \\
& = &
\andthen{p}{q\after\kappa_{1}} \ovee 
   (\instr_{p})^{*}([q \after \kappa_{2}, 0] \after [\kappa_{2}, \kappa_{1}]) \\
& = &
\andthen{p}{q\after\kappa_{1}} \ovee 
   ([\kappa_{2}, \kappa_{1}] \after \instr_{p})^{*}([q \after \kappa_{2}, 0]) \\
& = &
\andthen{p}{q\after\kappa_{1}} \ovee 
   (\instr_{p^{\perp}})^{*}([q \after \kappa_{2}, 0]) &
  \mbox{by Assumption~\ref{MeasurementAss}~\eqref{MeasurementAssInjection}} \\
& = &
\andthen{p}{q\after\kappa_{1}} \ovee
   \andthen{p^{\perp}}{q \after \kappa_{2}}.
\end{array}$$
For the second point we need to use
Assumption~\ref{MeasurementAss}~\eqref{MeasurementAssCollapse} and
also~\eqref{MeasurementAssInjection}. The proof involves some
bookkeeping with coprojections. We use the maps:
$$\begin{array}{rcl}
\rhd_{i}
& = &
\underbrace{[\kappa_{2}, \ldots, \kappa_{2}, \kappa_{1}, \kappa_{2}, \ldots,
   \kappa_{2}]}_{\text{$\kappa_1$ at the $i$-the position}} \;\colon\;
   n\cdot 1 \longrightarrow 1+1 \\
\swap_{i}
& = &
\underbrace{[\kappa_{i}, \kappa_{2}, \kappa_{3}, \ldots, \kappa_{i-1}, 
   \kappa_{1}, \kappa_{i+1}, \ldots, 
   \kappa_{n}]}_{\text{$\kappa_{i}$ at the first, $\kappa_1$ at the $i$-the position}}
  \;\colon\; n\cdot Y \conglongrightarrow n\cdot Y.
\end{array}$$

\noindent Then we can write $p_{i} = \rhd_{i} \after p = (\idmap+\,!)
\after \swap_{i} \after p$, where $\idmap+\,!  \colon n\cdot 1 =
1+(n-1)\cdot 1 \rightarrow 1+1$. Hence we get:
$$\begin{array}[b]{rcll}
\bigovee_{i}\, \andthen{p_i}{q_i}
& = &
\bigovee_{i}\, [q \after \kappa_{i}, \kappa_{2} \after\, !] \after \instr_{p_i} \\
& = &
\bigovee_{i}\, [q \after \kappa_{i}, \kappa_{2}] \after (\idmap+\,!) \after
   \instr_{(\idmap+\,!) \after \swap_{i} \after p} \quad \\
& = &
\bigovee_{i}\, [q \after \kappa_{i}, \kappa_{2}] \after (\idmap+\,!) \after
   \instr_{\swap_{i} \after p} 
   & \mbox{by Assumption~\ref{MeasurementAss}~\eqref{MeasurementAssCollapse}} \\
& = &
\bigovee_{i}\, [q \after \kappa_{i}, \kappa_{2} \after\, !] \after
   \swap_{i} \after \instr_{p} 
   & \mbox{by Assumption~\ref{MeasurementAss}~\eqref{MeasurementAssInjection}} \\
& = &
\bigovee_{i}\, [0, \ldots, 0, q \after \kappa_{i}, 0, \ldots, 0] \after
   \instr_{p} \\
& = &
\bigovee_{i}\, (\instr_{p})^{*}\big([0, \ldots, 0, q \after \kappa_{i}, 
   0, \ldots, 0]\big) \\
& = &
(\instr_{p})^{*}\big(\bigovee_{i}\, [0, \ldots, 0, q \after \kappa_{i}, 
   0, \ldots, 0]\big) \\
& = &
(\instr_{p})^{*}\big([q \after \kappa_{1}, \ldots, q \after \kappa_{n}]\big)
   & \mbox{by~\eqref{NaryCotupleOveeEqn}} \\
& = &
(\instr_{p})^{*}(q).
\end{array}\eqno{\qEd}$$
\end{myproof}\medskip

\noindent We obtain a consequence that highlights the role of side-effects in
formulas. In an informal reading the result below says: the sum of
`$p$ andthen $q$', and of `$p^{\perp}$ andthen $q$' is the same as
just having $q$. But in our present setting we have to take
side-effect $\nabla \after \instr_p$ of $p$ into account.

\begin{cor}
\label{TwoTestSefCor}
For predicates $p,q\colon X \rightarrow 1+1$ in a category
satisfying Assumption~\ref{MeasurementAss} one has:
\begin{equation}
\label{TwoTestSefAndThen}
\begin{array}{rcl}
\andthen{p}{q} \ovee \andthen{p^{\perp}}{q}
& = &
q \after \nabla \after \instr_{p}.
\end{array}
\end{equation}


\end{cor}\bigskip

\noindent We recall from Definition~\ref{SideEffectFreeDef} that the map $\nabla
\after \instr_{p}$ is the side-effect of the predicate $p$.  In case
$p$ is side-effect-free, that is, in case $\nabla \after \instr_{p} =
\idmap$, the above equation~\eqref{TwoTestSefAndThen} reduces to the
standard result $\andthen{p}{q} \ovee \andthen{p^{\perp}}{q} = q$.

\begin{myproof}
Directly from
Proposition~\ref{MeasSubstProp}~\eqref{MeasSubstPropTwo} we obtain:
$$\begin{array}[b]{rcl}
\andthen{p}{q} \ovee \andthen{p^{\perp}}{q}
& = &
\andthen{p}{[q,q] \after \kappa_{1}} \ovee 
   \andthen{p^{\perp}}{[q,q]\after\kappa_{2}} \\
& = &
(\instr_{p})^{*}([q,q]) \\
& = &
[q,q] \after \instr_{p} \\
& = &
q \after \nabla \after \instr_{p}.
\end{array}\eqno{\qEd}$$

\end{myproof}\medskip

\noindent We add another result illustrating that substitution $f^*(q) = q
\after f$ can be understood as weakest precondition $\wp(f)(q)$. We
apply this to the (binary case of the) test map
from~\eqref{TestMapEqn}. It yields a formula that looks very much like
the traditional weakest precondition formula for if-then-else:
$$\begin{array}{rcl}
\wp\big(\textsf{if }p \textsf{ then } f_{1} \textsf{ else } f_{2}\big)(q)
& = &
\big(p \wedge \wp(f_{1})(q)\big) \vee \big(\neg p \wedge \wp(f_{2})(q)\big).
\end{array}$$

\noindent In the present context the conjunctions are replaced by test
operators that deal with the side-effects involved.

\begin{lem}
\label{TestMapWPLem}
For parallel maps $f_{1},f_{2}\colon X \rightarrow Y$ in a category satisfying
Assumption~\ref{MeasurementAss} and a predicate $p\colon X \rightarrow
1+1$ on $X$ we have:
$$\begin{array}{rcl}
\wp\big(p?[f_{1},f_{2}]\big)(q)
& = &
\andthen{p}{\wp(f_{1})(q)} \ovee \andthen{p^\perp}{\wp(f_{2})(q)}.
\end{array}$$
\end{lem}

\begin{myproof}
We use Proposition~\ref{MeasSubstProp}~\eqref{MeasSubstPropTwo} to
obtain:
$$\begin{array}[b]{rcl}
\andthen{p}{\wp(f_{1})(q)} \ovee \andthen{p^\perp}{\wp(f_{2})(q)}
& = &
(\instr_{p})^{*}\big([\wp(f_{1})(q), \wp(f_{2})(q)]\big) \\
& = &
[q \after f_{1}, q \after f_{2}] \after \instr_{p} \\
& = &
q \after [f_{1}, f_{2}] \after \instr_{p} \\
& \smash{\stackrel{\eqref{TestMapEqn}}{=}} &
q \after p?[f_{1}, f_{2}] \\
& = &
\wp\big(p?[f_{1},f_{2}]\big)(q).
\end{array}\eqno{\qEd}$$
\end{myproof}

\auxproof{
\begin{lem}
Projections form a sub-effect algebra.\marginpar{Why is $p^\perp$ a
  projection if $p$ is?} In the commutative case $\andthen{p}{q}$ is
the meet of projections $p,q$.
\end{lem}
}\medskip

\noindent We continue with $C^*$-algebras as instance.  The definition of
commutativity $\andthen{p}{q} = \andthen{q}{p}$ used in
Definition~\ref{TestOperatorDef} turns out to coincide with
commutativity in $C^*$-algebras. This is a non-trivial result that
goes back to~\cite{GudderN02}. It uses the Fuglede-Putnam-Rosenblum
Theorem~\cite{Rudin87} for $C^*$-algebras. This theorem says: for
normal elements $a,b$ one has for each $x$,
\begin{equation}
\label{FPREqn}
\begin{array}{rcl}
a\cdot x = x \cdot b
& \Longrightarrow &
a^{*}\cdot x = x \cdot b^{*}.
\end{array}
\end{equation}

\noindent We recall that $a$ is normal if $a\cdot a^{*} = a^{*} \cdot
a$.

\begin{prop}
\label{CstarCommProp}
A $C^*$-algebra $A$ is commutative, in the sense that its
multiplication $\cdot$ is commutative, if and only if it is
commutative, in the sense of Definition~\ref{TestOperatorDef}:
$\andthen{e}{d} = \andthen{d}{e}$ holds for all effects $e,d\in
[0,1]_{A}$.
\end{prop}

%












\begin{myproof}
First assume that multiplication $\cdot$ of the $C^*$-algebra $A$ is
commutative. Then:
$$\begin{array}{rcccccccccccl}
\andthen{e}{d}
& = &
\sqrt{e} \cdot d \cdot \sqrt{e}
& = &
\sqrt{e} \cdot \sqrt{e} \cdot d 
& = &
e \cdot d
& = &
e\cdot \sqrt{d}\cdot \sqrt{d}
& = &
\sqrt{d} \cdot e \cdot \sqrt{d}
& = &
\andthen{d}{e}.
\end{array}$$

\noindent In the other direction, it suffices by
Lemma~\ref{EffectCommutationLem} to prove $e\cdot d = d\cdot e$ for
effects $e,d$. By assumption we have $\sqrt{e}\cdot d\cdot \sqrt{e} =
\andthen{e}{d} = \andthen{d}{e} = \sqrt{d} \cdot e \cdot \sqrt{d}$.
The product $\sqrt{e}\cdot \sqrt{d}$ is normal since:
$$\begin{array}{rcll}
\big(\sqrt{e}\cdot \sqrt{d}\big)^{*} \cdot \big(\sqrt{e}\cdot \sqrt{d}\big)
& = &
\sqrt{d}\cdot \sqrt{e} \cdot \sqrt{e}\cdot \sqrt{d} \\
& = &
\sqrt{d} \cdot e \cdot \sqrt{d} \\
& = &
\sqrt{e} \cdot d \cdot \sqrt{e} & \mbox{by assumption} \\
& = &
\sqrt{e}\cdot \sqrt{d} \cdot \sqrt{d}\cdot \sqrt{e} \\
& = &
\big(\sqrt{e}\cdot \sqrt{d}\big) \cdot \big(\sqrt{e}\cdot \sqrt{d}\big)^{*}.
\end{array}$$

\noindent We use Fuglede-Putnam-Rosenblum~\eqref{FPREqn} with $a =
\sqrt{e}\cdot \sqrt{d}$, $b=\sqrt{d}\cdot \sqrt{e}$ and
$x=\sqrt{e}$. The antecedent holds, since $a\cdot x = \sqrt{e} \cdot
\sqrt{d} \cdot \sqrt{e} = x \cdot b$. It yields as conclusion:
$$\begin{array}{rcccccl}
\sqrt{d} \cdot e
& = &
\big(\sqrt{e}\cdot \sqrt{d}\big)^{*}\cdot \sqrt{e}
& \smash{\stackrel{\eqref{FPREqn}}{=}} &
\sqrt{e}\cdot \big(\sqrt{d}\cdot \sqrt{e}\big)^{*}
& = &
e \cdot \sqrt{d}.
\end{array}$$

\noindent But then we are done since:
$$\begin{array}{rcccccccl}
e\cdot d
& = &
e\cdot \sqrt{d} \cdot \sqrt{d}
& = &
\sqrt{d} \cdot e \cdot \sqrt{d}
& = &
\sqrt{d} \cdot \sqrt{d} \cdot e 
& = &
d \cdot e.
\end{array}\eqno{\qEd}$$

\auxproof{
By exchanging $e,d$ we also get $\sqrt{e}\cdot d = d \cdot
\sqrt{e}$. In combination with the assumption $\sqrt{e}\cdot d\cdot
\sqrt{e} = \sqrt{d}\cdot e\cdot \sqrt{d}$ we are done:
$$\begin{array}{rcccccccccl}
e\cdot d
& = &
e\cdot \sqrt{d} \cdot \sqrt{d}
& = &
\sqrt{d} \cdot e \cdot \sqrt{d}
& = &
\sqrt{e} \cdot d \cdot \sqrt{e}
& = &
d \cdot \sqrt{e} \cdot \sqrt{e}
& = &
d \cdot e.
\end{array}\eqno{\qEd}$$
}
\end{myproof}

We conclude this section with an example that illustrates the use of
the logical test operators in an elementary probability calculation.

\begin{exa}
\label{PolarisationEx}
The idea in the famous polarisation experiment is to send photons
polarised as $\ket{\!\uparrow\!} = \ket{0} = \left(\begin{smallmatrix}
  1 \\ 0
\end{smallmatrix}\right)$ through filters.
\begin{enumerate}
\item[(1)] If these photons hit one filter with polarisation
  $\ket{\!\rightarrow\!} = \ket{1} = \left(\begin{smallmatrix} 0 \\ 1
\end{smallmatrix}\right)$, then nothing passes.

\item[(2)] Next assume there are two consecutive filters: first the
  $\ket{\!\uparrow\!}$-photons have to go through a filter with
  polarisation $\ket{\!\nearrow\!} =
  \frac{1}{\sqrt{2}}\left(\begin{smallmatrix} 1 \\ 1
\end{smallmatrix}\right)$, and then through a filter with polarisation
$\ket{\!\rightarrow\!} = \ket{1} = \left(\begin{smallmatrix} 0 \\ 1
\end{smallmatrix}\right)$. Surprisingly, now one quarter of the original
photons get through. We refer to \textit{e.g.}~\cite{RieffelP11} for
an explanation, and focus here on a calculation of this probability
$\frac{1}{4}$ in the current quantitative logic.
\end{enumerate}\medskip

\noindent For these experiments we use the $C^*$-algebra $A =
\B(\C^{2})$ with initial state $\omega_{\uparrow} \colon A \rightarrow
\C$ given by $\omega_{\uparrow}(M) = \bra{0}M\ket{0}$. The effect
$p_{\rightarrow} \in [0,1]_{A} = \Ef(\C^{2})$ corresponding to the
filter with polarisation $\ket{\!\rightarrow\!}$ is given by:
$$\begin{array}{rcccccl}
p_{\rightarrow}
& = &
\ket{\!\rightarrow\!}\bra{\!\rightarrow\!}
& = &
\left(\begin{smallmatrix} 0 & 1 \end{smallmatrix}\right)
   \left(\begin{smallmatrix} 0 \\ 1 \end{smallmatrix}\right) 
& = &
\left(\begin{smallmatrix} 0 & 0 \\ 0 & 1 \end{smallmatrix}\right).
\end{array}$$

\noindent It is clearly a projection. Evaluating this predicate in
the initial state $\omega_{\uparrow}$ yields the probability:
$$\begin{array}{rcccccccl}
\omega_{\uparrow} \models p_{\rightarrow}
& \smash{\stackrel{\eqref{ModelsEqn}}{=}} &
\omega_{\uparrow}(p_{\rightarrow})
& = &
\bra{\!\uparrow\!}p_{\rightarrow}\ket{\!\uparrow\!}
& = &
\left(\begin{smallmatrix} 1 & 0 \end{smallmatrix}\right)
\left(\begin{smallmatrix} 0 & 0 \\ 0 & 1 \end{smallmatrix}\right)
\left(\begin{smallmatrix} 1 \\ 0 \end{smallmatrix}\right) 
& = &
0.
\end{array}$$

\noindent This is the probability of the above first experiment:
nothing goes through.

Things become more interesting in the second experiment, where we
first send the $\ket{0}$-photon through the filter with polarisation
$\ket{\!\nearrow\!}$. This involves the effect:
$$\begin{array}{rcccccl}
p_{\nearrow}
& = &
\ket{\!\nearrow\!}\bra{\!\nearrow\!}
& = &
\frac{1}{\sqrt{2}}\left(\begin{smallmatrix} 1 \\ 1 \end{smallmatrix}\right)
\frac{1}{\sqrt{2}}\left(\begin{smallmatrix} 1 & 1 \end{smallmatrix}\right)
& = &
\frac{1}{2}\left(\begin{smallmatrix} 1 & 1 \\ 1 & 1 \end{smallmatrix}\right).
\end{array}$$

\noindent This $p_{\nearrow}$ is also a projection, so that
$\sqrt{p_{\nearrow}} = p_{\nearrow}$.

We now compute the probability via the test operators from
Section~\ref{TestSec}. First going through the filter $\nearrow$ and
then through $\rightarrow$ is expressed via the test-andthen operator.
The resulting probability is:
$$\begin{array}{rcl}
\omega_{\uparrow} \models \andthen{p_{\nearrow}}{p_{\rightarrow}}
& \smash{\stackrel{\eqref{CstarTestOperatorEqn}}{=}} &
\bra{\!\uparrow\!}\sqrt{p_{\nearrow}}\, p_{\rightarrow}\sqrt{p_{\nearrow}}
   \ket{\!\uparrow\!} \\
& = &
\left(\begin{smallmatrix} 1 & 0 \end{smallmatrix}\right)
\frac{1}{2}\left(\begin{smallmatrix} 1 & 1 \\ 1 & 1 \end{smallmatrix}\right)
\left(\begin{smallmatrix} 0 & 0 \\ 0 & 1 \end{smallmatrix}\right)
\frac{1}{2}\left(\begin{smallmatrix} 1 & 1 \\ 1 & 1 \end{smallmatrix}\right)
\left(\begin{smallmatrix} 1 \\ 0 \end{smallmatrix}\right) 
\hspace*{\arraycolsep} = \hspace*{\arraycolsep}
\frac{1}{4}\left(\begin{smallmatrix} 1 & 0 \end{smallmatrix}\right)
\left(\begin{smallmatrix} 1 & 1 \\ 1 & 1 \end{smallmatrix}\right)
\left(\begin{smallmatrix} 1 \\ 0 \end{smallmatrix}\right) 
\hspace*{\arraycolsep} = \hspace*{\arraycolsep}
\frac{1}{4}.
\end{array}$$

\noindent Instead of the andthen-test we can also use the weaker
then-test.  It adds the option that the photon does not get through
the first filter --- which happens with probability
$\frac{1}{2}$. Then:
$$\begin{array}{rcl}
\omega_{\uparrow} \models \implies{p_{\nearrow}}{p_{\rightarrow}}
& \smash{\stackrel{\eqref{CstarTestOperatorEqn}}{=}} &
\bra{\!\uparrow\!}\Big(\sqrt{p_{\nearrow}}\, p_{\rightarrow}\sqrt{p_{\nearrow}}\,
   + 1 - p_{\nearrow}\Big)\ket{\!\uparrow\!} \\
& = &
\left(\begin{smallmatrix} 1 & 0 \end{smallmatrix}\right)\Big(
\frac{1}{2}\left(\begin{smallmatrix} 1 & 1 \\ 1 & 1 \end{smallmatrix}\right)
\left(\begin{smallmatrix} 0 & 0 \\ 0 & 1 \end{smallmatrix}\right)
\frac{1}{2}\left(\begin{smallmatrix} 1 & 1 \\ 1 & 1 \end{smallmatrix}\right)
+ 
\left(\begin{smallmatrix} 1 & 0 \\ 0 & 1 \end{smallmatrix}\right)
-
\frac{1}{2}\left(\begin{smallmatrix} 1 & 1 \\ 1 & 1 \end{smallmatrix}\right)
\Big)\left(\begin{smallmatrix} 1 \\ 0 \end{smallmatrix}\right) \\
& = &
\left(\begin{smallmatrix} 1 & 0 \end{smallmatrix}\right)\Big(
\frac{1}{4}\left(\begin{smallmatrix} 1 & 1 \\ 1 & 1 \end{smallmatrix}\right)
+
\frac{1}{2}\left(\begin{smallmatrix} 1 & -1 \\ -1 & 1 \end{smallmatrix}\right)
\Big)\left(\begin{smallmatrix} 1 \\ 0 \end{smallmatrix}\right) \\
& = &
\left(\begin{smallmatrix} 1 & 0 \end{smallmatrix}\right)
\frac{1}{4}\left(\begin{smallmatrix} 3 & -1 \\ -1 & 3 \end{smallmatrix}\right)
\left(\begin{smallmatrix} 1 \\ 0 \end{smallmatrix}\right) \\
& = &
\frac{3}{4}.
\end{array}$$
\end{exa}


\section{Adding tensor products}\label{TensorSec}

Our third assumption involves tensor products, with some special
properties. We recall that a symmetric monoidal structure on a
category $\cat{B}$ involves a functor $\otimes\colon
\cat{B}\times\cat{B}\rightarrow\cat{B}$ and tensor unit $I\in\cat{B}$,
together with canonical isomorphisms:
\begin{equation}
\label{MacLaneEqn}
\xymatrix{
X\otimes (Y\otimes Z) \rto^-{\alpha}_-{\cong} & (X\otimes Y)\otimes Z
\qquad
X\otimes I\rto^-{\rho}_-{\cong} & X
\qquad
X\otimes Y \rto^-{\gamma}_-{\cong} & Y\otimes X
}
\end{equation}

\noindent One often writes $\lambda = \rho \after \gamma \colon
I\otimes X\conglongrightarrow X$. These isomorphisms make some obvious
diagrams commute, see~\cite{MacLane71}.

\begin{assumption}
\label{TensorAss} 
Our category satisfies not only Assumptions~\ref{CoprodAss}
and~\ref{MeasurementAss} but is also symmetric monoidal with the
additional requirements described in the next five points.
\begin{enumerate}
\item \label{TensorAssUnit} The final object $1$ is the tensor unit.
  As a result we have a `tensor with projections', since projections
  $\pi_{1}, \pi_{2}$ can be defined in:
\begin{equation}
\label{ProjTensorEqn}
\vcenter{\xymatrix{
X & X\otimes Y\ar[dl]^-{\idmap\otimes\,!}\ar[dr]_-{!\,\otimes\idmap}
   \ar[l]_-{\pi_1}\ar[r]^-{\pi_2} &
   Y \\
X\otimes 1\ar[u]_{\cong}^{\rho} & & 1\otimes Y\ar[u]^{\cong}_{\lambda}
}}
\end{equation}

\noindent These projections $\pi_{1},\pi_{2}$ are natural in $X$ and
$Y$.

\auxproof{
For $f\colon X \rightarrow X'$ and $g\colon Y \rightarrow Y'$ we have:
$$\begin{array}{rcl}
\pi_{2} \after (f\otimes g)
& = &
\lambda \after (!\otimes\idmap) \after (f\otimes g) \\
& = &
\lambda \after (!\otimes g) \\
& = &
\lambda \after (\idmap\otimes g) \after (!\otimes\idmap) \\
& = &
g \after \lambda \after (!\otimes\idmap) \\
& = &
g \after \pi_{2}.
\end{array}$$
}

\item \label{TensorAssDistr} The tensor product $\otimes$ distributes
  over finite coproducts. Explicitly, the following canonical maps are
  isomorphisms:
\begin{equation}
\label{DistrTensorEqn}
\xymatrix@C+.3pc{
(X\otimes A) + (Y\otimes A)
   \ar[rr]^-{[\kappa_{1}\otimes\idmap, \kappa_{2}\otimes\idmap]}_-{\cong} & &
   (X+Y) \otimes A
& 
0\ar[r]^-{!}_-{\cong} & A\otimes 0
}
\end{equation}

\noindent We shall write $\theta_{1} = [\kappa_{1}\otimes\idmap,
  \kappa_{2}\otimes\idmap]^{-1}$ for the inverse of this map on the
left. There is a similar distributivity map $\theta_{2} \colon
A\otimes (X+Y) \conglongrightarrow (A\otimes X) + (A\otimes Y)$ in the
other coordinate. These distributivity maps extend to $n$-ary
coproducts.

\item \label{TensorAssMonIsoPure} The monoidal
  isomorphisms~\eqref{MacLaneEqn} are pure --- see
  Definition~\ref{PureDef}.

\item \label{TensorAssTensorPure} If both $f,g$ are pure maps, then
  $f\otimes g$ is pure too.

\item \label{TensorAssInstr} For each $n$-test $p\colon X \rightarrow
  n\cdot 1$ and for each object $A$ the following diagrams commute.
\begin{equation}
\label{DistrMeasEqn}
\vcenter{\xymatrix@R-.5pc@C-.5pc{
X\otimes A\ar[rr]^-{\instr_{p}\otimes\idmap}\ar@/_2ex/[drr]_{\instr_{\pi_{1}^{*}(p)}} 
   & & (n\cdot X)\otimes A\ar[d]_{\cong}^{\theta_1}
& &
A\otimes X\ar[rr]^-{\idmap\otimes\instr_{p}}\ar@/_2ex/[drr]_{\instr_{\pi_{2}^{*}(p)}} 
   & & A\otimes (n\cdot X)\ar[d]_{\cong}^{\theta_2} \\
& & n\cdot (X\otimes A)
& &
& & n\cdot (A\otimes X)
}}
\end{equation}
\end{enumerate}
\noindent This concludes our assumptions about tensors.
\end{assumption}

The projections $\pi_{1}, \pi_{2}$ in~\eqref{ProjTensorEqn} for the
tensor $\otimes$, like in~\cite{Jacobs94a}, enable discarding of
resources. There are in general no associated diagonals $X \rightarrow
X\otimes X$ because of the no-cloning
Theorem~\cite{WootersZ82,Dieks82,Abramsky10b}.

The two distributivity isomorphisms $\theta_{1}, \theta_{2}$ interact
in the following way.
\begin{equation}
\label{DistrDistrDiag}
\hspace*{-1.1em}\vcenter{\xymatrix@R-.5pc@C-5pc{
& (A\!+\!B)\!\otimes\! (X\!+\!Y)\ar[dl]_{\theta_1}\ar[dr]^{\theta_2} & 
\\
(A\!\otimes\! (X\!+\!Y))\!+\!(B\!\otimes\!(X\!+\!Y))
      \ar[d]_{\theta_{2}+\theta_{2}} & &
   ((A\!+\!B)\!\otimes\! X)\!+\!((A\!+\!B)\!\otimes\! Y)
       \ar[d]^{\theta_{1}+\theta_{1}} 
\\
\big((A\!\otimes\! X)\!+\!(A\!\otimes\! Y)\big) \!+\!
    \big((B\!\otimes\! X)\!+\!(B\!\otimes\! Y)\big) 
    \ar[rr]_-{\widehat{\gamma}}^-{\cong} & &
    \big((A\!\otimes\! X)\!+\!(B\!\otimes\! X)\big) \!+\!
    \big((A\!\otimes\! Y)\!+\!(B\!\otimes\! Y\rlap{$)\big)$}
}}
\end{equation}

\noindent where $\widehat{\gamma} = [\kappa_{1}+\kappa_{1},
  \kappa_{2}+\kappa_{2}]$ is the map that swaps the inner two
occurrences; it is its own inverse.

\auxproof{
In the proof of commutation we use the inverses of the distribution
maps $\theta_{1}, \theta_{2}$, as in:
$$\begin{array}{rcl}
\lefteqn{\theta_{1}^{-1} \after (\theta_{2}^{-1}+\theta_{2}^{-1}) \after
   \widehat{\gamma}} \\
& = &
[\kappa_{1}\otimes\idmap, \kappa_{2}\otimes\idmap] \after
   ([\idmap\otimes\kappa_{1}, \idmap\otimes\kappa_{2}] +
    [\idmap\otimes\kappa_{1}, \idmap\otimes\kappa_{2}]) \after
   \widehat{\gamma} \\
& = &
\big[(\kappa_{1}\otimes\idmap) \after 
   [\idmap\otimes\kappa_{1}, \idmap\otimes\kappa_{2}], \;
   \kappa_{2}\otimes\idmap], \; (\kappa_{2}\otimes\idmap) \after
   [\idmap\otimes\kappa_{1}, \idmap\otimes\kappa_{2}]\big] \after
   \widehat{\gamma} \\
& = &
\big[[\kappa_{1}\otimes\kappa_{1}, \kappa_{1}\otimes\kappa_{2}], \;
   [\kappa_{2}\otimes\kappa_{1}, \kappa_{2}\otimes\kappa_{2}]\big] 
   \after [\kappa_{1}+\kappa_{1}, \kappa_{2}+\kappa_{2}] \\
& = &
\big[[\kappa_{1}\otimes\kappa_{1}, \kappa_{2}\otimes\kappa_{1}], \;
    [\kappa_{1}\otimes\kappa_{2}, \kappa_{2}\otimes\kappa_{2}]\big] \\
& = &
\big[(\idmap\otimes\kappa_{1}) \after
   [\kappa_{1}\otimes\idmap, \kappa_{2}\otimes\idmap], \;
   (\idmap\otimes\kappa_{2}) \after
   [\kappa_{1}\otimes\idmap, \kappa_{2}\otimes\idmap]\big] \\
& = &
[\idmap\otimes\kappa_{1}, \idmap\otimes\kappa_{2}] \after
   ([\kappa_{1}\otimes\idmap, \kappa_{2}\otimes\idmap] +
    [\kappa_{1}\otimes\idmap, \kappa_{2}\otimes\idmap]) \\
& = &
\theta_{2}^{-1} \after (\theta_{1}^{-1}+\theta_{1}^{-1}) 
\end{array}$$
}

Sometimes one writes $X^{\otimes n} = X\otimes \cdots \otimes X$ for
the $n$-fold tensor.  Because $\otimes$ distributes over $+$ and $1$
is the tensor unit, we have for the object $2 = 1+1$ the following
isomorphism:
$$\begin{array}{rcccccl}
2^{\otimes n}
& = &
\underbrace{2\otimes \cdots \otimes 2}_{\text{$n$ times}}
& \cong &
\underbrace{1 + \cdots + 1}_{\text{$2^{n}$ times}}
& = &
2^{n}.
\end{array}$$

\auxproof{
For $n=0$ we have, using that $1$ is the tensor unit:
$$\begin{array}{rcccl}
2^{\otimes 0}
& = &
1
& = &
2^{0}.
\end{array}$$

\noindent And for $n>0$ we have:
$$\begin{array}{rcl}
2^{\otimes n}
& = &
2^{\otimes (n-1)} \otimes (1+1) \\
& \cong &
2^{\otimes (n-1)} + 2^{\otimes (n-1)}  \\
& \cong &
2^{n-1} + 2^{n-1} \\
& \cong &
2^{n}.
\end{array}$$
}

\noindent Further we have $n\cdot X \cong X\otimes (n\cdot 1)$, since:
$$\begin{array}{rcccl}
n \cdot X
& \cong &
n \cdot (X\otimes 1)
& \cong &
X \otimes (n\cdot 1)
\end{array}$$

\noindent Hence, in presence of tensors we may describe instruments
also as maps $X \rightarrow X\otimes (n\cdot 1)$. The side-effect is
then obtained via the composite $X \rightarrow X\otimes (n\cdot 1)
\rightarrow X$, using the the first projection, as in
point~\eqref{ProjTensorEqn}.

\auxproof{
We use that the following diagram commutes:
$$\xymatrix@C+1pc{
n\cdot X\ar[r]^-{n\cdot \rho^{-1}}_-{\cong}\ar[ddrr]_{\nabla} &
   n\cdot (X\otimes 1)\ar[r]^-{[\idmap\otimes\kappa_{i}]_{i\leq n}}_-{\cong}
   \ar[rd]_{\nabla} &
   X \otimes (n\cdot 1)\ar[d]^-{\idmap\otimes\,! = \idmap\otimes\nabla}
\\
& & X\otimes 1\ar[d]^{\rho}_-{\cong}
\\
& & X
}$$
}

\begin{exas}
\label{TensorEx}
In the category $\Sets$, or in any extensive category, we can use
cartesian products $\times$ as tensors, see
Appendix~\ref{DiscProbSubsec}. They obviously distribute over
coproducts. In the Kleisli category $\Kl(\Dst)$ the projections
$\pi_{i} \colon X_{1}\times X_{2} \rightarrow \Dst(X_{i})$ are given
by $\pi_{i}(x_{1},x_{2}) = 1\ket{x_i}$. It is not hard to check that
instruments in $\Kl(\Dst)$ commute with distributivity maps, as in
diagrams~\eqref{DistrMeasEqn}.

\auxproof{
For an $n$-test $p\colon X \rightarrow \Dst(n)$ we have $\instr_{p}(x)
= \sum_{i}p(x)(i)\ket{\kappa_{i}x}$. There is the $n$-test
$\pi_{1}^{*}(p) = p \after \pi_{1} \colon X\times A \rightarrow
\Dst(n)$ with $\instr_{\pi_{1}^{*}(p)}(x,a) = \sum_{i}
p(x)(i)\ket{\kappa_{i}(x,a)}$. The map $\theta_{1}^{-1} \colon 
(n\cdot X)\otimes A \conglongrightarrow \Dst(n\cdot X\otimes A)$ is
given by:
$$\begin{array}{rcl}
\theta_{1}^{-1}
& = &
[\kappa_{1}\otimes\idmap, \ldots, \kappa_{n}\otimes\idmap] \\
& = &
[\dst \after ((\eta \after \kappa_{1})\times\eta), \ldots, 
   \dst \after ((\eta \after \kappa_{n})\times\eta)] \\
& = &
[\eta \after (\kappa_{1}\times\idmap), \ldots, 
   \eta \after (\kappa_{n}\times\idmap)] \\
& = &
\eta \after [\kappa_{1}\times\idmap, \ldots, \kappa_{n}\times\idmap].
\end{array}$$

\noindent Hence:
$$\begin{array}{rcl}
\big(\theta_{1}^{-1} \klafter \instr_{\pi_{1}^{*}(p)}\big)(x,a)
& = &
\Dst([\kappa_{1}\times\idmap, \ldots, \kappa_{n}\times\idmap])
   (\instr_{\pi_{1}^{*}(p)}(x,a)) \\
& = &
\sum_{i}p(x)(i)\ket{\kappa_{i}x, a} \\
& = &
\st(\instr_{p}(x), a) \\
& = &
\dst(\instr_{p}(x), \eta(a)) \\
& = &
\big(\instr_{p}\otimes\idmap\big)(x,a)
\end{array}$$
}

Tensors of $C^{*}$-algebras are a delicate matter. They exist in the
category of $C^*$-algebras with completely positive maps. Hence they
also exist on $\op{(\CstarCPU)}$, and the ``minimal'' tensor even
distributes over finite coproducts, as shown in~\cite{Cho14a}, see
Appendix~\ref{CstarSubsec}.

We need to check that our previous Assumptions~\ref{CoprodAss}
and~\ref{MeasurementAss} also hold in $\CstarCPU$. This is not a
problem. A useful fact is: when either $A$ or $B$ is a commutative
$C^*$ algebra, then a positive unital map $f\colon A \rightarrow B$
is automatically completely positive. Hence our predicates and tests
$\C^{n} \rightarrow A$ are completely positive. The instrument
$\instr_{\vec{e}}$ from~\eqref{MeasurementCstarEqn} is also completely
positive, due to its special form.

The projection $A\otimes B \rightarrow A$ in $\op{(\CstarCPU)}$,
described as a map $A \rightarrow A\otimes B$ in $\CstarCPU$, is given
by $a \mapsto a\sotimes 1$. When working in $\CstarCPU$ we often call
this map a coprojection, for obvious reasons, and write it as
$\kappa_1$.

We check the interaction between instruments and distributivity from
diagrams~\eqref{DistrMeasEqn} in the binary case (for $n=2$). So
assume we have a test $e \in [0,1]_{A}$, with associated instrument
$\instr_{e} \colon A\times A \rightarrow A$ given by
$\instr_{e}(a_{1}, a_{2}) = \sqrt{e}\cdot a_{1} \cdot \sqrt{e} +
\sqrt{1-e}\cdot a_{2} \cdot \sqrt{1-e}$. We get a predicate
$\kappa_{1}^{*}(e) = e\sotimes 1 \in [0,1]_{A\otimes B}$. Since
$(\sqrt{e}\sotimes 1)\cdot (\sqrt{e}\sotimes 1) = (e\sotimes 1)$ we
have $\sqrt{e\sotimes 1} = (\sqrt{e}\sotimes 1)$. Similarly for $1-e$
we have:
$$\begin{array}{rcccccl}
(\sqrt{1-e}\sotimes 1)\cdot (\sqrt{1-e}\sotimes 1)
& = &
(1-e)\sotimes 1 
& = &
1\sotimes 1 - e\sotimes 1 
& = &
1 - e\sotimes 1,
\end{array}$$

\noindent so that $\sqrt{1 - (e\sotimes 1)} = \sqrt{1-e}\sotimes 1$. Now
we can prove~\eqref{DistrMeasEqn}:
$$\begin{array}{rcl}
\lefteqn{\big(\instr_{\kappa_{1}^{*}(e)} \after 
   \theta_{1}\big)((a_{1},a_{2})\sotimes b)} \\
& = &
\instr_{\kappa_{1}^{*}(e)}(a_{1}\sotimes b, a_{2}\sotimes b) \\
& = &
\sqrt{e\otimes 1}\cdot (a_{1}\sotimes b) \cdot \sqrt{e\sotimes 1} +
\sqrt{1-(e\otimes 1)}\cdot (a_{2}\sotimes b) \cdot \sqrt{1-(e\sotimes 1)} \\
& = &
(\sqrt{e}\otimes 1)\cdot (a_{1}\sotimes b) \cdot (\sqrt{e}\sotimes 1) +
(\sqrt{1-e}\otimes 1)\cdot (a_{2}\sotimes b) \cdot (\sqrt{1-e}\sotimes 1) \\
& = &
(\sqrt{e}\cdot a_{1}\cdot \sqrt{e})\otimes b +
   (\sqrt{1-e}\cdot a_{2}\cdot \sqrt{1-e})\otimes b \\
& = &
(\sqrt{e}\cdot a_{1}\cdot \sqrt{e} +
   \sqrt{1-e}\cdot a_{2}\cdot \sqrt{1-e})\otimes b \\
& = &
\instr_{e}(a_{1}, a_{2})\sotimes b \\
& = &
(\instr_{e}\otimes\idmap)((a_{1},a_{2})\sotimes b).
\end{array}$$
\end{exas}

\noindent The object $2 = 1+1$ in a category satisfying
Assumption~\ref{TensorAss} carries a special multiplication map
$m\colon 2\otimes 2 \rightarrow 2$, namely $m =
[\idmap,\kappa_{2}\after\,!]  \after (\rho+\rho) \after \theta_{2}$,
obtained by going east-south-west in:
\begin{equation}
\label{TwoMultDiag}
\vcenter{\xymatrix@R-1pc@C+1pc{
(1+1)\otimes(1+1)\ar[d]_{m} 
   \ar[rr]^-{\theta_{2} = 
      [\idmap\otimes\kappa_{1}, \idmap\otimes\kappa_{2}]^{-1}}_-{\cong}
   & & ((1+1)\otimes 1)+((1+1)\otimes 1)
   \ar[d]_{\cong}^{\rho+\rho} \\
1+1 & & (1+1)+(1+1)\ar[ll]^-{[\idmap, \kappa_{2}\after\,!]} 
}}
\end{equation}

\begin{prop}
\label{TwoMonoidProp}
Let $\cat{B}$ be a category satisfying Assumption~\ref{TensorAss}.
\begin{enumerate}
\item \label{TwoMonoidPropMon} The map $m\colon 2\otimes 2 \rightarrow
  2$ defined in~\eqref{TwoMultDiag} makes $2 = 1+1\in\cat{B}$ into a
  commutative monoid in $\cat{B}$, with unit scalar $1 = \kappa_{1}
  \colon 1 \rightarrow 1+1=2$ as multiplicative unit, and zero scalar
  $0 = \kappa_{2} \colon 1 \rightarrow 2$ as zero element.

\item \label{TwoMonoidPropMult} Multiplication of scalars from
  Lemma~\ref{ScalarMonoidLem} can be expressed via this multiplication
  map $m$, namely as:
$$\xymatrix@C-0.5pc{
s\cdot r
= 
\big(1\ar[rr]^-{\lambda^{-1}=\rho^{-1}}_-{\cong} & &
   1\otimes 1\ar[r]^-{s\otimes r} & 2\otimes 2\ar[r]^-{m} & 2\big).
}$$

\item \label{TwoMonoidPropComm} This multiplication $\cdot$ of scalars
  is then also commutative: $s\cdot r = r\cdot s$. Thus, in presence
  of tensors, the scalars $\Pred(1)$ form a \emph{commutative} monoid.
\end{enumerate}
\end{prop}

\noindent Recall from Section~\ref{StatesSec} the distribution monad $\Dst_{M}
\colon \Sets \rightarrow \Sets$ associated with the effect monoid $M =
\Pred(1)$. This monad is commutative if $M=\Pred(1)$ is commutative,
see~\cite{Jacobs11c}.  By a general categorical result of
Kock~\cite{Kock71a,Kock71b} --- ``commutative theories have tensors''
--- this implies that the category $\Conv_{M} = \EM(\Dst_{M})$ of
Eilenberg-Moore algebras is symmetric monoidal closed.

\begin{myproof}
For the first point we concentrate on commutativity, which is the
most interesting part. We need to prove $m \after \gamma = m$. First
consider the diagram:
$$\xymatrix@R-1pc{
(1\!+\!1)\otimes (1\!+\!1)\ar[d]_{\gamma}^{\cong}\ar[r]^-{\theta_{1}}_-{\cong} & 
   (1\otimes (1\!+\!1))+(1\otimes(1\!+\!1))\ar[d]^{\gamma+\gamma}_{\cong} 
   \ar@/^2ex/[dr]^-{\theta_{2}+\theta_{2}}
\\
(1\!+\!1)\otimes (1\!+\!1)\ar[d]_{m}\ar[r]^-{\theta_{2}}_-{\cong}
   & ((1\!+\!1)\otimes 1)+((1\!+\!1)\otimes 1)
   \ar[d]_{\cong}^{\rho+\rho} &
   ((1\!\otimes\! 1)\!+\!(1\!\otimes\! 1))\!+\!
    ((1\!\otimes\! 1)\!+\!(1\!\otimes\! 1))
   \ar@/^2ex/[dl]^(0.5){\qquad\qquad(\rho+\rho)+(\rho+\rho) = 
        (\lambda+\lambda)+(\lambda+\lambda)} 
\\
1+1 & (1\!+\!1)+(1\!+\!1)\ar[l]^-{[\idmap, \kappa_{2}\after\,!]} 
}$$

\noindent Now we can use the interaction between $\theta_{2}$ and $\theta_{1}$
from~\eqref{DistrDistrDiag} in:
$$\begin{array}{rcl}
m \after \gamma
& = &
[\idmap, \kappa_{2} \after\, !] \after (\rho+\rho) \after \theta_{2} \after 
   \gamma \\
& = &
[\idmap, \kappa_{2} \after\, !] \after 
   ((\rho+\rho) + (\rho+\rho)) \after
   (\theta_{2}+\theta_{2}) \after \theta_{1} \\
& \smash{\stackrel{\eqref{DistrDistrDiag}}{=}} &
[\idmap, \kappa_{2} \after\, !] \after 
   ((\rho+\rho) + (\rho+\rho)) \after \widehat{\gamma} \after
   (\theta_{1}+\theta_{1}) \after \theta_{2} \\
& = &
[\idmap, \kappa_{2} \after\, !] \after \widehat{\gamma} \after 
   ((\rho+\rho) + (\rho+\rho)) \after 
   (\theta_{1}+\theta_{1}) \after \theta_{2} \\
& = &
[\idmap, \kappa_{2} \after\, !] \after [\kappa_{1}+\kappa_{1}, 
   \kappa_{2}+\kappa_{2}] \after (\rho+\rho) \after \theta_{2} \\
& = &
[\idmap, \kappa_{2} \after\, !] \after (\rho+\rho) \after \theta_{2} \\
& = &
m.
\end{array}$$

\noindent A crucial part of this proof is that $\lambda = \rho \colon
1\otimes 1 \rightarrow 1$. This is like in the Eckmann-Hilton style
argument in~\cite{KellyL80}. We also check that $\kappa_{1}\colon 1
\rightarrow 1+1 = 2$ is the unit for $m$:
\noindent Thus:
$$\begin{array}{rcl}
m \after (\idmap\otimes\kappa_{1}) 
& = &
[\idmap, \kappa_{2}\after\,!] \after (\rho+\rho) \after \theta_{2} \after
   (\idmap\otimes\kappa_{1}) \\
& = &
[\idmap, \kappa_{2}\after\,!] \after (\rho+\rho) \after \kappa_{1}
   \qquad \mbox{since } \theta_{2}^{-1} \after \kappa_{1} = \idmap\otimes\kappa_{1} \\
& = &
[\idmap, \kappa_{2}\after\,!] \after \kappa_{1} \after \rho \\
& = &
\rho.
\end{array}$$

\noindent Remaining steps to show that $m\colon 2\otimes 2 \rightarrow
2$ is a monoid are left to the reader.

\auxproof{
We still double-check that the upper and right parts of the above
diagram commute.
$$\begin{array}{rcl}
\theta_{1}^{-1} \after (\gamma+\gamma)
& = &
[\kappa_{1}\otimes\idmap, \kappa_{2}\otimes\idmap]  \after (\gamma+\gamma) \\
& = &
[(\kappa_{1}\otimes\idmap) \after \gamma, 
   (\kappa_{2}\otimes\idmap) \after \gamma] \\
& = &
[\gamma \after (\idmap\otimes\kappa_{1}), 
   \gamma \after (\idmap\otimes\kappa_{2})] \\
& = &
\gamma \after [\idmap\otimes\kappa_{1}, \idmap\otimes\kappa_{2}] \\
& = &
\gamma \after \theta_{2}^{-1} \\
(\rho+\rho) \after (\gamma+\gamma) \after (\theta_{2}^{-1}+\theta_{2}^{-1})
& = &
(\lambda+\lambda) \after ([\idmap\otimes\kappa_{1}, \idmap\otimes\kappa_{2}]
   + [\idmap\otimes\kappa_{1}, \idmap\otimes\kappa_{2}]) \\
& = &
[\lambda \after (\idmap\otimes\kappa_{1}), 
   \lambda \after (\idmap\otimes\kappa_{2})]
   + [\lambda \after (\idmap\otimes\kappa_{1}), 
    \lambda \after (\idmap\otimes\kappa_{2})] \\
& = &
[\kappa_{1} \after \lambda, \kappa_{2} \after \lambda] + 
   [\kappa_{1} \after \lambda, \kappa_{2} \after \lambda] \\
& = &
(\lambda+\lambda)+(\lambda+\lambda).
\end{array}$$
}

\auxproof{
$$\xymatrix@R-1pc@C+1pc{
1+1\ar[d]_{\rho^{-1}}^{\cong} \\
(1+1)\otimes 1\ar[d]_{\idmap\otimes\kappa_{1}}\ar@/^2ex/[drr]^-{\kappa_1} \\
(1+1)\otimes(1+1)\ar[d]_{m} 
   \ar[rr]^-{\theta_{2} = 
      [\idmap\otimes\kappa_{1}, \idmap\otimes\kappa_{2}]^{-1}}_-{\cong}
   & & ((1+1)\otimes 1)+((1+1)\otimes 1)
   \ar[d]_{\cong}^{\rho+\rho} \\
1+1 & & (1+1)+(1+1)\ar[ll]^-{[\idmap, \kappa_{2}\after\,!]} 
}$$

Similarly, we obtain that $\kappa_{2}$ is the zero element for $m$:
$$\begin{array}{rcl}
m \after (\idmap\otimes\kappa_{2})
& = &
[\idmap, \kappa_{2}\after\,!] \after (\rho+\rho) \after \theta_{2} \after
   (\idmap\otimes\kappa_{2}) \\
& = &
[\idmap, \kappa_{2}\after\,!] \after (\rho+\rho) \after \kappa_{2} \\
& = &
[\idmap, \kappa_{2}\after\,!] \after \kappa_{2} \after \rho  \\
& = &
\kappa_{2}\after\,! \after \rho \\
& = &
\kappa_{2}\after\,!.
\end{array}$$

For associativity of $m$ we need the following associativity result
for $\theta_{2}$.
$$\xymatrix{
U\otimes ((A+B)\otimes(X+Y))\ar[d]_{\idmap\otimes\theta_{2}}\ar[r]^-{\alpha} &
   (U\otimes (A+B))\otimes(X+Y)\ar[dd]^{\theta_{2}}
\\
U\otimes (((A+B)\otimes X)+((A+B)\otimes Y))\ar[d]_{\theta_{2}} &
\\
(U\otimes ((A+B)\otimes X))+(U\otimes ((A+B)\otimes Y))\ar[r]^-{\alpha+\alpha} &
   ((U\otimes (A+B))\otimes X)+((U\otimes (A+B))\otimes Y)
}$$

\noindent It is obtained via an easy computation:
$$\begin{array}{rcl}
\lefteqn{\big((\idmap\otimes\theta_{2}^{-1}) \after \theta_{2}^{-1} \after 
   (\alpha^{-1}+\alpha^{-1})\big)} \\
& = &
(\idmap\otimes\theta_{2}^{-1}) \after 
   [\idmap\otimes\kappa_{1}, \idmap\otimes\kappa_{2}] \after
   (\alpha^{-1}+\alpha^{-1})\big) \\
& = &
[(\idmap\otimes\theta_{2}^{-1}) \after (\idmap\otimes\kappa_{1}) \after
   \alpha^{-1}, (\idmap\otimes\theta_{2}^{-1}) \after 
   (\idmap\otimes\kappa_{2}) \after \alpha^{-1}] \\
& = &
[(\idmap\otimes(\idmap\otimes\kappa_{1})) \after \alpha^{-1}, 
   (\idmap\otimes(\idmap\otimes\kappa_{2})) \after \alpha^{-1}] \\
& = &
[\alpha^{-1} \after ((\idmap\otimes\idmap)\otimes\kappa_{1}),
   \alpha^{-1} \after ((\idmap\otimes\idmap)\otimes\kappa_{2})] \\
& = &
\alpha^{-1} \after [\idmap\otimes\kappa_{1}, \idmap\otimes\kappa_{2}] \\
& = &
\alpha^{-1} \after \theta_{2}^{-1}.
\end{array}$$

Finally we show associativity of $m$. The first few steps use
naturality of $\theta_{2}$. Then we can use the above diagram.
Interestingly, we have to use that $\kappa_2$ is zero element for $m$.
$$\begin{array}{rcl}
\lefteqn{m \after (m\otimes\idmap) \after \alpha} \\
& = &
[\idmap, \kappa_{2}\after\,!] \after (\rho+\rho) \after \theta_{2} \after
   (m\otimes\idmap) \after \alpha \\
& = &
[\idmap, \kappa_{2}\after\,!] \after (\rho+\rho) \after 
   ((m\otimes\idmap)+(m\otimes\idmap)) \after \theta_{2} \after \alpha \\
& = &
[\idmap, \kappa_{2}\after\,!] \after (m+m) \after (\rho+\rho) \after 
   (\alpha+\alpha) \after \theta_{2} \after (\idmap\otimes\theta_{2}) \\
& = &
[m, \kappa_{2}\after\,!] \after ((\idmap\otimes\rho)+(\idmap\otimes\rho)) 
   \after \theta_{2} \after (\idmap\otimes\theta_{2}) \\
& = &
[m, m \after (\idmap\otimes(\kappa_{2}\after\,!))] \after \theta_{2} \after
   (\idmap\otimes(\rho+\rho)) \after (\idmap\otimes\theta_{2}) \\
& & \qquad \mbox{since $\kappa_2$ is zero element for $m$} \\
& = &
m \after \nabla \after (\idmap + (\idmap\otimes(\kappa_{2}\after\,!))) 
   \after \theta_{2} \after
   (\idmap\otimes(\rho+\rho)) \after (\idmap\otimes\theta_{2}) \\
& = &
m \after \nabla \after \theta_{2} \after 
   (\idmap\otimes(\idmap+(\kappa_{2}\after\,!))) \after 
   (\idmap\otimes(\rho+\rho)) \after (\idmap\otimes\theta_{2}) \\
& = &
m \after (\idmap\otimes\nabla) \after 
   (\idmap\otimes(\idmap+(\kappa_{2}\after\,!))) \after 
   (\idmap\otimes(\rho+\rho)) \after (\idmap\otimes\theta_{2}) \\
& = &
m \after (\idmap\otimes[\idmap,\kappa_{2}\after\,!]) \after 
   (\idmap\otimes(\rho+\rho)) \after (\idmap\otimes\theta_{2}) \\
& = &
m \after (\idmap\otimes m).
\end{array}$$
}

For the second point we first prove $s\cdot r = m \after (s\otimes r)
\after \lambda^{-1}$ via the following diagram.
$$\xymatrix@C+1pc{
1\ar[d]_{r}\ar[r]^-{\lambda^{-1}=\rho^{-1}}_-{\cong} &
   1\otimes 1\ar[d]_{\idmap\otimes r}\ar[r]^-{s\otimes r} &
   2\otimes 2\ar[r]^-{\theta_{2}} &
   2\otimes 1 + 2\otimes 1\ar[r]^-{\rho+\rho} &
   2+2\ar[d]^{[\idmap, \kappa_{2} \after \,!]}
\\
2\ar[r]^-{\lambda^{-1}}\ar@{=}@/_4ex/[rrr] &
   1\otimes 2\ar[ur]^(0.5){s\otimes\idmap}\ar[r]^-{\theta_{2}} &
   1\otimes 1 + 1\otimes 1\ar[ur]^(0.4){s\otimes\idmap+s\otimes\idmap\quad}
      \ar[r]^-{\rho+\rho = \lambda+\lambda} &
   1+1\ar[ur]^{s+s}\ar[r]_-{[s,\kappa_{2}]} &
   2
}$$

\auxproof{
We briefly check the equality at the bottom of this rectangle.
$$\begin{array}{rcl}
\lambda \after \theta_{2}^{-1} \after (\rho+\rho)^{-1}
& = &
\lambda \after [\idmap\otimes\kappa_{1}, \idmap\otimes\kappa_{2}]
   \after (\lambda^{-1}+\lambda^{-1}) \\
& = &
[\lambda \after (\idmap\otimes\kappa_{1}) \after \lambda^{-1}, 
   \lambda \after (\idmap\otimes\kappa_{2}) \after \lambda^{-1}] \\
& = &
[\kappa_{1}, \kappa_{2}] \\
& = &
\idmap.
\end{array}$$
}

\bigskip

\noindent But now commutativity of the multiplication $\cdot$ of predicates
follows from commutativity of $m$ in:
$$\begin{array}[b]{rcl}
s\cdot r
\hspace*{\arraycolsep} = \hspace*{\arraycolsep}
m \after (s\otimes r) \after \lambda^{-1}
& = &
m \after \gamma \after (s\otimes r) \after \lambda^{-1} 
   \qquad \mbox{by the first point} \\
& = &
m \after (r\otimes s) \after \gamma \after \lambda^{-1} 
\hspace*{\arraycolsep} = \hspace*{\arraycolsep}
m \after (r\otimes s) \after \rho^{-1}
\hspace*{\arraycolsep} = \hspace*{\arraycolsep}
r\cdot s.
\end{array}\eqno{\qEd}$$
\end{myproof}

\begin{exas}
\label{TwoMonoidEx}
We briefly review what the multiplication map $m \colon 2\otimes 2
\rightarrow 2$ from~\eqref{TwoMultDiag} is in our main examples.  In
the category $\Sets$ we have $2 = \{0,1\}$ with $m\colon 2\times 2
\rightarrow 2$ given by multiplication, that is, by conjunction.
Similarly, in the Kleisli category $\Kl(\Dst)$ of the distribution
monad the map $m \colon 2\times 2 \rightarrow \Dst(2) \cong [0,1]$ is
multiplication. In the category of $C^*$-algebras this multiplication
is a positive unital map $m \colon \C^{2} \rightarrow
\C^{2}\otimes\C^{2}$, given on the two standard basis vectors $\ket{0}
= \left(\begin{smallmatrix} 1 \\ 0 \end{smallmatrix}\right)$ and
$\ket{1} = \left(\begin{smallmatrix} 0 \\ 1 \end{smallmatrix}\right)$
in $\C^2$ as:
\begin{equation}
\label{TwoMonoidExCstarEqn}
\begin{array}{rclcrcl}
m(\ket{0})
& = &
\ket{0}\sotimes\ket{0}
& \qquad\mbox{and}\qquad &
m(\ket{1})
& = &
\ket{0}\sotimes\ket{1} + \ket{1}\sotimes\ket{0} + \ket{1}\sotimes\ket{1}.
\end{array}
\end{equation}

\auxproof{
This multiplication map $m \colon \C^{2} \rightarrow
\C^{2}\otimes\C^{2}$ is given as composite:
$$\xymatrix@C+1pc{
\C^{2}\ar[d]_{m}\ar[rr]^-{\tuple{\idmap, !\,\after\pi_{2}}} & & 
   \C^{2}\oplus \C^{2}\ar[d]^{\rho^{-1}\oplus\rho^{-1}} \\
\C^{2}\otimes\C^{2} & & (\C^{2}\otimes\C)\oplus(\C^{2}\otimes\C)
   \ar[ll]^-{\tuple{\idmap\otimes\pi_{1}, \idmap\otimes\pi_{2}}^{-1}}_-{\cong}
}$$

\noindent The upper map $f = \tuple{\idmap, !\,\after\pi_{2}} \colon
\C^{2} \rightarrow \C^{2} \oplus \C^{2}$ is given by:
$$\begin{array}{rcl}
f\left(\begin{smallmatrix} z \\ w \end{smallmatrix}\right)
& = &
\left(\begin{smallmatrix} z \\ w \\ w \\ w \end{smallmatrix}\right).
\end{array}$$

\noindent The map $g = (\rho\oplus\rho) \after
\tuple{\idmap\otimes\pi_{1}, \idmap\otimes\pi_{2}} \colon
\C^{2}\otimes\C^{2} \rightarrow \C^{2}\oplus \C^{2}$ is:
$$\begin{array}{rcccl}
g\left(\left(\begin{smallmatrix} z_{1} \\ z_{2} \end{smallmatrix}\right) \sotimes
   \left(\begin{smallmatrix} w_{1} \\ w_{2} \end{smallmatrix}\right)\right)
& = &
\left(\begin{smallmatrix} 
w_{1}\left(\begin{smallmatrix} z_{1} \\ z_{2} \end{smallmatrix}\right) \\ 
w_{2}\left(\begin{smallmatrix} z_{1} \\ z_{2} \end{smallmatrix}\right)  
   \end{smallmatrix}\right) 
& = &
\left(\begin{smallmatrix} w_{1}z_{1} \\ w_{1}z_{2} \\
   w_{2}z_{1} \\ w_{2}z_{2} \end{smallmatrix}\right) 
\end{array}$$

\noindent Thus, with the multiplication map $m\colon \C^{2} \rightarrow
\C^{2}\otimes\C^{2}$ as defined above we have:
$$\begin{array}{rcl}
\big(g \after m\big)(\ket{0})
& = &
g\big(\ket{0}\sotimes\ket{0}\big) \\
& = &
g\big(\left(\begin{smallmatrix} 1 \\ 0 \end{smallmatrix}\right)
   \sotimes\left(\begin{smallmatrix} 1 \\ 0 \end{smallmatrix}\right)\big) \\
& = &
\left(\begin{smallmatrix} 1 \\ 0 \\ 0 \\ 0 \end{smallmatrix}\right) \\
& = &
f(\ket{0}) \\
\big(g \after m\big)(\ket{1})
& = &
g\big(\ket{0}\sotimes\ket{1} + \ket{1}\sotimes\ket{0} + 
   \ket{1}\sotimes\ket{1}\big) \\
& = &
g\big(\left(\begin{smallmatrix} 1 \\ 0 \end{smallmatrix}\right)
   \sotimes\left(\begin{smallmatrix} 0 \\ 1 \end{smallmatrix}\right)\big) +
g\big(\left(\begin{smallmatrix} 0 \\ 1 \end{smallmatrix}\right)
   \sotimes\left(\begin{smallmatrix} 1 \\ 0 \end{smallmatrix}\right)\big) +
g\big(\left(\begin{smallmatrix} 0 \\ 1 \end{smallmatrix}\right)
   \sotimes\left(\begin{smallmatrix} 0 \\ 1 \end{smallmatrix}\right)\big) \\
& = &
\left(\begin{smallmatrix} 0 \\ 0 \\ 1 \\ 0 \end{smallmatrix}\right) +
\left(\begin{smallmatrix} 0 \\ 1 \\ 0 \\ 0 \end{smallmatrix}\right) +
\left(\begin{smallmatrix} 0 \\ 0 \\ 0 \\ 1 \end{smallmatrix}\right) \\
& = &
\left(\begin{smallmatrix} 0 \\ 1 \\ 1 \\ 1 \end{smallmatrix}\right) \\
& = &
f(\ket{1}).
\end{array}$$
}
\end{exas}\medskip

\noindent In the remainder of this section we consider projections and pairings,
both for states and for predicates. We shall see that projections
after pairings on states return the original results, but not the
other way around. The latter precisely captures dependence and
entanglement, see Example~\ref{DependenceEntanglementEx}. We start
with the projections described in~\eqref{ProjTensorEqn}.

\begin{lem}
\label{ProjPureLem}
The projections $X \leftarrow X\otimes Y \rightarrow Y$
in~\eqref{ProjTensorEqn} that can be defined in a category satisfying
Assumption~\ref{TensorAss} are pure. As a consequence, also the unique
maps $X \rightarrow 1$ to the terminal object are pure.
\end{lem}

\begin{myproof}
We show that the first projection $\pi_{1}\colon X\otimes Y
\rightarrow X$ is pure. So let $p\colon X \rightarrow n\cdot 1$ be an
$n$-test. We have to prove $(n\cdot \pi_{1}) \after
\instr_{\pi_{1}^{*}(p)} = \instr_{p} \after \pi_{1}$, see
Definition~\ref{PureDef}.  We use requirement~\eqref{TensorAssInstr}
in Assumption~\ref{TensorAss} as the upper part of the diagram:
$$\xymatrix@R-.5pc@C-.5pc{
X\otimes Y\ar[dd]_{\pi_1}\ar[rr]^-{\instr_{p}\otimes\idmap}
   \ar@/_2ex/[drr]_{\instr_{\pi_{1}^{*}(p)}} 
   & & (n\cdot X)\otimes Y\ar[d]_{\cong}^{\theta_1}\ar@/^8ex/ [dd]^{\pi_1} \\
& & n\cdot (X\otimes Y)\ar[d]_{n\cdot\pi_{1}} \\
X\ar[rr]_-{\instr_p} & & n\cdot X
}$$

\noindent In this diagram we have $(n\cdot\pi_{1}) \after \theta_{1} = \pi_{1}$
since:
$$\begin{array}{rcccccccl}
\pi_{1} \after \theta_{1}^{-1}
& = &
\pi_{1} \after [\kappa_{1}\otimes\idmap]_{i\leq n}
& = &
[\pi_{1} \after (\kappa_{1}\otimes\idmap)]_{i\leq n}
& = &
[\kappa_{i} \after \pi_{1}]_{i\leq n}
& = &
n\cdot \pi_{1}.
\end{array}$$

\noindent But then we are done since by naturality of $\pi_1$ the outer
diagram commutes. Finally, the unique map $X \rightarrow 1$ is pure,
since it can be written as composite of pure maps:
$$\xymatrix{
X\ar[r]^-{\lambda^{-1}}_-{\cong} & 1\otimes X\ar[r]^-{\pi_{1}} & 1
}\eqno{\qEd}$$

\auxproof{
$$\begin{array}{rcccccl}
\pi_{1} \after \lambda^{-1} 
& = &
\rho \after (\idmap\otimes\,!) \after \lambda^{-1} 
& = &
\rho \after \lambda^{-1} \after \,! 
& = &
!.
\end{array}$$
}
\end{myproof}\medskip

\noindent In the sequel we restrict ourselves to the first projection $\pi_{1}
\colon X\otimes Y \rightarrow X$, since the second one can be handled
symmetrically. Via the predicate and state functors $\Pred$ and
$\Stat$ from the state-and-effect triangle~\eqref{GeneralTriangleDiag}
we obtain maps (of effect modules and of convex sets):
$$\xymatrix{
\Pred(X)\ar[rr]^-{\Pred(\pi_{1}) = (\pi_{1})^*} & & \Pred(X\otimes Y)
&
\Stat(X\otimes Y)\ar[rr]^-{\Stat(\pi_{1}) = (\pi_{1})_{*}} & & \Stat(X)
}$$

\noindent The first operation $(\pi_{1})^*$ sends a predicate on $X$
to a predicate on the extended type $X\otimes Y$. It corresponds to
\emph{weakening} in logic, that is, to addition of an unused
variable. The second operation $(\pi_{1})_{*}$ corresponds to taking
the \emph{marginal}, which, in the case of Hilbert spaces is the
partial trace, see Lemma~\ref{MarginalPartTraceLem} below.

In the category $\Sets$ the first projection applied to a state
$\omega = (x, y)\in X\times Y$ is simply the ordinary projection
$\pi_{1}\omega = x \in X$. The situation is slightly more interesting
in the Kleisli category $\Kl(\Dst)$ of the distribution monad. There,
for a distribution $\omega \in \Dst(X\times Y)$, the projection
$\pi_{1}\omega \in \Dst(X)$ is the marginal distribution, namely:
$$\begin{array}{rcccl}
\pi_{1}\omega(x)
& = &
\Dst(\pi_{1})(\omega)(x)
& = &
\sum_{y\in Y}\omega(x,y).
\end{array}$$

\noindent In $\op{(\CstarCPU)}$ the left marginal of a state $\omega
\colon A\otimes B \rightarrow \C$ is the map $\pi_{1}\omega \colon A
\rightarrow \C$ given by $\pi_{1}\omega(a) = \omega(a\sotimes 1)$.
When the $C^*$-algebra is of the form $A = \B(\H)$ for a
finite-dimensional Hilbert space $\H$, we obtain the so-called partial
trace.

\begin{lem}
\label{MarginalPartTraceLem}
Let $\H, \K$ be two finite-dimensional Hilbert space. The left
marginal, described on density matrices via the chain
$$\xymatrix{
\DM(\H\otimes \K) \smash{\stackrel{\eqref{HilbStatEqn}}{\cong}} 
   \Stat(\B(\H\otimes\K))
   \smash{\stackrel{\eqref{HilbCstarTensorEqn}}{\cong}} 
   \Stat(\B(\H)\otimes \B(\K))\ar[r]^-{(\pi_{1})_{*}} &
   \Stat(\B(\H)) \smash{\stackrel{\eqref{HilbStatEqn}}{\cong}}  \DM(\H)
}$$

\noindent is the partial trace operation $\rho \mapsto \tr_{\K}(\rho)$.
\end{lem}

\begin{myproof}
Each density matrix $\rho\in\DM(\H\otimes\K)$ gives rise to a state
$\omega_{\rho} \colon \B(\H) \otimes \B(\K) \rightarrow \C$ given by
$\omega_{\rho}(M\sotimes N) = \tr(\rho(M\otimes N))$. The resulting
marginal $\pi_{1}\omega_{\rho} \colon \B(\H) \rightarrow \C$ is
$\pi_{1}\omega_{\rho}(M) = \tr(\rho(M\otimes\idmap))$. The partial
trace $\tr_{\K}(\rho) \in \DM(\H)$ precisely satisfies
$\tr(\tr_{\K}(\rho)M) = \tr(\rho(M\otimes\idmap))$, see
also~\cite[Defn.~2.68]{HeinosaariZ12}. \QED
\end{myproof}

We continue with pairings.

\begin{defi}
\label{StatePredTensorDef}
Assume we are in category $\cat{B}$ satisfying
Assumption~\ref{TensorAss}.  For two states $\omega_{1} \colon 1
\rightarrow X_{1}$ and $\omega_{2} \colon 1 \rightarrow X_{2}$ we
write $\omega_{1}\sotimes\omega_{2} \colon 1 \rightarrow X_{1}\otimes
X_{2}$ for the state obtained by:
$$\xymatrix@C-.5pc{
\omega_{1}\sodot\omega_{2} = \big(1\ar[rr]^-{\lambda^{-1} = \rho^{-1}} & &
   1\otimes 1\ar[rr]^-{\omega_{1}\otimes\omega_{2}} & & 
   X_{1}\otimes X_{2}\big).
}$$

\noindent Similarly, for two predicates $p_{1} \colon X_{1}
\rightarrow 1+1$ and $p_{2} \colon X_{2} \rightarrow 1+1$ we define a
pairing predicate $p_{1}\sodot p_{2} \colon X_{1}\otimes X_{2}
\rightarrow 1+1$ as:
$$\xymatrix@C-1.5pc{
p_{1}\sodot p_{2} = (X_{1} \otimes X_{2}\ar[rrr]^-{p_{1}\otimes p_{2}} & & &
   2\otimes 2\ar[rr]^-{m} & & 2\big),
}$$

\noindent where $m\colon 2\otimes 2 \rightarrow 2$ is the
multiplication map from~\eqref{TwoMultDiag}.

One can freely extend these pairing definitions to multiple states and
predicates, like in $\omega_{1} \sodot \cdots \sodot \omega_{n}$ and
$p_{1} \sodot \cdots \sodot p_{n}$, where we implicitly use
associativity.
\end{defi}

The next result presents some basic results about these pairings.  The
first one, for instance, could be useful for modular reasoning, where
separate properties for separate states can be proved separately.

\begin{lem}
\label{StatePredTensorLem}
The pairings on states and predicates defined above satisfy:
\begin{enumerate}
\item \label{StatePredTensorLemModularity} $(\omega_{1}\sodot
  \omega_{2} \models p_{1} \sodot p_{2}) = m(\omega_{1}\models p_{1},
  \omega_{2} \models p_{2})$;

\item \label{StatePredTensorLemProj} $\pi_{i} \after
  (\omega_{1}\sodot\omega_{2}) = \omega_{i}$;

\item the pairing $\sodot$ on states is bi-affine,
  \textit{i.e.}~affine in each argument separately;

\item weakening is pairing with truth: $(\pi_{1})^{*}(p) = p \sodot 1$,
where $1 = \kappa_{1} \after\,!$ is the truth predicate;

\item the pairing $\sodot$ on predicates is a bihomomorphism of effect
  modules: the mappings $p\sodot (-)$ and $(-)\sodot p$ preserve
  finite sums $\ovee, 0$ and scalar multiplication, and satisfy
  $1\sodot 1 = 1$.
\end{enumerate}
\end{lem}

\begin{myproof}
The first point immediately holds, and for the second one we have:
$$\begin{array}{rcccccl}
\pi_{1} \after (\omega_{1}\sodot\omega_{2})
& = &
\rho \after (\idmap\otimes\,!) \after (\omega_{1}\otimes\omega_{2})
   \after \rho^{-1}
& = &
\rho \after (\omega_{1}\otimes\idmap) \after \rho^{-1}
& = &
\omega_{1}.
\end{array}$$

\auxproof{
For the first point we have:
$$\begin{array}{rcl}
\omega_{1}\sodot \omega_{2} \models p_{1} \sodot p_{2}
& = &
(p_{1} \sodot p_{2}) \after (\omega_{1}\sodot \omega_{2}) \\
& = &
m \after (p_{1}\otimes p_{2}) \after (\omega_{1}\otimes\omega_{2}) \after
   \rho^{-1} \\
& = &
m \after ((p_{1} \after \omega_{1})\otimes (p_{2} \after \omega_{2})) \after
   \rho^{-1} \\
& = &
m \after ((\omega_{1}\models p_{1})\otimes (\omega_{2}\models p_{2})) \after
   \rho^{-1} \\
& = &
m(\omega_{1}\models p_{1}, \omega_{2} \models p_{2}).
\end{array}$$
}

\noindent We show that $\sodot$ is bi-affine. Suppose we have an
$n$-test $p\colon 1 \rightarrow n\cdot 1$, and $n$ states $\omega_{i}
\colon 1 \rightarrow X$. Then for each state $\omega \colon 1
\rightarrow Y$,
$$\begin{array}{rcll}
\big(\bigovee_{i}p_{i}\omega_{i}\big)\sodot \omega
& = &
(([\omega_{1}, \ldots, \omega_{n}] \after p)\otimes \omega) \after \rho^{-1} \\
& = &
([\omega_{1}, \ldots, \omega_{n}]\otimes \omega) \after 
   (\instr_{p}\otimes\idmap) \after \rho^{-1}  
   & \mbox{by 
   Lemma~\ref{MeasurementAssSefLem}~\eqref{MeasurementAssSefLemOne}} \\
& = &
([\omega_{1}, \ldots, \omega_{n}]\otimes \omega) \after \theta_{1}^{-1} \after
   \instr_{\pi_{1}^{*}(p)} \after \rho^{-1} \quad
   & \mbox{by~\eqref{DistrMeasEqn}} \\
& = &
[\omega_{1}\otimes\omega, \ldots, \omega_{n}\otimes\omega] \after
   n\cdot \rho \after \instr_{p} \\
& & \qquad \mbox{since $\pi_{1} = 
      \rho\colon (n\cdot 1)\otimes 1 \rightarrow n\cdot 1$ is pure
    by \rlap{Assumption~\ref{TensorAss}~\eqref{TensorAssMonIsoPure}}} \\
& = &
[\omega_{1}\sodot\omega, \ldots, \omega_{n}\sodot\omega] \after p
   & \mbox{by 
   Lemma~\ref{MeasurementAssSefLem}~\eqref{MeasurementAssSefLemOne} 
   again} \\
& = &
\bigovee_{i} p_{i}(\omega_{i}\sodot\omega).
\end{array}$$\medskip

\noindent For the fourth point consider a `weakened' predicate
$(\pi_{1})^{*}(p) = p \after \pi_{1} \colon X\otimes Y \rightarrow X
\rightarrow 1+1$. We use that $\kappa_{1}\colon 1 \rightarrow 1+1$ is
the neutral element of $m$, see (the proof of)
Proposition~\ref{TwoMonoidProp}~\eqref{TwoMonoidPropMon} in:
$$\begin{array}[b]{rcl}
p \sodot 1
\hspace*{\arraycolsep} = \hspace*{\arraycolsep}
m \after (p \otimes 1) 
& = &
m \after (\idmap\otimes\kappa_{1}) \after (p \otimes\,!_{Y}) \\
& = &
\rho \after (p\otimes\idmap) \after (\idmap\otimes\,!_{Y}) \\
& = &
p \after \rho \after (\idmap\otimes\,!_{Y}) \\
& = &
p \after \pi_{1} \\
& = &
(\pi_{1})^{*}(p).
\end{array}$$
Finally, we show that $(-)\sodot p \colon \Pred(X) \rightarrow
\Pred(X\otimes Y)$ preserves sums and scalar
multiplication. Preservation of falsum is easy, since $0 = \kappa_{2}
\colon 1 \rightarrow 2$ is zero element for $m$, see
Proposition~\ref{TwoMonoidProp}~\eqref{TwoMonoidPropMon}:
$$\begin{array}{rcccccccl}
0\sodot p
& = &
m \after ((\kappa_{2} \after\,!_{Y})\otimes p)
& = &
\kappa_{2} \after\, !_{2\otimes 1} \after (!_{X}\otimes p) 
& = &
\kappa_{2} \after\, !_{X\otimes Y}
& = &
0.
\end{array}$$

\noindent Next, if $q_{1} \orthogonal q_{2}$, for $q_{1},q_{2} \colon
X \rightarrow 1+1$, via bound $b\colon X \rightarrow (1+1)+1$, then we
construct a bound $c$ on $X\otimes Y$ via the east-south-west
composite in:
$$\xymatrix@R-.5pc{
X\otimes Y\ar[d]_{c}\ar[r]^-{\idmap\otimes p} & 
  X\otimes(1+1)\ar[r]^-{\theta_{2}}_-{\cong} &
  X\otimes 1 + X\otimes 1\ar[d]^{\rho+!} \\
(1+1)+1 & & X+1\ar[ll]^-{[b,\kappa_{2}]}
}$$

\noindent This bound $c$ proves $q_{1}\sodot p \orthogonal q_{2}
\sodot p$ and $(q_{1}\sodot p) \ovee (q_{2} \sodot p) = (q_{1}\ovee
q_{2}) \sodot p$.  Preservation of scalar multiplication, and $1
\sodot 1 = 1$ are left to the reader. \QED

\auxproof{
$$\begin{array}{rcl}
[\idmap,\kappa_{2}] \after c
& = &
[\idmap,\kappa_{2}] \after [b,\kappa_{2}] \after (\rho+\,!) 
   \after \theta_{2} \after (\idmap\otimes p) \\
& = &
[[\idmap,\kappa_{2}] \after b,\kappa_{2}] \after (\rho+\,!) 
   \after \theta_{2} \after (\idmap\otimes p) \\
& = &
[q_{1},\kappa_{2}] \after (\rho+\,!) 
   \after \theta_{2} \after (\idmap\otimes p) \\
& = &
[\idmap,\kappa_{2}] \after (q_{1}+\idmap) \after (\rho+\,!) 
   \after \theta_{2} \after (\idmap\otimes p) \\
& = &
[\idmap,\kappa_{2}] \after (\rho+\,!) \after 
   ((q_{1}\otimes\idmap+\idmap)
   \after \theta_{2} \after (\idmap\otimes p) \\
& = &
[\idmap,\kappa_{2}] \after (\rho+\,!) \after 
   ((q_{1}\otimes\idmap+(q_{1}\otimes\idmap))
   \after \theta_{2} \after (\idmap\otimes p) \\
& = &
[\idmap,\kappa_{2} \after\,!] \after (\rho+\rho) \after \theta_{2} \after
   (q_{1}\otimes(\idmap+\idmap)) \after (\idmap\otimes p) \\
& = &
[\idmap,\kappa_{2} \after\,!] \after (\rho+\rho) \after \theta_{2} \after
   (q_{1}\otimes p) \\
& = &
m \after (q_{1}\otimes p) \\
& = &
q_{1}\sodot p
\\
{[[\kappa_{2},\kappa_{1}],\kappa_{2}]} \after c
& = &
[[\kappa_{2},\kappa_{1}],\kappa_{2}] \after [b,\kappa_{2}] \after (\rho+\,!) 
   \after \theta_{2} \after (\idmap\otimes p) \\
& = &
[[[\kappa_{2},\kappa_{1}],\kappa_{2}] \after b,\kappa_{2}] \after (\rho+\,!) 
   \after \theta_{2} \after (\idmap\otimes p) \\
& = &
[q_{2},\kappa_{2}] \after (\rho+\,!) 
   \after \theta_{2} \after (\idmap\otimes p) \\
& = &
\cdots \\
& = &
q_{2}\sodot p
\end{array}$$

\noindent But then:
$$\begin{array}{rcl}
(q_{1}\sodot p)\ovee(q_{2}\sodot p)
& = &
(\nabla+\idmap) \after c \\
& = &
(\nabla+\idmap) \after [b,\kappa_{2}] \after (\rho+\,!) 
   \after \theta_{2} \after (\idmap\otimes p) \\
& = &
[(\nabla+\idmap) \after b,\kappa_{2}] \after (\rho+\,!) 
   \after \theta_{2} \after (\idmap\otimes p) \\
& = &
[q_{1}\ovee q_{2},\kappa_{2}] \after (\rho+\,!) 
   \after \theta_{2} \after (\idmap\otimes p) \\
& = &
\cdots \\
& = &
(q_{1}\ovee q_{2})\sodot p.
\end{array}$$

Next, we have
$$\begin{array}{rcl}
(s\scalar q) \sodot p
& = &
m \after (([s,\kappa_{2}] \after q)\otimes p) \\
& = &
[\idmap,\kappa_{2}\after\,!] \after (\rho+\rho) \after \theta_{2} \after
   ([s,\kappa_{2}]\otimes (\idmap+\idmap)) \after (q\otimes p) \\
& = &
[\idmap,\kappa_{2}\after\,!] \after (\rho+\rho) \after 
   (([s,\kappa_{2}]\otimes\idmap)+([s,\kappa_{2}]\otimes\idmap)) \after 
   \theta_{2} \after (q\otimes p) \\
& = &
[\idmap,\kappa_{2}\after\,!] \after 
   ([s,\kappa_{2}]+[s,\kappa_{2}]) \after 
   (\rho+\rho) \after \theta_{2} \after (q\otimes p) \\
& = &
[[s,\kappa_{2}],\kappa_{2}\after\,!] \after 
   (\rho+\rho) \after \theta_{2} \after (q\otimes p) \\
& = &
[s,\kappa_{2}] \after [\idmap, \kappa_{2}\after\,!] \after 
   (\rho+\rho) \after \theta_{2} \after (q\otimes p) \\
& = &
[s,\kappa_{2}] \after (q\sodot p) \\
& = &
s\scalar (q\sodot p).
\end{array}$$

\noindent Finally,
$$\begin{array}{rcl}
1 \sodot 1
& = &
m \after ((\kappa_{1}\after\,!)\otimes(\kappa_{1}\after\,!)) \\
& = &
m \after (\idmap\otimes\kappa_{1}) \after 
   (\kappa_{1}\otimes\idmap) \after (!\otimes\,!) \\
& = &
\rho \after (\kappa_{1}\otimes\idmap) \after (!\otimes\,!) \\
& = &
\kappa_{1} \after \rho \after (!\otimes\,!) \\
& = &
\kappa_{1} \after\, ! \\
& = &
1.
\end{array}$$
}
\end{myproof}

\begin{exas}
\label{StatePredTensorEx}
We briefly review the above pairing operation $\sodot$ in the
probabilistic and quantum case. We will see that in both cases
$\sodot$ is multiplication.
\begin{enumerate}
\item In the Kleisli category $\Kl(\Dst)$ of the distribution monad
$\Dst$ on $\Sets$ the pairing $\varphi\sodot\psi \in \Dst(X\times Y)$
of two states/distributions $\varphi\in\Dst(X)$ and $\psi\in\Dst(Y)$
is given by their product:
$$\begin{array}{rcccl}
(\varphi\sodot\psi)(x,y)
& = &
\dst(\varphi,\psi)(x,y)
& = &
\varphi(x) \cdot \psi(y).
\end{array}$$

\noindent where $\dst$ is the `double-strength' map from
Appendix~\ref{DiscProbSubsec}. The pairing $p\sodot q \in
[0,1]^{X\times Y}$ of predicates $p\in [0,1]^{X}$ and $q\in [0,1]^{Y}$
is obtained similarly as:
$$\begin{array}{rcl}
(p\sodot q)(x,y)
& = &
p(x)\cdot q(y).
\end{array}$$

\auxproof{
We have $\widehat{p} \colon X \rightarrow \Dst(2)$ given by
$\widehat{p}(x)(1) = p(x)$ and $\widehat{p}(x)(0) = 1-p(x)$, and
similarly for $\widehat{q} \colon Y \rightarrow \Dst(2)$. Thus
we obtain:
$$\xymatrix{
\widehat{p\sodot q} = \big(X\times Y\ar[r]^-{\widehat{p}\otimes\widehat{q}}
 & \Dst(2\times 2)\ar[r]^-{\Dst(m)} & \Dst(2)\big)
}$$

\noindent where $m\colon 2\times 2 \rightarrow 2$ is multiplication.
This means:
$$\begin{array}{rcl}
\widehat{p\sodot q}(x,y)(1)
& = &
\sum_{b,b'\in 2} \widehat{p}(x)(b)\cdot \widehat{q}(y)(b') \cdot 
   (\eta \after m)(b,b')(1) \\
& = &
\sum_{b,b'\in 2, b\cdot b' = 1} \widehat{p}(x)(b)\cdot \widehat{q}(y)(b') \\
& = &
\widehat{p}(x)(1) \cdot \widehat{q}(y)(1) \\
& = &
p(x) \cdot q(y).
\end{array}$$
}

\noindent Thus, indeed, as in
Lemma~\ref{StatePredTensorLem}~\eqref{StatePredTensorLemModularity},
$$\begin{array}{rcl}
\varphi \sodot \psi \models p \sodot q
& = &
\sum_{x,y} (\varphi\sodot\psi)(x,y) \cdot (p\sodot q)(x,y) \\
& = &
\sum_{x,y} \varphi(x) \cdot \psi(y) \cdot p(x) \cdot q(y) \\
& = &
\big(\sum_{x} \varphi(x) \cdot p(x)\big) \cdot 
    \big(\sum_{y} \psi(y) \cdot q(y)\big) \\
& = &
(\varphi \models p) \cdot (\psi \models q).
\end{array}$$

\item For two $C^*$-algebras $A,B$ with effects $e\in [0,1]_{A}, d\in
  [0,1]_{B}$ we have the pairing effect $e \sodot d = e\sotimes d \in
  A\otimes B$. It is obtained as follows. The effect $e\in [0,1]_{A}$
  corresponds to the positive unital map $f_{e} \colon \C^{2}
  \rightarrow A$ given by $f_{e}(\ket{0}) = e$ and $f_{e}(\ket{1}) = 1
  - e$, and similarly for $d$. The function $f_{e\sodot d} \colon
  \C^{2} \rightarrow A\otimes B$ is then the composite $f_{e\sodot d}
  = (f_{e}\otimes f_{d}) \after m$, with multiplication $m$
  from~\eqref{TwoMonoidExCstarEqn}. Thus:
$$\qquad\begin{array}{rcl}
e\sodot d
\hspace*{\arraycolsep} = \hspace*{\arraycolsep}
f_{e\sodot d}(\ket{0})
\hspace*{\arraycolsep} = \hspace*{\arraycolsep}
(f_{e}\otimes f_{d})(m(\ket{0})) 
& = &
(f_{e}\otimes f_{d})(\ket{0}\sotimes\ket{0})
   \qquad \mbox{by~\eqref{TwoMonoidExCstarEqn}} \\
& = &
f_{e}(\ket{0}) \sotimes f_{d}(\ket{0}) \\
& = &
e \sotimes d.
\end{array}$$

\auxproof{
\noindent Correspondingly, as one would expect:
$$\begin{array}{rcl}
\widehat{e\sodot d}(\ket{1})
& = &
(\widehat{e}\otimes\widehat{d})(m(\ket{1})) \\
& = &
(\widehat{e}\otimes\widehat{d})(\ket{0}\sotimes\ket{1} +
   \ket{1}\sotimes\ket{0} + \ket{1}\sotimes\ket{1})
   \qquad \mbox{by~\eqref{TwoMonoidExCstarEqn}} \\
& = &
\widehat{e}(\ket{0}) \sotimes \widehat{d}(\ket{1}) + 
   \widehat{e}(\ket{1}) \sotimes \widehat{d}(\ket{0}) +
   \widehat{e}(\ket{1}) \sotimes \widehat{d}(\ket{1}) \\
& = &
(1-e)\sotimes d + e\sotimes (1-d) + (1-e)\sotimes (1-d) \\
& = &
1\sotimes d - e\sotimes d + e\sotimes 1 - e\sotimes d + 1\sotimes 1
   - 1\sotimes d - e\sotimes 1 + e\sotimes d \\
& = &
1\sotimes 1 - e\sotimes d \\
& = &
1 - e\sodot d.
\end{array}$$
}
\end{enumerate}
\end{exas}

\noindent Point~\eqref{StatePredTensorLemProj} in Lemma~\ref{StatePredTensorLem}
says that there is a retraction $\Stat(X)\times \Stat(Y)
\rightarrowtail \Stat(X\otimes Y)$, since first pairing and then
projecting gives the original output.  In the category of sets the
retraction is an isomorphism, but in general this is not the
case. There is no such isomorphism because of both entanglement and
dependence in the quantum world, and because of dependence in the
probabilistic world.  We briefly illustrate this fundamental
phenomenon in the current setting.

\begin{exas}[Dependence and entanglement]
\label{DependenceEntanglementEx}
In a probabilistic setting as given by the category $\Kl(\Dst)$,
consider a state on the tensor product $X\times Y$ of two sets
$X,Y$. It is a discrete probability distribution
$\varphi\in\Dst(X\times Y)$ on $X$ and $Y$. Because it describes
probabilities on two sets, such a $\varphi$ is often called a
\emph{joint} distribution. One can take its two marginals
$(\pi_{1})_{*}(\varphi) = \Dst(\pi_{1})(\varphi) \in \Dst(X)$ and
$(\pi_{2})_{*}(\varphi) = \Dst(\pi_{2})(\varphi) \in \Dst(Y)$.  If we
now pair these marginals we may ask if the original joint distribution
re-appears, that is if:
$$\begin{array}{rcl}
(\pi_{1})_{*}(\varphi) \sodot (\pi_{2})_{*}(\varphi)
& \smash{\stackrel{?}{=}} &
\varphi.
\end{array}$$

\noindent But this equation precisely expresses that $\varphi$ is a
\emph{product} distribution (\textit{i.e.}~is factorisable). Hence the
equation does not always hold.

Similarly, in the quantum world, a state $\omega \colon A\otimes B
\rightarrow \C$ can, in general, not be reconstructed from its two
marginals $(\pi_{1})_{*}(\omega) = \omega \after \kappa_{1} \colon A
\rightarrow \C$ and $(\pi_{2})_{*}(\omega) = \omega \after \kappa_{2}
\colon B \rightarrow \C$, since in general we do not have
$(\pi_{1})_{*}(\omega) \sodot (\pi_{2})_{*}(\omega) = \omega$.

Here is concrete qubit example, for $A = B = \B(\C^{2})$. We consider
the non-entangled vector $\ket{u} = \ket{00} = \ket{0}\sotimes\ket{0}
\in \C^{2} \otimes \C^{2}$ and also the entangled (EPR) vector
$\ket{v} = \frac{1}{\sqrt{2}}(\ket{00} + \ket{11}) \in
\C^{2}\otimes\C^{2}$. There are associated states $\omega_{u},
\omega_{v} \colon \B(\C^{2}\otimes \C^{2}) \rightarrow \C$, given by
$\omega_{u}(f) = \bra{u}f\ket{u}$ and $\omega_{v}(f) =
\bra{v}f\ket{v}$. The marginals are given by
$(\pi_{1})_{*}(\omega_{u})(M) = \omega_{u}(M\otimes\idmap)$ and
$(\pi_{2})_{*}(\omega_{u})(N) = \omega_{u}(\idmap\otimes N)$, and
similarly for $v$. The (non-entangled) state $\omega_{u}$ can be
reconstructed via pairing from its marginals:
$$\begin{array}{rcl}
\big((\pi_{1})_{*}(\omega_{u}) \sodot (\pi_{2})_{*}(\omega_{u})\big)(M\sotimes N)
& = &
(\pi_{1})_{*}(\omega_{u})(M) \cdot (\pi_{2})_{*}(\omega_{u})(N) \\
& = &
\omega_{u}(M\otimes\idmap) \cdot \omega_{u}(\idmap\otimes N) \\
& = &
\bra{00}(M\otimes\idmap)\ket{00}\cdot\bra{00}(\idmap\otimes N)\ket{00} \\
& = &
\big(\bra{0}M\ket{0}\cdot \inprod{0}{0}\big) \cdot
  \big(\inprod{0}{0} \cdot \bra{0}N\ket{0}\big) \\
& = &
\bra{0}M\ket{0}\cdot\bra{0}N\ket{0} \\
& = &
\bra{00}(M\sotimes N)\ket{00} \\
& = &
\omega_{u}(M\sotimes N).
\end{array}$$

\noindent Because $\ket{v} = \frac{1}{\sqrt{2}}(\ket{00} + \ket{11})$
is an entangled state, we do not have such an equation for $v$. We
step through the computation to indicate where this fails:
$$\begin{array}{rcl}
\lefteqn{\big((\pi_{1})_{*}(\omega_{v}) \sodot 
   (\pi_{2})_{*}(\omega_{v})\big)(M\sotimes N)} \\
& = &
\omega_{v}(M\otimes\idmap) \cdot \omega_{v}(\idmap\otimes N) \\
& = &
\frac{1}{2}(\bra{00}+\bra{11})(M\otimes\idmap)(\ket{00}+\ket{11}) \cdot
  \frac{1}{2}(\bra{00}+\bra{11})(\idmap\otimes N)(\ket{00}+\ket{11}) \\
& = &
\frac{1}{4}(\bra{0}M\ket{0} + \bra{1}M\ket{1})\cdot 
   (\bra{0}N\ket{0} + \bra{1}N\ket{1})
   \qquad \mbox{since } \inprod{0}{1} = \inprod{1}{0} = 0 \\
& = &
\frac{1}{4}\big(\bra{0}M\ket{0}\cdot \bra{0}N\ket{0} + 
   \bra{0}M\ket{0}\cdot \bra{1}N\ket{1} +
   \bra{1}M\ket{1}\cdot \bra{0}N\ket{0} + 
   \bra{1}M\ket{1}\cdot \bra{1}N\ket{1}\big) \\
& \neq &
\frac{1}{2}\big(\bra{0}M\ket{0}\cdot \bra{0}N\ket{0} + 
   \bra{0}M\ket{1}\cdot \bra{0}N\ket{1} +
   \bra{1}M\ket{0}\cdot \bra{1}N\ket{0} + 
   \bra{1}M\ket{1}\cdot \bra{1}N\ket{1}\big) \\
& = &
\frac{1}{2}\big(\bra{00}(M\otimes N)\ket{00} + 
  \bra{00}(M\otimes N)\ket{11} + \bra{11}(M\otimes N)\ket{00} + 
  \bra{11}(M\otimes N)\ket{11}\big) \\
& = &
\frac{1}{2}(\bra{00} + \bra{11})(M\sotimes N)(\ket{00} + \ket{11}) \\
& = &
\omega_{v}(M\sotimes N).
\end{array}$$

\auxproof{
Then there is an isomorphism $\beta \colon A\otimes B \cong
\B(\C^{4})$ via Kronecker product:
$$\begin{array}{rcl}
\beta\left(\left(\begin{smallmatrix} a & b \\ c & d \end{smallmatrix}\right)
\sotimes N\right)
& = &
\left(\begin{smallmatrix} aN & bN \\ cN & dN \end{smallmatrix}\right)
\\
\beta^{-1}\left(\begin{smallmatrix} M_{00} & M_{01} \\ M_{10} & M_{11} 
   \end{smallmatrix}\right)
& = &
\left(\begin{smallmatrix} 1 & 0 \\ 0 & 0 \end{smallmatrix}\right)
   \sotimes M_{00} +
\left(\begin{smallmatrix} 0 & 1 \\ 0 & 0 \end{smallmatrix}\right)
   \sotimes M_{01} +
\left(\begin{smallmatrix} 0 & 0 \\ 1 & 0 \end{smallmatrix}\right)
   \sotimes M_{10} +
\left(\begin{smallmatrix} 0 & 0 \\ 0 & 1 \end{smallmatrix}\right)
   \sotimes M_{11}.
\end{array}$$
}
\end{exas}\medskip

\noindent In the end we mention that the distribution of the tensor $\otimes$
over coproducts $+$ is useful to interpret operations in context. For
instance, the guarded command $p?[f_{1}, \ldots, f_{n}]$
from~\eqref{GuardedIfEqn} can now be interpreted more generally: given
an $n$-test $p\colon X \rightarrow n\cdot 1$, and $n$ maps $f_{i}
\colon X\otimes Z \rightarrow Y$ with an additional context parameter
$Z$ we can interpret the test map $p?[f_{1}, \ldots, f_{n}]$ 
formally via weakening as $\pi_{1}^{*}(p)?[f_{1}, \ldots, f_{n}]$ in:
\begin{equation}
\label{GuardedIfTensorEqn}
\vcenter{\xymatrix@R-.5pc{
X\otimes Z\ar[rr]^-{\instr_{p}\otimes\idmap}\ar@/_2ex/[drr]_{\instr_{\pi_{1}^{*}(p)}} & &
   (n\cdot X)\otimes Z\ar[d]^-{\theta_1}_-{\cong} & & \\
& & n\cdot (X\otimes Z)\ar[rr]^-{[f_{1}, \ldots, f_{n}]} & & Y.
}}
\end{equation}

\noindent The reasoning rule from Lemma~\ref{TestMapWPLem} still
applies, with weakened tests $\pi_{1}^{*}(p)$. As special case we have
an equation, for maps $h_{i} \colon X \rightarrow Y$,
$$\begin{array}{rcl}
p?[h_{1},\ldots, h_{n}]\otimes\idmap
& = &
\pi_{1}^{*}(p)?[h_{1}\otimes\idmap, \ldots, h_{n}\otimes\idmap].
\end{array}$$

\auxproof{
$$\begin{array}{rcl}
p?[h_{1},\ldots, h_{n}]\otimes\idmap
& = &
([h_{1},\ldots, h_{n}]\otimes\idmap) \after (\instr_{p}\otimes\idmap) \\
& = &
([h_{1},\ldots, h_{n}]\otimes\idmap) \after \theta_{1}^{-1} \after
   \instr_{\pi_{1}^{*}(p)} \\
& = &
([h_{1},\ldots, h_{n}]\otimes\idmap) \after 
   [\kappa_{1}\otimes\idmap, \ldots, \kappa_{n}\otimes\idmap] \after
   \instr_{\pi_{1}^{*}(p)} \\
& = &
[h_{1}\otimes\idmap, \ldots, h_{n}\otimes\idmap] \after 
   \instr_{\pi_{1}^{*}(p)} \\
& = &
\pi_{1}^{*}(p)?[h_{1}\otimes\idmap, \ldots, h_{n}\otimes\idmap].
\end{array}$$
}

\section{Quantum states}\label{QuantumStateSec}

In this section we add our final assumption about the existence of
a ``quantum object''. We use it in a subsequent quantum protocol.

\begin{assumption}
\label{QuantumStateAss}
Let $\cat{B}$ be a category satisfying Assumption~\ref{TensorAss}. This
category contains a special object $Q\in\cat{B}$ with two states and
a predicate:
$$\xymatrix{
1\ar[r]^-{\upstate} & Q
\qquad
1\ar[r]^-{\downstate} & Q
\qquad
Q\ar[r]^-{\isup} & 1+1
}$$

\noindent such that the following diagram commutes,
\begin{equation}
\label{QuantumStateDiag}
\vcenter{\xymatrix@R-.5pc@C+1pc{
1+1\ar[r]^-{[\upstate,\downstate]}\ar@{=}@/_2ex/[dr] & 
   Q\ar[d]^{\isup}\ar@/^2ex/[dr]^{\instr_{\isup}} \\
& 1+1\ar[r]_-{\upstate+\downstate} & Q+Q
}}
\end{equation}

\noindent and the predicate $\isup$ is \emph{not} side-effect-free, or
equivalently, the map $\isup\colon Q \rightarrow 1+1$ is \emph{not}
an isomorphism.

\auxproof{
If $\isup$ is an isomorphism, say with $f \after \isup = \idmap$ and
$\isup \after f = \idmap$, then $f = f \after \idmap = f \after \isup
\after [\upstate, \downstate] = [\upstate,\downstate]$. Thus
$[\upstate, \downstate]$ is the inverse. We show that $\isup$ is
side-effect-free:
$$\begin{array}{rcccccl}
\nabla \after \instr_{\isup}
& = &
\nabla \after (\upstate+\downstate) \after \instr_{\isup}
& = &
[\upstate, \downstate] \after \instr_{\isup}
& = &
\idmap.
\end{array}$$

\noindent Conversely, if $\isup$ is side-effect-free, then:
$$\begin{array}{rcccccl}
[\upstate,\downstate] \after \isup
& = &
\nabla \after (\upstate+\downstate) \after \instr_{\isup}
& = &
\nabla \after \instr_{\isup}
& = &
\idmap.
\end{array}$$
}
\end{assumption}

In general, the object $2 = 1+1$ comes equipped with two states,
namely $1 = \kappa_{1} \colon 1 \rightarrow 2$ and $0 = \kappa_{2}
\colon 1 \rightarrow 2$. The identity map $2 \rightarrow 1+1$ can be
understood as a `isup' predicate. In this way we obtain classical
states, a bit like the object $Q$ in the above assumption. But the
object $Q$ is different, since it is explicitly required that there is
\emph{not} an isomorphism $2 \cong Q$.

The states $\upstate, \downstate$ correspond to the operations
$\newstate_{0}, \newstate_{1}$ in~\cite{Staton15a}, and the predicate
$\isup$ is called $\measure$ there. The above description focuses on
the states and measurement operations, in line with the rest of the
paper. The unitary operations for state transformations are clearly
missing. They do play an important role in quantum programming
language, see~\cite{Selinger04,Staton15a,ValironRSSS15}.

\begin{exa}
\label{CstarQuantumStateEx}
In the opposite categories $\op{(\CstarPU)}$ and $\op{(\CstarCPU)}$ of
$C^*$-algebras with (completely) positive unital maps we can take $Q =
\B(\C^{2}) = \Mat_{2}(\C)$. Let's write an element $M\in Q$ as a
$2\times 2$ matrix:
$$\begin{array}{rcl}
M
& = &
\left(\begin{smallmatrix} M_{00} & M_{01} \\ M_{10} & M_{11}
   \end{smallmatrix}\right)
\qquad\mbox{with matrix entries}\qquad M_{ij} \in \C.
\end{array}$$

\noindent Then we can define two states $\upstate,\downstate$ as
(completely) positive unital maps $Q \rightarrow \C$ via:
$$\begin{array}{rccclcrcccl}
\upstate(M)
& = &
M_{00}
& = &
\bra{0}M\ket{0}
\qquad\mbox{and}\qquad
\downstate(M)
& = &
M_{11}
& = &
\bra{1}M\ket{1}.
\end{array}$$

\noindent As predicate $\isup\colon \C^{2} \rightarrow Q$ we define:
$$\begin{array}{rcccl}
\isup(z,w)
& = &
z\ket{0}\bra{0} + w\ket{1}\bra{1}
& = &
\left(\begin{smallmatrix} z & 0 \\ 0 & w \end{smallmatrix}\right).
\end{array}$$

\noindent It corresponds to the effect $\isup(1,0) =
\left(\begin{smallmatrix} 1 & 0 \\ 0 & 0 \end{smallmatrix}\right) =
\ket{0}\bra{0}$. This effect is obviously a projection, with
orthocomplement $(\ket{0}\bra{0})^{\perp} = \ket{1}\bra{1}$. The
associated instrument $\instr_{\isup} \colon Q \oplus Q
\rightarrow Q$ is thus:
$$\begin{array}{rcccl}
\instr_{\isup}(M,N)
& = &
\ket{0}\bra{0}M\ket{0}\bra{0} + \ket{1}\bra{1}N\ket{1}\bra{1}
& = &
\left(\begin{smallmatrix} M_{00} & 0 \\ 0 & N_{11} \end{smallmatrix}\right).
\end{array}$$

\noindent It is easy to see that Diagram~\eqref{QuantumStateDiag}
commutes:
$$\begin{array}{rcl}
\big(\tuple{\upstate, \downstate} \after \isup\big)(z,w)
& = &
(\upstate\left(\begin{smallmatrix} z & 0 \\ 0 & w \end{smallmatrix}\right),
\downstate\left(\begin{smallmatrix} z & 0 \\ 0 & w \end{smallmatrix}\right)) 
\hspace*{\arraycolsep} = \hspace*{\arraycolsep}
(z,w) \\
\big(\isup \after (\upstate\oplus\downstate)\big)(M,N)
& = &
\isup(M_{00}, N_{11}) 
\hspace*{\arraycolsep} = \hspace*{\arraycolsep}
\left(\begin{smallmatrix} M_{00} & 0 \\ 0 & N_{11} \end{smallmatrix}\right) 
\hspace*{\arraycolsep} = \hspace*{\arraycolsep}
\instr_{\isup}(M,N).
\end{array}$$
\end{exa}

\newpage
\begin{rem}
\label{QuantumStateRem}\hfill
\begin{enumerate}
\item In the opposite $\op{\Rng}$ of the category of rings --- not
  necessarily commutative --- we can use the object $Q = \Mat_{2}(\Z)$
  of $2\times 2$ matrices over the integers, with similar states and
  predicate. (In this ring example one does have the structure
  described Assumption~\ref{QuantumStateAss}, but
  Assumption~\ref{TensorAss} is not satisfied.)

\item Via the isomorphisms $\Stat(\B(\C^{2})) \cong \DM(\C^{2})$
  from~\eqref{HilbStatEqn} the state $\upstate$ in
  Example~\ref{CstarQuantumStateEx} corresponds to the density matrix
  $\ket{0}\bra{0} \in \DM(\C^{2})$, since for $M\in\B(\C^{2})$,
$$\begin{array}{rcccl}
\tr(\ket{0}\bra{0}M)
& = &
\bra{0}M\ket{0}
& = &
\upstate(M).
\end{array}$$
\end{enumerate}
\end{rem}\medskip

\noindent If we have a program $f\colon Q\otimes X \rightarrow Y$ we can insert
a state like $\upstate$ or $\downstate$ in the first
component. Category-theoretically this involves taking the composite:
$$\xymatrix@C+1pc{
X\cong 1\otimes X\ar[r]^-{\upstate\otimes\idmap} & Q\otimes X\ar[r]^-{f} & Y
}$$

\noindent In the language of~\cite{Staton15a} this operation would be
written as $\newstate(v.f(v,-))$. More informally we shall write it
below as: $\textsf{let } v = \upstate \textsf{ in } f(v\sotimes
-)$. We then have:
$$\begin{array}{rcl}
\big(\omega\models(\textsf{let }
v = \upstate \textsf{ in } f(v\sotimes -))^{*}(q)\big)
& = &
\big(\upstate\sodot\omega\models f^{*}(q)\big).
\end{array}$$\medskip

\noindent We conclude with a description of a familiar protocol in the
setting developed in this article.

\begin{exa}[Superdense coding]
\label{SuperDenseEx}
In the superdense coding protocol Alice sends two classical bits to
Bob by transferring her part of a shared, entangled quantum state.  In
a category with a quantum object $Q$ as in
Assumption~\ref{QuantumStateAss} this protocol can be described as a
map $\textsl{sdc} \colon 4 \rightarrow 4$ consisting of three
consecutive steps:
\begin{equation}
\label{SuperDenseShort}
\xymatrix@C+1pc{
\textsl{sdc} = \big(4\ar[r]^-{\textsl{init}} & 
   4 \otimes Q\otimes Q\ar[r]^-{t_A\otimes\idmap} &
   Q\otimes Q\ar[r]^-{t_B} & 4\big)
}
\end{equation}

\noindent The correctness of the protocol means that this composite is
the identity map.

We shall describe this \textsl{sdc} map~\eqref{SuperDenseShort} in
greater detail in the category $\op{(\CstarCPU)}$, with $Q =
\B(\C^{2})$ as in Example~\ref{CstarQuantumStateEx}. In pseudo-code
this protocol $\textsl{sdc} \colon \C^{4} \rightarrow \C^{4}$ is
described in Figure~\ref{SuperDenseFig}. It is not our aim to develop
the syntax of this code fragment in detail. Instead we explain its
interpretation in the category $\op{(\CstarCPU)}$. For the two test
operations we use the meaning described in~\eqref{TestMapEqn}.

We start by recalling some basic material. The Bell basis of
$\C^{4}$ is given by the four vectors:
$$\begin{array}{rclcrcl}
\ket{b_1} &=& \frac{1}{\sqrt{2}} (\ket{00} + \ket{11})
& \qquad &
\ket{b_2} &=& \frac{1}{\sqrt{2}} (\ket{01} + \ket{10}) \\
\ket{b_3} &=& \frac{1}{\sqrt{2}} (\ket{00} - \ket{11})
& &
\ket{b_4} &=& \frac{1}{\sqrt{2}} (\ket{01} - \ket{10})
\end{array}$$

\noindent The associated projections $e_{i} = \ket{b_{i}}\bra{b_{i}}
\in \B(\C^{4}) \cong \B(\C^{2})\otimes \B(\C^{2}) = Q\otimes Q$ are
described explicitly by the four matrices:
$$\begin{array}{rclcrclcrclcrcl}
e_1 &=& \frac{1}{2}\left(\begin{smallmatrix}
1 & 0 & 0 & 1 \\
0 & 0 & 0 & 0 \\
0 & 0 & 0 & 0 \\
1 & 0 & 0 & 1 \\
\end{smallmatrix}\right)
& \; &
e_2 &=& \frac{1}{2}\left(\begin{smallmatrix}
0 & 0 & 0 & 0 \\
0 & 1 & 1 & 0 \\
0 & 1 & 1 & 0 \\
0 & 0 & 0 & 0 \\
\end{smallmatrix}\right)
& \; &
e_3 &=& \frac{1}{2}\left(\begin{smallmatrix}
1 & 0 & 0 & -1 \\
0 & 0 & 0 & 0 \\
0 & 0 & 0 & 0 \\
-1 & 0 & 0 & 1 \\
\end{smallmatrix}\right)
& \; &
e_4 &=& \frac{1}{2}\left(\begin{smallmatrix}
0 & 0 & 0 & 0 \\
0 & 1 & -1 & 0 \\
0 & -1 & 1 & 0 \\
0 & 0 & 0 & 0 \\
\end{smallmatrix}\right)
\end{array}$$

\noindent They satisfy $e_{1} \ovee e_{2} \ovee e_{3} \ovee e_{4} =
\idmap$ and thus form a $4$-test $\C^{4} \rightarrow Q\otimes Q$, as
used in Bob's test in Figure~\ref{SuperDenseFig}. Since $e_{i}^{2} =
e_{i}$ we have $\sqrt{e_i} = e_{i}$. The associated instrument
$(Q\otimes Q)^{4} \rightarrow Q\otimes Q$ thus sends $(x_{1}, x_{2},
x_{3}, x_{4})$ to $\sum_{i}e_{i}\cdot x_{i} \cdot e_{i}$.

The first element $\ket{b_1} = \frac{1}{\sqrt{2}} (\ket{00} +
\ket{11}) \in Q\otimes Q$ of the Bell basis can be obtained from the
base vector $\ket{00} = \upstate \sotimes\upstate \in Q\otimes Q$ by
first applying $H\otimes\idmap$ --- where $H$ is Hadamard --- and then
conditional negation $\textsf{CNOT}$. This is done in the
initialisation fragment in Figure~\ref{SuperDenseFig}.  The resulting
initialisation map from~\eqref{SuperDenseShort} has type
$\textsf{init} \colon \C^{4}\otimes Q\otimes Q \rightarrow \C^{4}$ in
$\CstarCPU$ and is given by $\textsf{init}(\vec{z}\sotimes x) =
\bra{b_1}x\ket{b_1}\cdot \vec{z}$.

\auxproof{
$$\begin{array}{rcl}
\textsf{CNOT}(H\otimes\idmap)\ket{00}
& = &
\textsf{CNOT}(H\ket{0}\sotimes\ket{0}) \\
& = &
\frac{1}{\sqrt{2}}\textsf{CNOT}\left(
  \left(\begin{smallmatrix} 1 & 1 \\ 1 & -1 \end{smallmatrix}\right)
  \left(\begin{smallmatrix} 1 \\ 0 \end{smallmatrix}\right)
  \sotimes \left(\begin{smallmatrix} 1 \\ 0 \end{smallmatrix}\right)\right) \\
& = &
\frac{1}{\sqrt{2}}\textsf{CNOT}\left(
  \left(\begin{smallmatrix} 1 \\ 1 \end{smallmatrix}\right)
  \sotimes \left(\begin{smallmatrix} 1 \\ 0 \end{smallmatrix}\right)\right) \\
& = &
\frac{1}{\sqrt{2}}\left(\begin{smallmatrix}
1 & 0 & 0 & 0 \\
0 & 1 & 0 & 0 \\
0 & 0 & 0 & 1 \\
0 & 0 & 1 & 0 \\
\end{smallmatrix}\right)
\left(\begin{smallmatrix} 1 \\ 0 \\ 1 \\ 0 \end{smallmatrix}\right) \\
& = &
\frac{1}{\sqrt{2}}\left(\begin{smallmatrix}
   1 \\ 0 \\ 0 \\ 1 \end{smallmatrix}\right) \\
& = &
\frac{1}{\sqrt{2}}(\ket{00} + \ket{11}).
\end{array}$$
}


For Alice's test we need the four Pauli matrices in $\B(\C^{2})$:
$$\begin{array}{rclcrccclcrccclcrcccl}
\sigma_{1} 
& = & 
\left(\begin{smallmatrix} 1 & 0 \\ 0 & 1 \end{smallmatrix}\right)
& \quad &
\sigma_2 
& = &
X
& = &
\left(\begin{smallmatrix} 0 & 1 \\ 1 & 0 \end{smallmatrix}\right)
& \quad &
\sigma_3 
& = &
Z
& = &
\left(\begin{smallmatrix} 1 & 0 \\ 0 & -1 \end{smallmatrix}\right)
& \quad &
\sigma_4 
& = &
XZ
& = &
\left(\begin{smallmatrix} 0 & -1 \\ 1 & 0 \end{smallmatrix}\right).
\end{array}$$

\noindent When we apply $\sigma_{i}\otimes\idmap$ to the Bell basis
vectors $\ket{b_j}$ we see that:
\begin{equation}
\label{SuperDenseExchange}
\left\{{\renewcommand\arraystretch{1}\begin{array}{ll}
\sigma_{1}\otimes\idmap = \idmap \quad & \mbox{changes nothing} \\
\sigma_{2}\otimes\idmap & \mbox{exchanges $\ket{b_1}, \ket{b_2}$ and
   $\ket{b_3}, \ket{b_4}$} \\
\sigma_{3}\otimes\idmap & \mbox{exchanges $\ket{b_1}, \ket{b_3}$ and
   $\ket{b_2}, \ket{b_4}$} \\
\sigma_{4}\otimes\idmap & \mbox{exchanges $\ket{b_1}, \ket{b_4}$ and
   $\ket{b_2}, \ket{b_3}$.}
\end{array}}\right.
\end{equation}

\noindent Alice's test map $t_{A} \colon Q\otimes Q \rightarrow
\C^{4}\otimes Q\otimes Q$ in $\CstarCPU$ is given by $t_{A}(x) =
\sum_{i} \ket{i}\sotimes \B(\sigma_{i}\otimes\idmap)(x)$, where
$\ket{i}\in \C^{4}$ is a standard base vector.

Finally, the test operation performed by Bob in
Figure~\ref{SuperDenseFig} is the map $t_{B} \colon \C^{4} \rightarrow
Q\otimes Q$ given by $t_{B}(\vec{z}) = \sum_{i} z_{i}\cdot e_{i}$. Now
we are ready to calculate the composite~\eqref{SuperDenseShort}, in
reverse order:
$$\begin{array}{rcl}
\textsl{sdc}(z_{1}, z_{2}, z_{3}, z_{4})
& = &
\big(\textsf{init} \after t_{A} \after t_{B}\big)(z_{1}, z_{2}, z_{3}, z_{4}) \\
& = &
\big(\textsf{init} \after t_{A}\big)\big(\sum_{i} z_{i}\cdot e_{i}\big) \\
& = &
\textsf{init}\big(\sum_{j}\ket{j}\sotimes\B(\sigma_{j}\otimes\idmap)
   (\sum_{i} z_{i}\cdot e_{i})\big) \\
& = &
\textsf{init}\big(\sum_{i,j}z_{i}\cdot\ket{j}\sotimes
   (\sigma_{j}\otimes\idmap)^{\dag}\ket{b_i}\bra{b_i}
      (\sigma_{j}\otimes\idmap)\big) \\
& = &
\sum_{i,j}z_{i}\cdot \bra{b_1}(\sigma_{j}\otimes\idmap)^{\dag}\ket{b_i}\bra{b_i}
      (\sigma_{j}\otimes\idmap)\ket{b_1}\ket{j} \\
& = &
\sum_{j} z_{j}\cdot \ket{j} \qquad \mbox{by~\eqref{SuperDenseExchange}} \\
& = &
(z_{1}, z_{2}, z_{3}, z_{4}).
\end{array}$$

\noindent It remains a challenge to prove the correctness of protocols
like this via the dynamic logic that is sketched in this paper. Such a
challenge requires a further syntactical and logical development of
the logic, see~\cite{Adams14}, which is beyond the scope of the
current paper.
\end{exa}

\begin{figure}
\begin{center}
\fbox{$\begin{array}{rcl}
\textsl{sdc}(z_{1}, z_{2}, z_{3}, z_{4})
& = &
{\renewcommand\arraystretch{1}\begin{array}[t]{l}
\textsf{// initialisation steps} \\
\textsf{let }v_{1} = \upstate \textsf{ in} \\
\quad\textsf{let }v_{2} = \upstate \textsf{ in} \\
\quad\quad\textsf{let }b_{1} = 
   \textsf{CNOT}(H\otimes\idmap)(v_{1}\sotimes v_{2}) \textsf{ in} \\
\textsf{// Alice's test operation} \\
\textsf{let } t_{A} = \begin{array}[t]{l}
\textsf{begin test } (z_{1}, z_{2}, z_{3}, z_{4}) \\
\quad | \; 1 \longrightarrow b_{1} \\
\quad | \; 2 \longrightarrow (X\otimes\idmap)(b_{1}) \\
\quad | \; 3 \longrightarrow (Z\otimes\idmap)(b_{1}) \\
\quad | \; 4 \longrightarrow (XZ\otimes\idmap)(b_{1}) \\
\textsf{end test } \quad \textsf{in} 
\end{array} \\
\textsf{// Bob's test operation} \\
\textsf{let } t_{B} = \begin{array}[t]{l}
\textsf{begin test } t_{A} \\
\quad | \; e_{1} \longrightarrow 1 \\
\quad | \; e_{2} \longrightarrow 2 \\
\quad | \; e_{3} \longrightarrow 3 \\
\quad | \; e_{4} \longrightarrow 4 \\
\textsf{end test } \quad \textsf{in} 
\end{array} \\
t_{B}
\end{array}}
\end{array}\vspace*{-1em}$}
\end{center}
\caption{Pseudo code for the superdense coding function \textsl{sdc}
   from Example~\ref{SuperDenseEx}.}
\label{SuperDenseFig}
\end{figure}

\auxproof{

\section{Coherence space models}\label{CohSpSec}


We recall that a coherence space $A$ consists of a pair $A = (|A|,
\coh_{A})$, where $|A|$ is the underlying set, called the web of $A$,
and $\coh_{A}\, \subseteq |A|\times|A|$ is a reflexive and symmetric
relation, called coherence. The subscript $A$ in $\coh_A$ will be
omitted whenever it is clear from the context. A clique of $A$ is
subset $U \subseteq |A|$ satisfying $a\coh a'$ for all $a,a'\in U$.
We write $\Cliq(A) \subseteq \Pow(|A|)$ for the set of cliques of $A$,
which we order by inclusion $\subseteq$. 

One writes $a\incoh a'
$ for $a=a' \vee (a \not\coh a')$. We form a new
coherence space $A^{\perp} = (|A|, \incoh)$, so that $|A^{\perp}| =
|A|$. Notice that $A^{\perp\perp} = A$. A clique of $A^{\perp}$ may be 
called an anti-clique. 

\auxproof{
Clearly, $\incoh$ is reflexive. But it is also symmetric, since if 
$a\incoh a'$, then either:
\begin{itemize}
\item $a=a'$, but then $a'=a$ so that $a'\incoh a$

\item $a \not\coh a'$, but then also $a' \not\coh a$ since if $a'\coh a$
then $a\coh a'$ by symmetry of $\coh$, which is not the case.
\end{itemize}

Taking the incoherence of $\incoh$ is the original relation $\coh$ since:
$$\begin{array}{rcccccl}
a = a' \vee (a \not\incoh a')
& \Longleftrightarrow &
a = a' \vee (a \neq a' \wedge a \coh a') 
& \Longleftrightarrow &
a = a' \vee (a \coh a') 
& \Longleftrightarrow &
a \coh a'.
\end{array}$$
}

We list some basic properties
of $\Cliq(A)$, see \textit{e.g.}~\cite{Troelstra92}.
\begin{enumerate}
\item $a\coh a'$ iff $\{a,a'\} \in \Cliq(A)$.

\item $\emptyset\in\Cliq(A)$, and $\{a\}\in\Cliq(A)$ for each
  $a\in|A|$.

\item If $V \subseteq U \in \Cliq(A)$, then $V\in\Cliq(A)$, so that 
$\Cliq(X)$ is a downset. 

\item for a directed collection $U_{i} \in \Cliq(A)$, for $i\in I$, we
  have $\bigcup_{i\in I}U_{i} \in \Cliq(A)$.

\item If for an $I$-indexed collection $U_{i}\in\Cliq(A)$ we have
  $U_{i_1}\cup U_{i_2} \in \Cliq(A)$, for each pair $i_{1},i_{2}\in
  I$, then $\bigcup_{i\in I}U_{i} \in \Cliq(A)$.

\item Each $U\in\Cliq(A)$ is the directed join of the finite cliques
  contained in $U$.

\item For each set $X$, the pair $(X,=)$ forms a coherence space, with
  subsingletons as cliques: $\Cliq(X,=) = \{\emptyset\} \cup
  \set{\{x\}}{x\in X}$. We often simply write $X$ for $(X,=)$. In
  particular, we write $0$ for the coherence space $(\emptyset,=)$ and
  $1$ for $(\{0\}, =)$. Notice that $\Cliq(0) = \{\emptyset\}$ and
  $\Cliq(1) = \{\emptyset, \{0\}\}$.
\end{enumerate}

\auxproof{
\begin{enumerate}
\item Obvious.

\item Clearly, if $y,y' \in \{x\}$, then $y = x = y'$, so $y\coh y$
since $\coh$ is reflexive.

\item Obvious.

\item If $x_{1},x_{2}\in\bigcup_{i}C_{i}$, say $x_{1} \in C_{i_1}$, $x_{2}
\in C_{i_2}$, then, by directedness, there is an $i$ with $C_{i_1},
C_{i_2} \subseteq C_i$. Hence $x_{1}, x_{2} \in C_{i}$, and thus
$x_{1} \coh x_{2}$.

\item Take $J = \Powfin(I)$, with for $j\in J$, $D_{j} = \bigcup_{i\in j}C_{i}
\in \Cliq(X)$. These $D_j$ form a directed collection, so $\bigcup_{i\in I}
C_{i} = \bigcup_{j\in J}D_{j} \in \Cliq(X)$.

\end{enumerate}
}

For subsets $U,V \subseteq |A|$ one writes $U \orthogonal V$ if the
intersection $U\cap V$ is a subsingleton \textit{i.e.}~contains at
most one element. Formulated in terms of the associated characteristic
functions $\charac_{U}, \charac_{V} \colon |A| \rightarrow \{0,1\}$
this means $\sum_{a\in |A|} \charac_{U}(a)\cdot \charac_{V}(a) \leq
1$. Clearly, $\orthogonal$ is a symmetric relation.  For an arbitrary
subset $S \subseteq \Pow(|A|)$ one defines $\simop S \subseteq
\Pow(|A|)$ as:
$$\begin{array}{rcl}
\simop S
& = &
\set{V\subseteq |A|}{\allin{U}{S}{U\orthogonal V}}.
\end{array}$$

\noindent It is not hard to see that this operation $\simop{(-)}$
satisfies:
\begin{enumerate}
\item $\simop S$ is downclosed;

\item $S \subseteq T$ implies $\simop T \subseteq \simop S$;

\item $S \subseteq \simop\simop S$, and thus $\simop S =
  \simop\simop\simop S$;

\item $\Cliq(A^{\perp}) \subseteq \simop\Cliq(A)$;

\item $\simop\simop\Cliq(A) = \Cliq(A)$.
\end{enumerate}

\begin{enumerate}
\item Suppose $U' \subseteq U \in \simop S$, and let $V\in S$. We need
to prove that $U'\cap V$ is a subsingleton. So assume $a,a'\in U'\cap V$.
Then also $a,a' \in U\cap V$, so that $a=a'$, as required.

\item Assume $V\in \simop T$. This means that for each $U\in T$ we
  have $U\orthogonal V$. In particular, $U\orthogonal V$ for all $U\in
  S\subseteq T$. Hence $V\in \simop S$.

\item Let $U\in S$ and $V\in \simop S$. We have to prove $V\orthogonal
  U$. But we have $U\orthogonal V$, so we are done.

\item Let $U\in\Cliq(A)$ and $V\in\Cliq(A^{\perp})$. In order to show
  that $U\cap V$ is a subsingleton, assume $a,a'\in U\cap V$. Then
  $a,a'\in U$ so $a\coh a'$. But also $a,a'\in V$, so $a\incoh a' = (a
  = a' \vee a\not\coh a')$. This gives $a=a'$.

\item Assume $U \in \simop\simop\Cliq(A)$. In order to prove
  $U\in\Cliq(A)$, we assume $a,a'\in U$, and need to show $a\coh a'$.
  By the previous point we have $U\in\simop\Cliq(A^{\perp})$, and thus
  $\{a,a'\} \in \Cliq(A^{\perp})$, by the first point. If $a \not\coh
  a'$, then $a \incoh a'$ and thus $\{a,a'\} \in
  \Cliq(A^{\perp})$. But then $\{a,a'\} \cap \{a,a'\} = \{a,a'\}$ is a
  subsingleton, giving $a=a'$, and thus $a\coh a'$, which contradicts
  the assumption $a\not\coh a'$. Hence $a\coh a'$.
\end{enumerate}

The main observation is that there is an equivalence between:
$$\begin{prooftree}
\mbox{reflexive and symmetric relations } \coh\; \subseteq |A|\times|A| 
\Justifies
\mbox{subsets $S\subseteq \Pow(|A|)$ with $\simop\simop S = S$}
\end{prooftree}$$

\noindent In the direction downwards we send $\coh$ to
$\Cliq(A)$. Upwards we define $a\coh a'$ iff $\{a,a'\}\in S$. As shown
above, starting from $\coh$ we obtain the original relation, since:
$\{a,a'\} \in \Cliq(A)$ iff $a\coh a'$. Starting from $S$ we obtain
$a\coh a'$ iff $\{a,a'\} \in S$, so that the associated set of cliques
satisfies: $\Cliq(A) = S$.
\begin{itemize}
\item $(\subseteq)$ Let $U\in\Cliq(A)$; we prove $U\in \simop\simop S
  = S$. So let $V\in \simop S$ and $a,a'\in U\cap V$; we must prove
  $a=a'$. From $a,a'\in U$ we get $a\coh a'$ and thus $\{a,a'\} \in
  S$. But then $\{a,a'\} \cap V$ is a subsingleton. This means $a=a'$
  since $a,a'\in V$.

\item $(\supseteq)$ Assume $U \in S$ and $a,a'\in U$. We obtain
  $U\in\Cliq(A)$ if $\{a,a'\} \in S$. The latter follows if we prove
  $\{a,a'\} \in \simop\simop S = S$. So let $V\in \simop S$. We have
  to prove $\{a,a'\} \cap V$ is a subsingleton. But this follows since
  $\{a,a'\} \cap V \subseteq U \cap V$, and the latter is a
  subsingleton.
\end{itemize}

A \emph{linear} morphism $f\colon A \rightarrow B$ between two
coherence spaces $A,B$ is a monotone function $f \colon \Cliq(A)
\rightarrow \Cliq(B)$ between the corresponding sets of cliques,
satisfying: for each $U\in\Cliq(A)$ and $b \in f(U)$ there is an $a\in
U$ with $b\in f(\{a\})$. Such a function $f$ satisfies:
\begin{enumerate}
\item $f(\emptyset) = \emptyset$;
\item $f$ preserves directed joins;
\item $f(U\cup V) = f(U) \cup f(V)$, if $U,V$ and $U\cup V$ are cliques;
\end{enumerate}

{
\begin{enumerate}
\item Assume $b\in f(\emptyset)$; then there is an $a\in\emptyset$
  with $b\in f(\{a\})$. This is clearly impossible.

\item If $(U_{i})$ is a directed collection, then $\bigvee_{i}
  f(U_{i}) \subseteq f(\bigvee_{i}U_{i})$ by monotonicity. For the
  reverse direction, assume $b\in f(\bigvee_{i}U_{i})$. Then there is
  an $a \in \bigvee_{i}U_{i}$ with $b\in f(\{a\})$. Hence there is an
  $j$ with $a\in U_{j}$, and thus $b\in f(\{a\}) \subseteq f(U_{j})$.
  Thus $b\in \bigvee_{i} f(U_{i})$.

\item We have $(\supseteq)$ by monotonicity of $f$.  Hence assume
  $b\in f(U\cup V)$. Then there must be an $a\in U\cup V$ with $b\in
  f(\{a\})$. If $a\in U$, then $b \in f(\{a\}) \subseteq U$, and
  similarly, if $a\in V$, then $b \in f(\{a\}) \subseteq V$.
\end{enumerate}
}

\noindent We define a category $\CohL$ of coherence spaces with linear
maps. It is easy to see that the mapping $X \mapsto (X,=)$ yields a
functor $\Sets \rightarrow \CohL$. In fact, there is a full and
faithful functor from sets with partial functions to $\CohL$. Notice
that the object $0\in\CohL$ is both initial and final (and hence a
zero object).

\auxproof{
Composition of linear maps $A \stackrel{f}{\rightarrow} B
\stackrel{g}{\rightarrow} C$. Assume $U\in\Cliq(A)$ and $c\in
g(f(U))$.  Then there is a $b\in f(U)$ with $c \in g(\{b\})$.  But
then there is also an $a\in U$ with $b \in f(\{a\})$.  Since $\{b\}
\subseteq f(\{a\})$ and $g$ is monotone, we get $c \in g(\{b\})
\subseteq g(f(\{a\}))$, as required.

For a function $f\colon X \rightarrow Y$ we define $D(f) \colon
\Cliq(X,=) \rightarrow \Cliq(Y,=)$ by $D(f)(\emptyset) = \emptyset$
and $D(f)(\{x\}) = \{f(x)\}$. This $D(f)$ is clearly monotone. Also,
for each $U\in\Cliq(X,=)$ and $y\in D(f)(U)$, we have $U = \{x\}$ for
some $x\in X$, and thus $y \in D(f)(\{x\}) = \{f(x)\}$ means $y =
f(x)$. Hence there is indeed an $x\in U$ with $y\in D(f)(\{x\})$.

For an arbitrary coherence space $A$, there is clearly precisely one
function $f\colon \Cliq(A) \rightarrow \Cliq(1) = \{\emptyset\}$.  But
there is also precisely one morphism $0 \rightarrow A$, since any
linear $f\colon \Cliq(0) \rightarrow \Cliq(A)$ since
$f(\emptyset) = \emptyset$.
}

For two coherence space $A_{1},A_{2}$ we define their product
$A_{1}\times A_{2}$ coproduct $A_{1}+A_{2}$. The product is usually
written as $A_{1} \mathrel{\&} A_{2}$ and called `with', and the
coproduct is written as $A_{1}\oplus A_{2}$ and called `plus' in the
linear logic literature. Both product and coproduct have as underlying
web the disjoint union of webs: $|A_{1}\times A_{2}| = |A_{1}+A_{2}| =
|A_{1}|+|A_{2}|$, but the coherence relations differ. For product
one takes:
$$\begin{array}{rclcrclcc}
\kappa_{1}a \coh_{A_{1}\times A_{2}} \kappa_{1}a'
& \mbox{iff} &
a \coh_{A_1} a'
& \qquad &
\kappa_{2}a \coh_{A_{1}\times A_{2}} \kappa_{2}a'
& \mbox{iff} &
a \coh_{A_2} a'
& \qquad &
\kappa_{1}a \coh_{A_{1}\times A_{2}} \kappa_{2}a'.
\end{array}$$

\noindent And for coproduct only the third part differs:
$$\begin{array}{rclcrclcc}
\kappa_{1}a \coh_{A_{1} + A_{2}} \kappa_{1}a'
& \mbox{iff} &
a \coh_{A_1} a'
& \qquad &
\kappa_{2}a \coh_{A_{1}+ A_{2}} \kappa_{2}a'
& \mbox{iff} &
a \coh_{A_2} a'
& \qquad &
\kappa_{1}a \not\coh_{A_{1}+ A_{2}} \kappa_{2}a'.
\end{array}$$

\noindent The mapping $\Cliq(A_{1}) \times \Cliq(A_{2}) \rightarrow
\Cliq(A_{1}\times A_{2})$ given by $(U,V) \mapsto \kappa_{1}U \cup
\kappa_{2}V$ is then a bijection.

\auxproof{
We first check that $\kappa_{1}U \cup \kappa_{2}V$ is a clique. If
$x,x'\in\kappa_{1}U \cup \kappa_{2}V$ then $x \coh x'$ if
$x\in\kappa_{1}U,x'\in\kappa_{2}V$ or the other way around. If $x,x'$
come from the same component, then either
$x=\kappa_{1}a,x'=\kappa_{1}a'$, or $x=\kappa_{2}a,x'=\kappa_{2}a'$.
In the first case we have $a,a'\in U$, so $a\coh a'$ and thus
$x=\kappa_{1}a \coh \kappa_{1}a' = x'$. And in the second case we have
$a,a'\in V$, so $a\coh a'$ and thus also $x=\kappa_{2}a \coh
\kappa_{2}a' = x'$.
}

There are coprojection maps $\kappa_{i} \colon A_{i} \rightarrow
A_{1}+A_{2}$ in $\CohL$ given by $\kappa_{i}(U) =
\set{\kappa_{i}a}{a\in U}$, for $U\in\Cliq(A_{i})$. Clearly,
$\kappa_{i}(U)$ is again a clique. And if $x\in\kappa_{i}(U)$, then $x
= \kappa_{i}a$ for some $a\in U$, so that $x\in
\kappa_{i}(\{a\})$. Hence $\kappa_{i}$ is linear. If we have linear
maps $f_{i} \colon A_{i} \rightarrow B$, then we can define a function
$[f_{1}, f_{2}] \colon \Cliq(A_{1}+A_{2}) \rightarrow \Cliq(B)$ by:
$$\begin{array}{rcl}
[f_{1},f_{2}](W)
& = &
f_{1}\big(\kappa_{1}^{-1}(W)\big) \cup f_{2}\big(\kappa_{2}^{-1}(W)\big).
\end{array}$$

\auxproof{
Clearly, $\kappa_{1}^{-1}(W)$ is a clique: if $a,a'\in
\kappa_{1}^{-1}(W)$, then $\kappa_{1}a,\kappa_{1}a'\in W$, so
$\kappa_{1}a \coh \kappa_{1}a'$, and thus $a\coh a'$. Also, this
cotuple is linear: if $b \in [f_{1},f_{2}](W)$, then either $b \in
f_{1}\big(\kappa_{1}^{-1}(W)\big)$ or $b\in
f_{2}\big(\kappa_{2}^{-1}(W)\big)$. Without loss of generality we
handle the first case. There is then an $a\in \kappa_{1}^{-1}(W)$ with
$b\in f_{1}(\{a\})$. But this means $\kappa_{1}a \in W$ and thus $b
\in [f_{1}, f_{2}](\{\kappa_{1}a\}) = f_{1}(\{a\})$ since
$\kappa_{1}^{-1}(\{\kappa_{1}a\}) = \{a\}$ and
$\kappa_{2}^{-1}(\{\kappa_{1}a\}) = \emptyset$.
}

\noindent Then $[f_{1}, f_{2}] \after \kappa_{i} = f_{i}$, since
$\kappa_{j}^{-1}(\kappa_{i}U)$ is $U$ if $i=j$ and empty otherwise.
Here we use that linear maps preserve the emptyset. If also $g\colon
A_{1}+A_{2} \rightarrow B$ in $\CohL$ satisfies $g \after \kappa_{i} =
f_{i}$, then $g = [f_{1},f_{2}]$ because linear maps preserve $\cup$,
as described above.

\auxproof{
$$\begin{array}{rcl}
[f_{1},f_{2}](W)
& = &
f_{1}\big(\kappa_{1}^{-1}(W)\big) \cup f_{2}\big(\kappa_{2}^{-1}(W)\big) \\
& = &
g\big(\kappa_{1}\big(\kappa_{1}^{-1}(W)\big)\big) \cup 
   g\big(\kappa_{2}\big(\kappa_{2}^{-1}(W)\big)\big) \\
& = &
g\big(\set{\kappa_{1}a}{\kappa_{1}a\in W}\big) \cup
   g\big(\set{\kappa_{2}a}{\kappa_{2}a\in W}\big) \\
& = &
g\big(\set{\kappa_{1}a}{\kappa_{1}a\in W} \cup
   \set{\kappa_{2}a}{\kappa_{2}a\in W}\big) \\
& = &
g(W).
\end{array}$$
}

Notice that the coproduct $1+1$ has $2 = \{0,1\}$ as web, with
discrete coherence relation $i\coh j$ iff $i=j$. Thus there is a
3-element set of cliques $\Cliq(1+1) = \{\emptyset, \{0\}, \{1\}\}$.
Any linear map $f \colon A \rightarrow 1+1$

\bigskip

A \emph{probabilistic coherence space} (PCS) $A$ is given by a finite
web $|A|$ and a subset $P_{A} \subseteq \R_{\geq 0}^{|A|}$ of
functions satisfying $P_{A} = \simop\simop P_{A}$. First, two
functions $f,g\in \R_{\geq 0}^{|A|}$ are called polar, written as $f
\polar g$, if $\sum_{a\in|A|} f(a) \cdot g(a) \leq 1$. Then, for
$U\subseteq \R_{\geq 0}^{|A|}$ one write $\simop U =
\set{g}{\allin{f}{U}{f \polar g}}$. 

A map $M\colon A \rightarrow B$ between PCS's is a ``matrix'' 
$M\in \R_{\geq 0}^{|A|\times|B|}$ 
}

\section{Conclusions and further work}

This paper develops a quantitative categorical logic via four
assumptions, starting from predicates as maps $X \rightarrow 1+1$ and
states as maps $1 \rightarrow X$, in a dual relationship. It involves
abstract forms of validity, predicate- and state-transformers,
measurement instruments, and test programs \& predicates. The theory
applies to standard set-theoretic, probabilistic, and quantum models.

The first version of the article appeared online in 2012. Since then
the research has been extended in several other publications. In
particular, \cite{Adams14} introduces a predicate logic based on
effects, \cite{JacobsWW15a} elaborates the role of
coproduct-preservation by the states functor in a state-and-effect
triangle~\eqref{GeneralTriangleDiag}, \cite{ChoJWW15} describes how
instruments (from Section~\ref{MeasurementSec} arise in many examples
from quotients and comprehension in the logic, \cite{Cho15a} gives an
equivalent description of the notion of effectus as a category of
partial maps, and~\cite{Jacobs15b} describes many state-and-effect
triangle examples. There is more as yet unpublished work, especially
giving characterisations of `Boolean' and `classical' (probabilistic)
effectuses.



\subsection*{Acknowledgements} The work presented in this article
acquired its current shape during a period of roughly three years of
research, discussion, revision and extension. The main sparring
partner was Bas Westerbaan to whom I owe a lot. Other colleagues that
contributed with comments and feedback include Robin Adams, Kenta Cho,
Tobias Fritz, Robert Furber, Ichiro Hasuo, Chris Heunen, Jorik
Mandemaker, Frank Roumen, Sam Staton and Bram Westerbaan. I thank them
all, and also the anonymous referees for their constructive comments.

\bibliographystyle{plain}

\begin{thebibliography}{10}

\bibitem{Abramsky10b}
S.~Abramsky.
\newblock No-cloning in categorical quantum mechanics.
\newblock In S.~Gay and I.~Mackie, editors, {\em Semantical Techniques in
  Quantum Computation}, pages 1--28. Cambridge Univ. Press, 2010.

\bibitem{AbramskyC04}
S.~Abramsky and B.~Coecke.
\newblock A categorical semantics of quantum protocols.
\newblock In {\em Logic in Computer Science}, pages 415--425. IEEE, Computer
  Science Press, 2004.

\bibitem{AbramskyC09}
S.~Abramsky and B.~Coecke.
\newblock A categorical semantics of quantum protocols.
\newblock In K.~Engesser, {Dov}~M. Gabbai, and D.~Lehmann, editors, {\em
  Handbook of Quantum Logic and Quantum Structures: Quantum Logic}, pages
  261--323. North Holland, Elsevier, Computer Science Press, 2009.

\bibitem{Adams14}
R.~Adams.
\newblock {QPEL}: Quantum program and effect language.
\newblock In B.~Coecke, I.~Hasuo, and P.~Panangaden, editors, {\em Quantum
  Physics and Logic (QPL) 2014}, number 172 in Elect. Proc. in Theor. Comp.
  Sci., pages 133--153, 2014.

\bibitem{ArbibM86}
M.~Arbib and E.~Manes.
\newblock {\em Algebraic Approaches to Program Semantics}.
\newblock Texts and Monogr. in Comp. Sci. Springer, Berlin, 1986.

\bibitem{Arveson81}
W.~Arveson.
\newblock {\em An Invitation to $C^*$-Algebra}.
\newblock Springer-Verlag, 1981.

\bibitem{Awodey06}
S.~Awodey.
\newblock {\em Category Theory}.
\newblock Oxford Logic Guides. Oxford Univ. Press, 2006.

\bibitem{BirkhoffN36}
G.~Birkhoff and J.~von Neumann.
\newblock The logic of quantum mechanics.
\newblock {\em Ann. Math.}, 37:823--843, 1936.

\bibitem{Borceux94}
F.~Borceux.
\newblock {\em Handbook of Categorical Algebra}, volume 50, 51 and 52 of {\em
  Encyclopedia of Mathematics}.
\newblock Cambridge Univ. Press, 1994.

\bibitem{Busch09}
P.~Busch.
\newblock Effect.
\newblock In D.~Greenberger, K.~Hentschel, and F.~Weinert, editors, {\em
  Compendium of Quantum Physics}, pages 179--180. Springer, 2009.

\bibitem{BuschLM96}
P.~Busch, P.~Lahti, and P.~Mittelstaedt.
\newblock {\em The Quantum Theory of Measurement}.
\newblock Springer, $2^{\textrm{nd}}$ edition, 1996.

\bibitem{BuschS98}
P.~Busch and J.~Singh.
\newblock L{\"u}ders theorem for unsharp quantum measurements.
\newblock {\em Phys. Letters A}, 249:10–--12, 1998.

\bibitem{CarboniLW93}
A.~Carboni, S.~Lack, and R.~Walters.
\newblock Introduction to extensive and distributive categories.
\newblock {\em Journ. of Pure \& Appl. Algebra}, 84(2):145--158, 1993.

\bibitem{Chang58}
C.~Chang.
\newblock Algebraic analysis of many-valued logics.
\newblock {\em Trans. Amer. Math. Soc.}, 88:476--490, 1958.

\bibitem{Cho14a}
K.~Cho.
\newblock Semantics for a quantum programming language by operator algebras.
\newblock In B.~Coecke, I.~Hasuo, and P.~Panangaden, editors, {\em Quantum
  Physics and Logic (QPL) 2014}, number 172 in Elect. Proc. in Theor. Comp.
  Sci., pages 165--190, 2014.

\bibitem{Cho15a}
K.~Cho.
\newblock Total and partial computation in categorical quantum foundations.
\newblock In C.~Heunen en~J.~Vicary, editor, {\em Quantum Physics and Logic
  (QPL) 2015}, Elect. Proc. in Theor. Comp. Sci., 2015.

\bibitem{ChoJWW15}
K.~Cho, B.~Jacobs, A.~Westerbaan, and B.~Westerbaan.
\newblock Quotient comprehension chains.
\newblock In C.~Heunen en~J.~Vicary, editor, {\em Quantum Physics and Logic
  (QPL) 2015}, Elect. Proc. in Theor. Comp. Sci., 2015.

\bibitem{ChovanecK07}
F.~Chovanec and F.~K{\^o}pka.
\newblock D-posets.
\newblock In K.~Engesser, {Dov}~M. Gabbai, and D.~Lehmann, editors, {\em
  Handbook of Quantum Logic and Quantum Structures: Quantum Structures}, pages
  367–--428. North Holland, Elsevier, Computer Science Press, 2007.

\bibitem{CoeckeP08}
B.~Coecke and D.~Pavlovi{\'c}.
\newblock Quantum measurements without sums.
\newblock In G.~Chen, L.~Kauffman, and S.~Lamonaco, editors, {\em Mathematics
  of Quantum Computing and Technology}, pages 559--596. Taylor and Francis,
  2008.

\bibitem{DaviesL70}
E.~Davies and J.~Lewis.
\newblock An operational approach to quantum probability.
\newblock {\em Communic. Math. Physics}, 17:239--260, 1970.

\bibitem{dHondtP06a}
E.~D'Hondt and P.~Panangaden.
\newblock Quantum weakest preconditions.
\newblock {\em Math. Struct. in Comp. Sci.}, 16(3):429--451, 2006.

\bibitem{Dieks82}
D.~Dieks.
\newblock Communication by {EPR} devices.
\newblock {\em Phys. Letters A}, 92(6):271–--272, 1982.

\bibitem{Dijkstra75}
E.~Dijkstra.
\newblock Guarded commands, nondeterminacy and formal derivation of programs.
\newblock {\em Communications of the ACM}, 18(8):453--457, 1975.

\bibitem{DvurecenskijP00}
A.~Dvure\v{c}enskij and S.~Pulmannov{\'a}.
\newblock {\em New Trends in Quantum Structures}.
\newblock Kluwer Acad. Publ., Dordrecht, 2000.

\bibitem{FoulisB94}
D.~J. Foulis and M.K. Bennett.
\newblock Effect algebras and unsharp quantum logics.
\newblock {\em Found. Physics}, 24(10):1331--1352, 1994.

\bibitem{FurberJ13a}
R.~Furber and B.~Jacobs.
\newblock From {Kleisli} categories to commutative {$C^*$-algebras}:
  Probabilistic {Gelfand} duality.
\newblock In R.~Heckel and S.~Milius, editors, {\em Conference on Algebra and
  Coalgebra in Computer Science (CALCO 2013)}, number 8089 in Lect. Notes Comp.
  Sci., pages 141--157. Springer, Berlin, 2013.

\bibitem{FurberJ13b}
R.~Furber and B.~Jacobs.
\newblock Towards a categorical account of conditional probability, 2013.
\newblock QPL 2013, see arxiv.org/abs/1306.0831.

\bibitem{Giry82}
M.~Giry.
\newblock A categorical approach to probability theory.
\newblock In B.~Banaschewski, editor, {\em Categorical Aspects of Topology and
  Analysis}, number 915 in Lect. Notes Math., pages 68--85. Springer, Berlin,
  1982.

\bibitem{GiuntiniG94}
R.~Giuntini and H.~Greuling.
\newblock Toward a formal language for unsharp properties.
\newblock {\em Found. Physics}, 19:769–--780, 1994.

\bibitem{GoguenB92}
J.~Goguen and R.~Burstall.
\newblock Institutions: Abstract model theory for specification and
  programming.
\newblock {\em Journ. ACM}, 39(1):95--146, 1992.

\bibitem{GudderG02}
S.~Gudder and R.~Greechie.
\newblock Sequential products on effect algebras.
\newblock {\em Reports on Math. Physics}, 49(1):87--111, 2002.

\bibitem{GudderN02}
S.~Gudder and G.~Nagy.
\newblock Sequential quantum measurements.
\newblock {\em Journ. Math. Physics}, 42:5212--5222, 2001.

\bibitem{Haghverdi00b}
E.~Haghverdi.
\newblock Unique decomposition categories, geometry of interaction and
  combinatory logic.
\newblock {\em Math. Struct. in Comp. Sci.}, 10:205--231, 2000.

\bibitem{HarelKT00}
D~Harel, D.~Kozen, and J.~Tiuryn.
\newblock {\em Dynamic Logic}.
\newblock {MIT} Press, Cambridge, MA, 2000.

\bibitem{HeinosaariZ12}
T.~Heinosaari and M.~Ziman.
\newblock {\em The Mathematical Language of Quantum Theory. From Uncertainty to
  Entanglement}.
\newblock Cambridge Univ. Press, 2012.

\bibitem{Jacobs94a}
B.~Jacobs.
\newblock Semantics of weakening and contraction.
\newblock {\em Ann. Pure \& Appl. Logic}, 69(1):73--106, 1994.

\bibitem{Jacobs99a}
B.~Jacobs.
\newblock {\em Categorical Logic and Type Theory}.
\newblock North Holland, Amsterdam, 1999.

\bibitem{Jacobs10e}
B.~Jacobs.
\newblock Convexity, duality, and effects.
\newblock In C.~Calude and V.~Sassone, editors, {\em IFIP Theoretical Computer
  Science 2010}, number 82(1) in IFIP Adv. in Inf. and Comm. Techn., pages
  1--19. Springer, Boston, 2010.

\bibitem{Jacobs11c}
B.~Jacobs.
\newblock Probabilities, distribution monads, and convex categories.
\newblock {\em Theor. Comp. Sci.}, 412(28):3323--3336, 2011.

\bibitem{Jacobs13a}
B.~Jacobs.
\newblock Measurable spaces and their effect logic.
\newblock In {\em Logic in Computer Science}. IEEE, Computer Science Press,
  2013.

\bibitem{Jacobs13b}
B.~Jacobs.
\newblock On block structures in quantum computation.
\newblock In D.~Kozen and M.~Mislove, editors, {\em Math. Found. of Programming
  Semantics}, volume 298 of {\em Elect. Notes in Theor. Comp. Sci.}, pages
  233--255. Elsevier, Amsterdam, 2013.

\bibitem{Jacobs15a}
B.~Jacobs.
\newblock {Dijkstra} and {Hoare} monads in monadic computation.
\newblock {\em TCS}, 2015.
\newblock DOI \url{http://dx.doi.org/10.1016/j.tcs.2015.03.020}.

\bibitem{Jacobs15b}
B.~Jacobs.
\newblock A recipe for state and effect triangles.
\newblock In L.~Moss and P.~Sobocinski, editors, {\em Conference on Algebra and
  Coalgebra in Computer Science (CALCO 2015)}, LIPIcs, 2015.

\bibitem{JacobsM12a}
B.~Jacobs and J.~Mandemaker.
\newblock Coreflections in algebraic quantum logic.
\newblock {\em Found. of Physics}, 42(7):932--958, 2012.

\bibitem{JacobsM12b}
B.~Jacobs and J.~Mandemaker.
\newblock The expectation monad in quantum foundations.
\newblock In B.~Jacobs, P.~Selinger, and B.~Spitters, editors, {\em Quantum
  Physics and Logic (QPL) 2011}, volume~95 of {\em Elect. Proc. in Theor. Comp.
  Sci.}, pages 143--182, 2012.

\bibitem{JacobsW15a}
B.~Jacobs and A.~Westerbaan.
\newblock An effect-theoretic account of {Lebesgue} integration, 2015.
\newblock Math. Found. of Programming Semantics XXXI.

\bibitem{JacobsWW15a}
B.~Jacobs, B.~Westerbaan, and A.~Westerbaan.
\newblock States of convex sets.
\newblock In A.~Pitts, editor, {\em Foundations of Software Science and
  Computation Structures}, number 9034 in Lect. Notes Comp. Sci., pages
  87--101. Springer, Berlin, 2015.

\bibitem{Johnstone82}
P.~Johnstone.
\newblock {\em Stone Spaces}.
\newblock Number~3 in Cambridge Studies in Advanced Mathematics. Cambridge
  Univ. Press, 1982.

\bibitem{Kadison51}
R.~Kadison.
\newblock A representation theory for commutative topological algebra.
\newblock {\em Memoirs of the {AMS}}, 7, 1951.

\bibitem{Kalmbach83}
G.~Kalmbach.
\newblock {\em Orthomodular Lattices}.
\newblock Academic Press, London, 1983.

\bibitem{Katrnoska92}
F.~Katrno\v{s}ka.
\newblock Logics of idempotents of rings.
\newblock In C.~Bandt, J.~Flachsmeyer, and H.~Haase, editors, {\em Topology,
  Measure, and Fractals}, number~66 in Math. Research, pages 131--136. Academie
  Verlag, 1992.

\bibitem{KellyL80}
M.~Kelly and M.~Laplaza.
\newblock Coherence for compact closed categories.
\newblock {\em Journ. of Pure \& Appl. Algebra}, 19:193--213, 1980.

\bibitem{Kock71b}
A.~Kock.
\newblock Bilinearity and cartesian closed monads.
\newblock {\em Math. Scand.}, 29:161--174, 1971.

\bibitem{Kock71a}
A.~Kock.
\newblock Closed categories generated by commutative monads.
\newblock {\em Journ. Austr. Math. Soc.}, XII:405--424, 1971.

\bibitem{Kozen81}
D.~Kozen.
\newblock Semantics of probabilistic programs.
\newblock {\em Journ. Comp. Syst. Sci}, 22(3):328--350, 1981.

\bibitem{Kozen85}
D.~Kozen.
\newblock A probabilistic {PDL}.
\newblock {\em Journ. Comp. Syst. Sci}, 30(2):162--178, 1985.

\bibitem{Kraus83}
K.~Kraus.
\newblock {\em States, Effects, and Operations}.
\newblock Springer Verlag, Berlin, 1983.

\bibitem{Landsman09b}
N.P. Landsman.
\newblock Algebraic quantum mechanics.
\newblock In D.~Greenberger, K.~Hentschel, and F.~Weinert, editors, {\em
  Compendium of Quantum Physics}, pages 6--10. Springer, 2009.

\bibitem{MartinLof84}
P.~{Martin-L{\"o}f}.
\newblock {\em Intuitionistic Type Theory}.
\newblock Bibliopolis, Napoli, 1984.

\bibitem{MacLane71}
S.~\mbox{Mac Lane}.
\newblock {\em Categories for the Working Mathematician}.
\newblock Springer, Berlin, 1971.

\bibitem{MacLaneM92}
S.~\mbox{Mac Lane} and I.~Moerdijk.
\newblock {\em Sheaves in Geometry and Logic. A First Introduction to Topos
  Theory}.
\newblock Springer, New York, 1992.

\bibitem{Moggi91a}
E.~Moggi.
\newblock Notions of computation and monads.
\newblock {\em Inf. \& Comp.}, 93(1):55--92, 1991.

\bibitem{Mundici11}
D.~Mundici.
\newblock {\em Advanced {\L}ukasiewicz calculus and {MV}-algebras}, volume~35
  of {\em Trends in Logic}.
\newblock Springer, 2011.

\bibitem{Neumann32}
J.~von Neumann.
\newblock {\em Mathematische Grundlagen der Quantenmechanik}.
\newblock Springer, Berlin, 1932.

\bibitem{NordstromPS90}
B.~Nordstr{\"o}m, K.~Peterson, and J.M. Smith.
\newblock {\em Programming in Martin-L{\"o}f's Type Theory: an introduction}.
\newblock Number~7 in Logic Guides. Oxford Science Publ., 1990.

\bibitem{Ozawa84}
M.~Ozawa.
\newblock Quantum measuring processes of continuous observables.
\newblock {\em Journ. Math. Physics}, 25:79--87, 1984.

\bibitem{Panangaden09}
P.~Panangaden.
\newblock {\em Labelled {Markov} Processes}.
\newblock Imperial College Press, 2009.

\bibitem{PulmannovaG98}
S.~Pulmannov{\'a} and S.~Gudder.
\newblock Representation theorem for convex effect algebras.
\newblock {\em Commentationes Mathematicae Universitatis Carolinae},
  39(4):645--659, 1998.

\bibitem{Reichenbach49}
H.~Reichenbach.
\newblock {\em The theory of probability}.
\newblock Univ. California Press, 1949.

\bibitem{Rennela14a}
M.~Rennela.
\newblock Towards a quantum domain theory: Order-enrichment and fixpoints in
  \mbox{$W^*$-algebras}.
\newblock In B.~Jacobs, A.~Silva, and S.~Staton, editors, {\em Math. Found. of
  Programming Semantics}, number 308 in Elect. Notes in Theor. Comp. Sci.,
  pages 289–--307, 2014.

\bibitem{RennelaS15a}
M.~Rennela and S.~Staton.
\newblock Complete positivity and natural representation of quantum
  computations, 2015.
\newblock Math. Found. of Programming Semantics XXXI.

\bibitem{RieffelP11}
E.~Rieffel and W.~Polak.
\newblock {\em Quantum Computing. A Gentle Introduction}.
\newblock {MIT} Press, Cambridge, MA, 2011.

\bibitem{Rudin87}
W.~Rudin.
\newblock {\em Functional Analysis}.
\newblock McGraw-Hill Book Company, 1987.
\newblock Third, intern. edition.

\bibitem{Selinger04}
P.~Selinger.
\newblock Towards a quantum programming language.
\newblock {\em Math. Struct. in Comp. Sci.}, 14(4):527--586, 2004.

\bibitem{Staton15a}
S.~Staton.
\newblock Algebraic effects, linearity, and quantum programming languages.
\newblock In {\em Principles of Programming Languages}, pages 395--406. ACM
  SIGPLAN-SIGACT, 2015.

\bibitem{Stone49}
M.~Stone.
\newblock Postulates for the barycentric calculus.
\newblock {\em Ann. Math.}, 29:25--30, 1949.

\bibitem{Strocchi05}
F.~Strocchi.
\newblock {\em An Introduction to the Mathematical Structure of Quantum
  Mechanics: a short Course for Mathematicians}, volume~27 of {\em Adv. Series
  in Math. Physics}.
\newblock World Scientific, 2005.

\bibitem{Swirszcz74}
T.~Swirszcz.
\newblock Monadic functors and convexity.
\newblock {\em Bull. de l'Acad. Polonaise des Sciences. S\'er. des sciences
  math., astr. et phys.}, 22:39--42, 1974.

\bibitem{Takesaki01}
M.~Takesaki.
\newblock {\em Theory of Operator Algebras {I}}, volume 124 of {\em
  Encyclopedia of Mathematical Sciences}.
\newblock Springer, $2^{\textrm{nd}}$ edition, 2001.

\bibitem{Takesaki03b}
M.~Takesaki.
\newblock {\em Theory of Operator Algebras {III}}, volume 127 of {\em
  Encyclopedia of Mathematical Sciences}.
\newblock Springer, 2003.

\bibitem{Tomiyama57a}
J.~Tomiyama.
\newblock On the projection of norm one in {$W^*$-algebras}.
\newblock {\em Proc. Japan Acad.}, 10:608--612, 1957.

\bibitem{ValironRSSS15}
B.~Valiron, N.~Ross, D.~Scott Alexander, P.~Selinger, and J.~Smith.
\newblock Programming the quantum future.
\newblock {\em Communications of the ACM}, 58(8):52--61, 2015.

\bibitem{WootersZ82}
W.~Wootters and W.~Zurek.
\newblock A single quantum cannot be cloned.
\newblock {\em Nature}, 299:802--803, 1982.

\end{thebibliography}

\appendix
\section{}




This appendix describes some background information for our running
examples, involving discrete probability, continuous probability,
Hilbert spaces, and $C^*$-algebras.

\subsection{Discrete probability}\label{DiscProbSubsec}

To describe finite discrete probabilities categorically one uses the
distribution monad $\Dst\colon\Sets \rightarrow \Sets$. It maps a set
$X$ to the set $\Dst(X)$ of probability distributions over $X$, which we
describe as formal finite convex sums:
$$\textstyle \sum_{i} r_{i}\ket{x_i}
\qquad\mbox{where}\qquad
x_{i}\in X \mbox{ and } r_{i}\in [0,1] \mbox{ satisfy } \sum_{i}r_{i}=1.$$

\noindent We use the ``ket'' notation $\ket{-}$ to distinguish
elements $x\in X$ and their occurrences in formal sums. Alternatively
one may describe $\Dst(X)$ as the set of functions $\varphi \colon X
\rightarrow [0,1]$ with finite support and $\sum_{x}\varphi(x) =
1$. We freely switch between these two descriptions. Each function
$f\colon X \rightarrow Y$ gives a function $\Dst(f) \colon \Dst(X)
\rightarrow \Dst(Y)$, where:
\begin{equation}
\label{DstMapEqn}
\begin{array}{rcl}
\Dst(f)\big(\sum_{i}r_{i}\ket{x_i}\big)
& = &
\sum_{i}r_{i}\ket{f(x_{i})}.
\end{array}
\end{equation}

\noindent The unit $\eta \colon X \rightarrow \Dst(X)$ of this
distribution monad $\Dst$ sends $x\in X$ to the singleton/Dirac
distribution $\eta(x) = 1\ket{x}$. The Kleisli extension $f_{*} \colon
\Dst(X) \rightarrow \Dst(Y)$ of a function $f\colon X \rightarrow
\Dst(Y)$ is defined by:
\begin{equation}
\label{DstKleisliExtEqn}
\begin{array}{rcl}
f_{*}(\varphi)(y)
& = &
\sum_{x} \varphi(x)\cdot f(x)(y).
\end{array}
\end{equation}


Like for any monad, one can form the Kleisli category $\Kl(\Dst)$. In
this case we get the category of sets and stochastic matrices, as the
objects of $\Kl(\Dst)$ are sets, and its maps $X\rightarrow Y$ are
functions $X \rightarrow \Dst(Y)$. The unit function $\eta\colon X
\rightarrow \Dst(X)$ is then the identity map $X \rightarrow X$ in
$\Kl(\Dst)$. Composition of $f\colon X \rightarrow Y$ and $g\colon Y
\rightarrow Z$ in $\Kl(\Dst)$ yields a map $g \klafter f \colon X
\rightarrow Z$, which, as a function $X \rightarrow \Dst(Z)$ is given
by $g \klafter f = g_{*} \after f$, using Kleisli extension
from~\eqref{DstKleisliExtEqn}.


There is a forgetful functor $\Kl(\Dst) \rightarrow \Sets$, sending
$X$ to $\Dst(X)$ and $f$ to $\mu \after \Dst(f)$. It has a left
adjoint $J \colon \Sets \rightarrow \Kl(\Dst)$ which is the identity
on objects and sends $f$ to $\eta\after f$. The Kleisli category
$\Kl(\Dst)$ has coproducts $X_{1}+X_{2}$, on objects like in
$\Sets$. There are coprojection maps $J(\kappa_{i}) = \eta \after
\kappa_{i} \colon X_{i} \rightarrow \Dst(X_{1}+X_{2})$, where
$\kappa_{i} \colon X_{i} \rightarrow X_{1}+X_{2}$ are the
coprojections in $\Sets$. Cotuples in $\Kl(\Dst)$ are like in $\Sets$.
We have $\Dst(1) = 1$ --- making $\Dst$ an `affine' functor. As a
result, the singleton/final set $1 = \{0\}$ is final in the category
$\Kl(\Dst)$. The copower $n\cdot 1$ in $\Kl(\Dst)$ is thus the
$n$-element set $n$.

The distribution monad $\Dst$ is `commutative': there exists a natural
transformation $\dst \colon \Dst(X) \times \Dst(Y) \rightarrow
\Dst(X\times Y)$, given by $\dst(\varphi,\psi)(x,y) =
\varphi(x)\cdot\psi(y)$. It makes the monad `monoidal'. Via this map
$\dst$ the product $\times$ on $\Sets$ becomes a tensor on
$\Kl(\Dst)$. Explicitly, $\dst$ is used for functoriality: for
functions $f_{i}\colon X_{i} \rightarrow \Dst(Y_{i})$ we can define a
function $f_{1}\otimes f_{2} = \dst \after (f_{1}\times f_{2}) \colon
X_{1}\times X_{2} \rightarrow \Dst(Y_{1}\times Y_{2})$. The final
object $1\in\Kl(\Dst)$ is the unit for this tensor.

\subsection{Continuous probability}\label{ContProbSubsec}

We write $\Meas$ for the category with measurable spaces $X =
(X,\Sigma_{X})$ as objects, where $\Sigma_{X} \subseteq \Pow(X)$ is a
$\sigma$-algebra. A measurable space $X$ is \emph{discrete} if all
subsets of $X$ are measurable, \textit{i.e.}~if $\Sigma_{X} =
\Pow(X)$. A map $f\colon X \rightarrow Y$ in $\Meas$ is a
`measurable' function, satisfying $f^{-1}(M)\in \Sigma_X$ for each
measurable subset $M\in \Sigma_Y$. We use the unit interval $[0,1]$ as
measurable space, with `Borel' measurable subsets generated by the
intervals $[q,1]$, where $q$ is a rational number in $[0,1]$.

\auxproof{
We have $\neg[q,1] = [0,q)\in\Sigma_{[0,1]}$. If $r\in[0,1]$ is a real
number approximated by an ascending chain of rationals $0 < q_{1} < q_{2}
< \cdots < r$, then $[0,r) = \bigcup_{i}[0,q_{i})\in\Sigma_{[0,1]}$.
Hence $[r,1] = \neg[0,r) \in \Sigma_{[0,1]}$. Thus the open interval $(s,r) = 
\neg([0,s] \cup [r,1])\in\Sigma_{[0,1]}$, for $s<r$.
}

The (discrete) empty and singleton spaces $0$ and $1$ are initial and
final in the category $\Meas$. The (categorical) product $X_{1}\times
X_{2}$ of two measurable spaces $X_{i}$ carries the least
$\sigma$-algebra making both projections $\pi_{i}\colon X_{1}\times
X_{2} \rightarrow X_{i}$ measurable functions; equivalently, this
$\sigma$-algebra is generated by the rectangles $M_{1}\times M_{2}$
with $M_{i}\in\Sigma_{X_i}$. The coproduct $X_{1}+X_{2}$ in $\Meas$
involves the disjoint union of the underlying sets with the
$\sigma$-algebra given by the direct images $\kappa_{i}M =
\set{\kappa_{i}x}{x\in M}$ for $M\in\Sigma_{X_i}$, where
$\kappa_{i}\colon X_{i} \rightarrow X_{1}+X_{2}$ is the coprojection
map. Thus, if $X_{1}, X_{2}$ are discrete, then so is the coproduct
$X_{1}+X_{2}$. In particular, the copower $n\cdot 1$ is the
$n$-element discrete space.



A \emph{probability space} consists of a measurable space $X =
(X,\Sigma_{X})$ together with a function $\phi\colon \Sigma_{X}
\rightarrow [0,1]$ which satisfies $\phi(X) = 1$ and is countably
additive: $\phi\left(\bigovee_{i\in I}M_{i}\right) = \sum_{i\in
  I}\phi(M_{i})$, for each pairwise disjoint, countable collection of
measurable subsets $M_{i}\in\Sigma_{X}$. Here we use $\ovee$ for
\emph{disjoint} union, where $\Sigma_X$ is understood as effect
algebra. We write:
$$\begin{array}{rcl}
\Giry(X)
& = &
\set{\phi\colon\Sigma_{X}\rightarrow [0,1]}{\phi 
   \mbox{ is a probability measure}}.
\end{array}$$

\noindent For each finite discrete space $n$ we have $\Giry(n) \cong
\Dst(n)$. In particular, $\Giry(2) \cong [0,1]$ and $\Giry(1) \cong
1$.

Each measurable subset $M\in\Sigma_{X}$ yields a function $\ev_{M}
\colon \Giry(X) \rightarrow [0,1]$, namely $\ev_{M}(\phi) =
\phi(M)$. Thus one can equip the set $\Giry(X)$ with the least
$\sigma$-algebra making all these maps $\ev_{M}$ measurable. We obtain
the `Giry' functor $\Giry \colon \Meas\rightarrow\Meas$
from~\cite{Giry82} (see also~\cite{Kozen81}): for a map $f\colon X
\rightarrow Y$ in $\Meas$ we get a measurable function $\Giry(f)
\colon \Giry(X) \rightarrow \Giry(Y)$ given by:
$$\begin{array}{rcl}
\smash{\Giry(f)\big(\Sigma_{X} \stackrel{\phi}{\rightarrow} [0,1]\big)}
& = &
\smash{\big(\Sigma_{Y}\stackrel{f^{-1}}{\rightarrow} 
   \Sigma_{X} \stackrel{\phi}{\rightarrow} [0,1]\big).}
\end{array}$$

\auxproof{
For a probability measure $\phi$ on $X\times Y$ one gets a
probability measure $\Giry(\pi_{1})(\phi)$ on $X$, which is the
\emph{marginal} of $\phi$. It is given on $M\in\Sigma_{X}$ by:
$$\begin{array}{rcccl}
\Giry(\pi_{1})(\phi)(M)
& = &
\phi\big(\pi_{1}^{-1}(M)\big) 
& = &
\phi\big(M\times Y\big).
\end{array}$$

\noindent Probability measures are closed under convex sums, making
$\Giry(X)$ a convex set: for a finite collection $\phi_{i}\in\Giry(X)$
and $r_{i}\in[0,1]$ with $\sum_{i}r_{i}=1$ one has
$\sum_{i}r_{i}\phi_{i} \in\Giry(X)$.
}

For a probability measure $\phi\in\Giry(X)$ we describe integration,
but only for measurable functions (predicates) $X \rightarrow [0,1]$,
with the unit interval as codomain, and not for more general real- or
complex-valued functions (see~\cite{JacobsW15a} for a general
effect-theoretic account). For each $M\in\Sigma_{X}$ we write
$\indic{M} \colon X \rightarrow [0,1]$ for the \emph{indicator
  function} given by $\indic{M}(x) = 1$ for $x\in M$ and $\indic{M}(x)
= 0$ for $x\not\in M$. The mapping $M \mapsto \indic{M}$ is a
homomorphism of effect algebras form $\Sigma$ to the predicates on
$X$, \textit{i.e.}~to the measurable functions $X \rightarrow
[0,1]$. A \emph{step function} is a finite linear combination
$r_{1}\indic{M_1} + \cdots + r_{k}\indic{M_k} =
\sum_{i}r_{i}\scalar\indic{M_i}$ of indicator functions with $r_{i}\in
    [0,1]$ and $M_{i}\in\Sigma_X$ pairwise disjoint measurable
    subsets. It is a basic fact that each measurable predicate
    $p\colon X \rightarrow [0,1]$ can be approximated (from below) by
    step functions: $p = \bigvee_{n\in\NNO}p_{n}$. The (Lebesgue)
    integral is then defined as:
$$\begin{array}{rclcrcl}
\int p\intd\phi
& = &
\bigvee_{n\in\NNO} \int p_{n}\intd\phi \,\in\, [0,1]
& \qquad\mbox{where}\qquad &
\int \big(\textstyle\sum_{i}r_{i}\indic{M_i}\big)\intd\phi
& = &
\sum_{i}r_{i}\phi(M_{i}).
\end{array}$$

\noindent Unit $\eta\colon X \rightarrow \Giry(X)$ and multiplication
maps $\mu \colon \Giry^{2}(X) \rightarrow \Giry(X)$ can now be defined,
making $\Giry$ a monad on the category $\Meas$:
$$\begin{array}{rclcrcl}
\eta(x)(M)
& = &
\indic{M}(x)
& \qquad &
\mu(\Phi)(M)
& = &
\int \ev_{M} \intd\Phi.
\end{array}$$

\noindent For details, we refer to~\cite{Panangaden09}
or~\cite{Jacobs13a,JacobsW15a}. The Kleisli category $\Kl(\Giry)$ of
the Giry monad contains measurable spaces as objects; a morphism $X
\rightarrow Y$ in $\Kl(\Giry)$ is a measurable map $f\colon X
\rightarrow \Giry(Y)$, which can be understood as a Markov transition
map. There is an associated Kleisli extension map, which we write as
$f_{*} \colon \Giry(X) \rightarrow \Giry(Y)$, namely:
\begin{equation}
\label{GiryKleisliExtEqn}
\begin{array}{rclcl}
f_{*}(\phi)(N)
& = &
\int f(-)(N)\intd\phi
& \qquad\mbox{where}\qquad &
f(-)(N) \colon X \rightarrow [0,1].
\end{array}
\end{equation}

\noindent Coproducts $(+,0)$ are inherited from $\Meas$.  The
final/singleton space $1$ is also final in $\Kl(\Giry)$, since
$\Giry(1) = 1$. Further, the Giry monad is also commutative, like the
distribution monad $\Dst$, via a map $\dst \colon \Giry(X) \times
\Giry(Y) \rightarrow \Giry(X\times Y)$ determined by
$\dst(\phi,\psi)(M\times N) = \phi(M)\times\psi(N)$.  As a result, the
product $\times$ of measurable spaces becomes a tensor on the Kleisli
category $\Kl(\Giry)$, with the final space $1$ as tensor unit.  For
more information, see~\cite{Panangaden09}.

\auxproof{
On this tensor unit we have:
$$\begin{array}{rcccl}
\Giry(1)
& = &
\set{\phi\colon \Sigma_{1} \rightarrow [0,1]}{\phi \mbox{ is a
   probability measure}} 
& \cong &
1,
\end{array}$$

\noindent since $\phi(\emptyset) = 0$ and $\phi(1) = 1$. Hence there
is precisely one element in $\Giry(1)$. This makes $\Giry$ an affine
monad.
}

\auxproof{
\begin{lem}
\label{StepFunAproxLem}
For each map $p\colon X \rightarrow [0,1]$ in $\Meas$ there is a
sequence of step functions $p_{n} \leq p$ so that $p$ can be written
both as:
\begin{itemize}
\item pointwise join $p = \bigvee_{n\in\NNO}p_{n}$;

\item limit $p = \lim\limits_{n\rightarrow\infty}p_{n}$ of a 
uniformly convergent sequence.
\end{itemize}
\end{lem}

\begin{myproof}
Following~\cite{JacobsM12b} we define $p_{n}(x) = 0.d_{1}d_{2}\cdots
d_{n}$, where $d_{i}$ is the $i$-th decimal of $p(x)\in[0,1]$. This
$p_{n}$ takes at most $10^{n}$ different values, since
$d_{i}\in\{0,1,\ldots,9\}$. For each of these values $r_{i}\in[0,1]$
there is a measurable subset $M_{i} = p^{-1}(\{r_{i}\})\in\Sigma_{X}$,
since the singleton subset/interval $\{r_{i}\} = [r_{i}, r_{i}]
\subseteq [0,1]$ is measurable. Thus we can write $p_{n} =
\sum_{i}r_{i}\indic{M_i}$, so that it is a step function.

By construction, $p_{n} \leq p$. For each $\epsilon > 0$, take
$N\in\NNO$ such that for all decimals $d_i$ we have:
$$\begin{array}{rcl}
0.\underbrace{00\cdots 00}_{N\text{ times}}d_{1}d_{2}d_{3}\cdots
& < &
\epsilon.
\end{array}$$

\noindent Then for each $n\geq N$ we have $p(x)-p_{n}(x) < \epsilon$,
for all $x\in X$, and thus $d(p,p_{n}) \leq \epsilon$. Hence
$\bigvee_{n}p_{n} = p$ and $p = \lim\limits_{n\rightarrow\infty}p_{n}$. \QED
\end{myproof}

Next we summarise the main steps in defining the (Lebesgue) integral
for measurable predicates.

\begin{defi}
\label{LebIntDef}
Let $(X,\Sigma_{X},\phi)$ be a probability space. 
\begin{enumerate}
\item For $M\in\Sigma_{X}$ the integral of the associated indicator
  function is defined as:
$$\begin{array}{rcl}
\displaystyle\int \indic{M}\intd\phi
& = &
\phi(M) \;\in\; [0,1].
\end{array}$$

\item This definition is extended linearly to step functions:
$$\begin{array}{rcl}
\displaystyle\int \big(\textstyle\sum_{i}r_{i}\indic{M_i}\big)\intd\phi
& = &
\sum_{i}r_{i}\phi(M_{i})  \;\in\; [0,1].
\end{array}$$

\noindent (This sum is in $[0,1]$ since: $\sum_{i}r_{i}\phi(M_{i})
\leq \sum_{i}\phi(M_{i}) = \phi(\ovee_{i}M_{i}) \leq \phi(X) = 1$.)

\item Next, this integral is extended continuously to all measurable
  functions $p\colon X \rightarrow [0,1]$; after writing them as limit
  $\smash{p = \lim\limits_{n\rightarrow\infty}p_{n}}$ of step
  functions $p_{n}$ like in Lemma~\ref{StepFunAproxLem}, one defines:
$$\begin{array}{rcl}
\displaystyle\int p\intd\phi
& = &
\displaystyle\lim_{n\rightarrow\infty}
   \int p_{n}\intd\phi \;\in\; [0,1].
\end{array}$$

\noindent This integral $\int p\intd\phi$ is sometimes written as $E[p]$,
since it describes the \emph{expectation value} of the predicate $p$.
\end{enumerate}
\end{defi}

We list some basic properties of integration.

\begin{lem}
\label{LebIntLem}
Let $X$ be a measurable space.
\begin{enumerate}
\item For each $\phi\in\Giry(X)$ the operation $p\mapsto\int
  p\intd\phi$ is a map of effect modules $\Pred(X) \rightarrow [0,1]$
  that preserves pointwise limits.

\item For a map $f\colon X \rightarrow Y$ in $\Meas$ and predicate
$q\colon Y \rightarrow [0,1]$,
\begin{equation}
\label{LebIntCompEqn}
\begin{array}{rcl}
\displaystyle\int (q \after f) \intd \phi
& = &
\displaystyle\int q \intd \Giry(f)(\phi).
\end{array}
\end{equation}



\item For each $x\in X$ and $p\in\Pred(X)$ one has:
\begin{equation}
\label{LebIntUnitEqn}
\begin{array}{rcl}
\displaystyle\int p\intd\eta(x)
& = &
p(x),
\end{array}
\end{equation}

\noindent where $\eta_{X} \colon X \rightarrow \Giry(X)$ is the unit
map given by:
\begin{equation}
\label{UnitEqn}
\begin{array}{rcl}
\eta_{X}(x)(M)
& = &
\indic{M}(x).
\end{array}
\end{equation}
\end{enumerate}
\end{lem}

\noindent This unit $\eta$ yields a natural transformation $\eta
\colon \idmap \Rightarrow \Giry$. 

The unit map is well-defined: each $\eta(x) \colon \Sigma_{X} \rightarrow
[0,1]$ is a probability measure, since $\eta(x)(\emptyset) = 0$,
$\eta_{X}(x)(X) = 1$, and for a countably additive collection $(M_{i})_{i\in I}$
one has:
$$\begin{array}{rcl}
\eta(x)(\bigcup_{i}M_{i}) = 1
& \Longleftrightarrow &
x \in \bigcup_{i}M_{i} \\
& \Longleftrightarrow &
\ex{!i}{x\in M_{i}}  \\
& \Longleftrightarrow &
\ex{!i}{\eta(x)(M_{i}) = 1}  \\
& \Longleftrightarrow &
\sum_{i}\eta(x)(M_{i}) = 1.
\end{array}$$

\noindent Moreover, $\eta$ is natural: for $f\colon X \rightarrow Y$
in $\Meas$ one has:
$$\begin{array}{rcl}
\big(\Giry(f) \after \eta_{X}\big)(x)(M)
& = &
\Giry(f)(\eta(x))(M) \\
& = &
\big(\eta(x) \after f^{-1}\big)(M) \\
& = &
\eta(x)(f^{-1}(M)) \\
& = &
\left\{\begin{array}{ll}
1 \quad & \mbox{if }x\in f^{-1}(M) \\
0 & \mbox{otherwise}
\end{array}\right. \\
& = &
\left\{\begin{array}{ll}
1 \quad & \mbox{if }f(x)\in M \\
0 & \mbox{otherwise}
\end{array}\right. \\
& = &
\eta_{Y}(f(x))(M) \\
& = &
\big(\eta_{Y} \after f\big)(x)(M).
\end{array}$$

One may call the space $X$ separable if this unit map $\eta_X$ is
injective: $\allin{M}{\Sigma_X}{x\in M \Leftrightarrow x'\in M}$
implies $x=x'$.

\begin{enumerate}
\item Obviously:
$$\begin{array}{rcccccl}
\int 1\intd\phi
& = &
\int \indic{X}\intd\phi
& = &
\phi(X)
& = &
1 \\
\int 0\intd\phi
& = &
\int \indic{\emptyset}\intd\phi
& = &
\phi(\emptyset)
& = &
0.
\end{array}$$

\noindent The properties $\int (p\ovee q)\intd\phi = (\int p\intd\phi)
+ (\int q\intd\phi)$ and $\int (s \scalar p)\intd\phi = s\cdot (\int
p\intd\phi)$. The preservation of joins of ascending chains is
known as the Monotone Convergence Theorem~\cite{Panangaden09}.

\item For a step function $s = \sum_{i}r_{i}\indic{M_i}$ we have:
$$\begin{array}{rcccccl}
\displaystyle\int s\intd \eta(x)
& = &
\sum_{i} r_{i}\eta(x)(M_{i})
& = &
\sum_{i} r_{i}\indic{M_i}(x)
& = &
s(x).
\end{array}$$

\noindent Hence for a limit $\smash{p = \lim\limits_{n\rightarrow\infty}p_{n}}$
of step functions:
$$\begin{array}{rcccccl}
\displaystyle\int p\intd\eta(x)
& = &
\displaystyle\lim_{n\rightarrow\infty}\int p_{n}\intd\eta(x)
& = &
\displaystyle\lim_{n\rightarrow\infty} p_{n}(x)
& = &
p(x).
\end{array}$$

\item For a map $f\colon X \rightarrow Y$ and step function $s =
  \sum_{i}r_{i}\indic{N_i}$ with $N_{i}\in\Sigma_Y$ we write $t =
  \sum_{i}r_{i}\indic{f^{-1}(N_{i})}$. Then $t = s \after f$ since:
$$\begin{array}{rcccccl}
t(x) = r_{i}
& \Longleftrightarrow &
x\in f^{-1}(N_{i})
& \Longleftrightarrow &
f(x) \in N_{i}
& \Longleftrightarrow &
s(f(x)) = r_{i}.
\end{array}$$

\noindent Hence:
$$\begin{array}{rcl}
\displaystyle\int (s \after f) \intd \phi
\hspace*{\arraycolsep} = \hspace*{\arraycolsep}
\displaystyle\int t \intd \phi 
& = &
\sum_{i}r_{i}\phi\big(f^{-1}(N_{i})\big) \\[-.5em]
& = &
\sum_{i}r_{i}\Giry(f)(\phi)(N_{i}) 
\hspace*{\arraycolsep} = \hspace*{\arraycolsep}
\displaystyle\int s \intd \Giry(f)(\phi).
\end{array}$$

\noindent Hence for a limit $\smash{p = \lim\limits_{n\rightarrow\infty}p_{n}}$
of step functions:
$$\begin{array}{rcccccl}
\displaystyle\int (p \after f)\intd\phi
& = &
\displaystyle\lim_{n\rightarrow\infty}\int (p_{n} \after f)\intd\phi
& = &
\displaystyle\lim_{n\rightarrow\infty} \int p_{n} \intd\Giry(f)(\phi)
& = &
\displaystyle\int p \intd\Giry(f)(\phi).
\end{array}$$
\end{enumerate}

The next definition introduces two operations that are of fundamental
importance in this setting.

\begin{defi}
\label{SubstKleisliExtDef}
With an arbitrary measurable function $f\colon X
\rightarrow \Giry(Y)$ we associate two operations:
\begin{enumerate}
\item ``Kleisli extension'' $f^{\$} \colon \Giry(X) \rightarrow \Giry(Y)$,
given by:
\begin{equation}
\label{KleisliExtEqn}
\begin{array}{rcl}
f^{\$}(\phi)(N)
& = &
\displaystyle\int f(-)(N)\intd\phi \\
& = &
\displaystyle\int\big(\lamin{x}{X}{f(x)(N)}\big)\intd\phi.
\end{array}
\end{equation}

\noindent This uses that for $N\in\Sigma_Y$ one has a measurable
function $f(-)(N) \colon X \rightarrow [0,1]$.

\auxproof{
We check that $f(-)(N) \colon X \rightarrow [0,1]$ is measurable:
for $q\in [0,1]\cap\mathbb{Q}$ one has:
$$\begin{array}{rcl}
\big(f(-)(N)\big)^{-1}([q,1])
& = &
\setin{x}{X}{f(x)(N) \in [q,1]} \\
& = &
\setin{x}{X}{\ev_{N}(f(x)) \in [q,1]} \\
& = &
\setin{x}{X}{f(x) \in \ev_{N}^{-1}([q,1])} \\
& = &
f^{-1}\big(\ev_{N}^{-1}([q,1])\big).
\end{array}$$

\noindent The latter is a measurable subset of $X$ because $f \colon X
\rightarrow \Giry(Y)$ is a measurable function. 
}

\item ``Substitution'' $f^{*} \colon \Pred(Y) \rightarrow \Pred(X)$ given by:
\begin{equation}
\label{SubstEqn}
\begin{array}{rcccl}
f^{*}(q)
& = &
\displaystyle\int q\intd f(-)
& = &
\lamin{x}{X}{\displaystyle\int q\intd f(x)}.
\end{array}
\end{equation}
\end{enumerate}
\end{defi}

Since integration $\int (-)\intd\phi$ is a limit-preserving map of
effect modules (see Lemma~\ref{LebIntLem}), so is the substitution map
$f^*$, in a pointwise manner.

We check that substitution preserves the effect module structure
(pointwise):
$$\begin{array}{rcl}
f^{*}(p\ovee q)(x)
& = &
\int (p + q) \intd f(x) \\
& = &
\int p \intd f(x) + \int q \intd f(x) \\
& = &
f^{*}(p)(x) \ovee f^{*}(q)(x) \\
f^{*}(1)(x)
& = &
\int \indic{X}\intd f(x) \\
& = &
f(x)(X) \\
& = &
1 \\
f^{*}(s\scalar p)(x)
& = &
\int (s\cdot p) \intd f(x) \\
& = &
s\cdot \int p \intd f(x) \\
& = &
s\cdot f^{*}(p)(x) \\
& = &
(s \scalar f^{*}(p))(x) \\
f^{*}(\lim\limits_{n\rightarrow\infty}p_{n})(x)
& = &
\int (\lim\limits_{n\rightarrow\infty}p_{n}) \intd f(x) \\
& = &
\lim\limits_{n\rightarrow\infty}\int p_{n} \intd f(x) \\
& = &
\lim\limits_{n\rightarrow\infty} f^{*}(p_{n})(x) \\
& = &
(\lim\limits_{n\rightarrow\infty} f^{*}(p_{n}))(x) 
\end{array}$$

These operations of Kleisli extension $f^{\$}$ and substitution
$f^{*}$ are related in a basic manner, resembling a Galois connection.
This seemingly new observation gives a short proof of
Theorem~\ref{GiryMonadThm}.

\begin{prop}
\label{SubstKleisliExtProp}
For each map $f\colon X \rightarrow \Giry(Y)$ in $\Meas$, probability
measure $\phi\in\Giry(Y)$ and predicate $q\in\Pred(Y)$ one has:
$$\begin{array}{rcl}
\displaystyle\int f^{*}(q) \intd \phi 
& = &
\displaystyle\int q \intd f^{\$}(\phi).
\end{array}$$
\end{prop}

\begin{myproof}
Because of limit-preservation of substitution and integration it
suffices to prove the result for predicates given by step functions $s
= \sum_{i}r_{i}\indic{N_{i}} \in \Pred(Y)$. Then:
$$\begin{array}[b]{rcl}
\int f^{*}(s) \intd \phi 
& = &
\int f^{*}\big(\sum_{i}r_{i}\indic{N_{i}}\big) \intd \phi \\
& = &
\int \sum_{i}r_{i}f^{*}\big(\indic{N_{i}}\big) \intd \phi \\
& & \qquad
   \rlap{since $f^{*}$ is a map of effect modules} \\
& \smash{\stackrel{\eqref{SubstEqn}}{=}} &
\sum_{i}r_{i}\int \big(\lam{x}{\int \indic{N_{i}}\intd f(x)}\big) \intd \phi \\
& = &
\sum_{i}r_{i}\int f(-)(N_{i}) \intd \phi \\
& \smash{\stackrel{\eqref{KleisliExtEqn}}{=}} &
\sum_{i}r_{i}f^{\$}(\phi)(N_{i}) \\
& = &
\sum_{i}r_{i}\int \indic{N_{i}} \intd f^{\$}(\phi) \\
& = &
\int \sum_{i}r_{i}\indic{N_{i}} \intd f^{\$}(\phi) \\
& = &
\int s \intd f^{\$}(\phi).
\end{array}\hspace*{3em}\eqno{\qEd}$$
\end{myproof}

We are finally in a position to see that $\Giry$ is a monad. We do so
by following the formulation in terms of Kleisli extension.

\begin{theorem}[From~\cite{Giry82}]
\label{GiryMonadThm}
The functor $\Giry\colon\Meas\rightarrow\Meas$ is a monad, with
unit $\eta$ from~\eqref{UnitEqn} and Kleisli extension $(-)^{\$}$ 
from~\eqref{KleisliExtEqn}.
\end{theorem}

\begin{myproof}
We check the equations for Kleisli extension: the unit equations
$\eta^{\$} = \idmap$ and $f^{\$} \after \eta = f$ are obtained as
follows.
$$\begin{array}{rcl}
\eta^{\$}(\phi)(M)
& = &
\int \eta(-)(M)\intd\phi \\
& \smash{\stackrel{\eqref{UnitEqn}}{=}} &
\int \indic{M}\intd\phi \\
& = &
\phi(M) \\
\big(f^{\$} \after \eta\big)(x)(N) 
& = &
f^{\$}(\eta(x))(N) \\
& \smash{\stackrel{\eqref{KleisliExtEqn}}{=}} &
\int f(-)(N)\intd\eta(x) \\
& \smash{\stackrel{\eqref{LebIntUnitEqn}}{=}} &
f(x)(N).
\end{array}$$

\noindent The composition equation $g^{\$} \after f^{\$} = (g^{\$} \after f)^{\$}$
requires a bit more care:
$$\begin{array}[b]{rcl}
\big(g^{\$} \after f^{\$}\big)(\phi)(K)
& = &
g^{\$}(f^{\$}(\phi))(K) \\
& \smash{\stackrel{\eqref{KleisliExtEqn}}{=}} &
\int g(-)(K) \intd f^{\$}(\phi) \\
& = &
\int f^{*}\big(g(-)(K)\big) \intd \phi \\
& & \qquad \mbox{by Proposition~\ref{SubstKleisliExtProp}} \\
& \smash{\stackrel{\eqref{SubstEqn}}{=}} &
\int \big(\lam{x}{\int g(-)(K)\intd f(x)}\big) \intd \phi \\
& \smash{\stackrel{\eqref{KleisliExtEqn}}{=}} &
\int \big(\lam{x}{g^{\$}(f(x))(K)}\big) \intd \phi \\
& = &
\int (g^{\$} \after f)(-)(K) \intd \phi \\
& = &
(g^{\$} \after f)^{\$}(\phi)(K).
\end{array}\eqno{\qEd}$$
\end{myproof}

As a result, composition in the Kleisli category $\Kl(\Giry)$ is given
as follows. For $f\colon X \rightarrow \Giry(Y)$ and $g\colon Y
\rightarrow \Giry(Z)$ we have:
\begin{equation}
\label{KleisliCompEqn}
\begin{array}{rcl}
(g \after f)(x)(K)
& = &
\displaystyle\int g(-)(K)\intd f(x)
\end{array}
\end{equation}

\noindent where $x\in X$ and $K\in\Sigma_{Z}$.

The multiplication $\mu \colon \Giry^{2}(X) \rightarrow \Giry(X)$
of the monad is given on $\Phi\in\Giry^{2}(X)$ and $M\in\Sigma_{X}$ by:
\begin{equation}
\label{MuEqn}
\begin{array}{rcl}
\mu(\Phi)(M)
\hspace*{\arraycolsep} = \hspace*{\arraycolsep}
\big(\idmap[\Giry(X)]\big)^{\$}(\Phi)(M)
& = &
\displaystyle\int\idmap(-)(M)\intd\Phi \\
& = &
\displaystyle\int\ev_{M}\intd\Phi.
\end{array}
\end{equation}
}

\subsection{Hilbert spaces}\label{HilbSubsec}

A Hilbert space is a vector space (over $\C$) with an inner product
such that the induced norm makes the space complete. We write $\Hilb$
for the category of Hilbert spaces, with linear and bounded
(equivalently: continuous) maps between them. This category is
symmetric monoidal, with the complex numbers $\C$ as unit for the
tensor $\otimes$.  There are also finite biproducts $(\oplus, 0)$,
given by the cartesian product of the underlying space. Tensors
$\otimes$ distribute over finite biproducts. Further, $\Hilb$ is a
`dagger category', where $f^{\dag}$ is the conjugate transpose of a
map $f$. A map $f$ is an isometry (aka.\ dagger monic) if $f^{\dag}
\after f = \idmap$, and a unitary if not only $f^{\dag} \after f =
\idmap$ but also $f \after f^{\dag} = \idmap$. In that case $f^{\dag}$
is the inverse of $f$. For a Hilbert space $\H$ we write $\B(\H)$ for
the set of endomaps $\H \rightarrow \H$ in $\Hilb$. There is a subset
$\Ef(\H) = \setin{A}{\B(\H)}{0 \leq A \leq \idmap}$ of
effects. Hilbert spaces will be used mostly as source of examples for
$C^*$-algebras, via this $\B(-)$ construction.

\subsection{$C^*$-algebras}\label{CstarSubsec}

In the current context we use $C^*$-algebra over the complex numbers
$\C$, with a multiplicative unit element $1$.  Thus a $C^*$-algebra is
a vector space over $\C$ which is *-algebra --- with multiplication
$\cdot$, involution $(-)^*$, and unit $1$ --- and carries a norm
$\|-\|$ that makes the space complete and interacts with
multiplication via $\|a\cdot b\| \leq \|a\| \cdot \|b\|$ and
$\|a^{*}\cdot a\| = \|a\|^{2}$. A $C^*$-algebra is called
\emph{commutative} if its multiplication is commutative, and
\emph{finite-dimensional} if it has finite dimension when considered
as a vector space.

An element $a\in A$ of a $C^*$-algebra $A$ is called self-adjoint if
$a^{*} = a$ and positive if it is of the form $a = x^{*}x$ for some
element $x\in A$. There is a partial order $a\leq b$ iff $b - a$ is
positive.

A linear map $f\colon A \rightarrow B$ between $C^*$-algebras $A,B$ is
called:
\begin{enumerate}
\item multiplicative (M) if $f(a\cdot a') = f(a) \cdot f(a')$ for all
  $a,a'\in A$;

\item involutive (I) if $f(a^{*}) = f(a)^{*}$ for all $a\in A$;

\item unital (U) if $f(1) = 1$, and subunital (sU) if $0 \leq f(1)
  \leq 1$;

\item positive (P) if $f(a) \geq 0$ for each $a\geq 0$;

\item completely positive (CP) if for each $n\in\NNO$ the map
  $\Mat_{n}(f) \colon \Mat_{n}(A) \rightarrow \Mat_{n}(B)$ is
  positive, where $\Mat_{n}(A)$ is the $C^*$-algebra of $n\times n$
  matrices with entries from $A$, with matrix multiplication.
\end{enumerate}

\noindent We call $f$ a MIU-map if it is multiplicative, involutive,
and unital, and a (C)PU-map if it (completely) positive and
unital. There are implications MIU $\Rightarrow$ CPU $\Rightarrow$ PU
$\Rightarrow$ I. A MIU-map is traditionally called a
*-homomorphism. If either $A$ or $B$ is commutative, then a positive
map $f\colon A \rightarrow B$ is automatically completely positive. We
use categories of $C^*$-algebras $\CstarMIU \hookrightarrow \CstarCPU
\hookrightarrow \CstarPU$ with maps as indicated by the
subscripts. There are obvious subcategories $\CCstarPU$ and
$\CCstarMIU$ of \emph{commutative} $C^*$-algebras with PU/MIU
maps. The famous Gelfand duality involves an equivalence $\CH \simeq
\op{(\CCstarMIU)}$ where $\CH$ is the category of compact Hausdorff
spaces and continuous maps.  For PU-maps an equivalence $\Kl(\Rad)
\simeq \op{(\CCstarPU)}$ is shown in~\cite{FurberJ13a}, where $\Rad$
is the Radon monad on $\CH$.

The collection of endomaps $\B(\H)$ of a Hilbert space $\H$ is a
$C^*$-algebra, with multiplication given by composition. This mapping
$\H \mapsto \B(\H)$ gives rise to functors:
$$\xymatrix{
\Hilb_{\textrm{unit}}\ar[r]^-{\B} & {\op{(\CstarMIU)}}
& &
\Hilb_{\textrm{isom}}\ar[r]^-{\B} & {\op{(\CstarCPU)}}
}$$

\noindent where $\Hilb_{\textrm{unit}}$,
resp.\ $\Hilb_{\textrm{isom}}$, is the category of Hilbert spaces with
unitaries, resp.\ isometries, as morphisms. A map $f\colon \H
\rightarrow \K$ is sent to $\B(f) = f^{\dag}(-)f \colon \B(\K)
\rightarrow \B(\H)$.

The algebra $\C$ of complex numbers is initial among $C^*$-algebras,
in all categories $\CstarMIU$, $\CstarCPU$, $\CstarPU$. Similarly,
there are finite products of $C^*$-algebras, via the trivial singleton
space $\{0\}$ as final object, and via direct sums $\oplus$ of vector
spaces (\textit{i.e.}~products of the underlying sets) as cartesian
products. The operations are used pointwise. There are also tensor
products of $C^*$-algebras. These are described in more detail
in~\cite[section IV.4]{Takesaki01}, but we outline them here. The
$C^*$-tensors for two $C^*$-algebras $A$ and $B$ are obtained by
taking the usual tensor of underlying vector spaces $A \otimes B$,
defining a *-algebra structure as follows:
$$\begin{array}{rclcrcl}
(a_1 \sotimes b_1)(a_2 \sotimes b_2) 
& = & 
(a_1a_2) \sotimes (b_1b_2)
& \qquad &
(a \sotimes b)^* 
& = &
a^* \sotimes b^*.
\end{array}$$

\noindent One then obtains a $C^*$-algebra by introducing a $C^*$-norm
compatible with the *-algebra structure, and taking the
completion. There are minimal and maximal, or injective and projective
$C^*$-norms, but if $A$ or $B$ is finite dimensional these coincide,
see~\cite[chapter XV, 1.4-1.6]{Takesaki03b}. Since $A \otimes B$ is
the completion of the algebraic tensor of $A$ and $B$, the span of
elements of the form $a \sotimes b$ is dense, and in fact in the
finite dimensional case, $A \otimes B$ is just the linear span of such
elements.

The positive cone of $A \otimes B$ contains the positive elements
according to the multiplication and involution, \textit{i.e.}~$a \in A
\otimes B$ is positive if $a = x^*x$ for some other element $x$ in the
tensor product. We note at this point that this cone is larger than
the cone obtained by taking sums of elements $a \sotimes b$ with $a
\in A$ and $b \in B$ both positive. The effect of this is that no
$C^*$-tensor is a functor on $\CstarPU$. Instead, the maps that can be
tensored are the completely positive ones (see \cite[\S IV.3, and \S
  IV.4 Prop.~23]{Takesaki01}). In short, it is the category
$\CstarCPU$ that is symmetric monoidal, with the (initial) object $\C$
as tensor unit. In~\cite{Cho14a} it is shown that the ``minimal''
tensor $\otimes$ on $\CstarCPU$ distributes over finite products
$(\oplus, \{0\})$. We use this result in the form: the opposite
category $\op{(\CstarCPU)}$ is symmetric monoidal with the (final)
object $\C$ as tensor unit, and $\otimes$ distributing over finite
coproducts.

For $n\in\NNO$ one can write the matrix algebra $\Mat_{n}(A)$ as
tensor product $\Mat_{n}(A) \cong \Mat_{n}(\C) \otimes A = \B(\C^{n})
\otimes A$.  In particular, when $A$ is of the form $\B(\H)$ we can
move the tensor inside $\B(-)$ as in:
\begin{equation}
\label{HilbCstarTensorEqn}
\begin{array}{rcccccl}
\B(\C^{n}) \otimes \B(\H)
& \cong &
\Mat_{n}(\B(\H))
& \cong &
\B(\H^{n})
& \cong &
\B(\C^{n}\otimes \H),
\end{array}
\end{equation}

\noindent where $\H^{n}$ is the $n$-fold biproduct $\H \oplus \cdots
\oplus \H \cong \C^{n}\otimes \H$.

A $W^*$-algebra (also called a von Neumann algebra) is a special kind
of $C^*$ algebra that satisfies additional closure properties,
see \textit{e.g.}~\cite{Takesaki01} for details.

\end{document}